\documentclass[11pt]{article}
\usepackage[T1]{fontenc}
\usepackage{lmodern,microtype,amsmath,amssymb,amsthm,mathtools,esint}
\usepackage[a4paper,margin=0.9in]{geometry}
\usepackage[colorlinks=true,linkcolor=blue]{hyperref}
\hypersetup{pdftitle={Existence of complete singular solutions of the Yamabe equation at prescribed points},pdfauthor={Zheng-Chao Han, Qinfeng Jiang, Hua-Yang Wang, Jingang Xiong, Lei Zhang},pdfsubject={Fixed neck variational construction in a positive conformal class}}
\usepackage{enumitem,array,longtable,booktabs,mathrsfs,needspace}
\setlist{nosep,leftmargin=*}
\allowdisplaybreaks[2]
\numberwithin{equation}{section}
\numberwithin{table}{section}
\numberwithin{figure}{section}
\newtheorem{theorem}{Theorem}[section]
\newtheorem{lemma}[theorem]{Lemma}
\newtheorem{proposition}[theorem]{Proposition}
\newtheorem{corollary}[theorem]{Corollary}
\theoremstyle{definition}
\theoremstyle{remark}

\theoremstyle{plain}
\newtheorem*{proposition*}{Proposition}
\newcommand{\thistheoremname}{}
\newtheorem*{genericthm*}{\thistheoremname}
\newenvironment{namedthm}[1]
  {\renewcommand{\thistheoremname}{#1}\begin{genericthm*}}
  {\end{genericthm*}}
  
\newcommand{\R}{\mathbb R}
\newcommand{\Sph}{\mathbb S}
\newcommand{\eps}{\varepsilon}
\newcommand{\dd}{\,\mathrm{d}}
\newcommand{\vol}{\dd V_g}
\newcommand{\norm}[1]{\left\|#1\right\|}
\newcommand{\K}{\mathcal K}
\newcommand{\Ff}{\mathcal F}
\newcommand{\amp}{a}
\newcommand{\Om}{\Omega_\rho}
\newcommand{\alphaH}{\alpha_{\rm H}}
\newcommand{\bgpar}{{\boldsymbol{\eta}}}
\newcommand{\bgshape}{{\mathbf s}}
\newcommand{\bgcoord}{{\boldsymbol{\pi}_E}}
\newcommand{\radsrc}{{\mathfrak f}}
\newcommand{\energyvar}{{\mathfrak e}}
\newcommand{\kelpar}{{\mathbf k}}
\newcommand{\trpar}{{\mathfrak t}}
\newcommand{\Xtrace}{{\mathcal X_{\rm tr}}}
\newcommand{\Ytrace}{{\mathcal Y_{\rm tr}}}
\newcommand{\Zcollar}{{\mathcal Z_{\rm col}}}
\newcommand{\trlow}{{\mathbf b_{\rm tr}}}
\newcommand{\uF}{{u^{\mathrm F}}}

\title{\bfseries Existence of complete singular solutions\\
of the Yamabe equation at prescribed points}
\author{Zheng-Chao Han,\quad Qinfeng Jiang,\quad
Hua-Yang Wang\thanks{H.-Y. Wang was partially supported by NSFC grant 12601183.}\\[2mm]
Jingang Xiong\thanks{J. Xiong was partially supported by NSFC grant 12325104.},
\quad Lei Zhang\thanks{Lei Zhang was partially supported by Simons Foundation
Grant SFI-MPS-TSM-00013752.}}
\date{}
\begin{document}
\maketitle
\begin{abstract}
Let $(M^n,g)$ be smooth, connected and closed, with $n\ge3$, and let
$P$ be a prescribed nonempty finite set. We characterize the conformal
classes $[g]$ that admit a complete metric of positive constant scalar
curvature on $M\setminus P$. Such a metric exists precisely when
$Y(M,[g])>0$ and either $[g]$ is nonround or $|P|\ge2$.
In a nonround positive conformal class, the metrics can be chosen with
Fowler asymptotics and arbitrary fixed ratios of their small asymptotic
neck sizes.
As the neck sizes tend to zero, every sequence of these metrics has a
subsequence converging smoothly away from the punctures to a smooth
Yamabe minimizing metric.
The construction allows degenerate Yamabe minimizers and arbitrary
curvature jets at the prescribed points. For nonempty closed singular sets $\Lambda\subset M$ with empty
    interior, we also announce the same existence criterion with the
    additional requirement that the Bessel capacity
    $C_{1+2/n,n/2}(\Lambda)$ vanish.
\end{abstract}

\section{Introduction}\label{sec:introduction}
The positive singular Yamabe problem concerns conformal metrics on
punctured compact manifolds. One seeks a complete metric of positive
constant scalar curvature in a prescribed conformal class. We give a
characterization for finite puncture sets. Throughout, conformal
equivalence means equivalence by a conformal diffeomorphism. The
conformal class on a punctured manifold is the restriction of the
given smooth conformal class.
Write $Y(M,[g])$ for the Yamabe constant. A conformal class is
nonround if it is not conformally equivalent to $(\Sph^n,[g_{\rm rd}])$.

\begin{theorem}\label{thm:main}
Let $(M^n,g)$ be smooth, connected and closed, with $n\ge3$.
Let $P=\{p_1,\ldots,p_m\}$ be a nonempty finite set. Then
$M\setminus P$ carries a complete metric in $[g]|_{M\setminus P}$
of positive constant scalar curvature if and only if
\[
 Y(M,[g])>0\quad\text{and}\quad
 \bigl((M,[g])\not\simeq(\Sph^n,[g_{\rm rd}])\ \text{or}\ m\ge2\bigr).
\]
\end{theorem}

We use \emph{Fowler asymptotics} in the following sense. In local normal
coordinates centered at a puncture, after a smooth conformal normalization
of the background, write $r$ for the radial coordinate and $u$ for the
corresponding conformal factor. There are a positive periodic Fowler
profile $v_\eps$ and a phase $\tau$ such that
\[
 r^{(n-2)/2}u(r,\theta)-v_\eps(-\log r+\tau)\longrightarrow0
 \quad\text{as }r\to0,
\]
uniformly in the angular variable. The profile satisfies
\begin{equation}
 v_\eps''-\frac{(n-2)^2}{4}v_\eps
       +\frac{n(n-2)}4v_\eps^{(n+2)/(n-2)}=0,
 \qquad \min v_\eps=\eps.
 \label{eq:fowler-ode}
\end{equation}
Its minimum $\eps$ is the \emph{neck parameter} for scalar curvature
$n(n-1)$. The associated metric on a half-cylinder is
asymptotic to $v_\eps(t+\tau)^{4/(n-2)}(\dd t^2+g_{\Sph^{n-1}})$.
Such an asymptotic geometric end is also called a Delaunay end.
In a nonround positive conformal class, we can prescribe arbitrary fixed
positive ratios of these small neck parameters.
\begin{corollary}\label{cor:neck-ratios}
Let $(M^n,g)$ be smooth, connected and closed, with $n\ge3$.
Assume $Y(M,[g])>0$ and that $[g]$ is nonround.
Let $P=\{p_1,\ldots,p_m\}$ be a nonempty finite set.
Fix $t_1,\ldots,t_m>0$. There is $\eps_0>0$ such that, for
each $0<\eps<\eps_0$, there is a complete metric
$g_\eps\in[g]|_{M\setminus P}$ of scalar curvature $n(n-1)$.
At $p_i$ its Fowler asymptotic profile has minimum exactly $t_i\eps$,
for every $1\le i\le m$.

For a suitable choice of $g_\eps$, every sequence $\eps_j\downarrow0$
has a subsequence along which $g_{\eps_j}$ converges smoothly on
compact subsets of $M\setminus P$. The limit extends to a smooth
Yamabe minimizing metric on $M$ of scalar curvature $n(n-1)$.
\end{corollary}

Schoen \cite{schoen} constructed complete conformal metrics of positive
constant scalar curvature on the round sphere with any prescribed finite
set of at least two punctures. Mazzeo and Pacard \cite{mp}
gave a gluing construction for these prescribed sets.
The obstruction for a single puncture follows from the classification of
positive entire solutions of the critical Euclidean equation by
Caffarelli, Gidas and Spruck \cite{cgs}.

Byde \cite{byde} attached a Delaunay end at a prescribed point of a
nondegenerate compact background of positive constant scalar curvature
that is locally conformally flat near that point.
For a compact background of scalar curvature $n(n-1)$ with invertible
Yamabe linearization, Santos \cite[Theorems~1.5 and~6.1]{santos}
constructed Delaunay ends at prescribed finite point sets.
His construction in dimensions $3\le n\le5$ applies to arbitrary such
sets. In dimensions $n\ge6$, his hypotheses include
\[
 \nabla^jW_g(p_i)=0,\qquad
 0\le j\le\left\lfloor\frac{n-6}{2}\right\rfloor,
 \qquad 1\le i\le m,
\]
where $W_g$ is the Weyl tensor.

Part I \cite{part1} establishes the classification of isolated singularities
and the sharp critical upper bound for smooth background metrics in
dimensions $7\le n\le29$, and reviews the previously known classification
for $3\le n\le6$. These results cover all dimensions $3\le n\le29$.
In every dimension $n\ge30$, Part~II~\cite{part2} constructs a complete singular solution whose conformal factor grows faster than the Fowler rate along a sequence approaching the singular  point.

Our construction allows arbitrary curvature jets at the prescribed points
and Yamabe minimizers that need not be isolated or nondegenerate. 

Subsequent  work related to Theorem~\ref{thm:main} is discussed in the \hyperref[sec:addednote]{Added note} toward the end of the paper. The \hyperref[sec:closed-set-announcement]{Appendix at the end of the paper} announces the Bessel-capacity criterion for nonempty closed singular sets with empty interior. Its proof will appear in \cite{part4}.

\subsection*{Overview of the proofs}

Proposition~\ref{prop:necessity-round} proves the necessary positivity
of the Yamabe constant and settles the round case. In the nonround
case with positive Yamabe constant, choose a Yamabe minimizing
representative $g$, normalized to scalar curvature $n(n-1)$, as a
fixed reference metric. We write the metrics to be constructed as
conformal multiples of this fixed $g$.

Let $\K$ denote the set of positive smooth functions $v$ on $M$
for which $v^{4/(n-2)}g$ is a Yamabe minimizing metric of scalar
curvature $n(n-1)$. In particular, $1\in\K$. Choosing $g$ as the
reference metric does not prescribe the limiting conformal factor
to be $1$. As the neck sizes tend to zero, the construction gives
subsequential convergence, smoothly on compact subsets of
$M\setminus P$, to a function $v_*\in\K$. The corresponding limiting
metric is $v_*^{4/(n-2)}g$, which need not equal $g$.

Fix positive constants $t_1,\ldots,t_m$. For each sufficiently small
$\eps>0$, prescribe the Fowler neck parameter
$\eps_i=t_i\eps$ at $p_i$; thus the ratios of the prescribed neck
sizes are independent of $\eps$. We seek a positive conformal
factor $u$ on $M\setminus P$ satisfying
\[
 \Ff(u)=-\Delta_gu+\frac{n(n-2)}4u
              -\frac{n(n-2)}4u^{(n+2)/(n-2)}=0.
\]
We first construct a finite-dimensional family of smooth conformal
factors near $\K$, each satisfying the  Yamabe equation after projection to
the orthogonal complement of 
a  finite-dimensional space. We then attach the singular ends to members of
this family and choose the remaining parameters so that all finite-dimensional
residual terms in the equation vanish.

\medskip\noindent
\textit{A uniform reduction near the minimizing set.}
The set $\K$ is compact. For a smooth positive function $u$ on $M$, the
linearization of the Yamabe equation is
\[
 L_u=D\Ff(u)=-\Delta_g+\frac{n(n-2)}4
                   -\frac{n(n+2)}4u^{4/(n-2)}.
\]
The uniform $L^\infty$ bounds on $\K$ allow us to choose a fixed
finite-dimensional space $E$, spanned by sufficiently many low
eigenfunctions of $-\Delta_g$ and containing the constants, such
that the quadratic form of $L_u$ is uniformly coercive on
$E^\perp$ for positive $u$ in a fixed neighborhood of $\K$.
Here orthogonality is taken in $L^2(M,g)$, and $\Pi_E$ denotes the
corresponding orthogonal projection onto $E$. We invert the
projected operator
\[
 (I-\Pi_E)L_u\big|_{E^\perp},
\]
rather than requiring the full linearization to be invertible.
This allows us to solve the $E^\perp$ component of the equation
while retaining the component in $E$ as a finite-dimensional
problem.

Fix an orthonormal basis $\Psi_1,\ldots,\Psi_{N_E}$ of $E$. A numerical
coefficient vector $\bgpar=(\eta_1,\ldots,\eta_{N_E})$ parametrizes the family
\[
 U_\bgpar=\sum_{j=1}^{N_E}\eta_j\Psi_j+h(\bgpar),\qquad
 h(\bgpar)\in E^\perp,\qquad (I-\Pi_E)\Ff(U_\bgpar)=0.
\]
Lemma~\ref{lem:compact-graph} constructs this smooth family and shows
that it contains every member of $\K$. Its remaining equation error has
the form $\Ff(U_\bgpar)=\sum_jd_{0,j}(\bgpar)\Psi_j$.
Let $J_g$ be the energy with first variation
$DJ_g(u)[\phi]=\int_M\Ff(u)\phi\vol$. For the reduced energy
$j_0(\bgpar)=J_g(U_\bgpar;M)$, differentiation gives
$\partial_{\eta_j}j_0=d_{0,j}$. Thus its critical points solve the full equation.

The minimizing property of $\K$ is used after separating the amplitude
$\int_M U_\bgpar\vol$ from  the normalized direction vector
$\bgpar/\int_M U_\bgpar\vol$.
The $E^\perp$ correction has zero integral because constants belong to $E$.
By Lemma~\ref{lem:negative-scale}, the reduced energy $j_0(\bgpar)$ is strictly concave in the amplitude direction near
$\K$, so the amplitude critical-point equation determines the amplitude
as a smooth function of the normalized direction. After this substitution,
the remaining reduced energy has its minimum precisely at the directions
corresponding to $\K$. We choose a compact neighborhood of that minimum
set whose boundary has strictly larger energy. This boundary gap, rather
than nondegeneracy of an individual minimizer, will select the background
after the singular ends are inserted.

\medskip\noindent\textit{Joining local ends to the exterior solution.}
For $3\le n\le7$, choose disjoint small coordinate balls
$B_{\rho_i}(p_i)$ and write
$\Om=M\setminus\bigcup_iB_{\rho_i}(p_i)$ for the exterior region.
On each punctured ball, we begin with a family of Euclidean Fowler
models with prescribed neck $\eps_i$, an adjustable phase, and $n$
parameters controlling the first spherical-harmonic boundary data.
We also allow the boundary values in the higher spherical harmonics
to vary.

To obtain local solutions for the curved metric, we solve for a
correction to these models. Lemma~\ref{lem:model-inverse} provides
an inverse for the Euclidean Yamabe linearization at each model,
and Lemma~\ref{lem:geometry} estimates the model's error in the
curved metric. Applying this inverse to that error gives a linear
correction. On $S_{\rho_i}=\partial B_{\rho_i}(p_i)$, its boundary
value and its normal derivative multiplied by $\rho_i$ are
$O(\rho_i^{(10-n)/2})=o(\rho_i)$ for $3\le n\le7$.
Together with the estimates for the adjustable boundary data and
the nonlinear terms, this estimate allows us to solve the correction
equation by contraction. Proposition~\ref{prop:local} gives the
resulting family of exact local singular solutions, whose parameters
and boundary data remain available for matching.

The exterior construction in Subsection~\ref{sec:exterior-inverse}
starts from $U_\bgpar$ and adds the terms needed to match the
local singular solutions.
Near each $p_i$, we introduce a term with an adjustable coefficient
and leading radial behavior $r^{2-n}$, where $r$ is the local
radial coordinate. These terms supply the decaying radial component
of the Fowler models near the matching spheres.
We also add Euclidean harmonic extensions of adjustable nonconstant
boundary values on each $S_{\rho_i}$, cut off away from those spheres.

For each choice of these coefficients and boundary data, we use
the linearization at $U_\bgpar$ and a contraction argument to solve
for a further small correction, together with the coefficients
$d_{\eps,j}$. The resulting exterior function $u^{\rm out}$
satisfies
\[
 \Ff(u^{\rm out})=\sum_jd_{\eps,j}\Psi_j,
 \qquad
 \int_{\Om}u^{\rm out}\Psi_j\vol=\eta_j.
\]
The integrals in the second condition are the \emph{exterior moments}.
They are prescribed on $\Om$, whereas the normalization of
$U_\bgpar$ uses integrals over the whole manifold $M$.

In Subsection~\ref{subsec:low-boundary-matching}, we
 choose the local parameters and the boundary data so that
the local and exterior functions, and their normal derivatives,
agree on each $S_{\rho_i}$. At $p_i$ we express both factors relative to the same local metric
    $h_i=\Theta_i^{4/(n-2)}g$ and take both normal derivatives in the
    direction of increasing radius. Separating the constant, degree-one, and higher
spherical-harmonic components lets us solve these matching equations
for the finite parameters and the remaining boundary functions.
For fixed neck sizes, matching determines these data as functions
of $\bgpar$. The joined family solves the Yamabe equation on the
punctured balls, while its exterior error has coefficients
$d_{\eps,j}(\bgpar)$.

\medskip\noindent\textit{The local construction in higher dimensions.}
For $n\ge8$, we construct a local family on each punctured ball by
joining three pieces: an inner end asymptotic to a Fowler solution
with prescribed neck $\eps_i$, a middle annulus modeled on a smooth
Euclidean bubble, and an outer annulus connecting the bubble to the
boundary data at $S_{\rho_i}$. The two spheres separating these pieces
lie inside $B_{\rho_i}(p_i)$; the local family will subsequently be
joined to the global exterior function at $S_{\rho_i}$.

The smooth bubble in the middle annulus depends on the scale $b_i$
and the $n$-component center parameter $\zeta_i$.
It solves the Euclidean Yamabe equation, but generally has a nonzero
error when substituted into the equation for the curved metric.
If the leading nonzero metric Taylor term has degree $q_i$ and
$n\ge4q_i+2$, we add the first metric correction constructed in
\cite{part1} to cancel the error linear in that term,  as specified
in Lemma~\ref{lem:first-metric-correction}.
We then solve the middle equation allowing a remaining error in
the $n+1$ directions associated with variations of the bubble's
scale and center.  

Proposition~\ref{prop:inner-graph} constructs the exact inner
solution from a radial Fowler reference profile with prescribed
neck $\eps_i$. Proposition~\ref{prop:inner-dn} allows its full
boundary value to be prescribed. The inner phase and angular
adjustments are then determined by the boundary values supplied
through the internal matching; they are not retained as additional
free parameters of the assembled local family.
Throughout the construction, the singularity remains at $p_i$.

On the outer annulus inside $B_{\rho_i}(p_i)$, the construction
in Subsection~\ref{subsec:outer-annular-problem} solves the
Yamabe equation exactly with prescribed values on its two
boundary spheres.
These boundary values remain adjustable so that the annular solution
can be joined to the middle-region solution at its inner boundary
and to the global exterior function at $S_{\rho_i}$.

Lemma~\ref{lem:bridge-transmission} determines the boundary data
on the two internal spheres so that the functions and their
normal derivatives agree there.
Proposition~\ref{prop:bridge} gives the resulting positive local
family, with the prescribed Fowler neck and outer boundary values.
Its only equation error is supported in the middle annulus and
lies in the $n+1$ bubble directions just described.
We then join this local family to the global exterior function
on $S_{\rho_i}$, imposing the same exterior moment conditions
as in the construction for $3\le n\le7$.

\medskip\noindent\textit{A finite action and the remaining equations.}
The ordinary energy diverges on a Fowler end. We truncate the end near
the puncture and subtract an explicit expression determined by its
asymptotic Fowler profile. The subtraction uses the actual asymptotic
phase, which varies with the constructed solution. Its derivative cancels
the boundary term at the shrinking cutoff sphere. The remaining first
variation consists of the equation residual and the boundary term at the
outer sphere; Lemma~\ref{lem:end-first-variation} gives the formula,
adapted from \cite{part2}. The cutoff limit is taken for each fixed neck size;
only afterward do we study the family as $\eps\to0$.

For $3\le n\le7$, let $j_\eps(\bgpar)$ denote the sum of the
exterior energy and the renormalized energies on the punctured balls.
Proposition~\ref{prop:low-global-action} establishes the
first-variation identity and convergence stated below.
In its first variation, the boundary terms cancel because the
functions and their normal derivatives agree on each common sphere.
The equation holds throughout each punctured ball, so only the
exterior error contributes. Differentiating the exterior moment
conditions therefore gives
\[
 \partial_{\eta_l}j_\eps
 =\sum_j d_{\eps,j}
      \int_{\Om}\Psi_j\,\partial_{\eta_l}u^{\rm out}\vol
 =d_{\eps,l}.
\]
Moreover,
\[
 j_\eps-C_\eps\longrightarrow j_0
 \quad\text{in }C^2,
\]
where $C_\eps$ is independent of $\bgpar$.
For each normalized spectral direction, the amplitude critical-point
equation has a solution close to that for $j_0$. Substituting this
amplitude and minimizing in the remaining variables gives an interior
critical point by the boundary energy gap. The displayed
first-variation identity then eliminates the exterior error.

For $n\ge8$, the joined family also depends on the bubble scales
$b_i$ and center parameters $\zeta_i$. Denote its renormalized
energy by $I_\eps$. 
Proposition~\ref{prop:global-high-action} gives its first
variation and its expansion in the remaining parameters.
Boundary terms again cancel after matching,
but equation errors remain both on the global exterior region
and in the middle annuli. The corrections in the middle annuli are normalized to have zero
pairing with the functions spanning the local equation errors.
At fixed bubble scales and centers, the compact background affects
the local families through their matched outer boundary values.
The functions used in the normalization remain fixed under these
variations. Differentiating the zero-pairing conditions therefore
shows that these variations have zero pairing with the local errors.
 The exterior moment conditions
therefore identify the derivatives of $I_\eps$ in $\bgpar$ with
the exterior error coefficients. The derivatives in the bubble
scales and centers give an invertible linear system for the
$n+1$ local error coefficients at each puncture. Thus a critical
point of $I_\eps$ eliminates both the exterior and the local errors.

Up to an additive constant independent of the remaining parameters,
the expansion of $I_\eps$ is $j_0(\bgpar)$ plus contributions from
the punctures. Each end contributes a positive multiple of
$\eps_i^2\log b_i$. For a nonzero metric Taylor term of degree $q_i$ retained in
the leading expansion, Lemma~\ref{lem:tt-coefficient} identifies
the quadratic curvature contribution from \cite{part1}: it has order
$b_i^{2q_i}$ and a strictly negative coefficient.
The interaction with the positive background value $U_\bgpar(p_i)$
also has a negative coefficient, but has order $b_i^{(n-2)/2}$.
Comparing these powers explains the division at $n=4q_i+2$:
curvature dominates above this threshold, the background interaction
dominates below it, and both contribute at equality. When no
curvature term of sufficiently low degree is present, the background
interaction supplies the leading negative contribution.
The derivatives of the logarithmic term and the leading negative
contribution with respect to $\log b_i$ have opposite signs;
balancing them selects the bubble scale.

The leading action is strictly concave in the amplitude and
logarithmic scale variables. In the normalized background and
bubble-center variables, the compact minimizing set and the decay
of the center-dependent coefficients give a positive boundary
energy gap. We control the remainder in this expansion and its
first parameter derivatives at order $o(\eps^2)$.
 Lemma~\ref{lem:c1-minimax} uses this $C^1$ control, without
requiring second derivatives of the remainder, to obtain the
required interior critical point.

Finally, the constructed ends satisfy
$u\ge c\eps_i r^{-(n-2)/2}$. A curve approaching $p_i$ therefore has
conformal length bounded below by a positive multiple of
$\int |\dd r|/r$, which diverges. The selected compact backgrounds
approach $\K$, and all local and exterior corrections vanish away from
the punctures. Compactness of $\K$ then gives the subsequential convergence
asserted in Corollary~\ref{cor:neck-ratios}.

\section{The compact reduction and the renormalized action}\label{sec:compact}
The variational construction requires a finite-dimensional family of
smooth backgrounds whose equation residual is the gradient of its energy,
and a finite action for singular ends with a controlled first variation.
We establish these ingredients after treating positivity and the round case.

\subsection{Equation, normalization, and the round case}

Put
\[
 p:=\frac{n+2}{n-2},\qquad
 c(n):=\frac{n-2}{4(n-1)}.
\]
For a smooth metric $g$, let $\operatorname{Scal}_g$ denote its scalar curvature.
We use the Laplacian with nonpositive spectrum. Define
\[
 \Ff_g(u):=-\Delta_gu+c(n)\operatorname{Scal}_gu-\frac{n(n-2)}4u^p.
\]
Its energy on a smooth domain $D$ is
\begin{equation}
 J_g(u;D):=\frac12\int_D(|\nabla u|^2+c(n)\operatorname{Scal}_gu^2)\vol
       -\frac{(n-2)^2}{8}\int_Du^{p+1}\vol.
 \label{eq:energy}
\end{equation}

On the whole manifold we abbreviate $J(u):=J_g(u;M)$ and fix
$0<\alphaH<1$ as the H\"older exponent.
Our equation gives scalar curvature $n(n-1)$; Parts I and II
use coefficient $n(n-2)$ and scalar curvature $4n(n-1)$.
Their conformal factors are related by $u=2^{(n-2)/2}u^{\rm I}$.
Write $g_{\rm E}=(\delta_{ij})$, and denote volume and hypersurface
measures by $\dd V_g$ and $\dd S_g$, respectively.

\begin{proposition}\label{prop:necessity-round}\label{sec:minimal}
If $[g]|_{M\setminus P}$ contains a metric of constant positive scalar
curvature, then $Y(M,[g])>0$. If $(M,[g])$ is conformally equivalent to
the round sphere, then $[g]|_{M\setminus P}$ contains a complete metric
of constant positive scalar curvature if and only if $|P|\ge2$.
\end{proposition}
\begin{proof}
After a constant rescaling, its positive conformal factor solves
$(-\Delta_g+c(n)\operatorname{Scal}_g)u=\frac{n(n-2)}4u^p$ on $M\setminus P$.
For every $\phi\in C_c^\infty(M\setminus P)$, the ground-state
identity gives
\[
 \int_M(|\nabla\phi|^2+c(n)\operatorname{Scal}_g\phi^2)\vol
 =\int_Mu^2|\nabla(\phi/u)|^2\vol
       +\frac{n(n-2)}4\int_Mu^{p-1}\phi^2\vol.
\]
Let $\psi>0$ be a first eigenfunction of the compact conformal
Laplacian, with eigenvalue $\lambda_1$. Since finite point sets have
zero $H^1$ capacity for $n\ge3$, choose smooth cutoffs $0\le\chi_j^{\rm cap}\le1$
vanishing near $P$, with $\chi_j^{\rm cap}\psi\to\psi$ in $H^1$ and
$\chi_j^{\rm cap}=1$ eventually on each fixed compact subset of $M\setminus P$.
Testing with $\chi_j^{\rm cap}\psi$ and passing to the limit gives, on any
such subset $K$ of positive volume,
\[
 \lambda_1\int_M\psi^2\vol\ge \frac{n(n-2)}4\int_Ku^{p-1}\psi^2\vol>0.
\]
Thus $\lambda_1>0$, equivalently $Y(M,[g])>0$.

Suppose the conformal class is round and $|P|=1$. Stereographic
projection gives a positive entire solution of the critical Euclidean
equation. By \cite[Corollary 8.2 and Remark 8.3]{cgs}, this solution
is a standard bubble and decays as $|x|^{2-n}$. A radial curve to
infinity has finite conformal length, so the metric is incomplete.
For $|P|\ge2$, the round construction gives a complete metric of
positive constant scalar curvature at the prescribed points
\cite[Theorem 1]{mp}.
\end{proof}

\subsection{A uniform family near the minimizing set}
\label{subsec:compact-family}

We seek a family of smooth positive conformal factors containing all
normalized minimizers, with the remaining equation components equal to
the derivatives of a finite-dimensional energy.

We henceforth work in a nonround positive conformal class, where
$0<Y(M,[g])<Y(\Sph^n)$ \cite{lp}.
By \cite[Theorem 4.5]{lp}, choose a Yamabe minimizer and rescale its
metric to scalar curvature $n(n-1)$. Denote this fixed reference metric
by $g$ and abbreviate $\Ff_g$ to $\Ff$.
The equation becomes
\begin{equation}
 \Ff(u)=-\Delta_gu+\frac{n(n-2)}4u-\frac{n(n-2)}4u^p=0.\label{eq:equation}
\end{equation}
Let $\K$ be the set of all smooth positive conformal factors $u$ for
which $u^{4/(n-2)}g$ is a Yamabe minimizing metric of scalar curvature
$n(n-1)$. Its members satisfy \eqref{eq:equation}. Testing that equation
and using the Yamabe quotient gives
\[
 M_{\K}:=\int_Mu^{p+1}\vol
 =\left(\frac{Y(M,[g])}{n(n-1)}\right)^{n/2}
 <\operatorname{Vol}(\Sph^n,g_{\rm rd}),\qquad u\in\K.
\]
We use the strict inequality $Y(M,[g])<Y(\Sph^n)$ to obtain
uniform $C^{j,\alphaH}$ bounds on $\K$ for every fixed integer $j\ge0$. For $u\in\K$, test \eqref{eq:equation} against
$u^{1+2\delta_1}$ and apply the sharp Sobolev inequality
\cite[Theorems 2.3 and 3.3]{lp} to $v=u^{1+\delta_1}$.
Since $\frac{n(n-2)}4M_{\K}^{2/n}=c(n)Y(M,[g])$, H\"older's inequality gives
\[
 \|v\|_{p+1}^2
 \le (1+\eta)\frac{(1+\delta_1)^2}{1+2\delta_1}
       \frac{Y(M,[g])}{Y(\Sph^n)}\|v\|_{p+1}^2
       +C_{\eta,\delta_1}\|v\|_2^2.
\]
Choose $\eta,\delta_1>0$ so that the coefficient of
$\|v\|_{p+1}^2$ on the right is less than one. Also require
$2(1+\delta_1)\le p+1$. Absorption gives a uniform
$L^{(p+1)(1+\delta_1)}$ bound. Since
\[
 (p+1)(1+\delta_1)>\frac n2(p-1)=p+1,
\]
\cite[Theorem 4.1]{lp} gives a uniform $C^{2,\alphaH}$ bound.
Limits retain the common mass and equation, and are positive by the
strong maximum principle. Compactness and elliptic bootstrapping give
\[
 \inf_{u\in\K}\min_Mu>0,\qquad
 \sup_{u\in\K}\|u\|_{C^{j,\alphaH}(M)}<\infty
 \quad\text{for every fixed }j,
\]
with compactness in every $C^j$ norm.
The bounds at higher orders also give compactness in
$C^{2,\alphaH}(M)$, the space used for the implicit function theorem.
Constants below may depend on the metric, punctures, positive neck
ratios and compact parameter boxes; uniformity in $\eps$ is stated
in the respective estimates.

We adapt the projected implicit function construction of
\cite[Lemma 2.2 and Appendix A]{ens} using one spectral space over
all of $\K$. For positive $w$, put
\[
 A(w)=\int_M(|\nabla w|^2+\frac{n(n-2)}4w^2)\vol,\qquad
 B(w)=\int_Mw^{p+1}\vol.
\]
The Yamabe inequality is
$A(w)\ge c(n)Y(M,[g])B(w)^{2/(p+1)}$.
The linearization at a positive function $u$ is
\[
 L_u=D\Ff(u)=-\Delta_g+\frac{n(n-2)}4-\frac{n(n+2)}4u^{p-1}.
\]
The uniform bounds on $\K$ allow us to choose a bounded open convex
neighborhood $\mathcal U\subset C^{2,\alphaH}(M)$ such that
$\inf_{u\in\mathcal U}\min_Mu>0$.
Set $C_0=\frac{n(n+2)}4\sup_{u\in\mathcal U}\|u\|_\infty^{p-1}$.
Choose $E$ to be the span of sufficiently many low eigenfunctions of
$-\Delta_g+\frac{n(n-2)}4$, including the constants, and let $\Lambda$
be the first excluded eigenvalue. Choose the cutoff so that
$\Lambda\ge2C_0$.
Write $N_E=\dim E$, and let $\Pi_E$ be the orthogonal projection
onto $E$. We use $I-\Pi_E$ for the projection onto $E^\perp$.
The inequality $\Lambda\ge2C_0$ also gives, for $z\in E^\perp\cap H^1(M)$,
\begin{equation}
 \langle L_uz,z\rangle=J''(u)[z,z]
 \ge(1-C_0/\Lambda)A(z)\ge c\norm z_{H^1}^2,
 \qquad u\in\mathcal U,\quad z\in E^\perp.
 \label{eq:high-coercive}
\end{equation}

The subspace $E^\perp$ need not be invariant under multiplication by
$u^{p-1}$. It is the projected operator $(I-\Pi_E)L_u$ on $E^\perp$
that is uniformly invertible.

Fix an $L^2(M,g)$-orthonormal basis $\Psi_1,\ldots,\Psi_{N_E}$ of $E$
with $\Psi_1=\operatorname{Vol}_g(M)^{-1/2}$. Define the coefficient map by
\[
 \bgcoord u=\left(\int_Mu\Psi_j\vol\right)_{j=1}^{N_E},\qquad
 \Pi_Eu=\sum_{j=1}^{N_E}(\bgcoord u)_j\Psi_j.
\]
The coefficient vector $\bgpar\in\R^{N_E}$ parametrizes the following family.

\begin{lemma}\label{lem:compact-graph}
There is a smooth family of positive functions on a neighborhood
$O\subset\R^{N_E}$ of the entire set $\bgcoord\K$, given by
\[
 U_\bgpar=\sum_{j=1}^{N_E}\eta_j\Psi_j+h(\bgpar),\qquad
 h(\bgpar)\in E^\perp,\qquad (I-\Pi_E)\Ff(U_\bgpar)=0.
\]
The cutoff can be chosen so that $\bgcoord$ is also injective on
$\{u/\int_Mu\vol:u\in\K\}$.
The family satisfies $\bgcoord U_\bgpar=\bgpar$. Write
$\Ff(U_\bgpar)=\sum_jd_{0,j}(\bgpar)\Psi_j$. Then
\begin{equation}
 j_0(\bgpar):=J(U_\bgpar),\qquad Dj_0(\bgpar)=d_0(\bgpar).\label{eq:compact-gradient}
\end{equation}
For every fixed integer $k\ge0$, the map $\bgpar\mapsto U_\bgpar$ is smooth
with values in $C^{k,\alphaH}(M)$. Its bounds are uniform on fixed
relatively compact neighborhoods of $\bgcoord\K$ in $O$, with constants depending on $k$.
\end{lemma}
\begin{proof}
We first prove that $\bgcoord$ is injective on
$\{u/\int_Mu\vol:u\in\K\}$. If
$w=u/\int_Mu\vol$ with $u\in\K$, then
\[
 (-\Delta_g+\tfrac{n(n-2)}4)w=\lambda(w)w^p,
 \qquad
 \lambda(w)=\frac{n(n-2)}4
      \left(\frac{M_\K}{B(w)}\right)^{2/n}.
\]
For two such functions $w_0,w_1$, set $z=w_1-w_0$.
The uniform positive bounds on the minimizers give
$|B(w_1)-B(w_0)|+|\lambda(w_1)-\lambda(w_0)|\le C\|z\|_2$.
Testing the difference equation against $z$ yields
$A(z)\le C_\K\|z\|_2^2$.
If $\Pi_Ew_1=\Pi_Ew_0$, then
$\Lambda\|z\|_2^2\le A(z)\le C_\K\|z\|_2^2$.
Thus choosing $\Lambda>C_\K$ as well gives $z=0$.
We make both cutoff choices once and keep $E$ fixed below.

For the linearization $L_u$ defined above, \eqref{eq:high-coercive} shows that
$(I-\Pi_E)L_u:E^\perp\cap H^1\to E^\perp\cap H^{-1}$ is an isomorphism.
Here $f\in E^\perp\cap H^{-1}$ means
$\langle f,\Psi_j\rangle_{H^{-1},H^1}=0$ for every $j$.
For $(I-\Pi_E)L_uz=f$, the identity $L_uz=f+\Pi_E(L_uz)$ and the fixed smooth
basis of $E$ give
\[
 \|z\|_{H^1}\le C\|f\|_{H^{-1}},\qquad
 \|\Pi_E(L_uz)\|_{C^{0,\alphaH}}\le C\|z\|_{H^1}.
\]
Schauder estimates \cite[Theorems 2.4--2.5]{lp} therefore give the
uniform inverse bound $\|z\|_{C^{2,\alphaH}}\le C\|f\|_{C^{0,\alphaH}}$.
The implicit function theorem applied to
$(I-\Pi_E)\Ff(\sum_j\eta_j\Psi_j+h)$ now yields a
smooth graph near each $u\in\K$. Since $\Ff(U_\bgpar)\in E$,
elliptic regularity gives smooth dependence in every fixed spatial norm.

The local graphs agree whenever their $\Pi_E$ projections coincide. Indeed,
for two graph values $u_0,u_1\in\mathcal U$ with
$z=u_1-u_0\in E^\perp$, their residuals lie in $E$, so
\[
 0=DJ(u_1)[z]-DJ(u_0)[z]
   =\int_0^1J''(u_0+tz)[z,z]\dd t
   \ge c\|z\|_{H^1}^2.
\]
Convexity of $\mathcal U$ justifies the segment. A finite cover of
$\K$ gives the common graph, and $\Pi_E\partial_{\eta_l}U_\bgpar=\Psi_l$ gives
\[
 \partial_{\eta_l}j_0=\int_M\Ff(U_\bgpar)\partial_{\eta_l}U_\bgpar\vol=d_{0,l}.
\]
\end{proof}

Fix a relatively compact neighborhood $O_0$ of $\bgcoord\K$ with
$\overline{O_0}\subset O$. We take all uniform background estimates on
$\overline{O_0}$.

\medskip\noindent\textit{Solving a linear equation with prescribed moments.}
For $L_\bgpar=D\Ff(U_\bgpar)$ and prescribed
$(f,\mathbf m)\in C^{0,\alphaH}(M)\times\R^{N_E}$, the unknowns are a
function $z$ and coefficients $d_j$ satisfying
\[
 L_\bgpar z-\sum_{j=1}^{N_E}d_j\Psi_j=f,
 \qquad \int_Mz\Psi_j\vol=m_j.
\]
The moments prescribe $z_E:=\Pi_Ez=\sum_jm_j\Psi_j$; the projected
inverse then determines the remaining part and the residual coefficients:
\begin{equation}
\begin{aligned}
 z&=z_E+((I-\Pi_E)L_\bgpar|_{E^\perp})^{-1}
                 (I-\Pi_E)(f-L_\bgpar z_E),\\
 d_j&=\int_M(L_\bgpar z-f)\Psi_j\vol.
\end{aligned}
 \label{eq:augmented-compact}
\end{equation}
This gives the unique augmented solution on the compact manifold.
The map $(f,\mathbf m)\mapsto(z,d)$ and its first $\bgpar$ derivatives
are bounded from $C^{0,\alphaH}\times\R^{N_E}$ to
$C^{2,\alphaH}\times\R^{N_E}$ on compact subsets of $O$.
Differentiating the product of $(I-\Pi_E)L_\bgpar|_{E^\perp}$ and its
inverse gives the first $\bgpar$-derivative bound.

\subsection{Nehari reduction and the minimum set}
\label{subsec:amplitude}

The projected equation has reduced the compact problem to finding a
critical point of $j_0$. However, Yamabe minimizers minimize the
scale-invariant Yamabe quotient, not the unconstrained action $J$.
We first remove its negative amplitude direction by the Nehari
condition $DJ(u)[u]=0$.
For a fixed positive function $w$, direct differentiation of
\[
 J(tw)=\frac{t^2}{2}A(w)
       -\frac{n(n-2)}{4(p+1)}t^{p+1}B(w)
\]
shows that its unique positive stationary point is a strict maximum.
The resulting value is
\[
 \max_{t>0}J(tw)
 =\frac1n\left(\frac{n(n-2)}4\right)^{-(n-2)/2}
       \left(\frac{A(w)}{B(w)^{2/(p+1)}}\right)^{n/2}.
\]
Thus the ray reduction and the Yamabe quotient encode the same
minimization problem. We implement the amplitude condition within the
projected family $U_\bgpar$, rather than multiply its members by an
arbitrary constant, which would generally destroy the projected equation.

Since $1\in E$, the function $h(\bgpar)\in E^\perp$ has zero integral.
We therefore use the integral of $U_\bgpar$ as the amplitude. For
$\bgpar$ near $\bgcoord\K$, write $\bgpar=\amp\bgshape$, where
\[
 \amp=\int_MU_\bgpar\vol
      =\operatorname{Vol}_g(M)^{1/2}\eta_1>0,
 \qquad
 \bgshape\in\mathcal H_E:=
 \{\mathbf s\in\R^{N_E}:\operatorname{Vol}_g(M)^{1/2}s_1=1\}.
\]
For the minimizers themselves, let
\[
 \widehat\K=\{u/\!\int_Mu\vol:u\in\K\},
 \qquad \mathcal S=\bgcoord\widehat\K.
\]
For $w\in\widehat\K$, the unique multiplier $\amp>0$ for which
$\amp w\in\K$ is
$\amp=(M_\K/B(w))^{1/(p+1)}$. The identity
$\int_Mu^{p+1}\vol=M_{\K}$ for $u\in\K$ implies that
$u\mapsto u/\int_Mu\vol$ is injective on $\K$.
This map and $\Pi_E|_{\widehat\K}$ are continuous injections of
compact sets, so
\[
 \Gamma=\left\{\left(\int_Mu\vol,
        \bgcoord\left(u\Big/\int_Mu\vol\right)\right):u\in\K\right\}
       =\{(\amp_*(\bgshape),\bgshape):\bgshape\in\mathcal S\}
\]
is a compact continuous graph with amplitudes bounded away from zero.
On $\Gamma$, $U_{\amp \bgshape}\in\K$ and $\partial_\amp j_0(\amp \bgshape)=0$.

Since $\bgcoord U_{\amp\bgshape}=\amp\bgshape$, 
differentiating in $\amp$, with $\bgshape$ fixed, gives
$\bgcoord(\partial_\amp U_{\amp\bgshape})=\bgshape$.
Consequently,
\[
 \partial_\amp U_{\amp\bgshape}
       -\amp^{-1}U_{\amp\bgshape}\in E^\perp.
\]
Since $\Ff(U_{\amp\bgshape})\in E$, 
\[
DJ(U_{\amp\bgshape})[\partial_\amp U_{\amp\bgshape}- \amp^{-1}U_{\amp\bgshape}]
=\int_M\mathcal F(U_{\amp\bgshape})\,[\partial_\amp U_{\amp\bgshape}-\amp^{-1}U_{\amp\bgshape}]\,dV_g
=0,
\]
so
\begin{equation}
 \partial_\amp j_0(\amp\bgshape)
 =DJ(U_{\amp\bgshape})[\partial_\amp U_{\amp\bgshape}]
 =\amp^{-1}DJ(U_{\amp\bgshape})[U_{\amp\bgshape}].
 \label{eq:projected-nehari-identity}
\end{equation}
We call
\begin{equation}
 \frac{\partial}{\partial\amp}j_0(\amp\bgshape)=0
 \label{eq:amplitude-equation}
\end{equation}
the amplitude equation: $\amp$ is the unknown and $\bgshape$ is held
fixed. By \eqref{eq:projected-nehari-identity}, it is exactly the
Nehari condition on this family, namely
$A(U_{\amp\bgshape})=\frac{n(n-2)}4 B(U_{\amp\bgshape})$.

\begin{lemma}\label{lem:negative-scale}
Near $\Gamma$, $\partial_\amp^2j_0(\amp\bgshape)$ is uniformly negative.
There is a smooth function $\amp_0$ on a neighborhood of $\mathcal S$
in $\mathcal H_E$, with $\amp_0|_{\mathcal S}=\amp_*$, such that
\[
 \left.\frac{\partial}{\partial\amp}j_0(\amp\bgshape)
       \right|_{\amp=\amp_0(\bgshape)}=0.
\]
For every fixed compact subset $D_0$ of this neighborhood, a sufficiently
small $C^2$ perturbation $j$ of $j_0$, on a fixed neighborhood of
$\{\amp_0(\bgshape)\bgshape:\bgshape\in D_0\}$, has a corresponding
$C^1$ function $\amp_j$, defined near $D_0$, satisfying
\[
 \left.\frac{\partial}{\partial\amp}j(\amp\bgshape)
       \right|_{\amp=\amp_j(\bgshape)}=0.
\]
For $\bgshape\in D_0$ this solution is unique in a fixed small interval
about $\amp_0(\bgshape)$, and $\amp_j\to\amp_0$ in $C^1(D_0)$ as
$j\to j_0$ in $C^2$.
\end{lemma}
\begin{proof}
At $u\in\K$, put $\bgpar=\bgcoord u$, $W=D_\bgpar U_\bgpar[\bgpar]$, and $z=u-W\in E^\perp$.
The differentiated graph equation gives $(I-\Pi_E)L_uW=0$, hence
$J''(u)[W,z]=0$. Coercivity gives
\[
\begin{aligned}
 J''(u)[u,u]&=J''(u)[W,W]+J''(u)[z,z]\ge J''(u)[W,W],\\
 J''(u)[u,u]&=\frac{n(n-2)}4(1-p)\int_Mu^{p+1}\vol=-nM_{\K}.
\end{aligned}
\]
Since $DJ(u)=0$, we obtain $D^2j_0(\bgpar)[\bgpar,\bgpar]=J''(u)[W,W]\le-nM_{\K}$.
For $\bgpar=\amp \bgshape$ on $\Gamma$,
\[
 \partial_\amp^2j_0(\amp \bgshape)=\amp^{-2}D^2j_0(\bgpar)[\bgpar,\bgpar]\le-2\eta_0<0,
\]
because the amplitudes are bounded. By continuity,
$\partial_\amp^2j_0\le-\eta_0$ on a neighborhood.

The scalar implicit function theorem gives local roots agreeing
with $\amp_*$ on $\mathcal S$. Compactness and continuity of
$\amp_*$ place overlapping roots in a common interval of strict
monotonicity, so they agree. On every compact subset of the domain
of $\amp_0$, choose
$r_0>0$ with $\partial_\amp^2j_0\le-\eta_0$ for
$|\amp-\amp_0(\bgshape)|\le r_0$. Sufficiently small $C^2$ perturbations preserve the strict negativity
and the opposite endpoint signs, giving a unique nearby root
$\amp_j(\bgshape)$ for each direction. The implicit function theorem
makes this root $C^1$ in $\bgshape$; its derivative formula and the
uniform lower bound for $-\partial_\amp^2j$ give
$\amp_j\to\amp_0$ in $C^1$ on the chosen compact subset.
\end{proof}

We now identify the minimum set after solving the amplitude equation.
By \eqref{eq:projected-nehari-identity}, $U=U_{\amp_0(\bgshape)\bgshape}$ satisfies $A(U)=\frac{n(n-2)}4B(U)$ and $J(U)=A(U)/n$.
The Yamabe inequality gives
\begin{equation}
 \varphi_0(\bgshape):=j_0(\amp_0(\bgshape)\bgshape)=\frac{A(U)}n
 \ge\frac1n \left(\frac{n(n-2)}4\right)^{-(n-2)/2}\bigl(c(n)Y(M,[g])\bigr)^{n/2}.
 \label{eq:reduced-minimum}
\end{equation}
Equality holds precisely when $U$ is a Yamabe minimizer;
$A(U)=\frac{n(n-2)}4B(U)$ then normalizes its equation to \eqref{eq:equation},
so $U\in\K$. In this equality case, constants belonging to $E$ give
\[
 \int_MU\vol=\amp_0(\bgshape),\qquad
 \bgshape=\bgcoord\left(U/\!\int_MU\vol\right)\in\mathcal S.
\]
Conversely, if $\bgshape=\bgcoord(u/\int_Mu\vol)$ for $u\in\K$,
uniqueness in Lemmas~\ref{lem:negative-scale} and~\ref{lem:compact-graph}
gives
$\amp_0(\bgshape)=\int_Mu\vol$ and $U_{\amp_0(\bgshape)\bgshape}=u$.
Hence the minimum set of $\varphi_0$ is exactly $\mathcal S$.

\subsection{A finite action for a given singular family}
\label{subsec:end-action-criterion}

For a Euclidean Fowler solution, the functional $J_{g_{\rm E}}$
in \eqref{eq:energy}, evaluated on an annulus about the singular
point, tends to $+\infty$ as the inner radius tends to zero.
We compute the terms depending on that radius to determine
the subtraction used to define a finite action.
Throughout this calculation the neck parameter $\eps$ is fixed.

For sufficiently small $\eps>0$, choose the phase of the Fowler profile
in \eqref{eq:fowler-ode} so that
\[
 v_\eps(0)=\eps,\qquad v_\eps'(0)=0.
\]
The radial Euclidean conformal factor is
$r^{-(n-2)/2}v_\eps(-\log r+\tau)$. Its phase is defined modulo the
period; a real lift of the phase will be chosen on each parameter
neighborhood used for the action.
The profile $2^{-(n-2)/2}v_\eps$ uses the normalization of Parts I and II,
with minimum $2^{-(n-2)/2}\eps$.
The Hamiltonian is 
\begin{equation}
 h_\eps=\frac12(v_\eps')^2-\frac{(n-2)^2}{8}v_\eps^2
               +\frac{(n-2)^2}{8}v_\eps^{p+1}<0.
 \label{eq:primitive}
\end{equation}
In particular,
$h_\eps=-\frac{(n-2)^2}{8}\eps^2(1-\eps^{4/(n-2)})$.
On the product cylinder the radial action density of this profile is
\[
 |\Sph^{n-1}|\left[
  \tfrac12(v_\eps')^2+\tfrac{(n-2)^2}{8}v_\eps^2
                      -\tfrac{(n-2)^2}{8}v_\eps^{p+1}\right]
 =|\Sph^{n-1}|\bigl[(v_\eps')^2-h_\eps\bigr].
\]

Fix a smooth local metric $h$ and geodesic normal coordinates $X$ at
the puncture. For $r=|X|$, write
\[
 \gamma=r^{-2}h,\qquad s=-\log r,\qquad
 u=r^{-(n-2)/2}v,\qquad s_0=-\log\rho.
\]
The metric $\gamma$ has the form $\dd s^2+g_s$; hence
$\dd V_\gamma=\dd s\,\dd S_\gamma(s)$, with $\dd S_\gamma(s)$ the
induced measure on the angular slice. Conformal covariance on a finite
annulus gives
\begin{equation}
 J_h(u;e^{-T}<r<\rho)=J_\gamma(v;s_0<s<T)
 +\frac{n-2}{4}\int_{\{s=T\}}v^2\dd S_\gamma
 -\frac{n-2}{4}\int_{\{s=s_0\}}v^2\dd S_\gamma.
 \label{eq:finite-cylinder-energy}
\end{equation}
 For a given family of positive functions,
we use the fixed-outer-boundary version of the subtraction in
\cite[Lemma 3.3 and equation (3.10)]{part2}, in the present normalization.

\begin{lemma}[First variation of the renormalized end action]
\label{lem:end-first-variation}
Let $\mathcal V$ be an open subset of a Banach space $\mathscr B$.
Suppose a family of positive functions
$u(\,\cdot\,;\vartheta)$, $\vartheta\in\mathcal V$, is given on
$\{X:0<|X|\le\rho\}$. Its members are smooth in space away from the
puncture and depend $C^1$ on $\vartheta$ in $C^{2,\alphaH}$ on compact
annuli. 
On a fixed closed annulus immediately inside $S_\rho$, including
$S_\rho$ itself, assume that the family depends $C^1$ on
$\vartheta$ with values in $C^{2,\alphaH}$.
Parameter derivatives are Fr\'echet derivatives in $\mathscr B$ and keep
$h$, $\eps$, $\rho$, and the physical coordinates fixed.
Let $\tau\in C^1(\mathcal V;\R)$ be a chosen real lift of the family's
asymptotic Fowler phase, and set
\[
 v(s,\theta;\vartheta)
 =e^{-(n-2)s/2}u(e^{-s}\theta;\vartheta).
\]
For this family, define its angularly integrated cylindrical action
 density by
\[
 \mathcal L(s,\vartheta)=
 \int_{\{s\}\times\Sph^{n-1}}
 \left[\frac12|\nabla_\gamma v|_\gamma^2
       +\frac{c(n)}2\operatorname{Scal}_\gamma v^2
       -\frac{(n-2)^2}{8}v^{p+1}\right]\dd S_\gamma.
\]
Thus $J_\gamma(v;s_0<s<T)=\int_{s_0}^T\mathcal L(s,\vartheta)\dd s$.

All bounds and limits below are locally uniform in $\vartheta$,
after restriction to a parameter neighborhood. Suppose there is a
nonnegative $m\in L^1([s_0,\infty))$, independent of $\vartheta$ on
that neighborhood, such that
\begin{equation}
 \left\|D_\vartheta^j
  \left\{\mathcal L(s,\vartheta)
   -|\Sph^{n-1}|\bigl[(v_\eps')^2-h_\eps\bigr]
                           (s+\tau(\vartheta))\right\}\right\|
 \le m(s),\qquad j=0,1.
 \label{eq:terminal-dominating-density}
\end{equation}
For $j=0$ the norm is absolute value, and for $j=1$ it is the
operator norm on $\mathscr B$. The dominating function $m$ may depend on the fixed neck and the parameter
neighborhood; this lemma requires no uniform bound as $\eps\to0$.
Suppose also that
\begin{equation}
 \int_{\{s=T\}}v^2\dd S_\gamma
       -|\Sph^{n-1}|v_\eps(T+\tau)^2\longrightarrow0
       \quad\text{in }C^1\text{ as }T\to\infty,
 \label{eq:terminal-quadratic-error}
\end{equation}
and, writing $\dot u=D_\vartheta u[\dot\vartheta]$ and
$\dot\tau=D_\vartheta\tau[\dot\vartheta]$, that
\begin{equation}
 \int_{S_{e^{-T}}}\partial_{\nu_h}u\,\dot u\,\dd S_h
 -|\Sph^{n-1}|\bigl[(v_\eps')^2+\tfrac{n-2}{2}v_\eps v_\eps'\bigr]
                    (T+\tau)\,\dot\tau\longrightarrow0
 \label{eq:terminal-conormal-error}
\end{equation}
uniformly for $\|\dot\vartheta\|_{\mathscr B}\le1$.
The normal at this inner cutoff points towards the puncture.

Then the following limit exists locally in $C^1$ with respect to
$\vartheta$ and defines the renormalized end action:
\begin{equation}
 \begin{aligned}
 S^h_{\eps,\rho}(u):=\lim_{T\to\infty}\Bigl\{
 &J_h(u;e^{-T}<r<\rho)\\
 &-|\Sph^{n-1}|\left[\int_0^{(T+\tau)}(v_\varepsilon'(s))^2\,ds
     +\frac{n-2}{4}v_\eps(T+\tau)^2-h_\eps T\right]\Bigr\}.
 \end{aligned}
 \label{eq:end-action}
\end{equation}
Here $u=u(\,\cdot\,;\vartheta)$ and $\tau=\tau(\vartheta)$; the action
uses the chosen real phase lift. Its differential is
\begin{equation}
 DS^h_{\eps,\rho}(u)[\dot u]
 =\lim_{T\to\infty}\int_{\{e^{-T}<r<\rho\}}
            \Ff_h(u)\dot u\,\dd V_h
       +\int_{S_\rho}\partial_{\nu_h}u\,\dot u\,\dd S_h.
 \label{eq:renormalized-first-variation}
\end{equation}
The normal on $S_\rho$ points out of the punctured ball.
\end{lemma}
\begin{proof}
We first establish existence of the limit. 
The expression
inside the limit in \eqref{eq:end-action} equals
\[
 \begin{aligned}
 &\int_{s_0}^T\left\{\mathcal L(s,\vartheta)
      -|\Sph^{n-1}|[(v_\eps')^2-h_\eps](s+\tau)\right\}\dd s\\
 &\quad+\frac{n-2}{4}\left[
        \int_{\{s=T\}}v^2\dd S_\gamma
        -|\Sph^{n-1}|v_\eps(T+\tau)^2\right]\\
 &\quad-\frac{n-2}{4}\int_{\{s=s_0\}}v^2\dd S_\gamma
-|\mathbb S^{n-1}|
\left[
\int_0^{s_0+\tau}(v_\varepsilon'(s))^2\,ds
-h_\varepsilon s_0
\right].
 \end{aligned}
\]
The first line converges with one parameter derivative by
\eqref{eq:terminal-dominating-density}; the second tends to zero in
$C^1$ by \eqref{eq:terminal-quadratic-error}. 
The final line is $C^1$ in $\vartheta$: the trace
$v(s_0,\cdot;\vartheta)$ depends $C^1$ on $\vartheta$ by the
assumed parameter regularity up to $S_\rho$, and $\tau$ is
$C^1$ by hypothesis. Thus the limit exists locally in $C^1$
with respect to $\vartheta$ and defines the renormalized end action.

We next differentiate the finite-domain action, keeping $T$ fixed,
before passing to infinity. The derivative of the subtraction is
\begin{equation}
 \begin{aligned}
 &D_\vartheta\!\left(|\Sph^{n-1}|\left[\int_0^{(T+\tau)}(v_\varepsilon'(s))^2\,ds
       +\tfrac{n-2}{4}v_\eps(T+\tau)^2-h_\eps T\right]\right)
       [\dot\vartheta]\\
 &\qquad=|\Sph^{n-1}|\bigl[(v_\eps')^2
       +\tfrac{n-2}{2}v_\eps v_\eps'\bigr](T+\tau)\,\dot\tau.
 \end{aligned}
 \label{eq:terminal-parameter-cancellation}
\end{equation}
\[
\begin{aligned}
D_\vartheta J_h(u;e^{-T}<r<\rho)[\dot\vartheta]
={}&
\int_{\{e^{-T}<r<\rho\}}
\mathcal F_h(u)\,\dot u\,dV_h\\
&+\int_{S_\rho}\partial_{\nu_h}u\,\dot u\,dS_h
+\int_{S_{e^{-T}}}\partial_{\nu_h}u\,\dot u\,dS_h.
\end{aligned}
\]
By \eqref{eq:terminal-conormal-error}, the last integral minus
the derivative of the subtracted expression in
\eqref{eq:terminal-parameter-cancellation} tends to zero.
The remaining two terms give
\eqref{eq:renormalized-first-variation}.

\end{proof}
The first-variation identity is the version of
\cite[equation (3.11)]{part2} with the outer sphere $S_\rho$ fixed. 

For a family satisfying Lemma~\ref{lem:end-first-variation}, we can
now express its action in the fixed global metric. For the change $h=\Theta^{4/(n-2)}g$, $u_g=\Theta u_h$,
express both energy densities with respect to $\dd V_g$.
The quadratic density of $J_g(u_g;D)$ minus that of $J_h(u_h;D)$
is $\frac12\operatorname{div}_g(\Theta u_h^2\nabla_g\Theta)$.
This follows from conformal covariance and expansion of
$\nabla_g(\Theta u_h)$.
The nonlinear densities agree. Integration over $D$ gives
\begin{equation}
 J_h(u_h;D)=J_g(u_g;D)
 -\frac12\int_{\partial D}u_g^2\partial_{\nu_g}\log\Theta\dd S_g.
 \label{eq:smooth-conformal-boundary}
\end{equation}
Applying this identity with base metric $h$, transformed metric
$\gamma=r^{-2}h$, and $u=r^{-(n-2)/2}v$ gives
\eqref{eq:finite-cylinder-energy}.
For a smooth conformal change, define
\begin{equation}
 S^g_{\eps,\rho}(u_g)=S^h_{\eps,\rho}(u_h)
 +\frac12\int_{S_\rho}u_g^2\partial_{\nu_g}\log\Theta\dd S_g,
 \qquad u_g=\Theta u_h.
 \label{eq:action-metric-conversion}
\end{equation}
Here $\nu_g$ points from the end into the exterior. For the local families used below,
\[
 \int_{S_{e^{-T}}}u_g^2\partial_{\nu_g}\log\Theta\,\dd S_g
 \longrightarrow0
\]
with one parameter derivative as $T\to\infty$. Hence this definition
agrees with the finite-domain conformal identity after renormalization. The metric
superscript is omitted when the metric is fixed.

\section{The construction in dimensions three through seven}\label{sec:low-construction}
We construct a singular solution in each punctured coordinate ball
and an exterior function whose equation error is a linear combination
of the restrictions of $\Psi_1,\ldots,\Psi_{N_E}$. We then choose
the boundary data so that the corresponding conformal factors in
the metric $g$, and their $g$-normal derivatives, agree. The remaining
exterior residual is detected by the gradient of the joined action.

Following Parts I and II, at each prescribed point choose G\"unther's
exact conformal normal coordinates \cite{gunther,part1,part2} for a local
metric $h=\Theta^{4/(n-2)}g$, with $\Theta(0)=1$ and $d\Theta(0)=0$.
In these fixed geodesic normal coordinates $X$,
\[
 \det h\equiv1,\qquad
 \sum_j h^{ij}X_j=X_i,\qquad
 h-g_{\rm E}=O(r^2),\qquad \operatorname{Scal}_h=O(r^2).
\]
The vanishing quadratic volume coefficient gives
$\operatorname{Ric}_h(0)=0$, and hence
\[
 \Delta_hX_j=-\frac23\sum_k(\operatorname{Ric}_h)_{jk}(0)X_k
                  +O(r^2)=O(r^2).
\]
These curvature consequences also follow from
\cite[Theorem 5.1 and Corollary 5.4]{lp}.
The normalization imposes no restriction on the original Weyl jets.
All local metrics and charts are chosen before $\rho$ and before any
background or matching parameter is varied. A global conformal factor
$W$ is represented in this local metric by $\Theta^{-1}W$.

Fix the prescribed positive neck ratios $t_1,\ldots,t_m$.
In this section we use the matching radius $\rho$ as the small
construction parameter. For each sufficiently small $\rho>0$, set
\begin{equation}
 \eps=\rho^{(2n-1)/4},\qquad
 \rho_i=t_i^{4/(2n-1)}\rho,\qquad
 \eps_i=\rho_i^{(2n-1)/4}=t_i\eps,
 \qquad
 \Om=M\setminus\bigcup_{i=1}^m B_{\rho_i}(p_i).
 \label{eq:scales}
\end{equation}
The balls are disjoint for small $\rho$, their radii are comparable,
and the neck ratios are exactly the prescribed ratios.
Parameter derivatives in this section vary the background or matching
data while keeping $\rho$ and the $t_i$ fixed; hence the necks and all
matching spheres are fixed as well. The numerical background vector $\bgpar$ ranges over the fixed compact
neighborhood of $\bgcoord\K$ selected in Section~\ref{sec:compact}.

\subsection{Local singular solutions and their boundary data}
\label{sec:local-low}

We first work at one puncture and suppress its index, writing $\rho$
and $\eps$ for its radius and neck.

The local functions are defined on $B_\rho(0)\setminus\{0\}$ and have
regularity up to the outer sphere $S_\rho$.

\medskip\noindent\textit{Interior and boundary norms.}
For $0<r\le\rho/2$, put $A_r=\{X:r\le|X|\le2r\}$ and define
\[
 \|z\|_{C^{j,\alphaH}_{\rm sc}(A_r)}
 :=\|z(r\,\cdot)\|_{C^{j,\alphaH}(\{1\le|y|\le2\})},
 \qquad
 \|z\|_{j,\alphaH;\gamma;\rho}
 :=\sup_{0<r\le\rho/2}r^{-\gamma}
           \|z\|_{C^{j,\alphaH}_{\rm sc}(A_r)}.
\]
The scaled annular norm includes factors $r^k$ on derivatives of order
$k$ and $r^{j+\alphaH}$ on the highest H\"older seminorm.
The factor $r^{-\gamma}$ controls growth or decay towards the puncture;
$\rho$ specifies the domain and its outer scale. In particular,
\[
 \|z\|_{j,\alphaH;\gamma;\rho}
 =\rho^{-\gamma}\|z(\rho\,\cdot)\|_{j,\alphaH;\gamma;1}.
\]
We write $C^{j,\alphaH}_\gamma(\overline{B_\rho(0)}\setminus\{0\})$
for functions with finite weighted norm and $C^{j,\alphaH}$ regularity
on closed annuli away from zero, including the outer boundary.
This notation imposes no ordinary regularity at the puncture.

For a function $\phi$ on $S_\rho$, define
\[
 \|\phi\|_{C^{j,\alphaH}_{\rm sc}(S_\rho)}
 :=\|\phi(\rho\,\cdot)\|_{C^{j,\alphaH}(\Sph^{n-1})}.
\]
This norm includes the factors $\rho^k$ on tangential derivatives of
order $k$ and $\rho^{j+\alphaH}$ on the highest H\"older seminorm.
There is no additional factor $\rho^{-\gamma}$ multiplying the whole norm.
Boundary data written directly as functions of $\theta$ are already
on the unit sphere and use its ordinary H\"older norm.
The outer annulus gives, for $j\ge1$,
\[
 \begin{aligned}
 \|z|_{S_\rho}\|_{C^{j,\alphaH}_{\rm sc}(S_\rho)}
 &\le C\rho^\gamma\|z\|_{j,\alphaH;\gamma;\rho},\\
 \|\rho\partial_rz|_{S_\rho}\|_{C^{j-1,\alphaH}_{\rm sc}(S_\rho)}
 &\le C\rho^\gamma\|z\|_{j,\alphaH;\gamma;\rho}.
 \end{aligned}
\]
At derivative order $j$, we write $z=O_{\rm sc}(r^\gamma)$ when
$\|z\|_{C^{j,\alphaH}_{\rm sc}(A_r)}\le Cr^\gamma$ for
$0<r\le\rho/2$, equivalently $\|z\|_{j,\alphaH;\gamma;\rho}\le C$.
The constant is uniform in $r$ and in the small scales and parameter
ranges specified in the estimate. Unless otherwise stated, estimates
for a function use $j=2$. Parameter derivatives are included only when stated.

We use $\dd\theta$ for the round
measure on $\Sph^{n-1}$. Euclidean indices range over $1,\ldots,n$, and angular indices over
$1,\ldots,n-1$. Cylinder indices range over $0,\ldots,n-1$;
the index $0$, also denoted by $t$, is axial.
For functions on $\Sph^{n-1}$, set
\[
 \fint  f=|\Sph^{n-1}|^{-1}\int_{\Sph^{n-1}}f\dd\theta.
\] We use the orthogonal
spherical harmonic projections
\[
 \Pi_0f:=\fint f,\qquad
 \Pi_1f:=n\sum_{i=1}^n\fint (\theta_i f)\theta_i,\qquad
 \Pi_{\ge2}f:=f-\Pi_0f-\Pi_1f.
\]
The coefficient vector is $[\Pi_1f]=(n\fint (\theta_i f))_{i=1}^n$,
so $\Pi_1f=[\Pi_1f]\cdot\theta$.

\medskip\noindent\textit{An exact two-ended Fowler model.}
We use the Fowler and Kelvin formulas in
\cite[Proposition 2.2 and Corollary 2.1]{santos}
and \cite[Section 2, equations (7)--(12)]{kmps}.
The coefficient of the Euclidean equation in these model formulas
is $\frac{n(n-2)}4$, as here.
Use the Fowler profile with $v_\eps(0)=\eps$ and $v_\eps'(0)=0$.
Let $c$ vary in a fixed compact subinterval of $(0,\infty)$, and set
\[
 R=R(\eps,c)=\left(\frac{\eps}{2c}\right)^{2/(n-2)},
 \qquad
 \uF_{\eps,c}(r)=r^{-(n-2)/2}
                  v_\eps\!\left(\log\frac Rr\right).
\]
For fixed $\eps$, the parameter $c$ determines the phase $\log R$.

Denote the special conformal deformation fixing the puncture by
$\uF_{\eps,c,\beta}$, with $\beta$ in a fixed bounded subset of $\R^n$.
For $\beta\ne0$, it is an exact Euclidean Yamabe solution with
singularities at $0$ and $\beta/|\beta|^2$; for $\beta=0$ the second
end is at infinity and the solution is radial. We use only its
restriction to a ball with $\rho|\beta|<1/2$, excluding the second
singularity. The explicit formula is
\[
 \uF_{\eps,c,\beta}(X)
 =\bigl|X-\beta|X|^2\bigr|^{-(n-2)/2}
 v_\eps\bigl(-2\log|X|+\log|X-\beta|X|^2|+\log R\bigr).
\]
In calculations with fixed parameters, write
$\uF=\uF_{\eps,c,\beta}$ and $\uF_0=\uF_{\eps,c,0}$ for its radial
member. Its first Kelvin variation is
\begin{equation}
 \uF(X)=\uF_0(r)
 + (\beta\cdot X)\bigl[(n-2)\uF_0(r)+r\uF_0'(r)\bigr]+E_\beta(X),
 \qquad
 \|E_\beta\|_{C^{2,\alphaH}_{\rm sc}(A_r)}
       \le C|\beta|^2r^{(6-n)/2}.
 \label{eq:fowler-kelvin-expansion}
\end{equation}
This is the first-order term in the refined Fowler asymptotics of
\cite[Theorem 1.1]{part1}, in our normalization. The parameter $\beta$
adjusts the degree-one boundary data while preserving the puncture.
The exact model $\uF_{\eps,c,\beta}$ solves the Euclidean equation.
Finally,
\[
 \frac R\rho=(2c)^{-2/(n-2)}\rho^{3/[2(n-2)]}\longrightarrow0,
\]
so the matching sphere lies well outside the reference Fowler scale.
On the fixed parameter boxes, positivity of $v_\eps$ gives
\begin{equation}
 C^{-1}\eps r^{-(n-2)/2}\le \uF\le Cr^{-(n-2)/2},\qquad
 R/\rho\asymp\rho^{3/[2(n-2)]}.\label{eq:model-bounds}
\end{equation}

The Hamiltonian and differentiated Fowler ODE give
$|v_\eps'|\le \frac{n-2}{2}v_\eps$ and $|v_\eps^{(j)}|\le C_jv_\eps$.
With $s=-2\log r+\log|X-\beta r^2|+\log R$, differentiation gives
\[
\begin{aligned}
 \partial_c\uF&=-\frac{2\uF}{(n-2)c}\frac{v_\eps'(s)}{v_\eps(s)},\\
 \partial_{\beta_j}\uF
 &=-\uF\left(\frac{v_\eps'(s)}{v_\eps(s)}-\frac{n-2}{2}\right)
       \frac{r^2(X_j-\beta_jr^2)}{|X-\beta r^2|^2},\\
 |\partial_c\uF|&\le C\uF,\qquad |\partial_\beta \uF|\le Cr\uF,\qquad
 (\uF)^{p-1}\le Cr^{-2}.
\end{aligned}
\]
These bounds hold in scaled annular norms. Moreover,
\[
 \log\frac{\uF}{\uF_{\eps,c,0}}
 =-\frac{n-2}{2}\log\frac{|X-\beta r^2|}{r}
   +\int_{-\log r+\log R}^{s}\frac{v_\eps'(t)}{v_\eps(t)}\dd t
 =O(r),
\]
so $\uF\asymp \uF_{\eps,c,0}$ on the fixed set of Kelvin parameters.

\begin{lemma}\label{lem:model-inverse}
Fix $3\le n\le7$, $\mu=9/8$, $0<c_*<c^*$, and $B_*>0$. Let
$c\in[c_*,c^*]$, $|\beta|\le B_*$, and
$\eps=\rho^{(2n-1)/4}$, with $R=R(\eps,c)$ as above.
For $\uF=\uF_{\eps,c,\beta}$ on $B_\rho(0)\setminus\{0\}$, define
\[
 \mathscr L_{\uF}=-D\Ff_{g_{\rm E}}(\uF)
 =\Delta+\frac{n(n+2)}4(\uF)^{4/(n-2)}.
\]
For all sufficiently small $\rho>0$, this Euclidean operator is an
isomorphism
\[
 \begin{aligned}
 \mathscr L_{\uF}:\quad
 &\{w\in C^{2,\alphaH}_\mu(\overline{B_\rho(0)}\setminus\{0\}):
                 \Pi_{\ge2}(w|_{S_\rho})=0\}\\
 &\hspace{3em}\longrightarrow
 C^{0,\alphaH}_{\mu-2}(\overline{B_\rho(0)}\setminus\{0\}).
 \end{aligned}
\]
Its inverse $G_{\uF}$ satisfies
\[
 \|G_{\uF}f\|_{2,\alphaH;\mu;\rho}
 \le C\|f\|_{0,\alphaH;\mu-2;\rho},
\]
with $C$ independent of $\rho,c,\beta$ in these ranges.
The boundary condition prescribes the harmonics of degree at least two.
The weighted condition determines the mean and degree-one components of
$G_{\uF}f$.
\end{lemma}
\begin{proof}
For $\beta=0$, \cite[Proposition 2.3 and Corollaries 2.3--2.4]{santos}
provides a right inverse with the uniform bound in the statement.
To prove injectivity, write a solution of $\mathscr L_{\uF}w=0$ as
$w=r^{-(n-2)/2}z(s,\theta)$, where $s=-\log r$. Its coefficient in any
spherical harmonic of degree $l$ satisfies
\[
 -z_l''+\left(\frac{(2l+n-2)^2}{4}
       -\frac{n(n+2)}4v_\eps(s+\log R)^{p-1}\right)z_l=0.
\]
Membership in $C^{2,\alphaH}_\mu$ requires
$|z_l|+|z_l'|=O(e^{-(n-2+2\mu)s/2})$ as $s\to\infty$.
For $l=0$, differentiating the Fowler profile in its phase and neck
parameter gives the two independent solutions
\[
 v_\eps'(s+\log R),\qquad
 (\partial_\eps v_\eps)(s+\log R),
\]
where $\partial_\eps$ holds the argument of $v_\eps$ fixed.
The first is periodic; the second is a sum of a periodic function and
$s$ times a periodic function, as follows by differentiating the
periodicity relation. No nonzero combination has the required decay.
For $l=1$, Euclidean translations and the Kelvin deformation give
\cite[equations (11)--(12)]{kmps}
\[
 e^{s}\left(\frac{n-2}{2}v_\eps+v_\eps'\right)(s+\log R),
 \qquad
 e^{-s}\left(\frac{n-2}{2}v_\eps-v_\eps'\right)(s+\log R).
\]
Their corresponding functions $w$ have orders $r^{-n/2}$ and
$r^{(4-n)/2}$, respectively. Both are excluded because
$\mu>(4-n)/2$.
For $l\ge2$, the potential in the displayed ODE is at least
$(n+2)/2$, since $v_\eps\le1$. The boundary condition gives
$z_l(-\log\rho)=0$, and the weight gives $z_l(s)\to0$.
The maximum principle on the half-line therefore implies $z_l=0$.
Choose $\rho$ so that $\rho B_*<1/2$, which keeps the second
singularity outside the ball. For bounded $\beta$, the model estimates give
$\uF_{\eps,c,\beta}/\uF_{\eps,c,0}=1+O_{\rm sc}(|\beta|r)$ and
$(\uF_{\eps,c,0})^{p-1}\le Cr^{-2}$. Hence
\[
 \|(\mathscr L_{\uF_{\eps,c,\beta}}-\mathscr L_{\uF_{\eps,c,0}})w\|_{0,\alphaH;\mu-2;\rho}
 \le C|\beta|\rho\|w\|_{2,\alphaH;\mu;\rho}.
\]
For small $\rho$, composing this operator difference with the
radial inverse has norm less than one. A Neumann series therefore
gives the inverse of $\mathscr L_{\uF}$.
\end{proof}

We estimate the residual $\Ff_h(\uF)$ of the Euclidean model in the
metric $h$.

\begin{lemma}\label{lem:geometry}
For a radial function $f$,
\begin{align}
 (\Delta_h-\Delta)f&=0,\label{eq:radial-cancel}\\
 (\Delta_h-\Delta)(X_jf)&=f\Delta_hX_j.\label{eq:linear-cancel}
\end{align}
For the model $\uF$,
\begin{equation}
 \norm{(\Delta_h-c(n)\operatorname{Scal}_h)\uF+\frac{n(n-2)}4(\uF)^p}_{0,\alphaH;\mu-2;\rho}
       \le C\rho^{(10-n)/2-\mu},\label{eq:geometric-source}
\end{equation}
also after one $c$ or $\beta$ derivative.
\end{lemma}
\begin{proof}
Using $\det h=1$ and $\sum_jh^{ij}X_j=X_i$ gives
\[
 \Delta_hf=f''+\frac{n-1}{r}f',
 \qquad \langle\nabla_hX_j,\nabla_hf\rangle_h=X_jf'/r.
\]
The first formula proves \eqref{eq:radial-cancel}; the product rule
for $\Delta_h(X_jf)$ proves \eqref{eq:linear-cancel}.
Apply the expansion \eqref{eq:fowler-kelvin-expansion}, with
$\uF_0=\uF_{\eps,c,0}$ its radial member.
The Euclidean equation for $\uF$ makes the residual in
\eqref{eq:geometric-source} equal to
$(\Delta_h-\Delta)\uF-c(n)\operatorname{Scal}_h\uF$.
The radial profiles have scaled derivatives $O(r^{-(n-2)/2})$.
The four geometric contributions are
\[
 \begin{aligned}
 (\Delta_h-\Delta)\uF_0&=0,\\
 (\Delta_h-\Delta)\sum_j\beta_jX_j((n-2)\uF_0+r\uF_0')&=O(|\beta|r^{(6-n)/2}),\\
 (\Delta_h-\Delta)E_\beta&=O(|\beta|^2r^{(6-n)/2}),\qquad
 \operatorname{Scal}_h\uF=O(r^{(6-n)/2}).
 \end{aligned}
\]
The second line uses $\Delta_hX_j=O(r^2)$ and the radial cancellation;
the bound for $(\Delta_h-\Delta)E_\beta$ uses
$h-g_{\rm E}=O(r^2)$. The scalar curvature bound gives the estimate
for $\operatorname{Scal}_h\uF$. These estimates hold in the scaled $C^{0,\alphaH}$ norms and retain
the displayed bounds after one $c$ derivative. By
\eqref{eq:linear-cancel}, differentiating the linear Kelvin
contribution to the residual with respect to $\beta_k$ gives
\[
 \bigl[(n-2)\uF_0+r\uF_0'\bigr]\Delta_hX_k
 =O(r^{(6-n)/2}).
\]
Taylor's formula applied to the explicit Kelvin model gives
$(\Delta_h-\Delta)\partial_{\beta_k}E_\beta
=O(|\beta|r^{(6-n)/2})$.
The scalar-curvature contribution satisfies
$\operatorname{Scal}_h\,\partial_{\beta_k}\uF
=O(r^{(8-n)/2})$.
On the fixed parameter box, the full $\beta_k$ derivative of the
residual is therefore $O(r^{(6-n)/2})$.
Multiplication by $r^{2-\mu}$ gives
\eqref{eq:geometric-source}, since $\frac{10-n}{2}-\mu\ge3/8$.
\end{proof}

For $\phi\in C^{2,\alphaH}(\Sph^{n-1})$ with $\Pi_{\ge2}\phi=\phi$,
let $P_\rho\phi$ be its harmonic extension to $B_\rho$, with boundary
value $\phi(\theta)$ at $X=\rho\theta$. By \cite[Proposition 2.4 and equations
(2.17)--(2.18)]{santos}, it satisfies
\[
 \norm{P_\rho\phi}_{C^{2,\alphaH}_{\rm sc}(A_r)}
 \le C\norm\phi(r/\rho)^2,\qquad
 \norm{P_\rho\phi}_{2,\alphaH;\mu;\rho}\le C\rho^{-\mu}\norm\phi.
\]
Here $\|\phi\|=\|\phi\|_{C^{2,\alphaH}(\Sph^{n-1})}$.

We seek a solution of the form $V=\uF+P_\rho\phi+v$.
Set $\delta_0=1/100$, fix $C_*$ large before choosing $\rho$ small,
and require $v$ to belong to the ball
\[
 \left\{v\in C^{2,\alphaH}_\mu(\overline{B_\rho}\setminus\{0\}):
   \Pi_{\ge2}v|_{S_\rho}=0,\quad
   \|v\|_{2,\alphaH;\mu;\rho}\le C_*\rho^{1-\mu}\right\}.
\]
The trace estimates give
\[
 \|v|_{S_\rho}\|_{C^{2,\alphaH}_{\rm sc}(S_\rho)}
 +\|\rho\partial_rv|_{S_\rho}\|_{C^{1,\alphaH}_{\rm sc}(S_\rho)}
 \le CC_*\rho.
\]
Thus the boundary data are $O(\rho)$ even though the radius
$C_*\rho^{1-\mu}$ in the weighted norm grows as $\rho\to0$.

\begin{proposition}\label{prop:local}
For all sufficiently small $\rho$, set $\eps=\rho^{(2n-1)/4}$ as in
\eqref{eq:scales}. For $(c,\beta)$ in the fixed parameter boxes and
$\norm\phi\le C\rho^{1-\delta_0}$, there is a unique $v$ in the
displayed ball such that
\[
 V=\uF+P_\rho\phi+v>0,\qquad
 (\Delta_h-c(n)\operatorname{Scal}_h)V+\frac{n(n-2)}4V^p=0.
\]
It depends $C^1$ on $(c,\beta,\phi)$, with
\begin{align}
 \norm v_{2,\alphaH;\mu;\rho}
 +\norm{\partial_cv}_{2,\alphaH;\mu;\rho}
 +\norm{\partial_\beta v}_{2,\alphaH;\mu;\rho}
 &=o(\rho^{1-\mu}),\label{eq:local-small}\\
 \norm{D_\phi v}_{\mathcal L(C^{2,\alphaH}_{\ge2}(\Sph^{n-1}),C^{2,\alphaH}_\mu(\overline{B_\rho}\setminus\{0\}))}
 &\le o(1)\rho^{-\mu}.\label{eq:local-lip}
\end{align}
The traces $v|_{S_\rho}$ and $\rho\partial_rv|_{S_\rho}$ are
$o(\rho)$ in scaled $C^{2,\alphaH}$ and $C^{1,\alphaH}$ norms,
respectively, also after one derivative in $c$ or any $\beta_j$.
\end{proposition}
\begin{proof}
Put $w=P_\rho\phi$. Within this proof abbreviate the interior norms
$\|\cdot\|_{2,\alphaH;\mu;\rho}$ and
$\|\cdot\|_{0,\alphaH;\mu-2;\rho}$ by
$\|\cdot\|_\mu$ and $\|\cdot\|_{\mu-2}$, respectively. For $v$ in the
displayed ball, the lower bound for $\uF$ gives
\[
 \|(w+v)/\uF\|_{C^{2,\alphaH}_{\rm sc}(A_r)}
 \le\eta_\rho:=C\eps^{-1}\rho^{n/2-\delta_0}
 =C\rho^{1/4-\delta_0}\longrightarrow0.
\]
For $\eta_\rho\le1/2$, the nonlinear remainder is
\[
 Q_{\uF}(z)=\frac{n(n-2)}4
       \bigl[(\uF+z)^p-(\uF)^p-p(\uF)^{p-1}z\bigr].
\]
The equation for $v$ is equivalent to
\[
 v=G_{\uF}\bigl[(\Delta-\Delta_h)(\uF+w+v)
 +c(n)\operatorname{Scal}_h(\uF+w+v)-Q_{\uF}(w+v)-\frac{n(n+2)}4(\uF)^{p-1}w\bigr].
\]
Taylor's formula and the weighted product estimate
\eqref{eq:weighted-product-calculus} give
\[
 \frac{\|Q_{\uF}(w)\|_{\mu-2}}{\rho^{1-\mu}}
 \le C\rho^{1/4-2\delta_0}=o(1),\qquad
 \|D_vQ_{\uF}(w+v)\|_{\mu\to\mu-2}\le C\eta_\rho.
\]
The operator $\Delta_h-\Delta-c(n)\operatorname{Scal}_h$
has norm $O(\rho^2)$ from the space with norm $\|\cdot\|_\mu$
to the space with norm $\|\cdot\|_{\mu-2}$, and
\[
 \|(\Delta_h-\Delta)w\|_{\mu-2}
 \le C\rho^{3-\delta_0-\mu};
\]
the scalar curvature term gains two further powers of $\rho$.
To estimate $(\uF)^{p-1}P_\rho\phi$, use $\uF\asymp\uF_0:=\uF_{\eps,c,0}$,
$\uF_0\le2c$ on $R\le r\le\rho$, and $(\uF_0)^{p-1}\le Cr^{-2}$ on $r\le R$.
The two radial ranges give
\[
 r^{2-\mu}\|(\uF)^{p-1}P_\rho\phi\|_{C^{0,\alphaH}_{\rm sc}(A_r)}
 \le C\|\phi\|\rho^{-2}
 \begin{cases}
 r^{2-\mu},&r\le R,\\
 r^{4-\mu},&R\le r\le\rho.
 \end{cases}
\]
Taking the supremum yields
\begin{equation}
 \|(\uF)^{p-1}P_\rho\phi\|_{\mu-2}
 \le C\rho^{-\mu}\bigl[\rho^2+(R/\rho)^{2-\mu}\bigr]\|\phi\|.
 \label{eq:local-potential-split}
\end{equation}
For $\|\phi\|\le C\rho^{1-\delta_0}$, this bound for
$\|(\uF)^{p-1}P_\rho\phi\|_{\mu-2}$ is $o(\rho^{1-\mu})$, since
$-\delta_0+3(2-\mu)/[2(n-2)]\ge101/400$.
The geometric residual in Lemma~\ref{lem:geometry}, after applying
$G_{\uF}$, has norm at most $C\rho^{(10-n)/2-\mu}$. Relative to the
radius of the ball, its bound is
\[
 \frac{\rho^{(10-n)/2-\mu}}{\rho^{1-\mu}}
 =\rho^{(8-n)/2}\longrightarrow0,\qquad 3\le n\le7.
\]
Together these bounds show that the map on the right is a contraction with derivative norm less than $1/2$ and value
$o(\rho^{1-\mu})$ at zero. It has a unique fixed point in the
stated ball, with $\|v\|_\mu=o(\rho^{1-\mu})$ and
$V\ge(1-\eta_\rho)\uF>0$.

The linearization of $-\Ff_h(V)=0$ is
$L_V^h=-D\Ff_h(V)
=\Delta_h-c(n)\operatorname{Scal}_h+\frac{n(n+2)}4V^{p-1}$.
It differs from $\mathscr L_{\uF}$ by $O(\rho^2+\eta_\rho)$ between the spaces
in Lemma~\ref{lem:model-inverse}. Its inverse is uniformly
bounded, and the implicit function theorem gives $C^1$ dependence.
Differentiating the equation, using $\mathscr L_{\uF}(\uF_s)=0$ for $s=c,\beta_j$, gives
\[
\begin{aligned}
 L_V^h v_s
 &=-(\Delta_h-\Delta-c(n)\operatorname{Scal}_h)\uF_s
   -\frac{n(n+2)}4(V^{p-1}-(\uF)^{p-1})\uF_s,\\
 L_V^hD_\phi v[\dot\phi]
 &=-(\Delta_h-\Delta-c(n)\operatorname{Scal}_h)P_\rho\dot\phi
   -\frac{n(n+2)}4V^{p-1}P_\rho\dot\phi.
\end{aligned}
\]
Apply the uniform inverse of $L_V^h$, Lemma~\ref{lem:geometry}
and \eqref{eq:local-potential-split}:
\[
\begin{aligned}
 \|v_s\|_\mu
 &\le C\rho^{(10-n)/2-\mu}
   +C\rho^{-\mu}\bigl[\rho^2+(R/\rho)^{2-\mu}\bigr]\|\phi\|
   +C\|v\|_\mu
   =o(\rho^{1-\mu}),\\
 \|D_\phi v[\dot\phi]\|_\mu
 &\le C\rho^{-\mu}\bigl[\rho^2+(R/\rho)^{2-\mu}\bigr]\|\dot\phi\|
   =o(1)\rho^{-\mu}\|\dot\phi\|.
\end{aligned}
\]
Here $|(V^{p-1}-(\uF)^{p-1})\uF_s|\le C(\uF)^{p-1}(|w|+|v|)$ and
$V\asymp\uF$. The corresponding scaled H\"older estimates follow
from the already established relative bound
$\|(V-\uF)/\uF\|_{C^{2,\alphaH}_{\rm sc}(A_r)}\le\eta_\rho$.
This proves \eqref{eq:local-small}--\eqref{eq:local-lip}. The trace estimates give
\[
 \begin{aligned}
 &\|v|_{S_\rho}\|_{C^{2,\alphaH}_{\rm sc}(S_\rho)}
 +\|\rho\partial_rv|_{S_\rho}\|_{C^{1,\alphaH}_{\rm sc}(S_\rho)}\\
 &\qquad\le C\rho^\mu\|v\|_{2,\alphaH;\mu;\rho}=o(\rho).
 \end{aligned}
\]
The same argument applied to $\partial_cv$ and $\partial_{\beta_j}v$
gives their boundary estimates.
\end{proof}

To compute the leading boundary data, return to the radial Euclidean
model before the Kelvin deformation and the curved-metric correction:
\[
 \uF_0(r):=\uF_{\eps,c,0}(r)
 =r^{-(n-2)/2}v_\eps\!\left(\log\frac Rr\right),
 \qquad R=\left(\frac{\eps}{2c}\right)^{2/(n-2)}.
\]
Its data at $r=R$ are
$\uF_0(R)=2c$ and $R\uF_0'(R)=-(n-2)c$.
Integrating
$\uF_0''+(n-1)\uF_0'/r+\frac{n(n-2)}4(\uF_0)^p=0$ with these data gives
\[
 \uF_0(r)=c+\frac{\eps^2}{4c}r^{2-n}
 -\frac n4\int_R^r s\uF_0(s)^p
      \bigl[1-(s/r)^{n-2}\bigr]\dd s,
 \qquad R\le r\le\rho.
\]
Since $\uF_0\le2c$ on this annulus, the integral term and its scaled
radial derivative are $O(r^2)$. Thus
\[
 \uF_0(r)=c+\frac{\eps^2}{4c}r^{2-n}+O(r^2),\qquad R\le r\le\rho,
\]
so $c$ is the constant coefficient in this annular expansion.
The same bounds hold after one $c$ derivative with $\eps,r$ fixed:
$|\partial_c\uF_0|\le C\uF_0$ controls the differentiated integrand,
and $\partial_cR=-2R/[(n-2)c]$ makes the lower-endpoint contribution
$O(R^2)$. 

Since $R/\rho\to0$, the radial equation and the preceding expansion
  give uniformly bounded scaled derivatives of $\uF_0$ and
  $\partial_c\uF_0$ on $\rho/2\le r\le2\rho$.
  Expanding the exact Kelvin formula on this annulus therefore gives
  \[
   \|E_\beta\|_{C^{2,\alphaH}_{\rm sc}(A_{\rho/2})}
   +\|\partial_cE_\beta\|_{C^{2,\alphaH}_{\rm sc}(A_{\rho/2})}
   \le C|\beta|^2\rho^2,
   \qquad
   \|D_\beta E_\beta\|_{C^{2,\alphaH}_{\rm sc}(A_{\rho/2})}
   \le C|\beta|\rho^2.
  \]
  These bounds also control the scaled normal traces on $S_\rho$.
  For bounded $\beta$, the quadratic Kelvin remainder and its
  $c$ derivative are therefore $O(\rho^2)$ at the matching sphere,
  while its derivative with respect to $B=\rho\beta$ is $O(\rho)$.

Put $B=\rho\beta$ and
\[
 b_\eps(c)=\frac{\eps^2\rho^{2-n}}{4c}
           =\frac{\rho^{3/2}}{4c}.
\]
For the full local solution $V$, the mean and first harmonic
Dirichlet and scaled normal data on $S_\rho$ are
\begin{equation}
 \begin{aligned}
 \Pi_0V|_{S_\rho}&=c+b_\eps(c)+o(\rho),&
 \Pi_0(\rho\partial_rV)|_{S_\rho}&=(2-n)b_\eps(c)+o(\rho),\\
 [\Pi_1V|_{S_\rho}]&=(n-2)cB+o(\rho),&
 [\Pi_1(\rho\partial_rV)|_{S_\rho}]&=(n-2)cB+o(\rho).
 \end{aligned}
 \label{eq:low-model-traces}
\end{equation}
The pole cancels from $(n-2)\uF_0+r\uF_0'$, leaving $(n-2)c+O(r^2)$
in the Kelvin first variation; $P_\rho\phi$ has zero low projection.
By the model expansion and Proposition~\ref{prop:local}, each
remainder in \eqref{eq:low-model-traces} and its $c$ derivative
are $o(\rho)$. Its $B$ and $\phi$ derivatives are $o(1)$, with the
other inputs fixed.
For the $B$ derivative, use $\partial_Bv=\rho^{-1}\partial_\beta v$.
Its scaled Dirichlet and normal trace norms on $S_\rho$ are $o(1)$.
The derivatives in $(c,B,\phi)$ of
\[
 \Pi_{\ge2}(V-P_\rho\phi)|_{S_\rho},
 \qquad
 \Pi_{\ge2}\bigl(\rho\partial_r(V-P_\rho\phi)\bigr)|_{S_\rho}
\]
are $o(1)$ in the scaled $C^{2,\alphaH}$ and $C^{1,\alphaH}$
boundary norms, respectively.

\subsection{The exterior equation with prescribed spectral moments}
\label{sec:exterior-inverse}

The local necks are now fixed, but their outer boundary traces remain
adjustable. On the complementary domain we solve the Yamabe equation
up to an error in $E$, while prescribing the volume integrals against
$\Psi_1,\ldots,\Psi_{N_E}$.

The exterior domain is $\Om=M\setminus\bigcup_iB_{\rho_i}(p_i)$,
with the radii fixed at the start of this section.
Choose a fixed small chart radius $r_0$, independent of $\rho$, with
disjoint coordinate balls $B_{r_0}(p_i)$, and take $\rho$ small enough
that $\rho_i<r_0/2$ for every $i$. Let
$K_0=M\setminus\bigcup_iB_{r_0/2}(p_i)$. For $\rho_i\le r\le r_0/2$
let $A_{i,r}$ denote the coordinate annulus $r\le|X_i|\le2r$.
Define
\[
 \|z\|_{j,\alphaH;\gamma;\Om}
 :=\|z\|_{C^{j,\alphaH}(K_0)}
   +\max_i\sup_{\rho_i\le r\le r_0/2}
      r^{-\gamma}\|z\|_{C^{j,\alphaH}_{\rm sc}(A_{i,r})}.
\]
The scaled norm uses the fixed chart at $p_i$ and the map $X_i=ry$.
The space $C^{j,\alphaH}_\gamma(\overline\Om)$ includes regularity up
to each inner boundary sphere. Boundary norms on $S_{\rho_i}$ use the
pullback $X_i=\rho_i\theta$ to $\Sph^{n-1}$.

Set $\nu=7/4-n$ and write $\Psi=(\Psi_1,\ldots,\Psi_{N_E})$.
\begin{lemma}\label{lem:exterior-inverse}
For all sufficiently small $\rho$, and for prescribed
$f\in C^{0,\alphaH}_{\nu-2}(\overline\Om)$ and
$\mathbf m\in\R^{N_E}$, the augmented problem
\begin{equation}
 L_\bgpar z-\sum_{j=1}^{N_E}d_j\Psi_j=f\quad\hbox{on }\Om,\qquad
 (I-\Pi_0)(z|_{S_{\rho_i}})=0\quad(1\le i\le m),\qquad
 \int_{\Om}z\Psi\vol=\mathbf m,
 \label{eq:augmented-exterior}
\end{equation}
has a linear solution operator, $C^1$ in $\bgpar$, returning
$z\in C^{2,\alphaH}_\nu(\overline\Om)$ and $d\in\R^{N_E}$.
The boundary condition requires a constant on each sphere, but these
constants may differ. The solution operator chooses them as part of
its output and satisfies
\begin{equation}
 \norm z_{2,\alphaH;\nu;\Om}+|d|
 \le C(\norm f_{0,\alphaH;\nu-2;\Om}+|\mathbf m|),
 \label{eq:augmented-bound}
\end{equation}
uniformly, also after one $\bgpar$ derivative at fixed $(f,\mathbf m)$.
\end{lemma}
\begin{proof}
On each Euclidean annulus $\rho_i<|X_i|<r_0$, use the right
inverse $G_{\rm E}$ from \cite[Lemma 4.2]{santos}. It maps the
weighted $C^{0,\alphaH}_{\nu-2}$ source space to
$C^{2,\alphaH}_{\nu}$, satisfies $-\Delta G_{\rm E}f=f$ and
$(I-\Pi_0)(G_{\rm E}f)|_{|X_i|=\rho_i}=0$, and has norm bounded
independently of $\rho_i$.
On a puncture ball of fixed small radius $r_0$, the coefficient errors
of orders two and one are $O(r^2)$ and $O(r)$; the potential is bounded.
The weighted operator bound is $\|(L_\bgpar+\Delta)G_{\rm E}\|\le Cr_0^2$, uniformly in the hole
radii and the compact parameter. The Neumann series defines
\[
 G_\bgpar=G_{\rm E}[I+(L_\bgpar+\Delta)G_{\rm E}]^{-1},\qquad L_\bgpar G_\bgpar=I.
\]

Cut off these local solutions to obtain $z_0$ with constant inner
traces and $L_\bgpar z_0=f$ near the holes. Its moments are bounded since
$\int_{\rho_i}^{r_0}r^{\nu+n-1}\dd r\le C$ and $\nu+n=7/4>0$.
The corrected source $f-L_\bgpar z_0$ vanishes near every puncture,
so its zero extension is defined on $M$.
Apply \eqref{eq:augmented-compact} to
\[
 \left(f-L_\bgpar z_0,\ \mathbf m-\int_{\Om}z_0\Psi\vol\right)
\]
and denote the result by $(z_1,d)$.
Choose a fixed $r_1\ll r_0$ and then $\rho$ small enough that
$\rho_i<r_1$. Take $\chi=0$ for $r\le r_1$ and $\chi=1$ for $r\ge2r_1$,
with $f-L_\bgpar z_0=0$ on the transition region.
For $(z_0+\chi z_1,d)$, the equation and moment errors are
\[
 (\chi-1)\sum_jd_j\Psi_j+[L_\bgpar,\chi]z_1,
 \qquad -\int_M(1-\chi)z_1\Psi\vol.
\]
Measure the equation error in $C^{0,\alphaH}_{\nu-2}$ and the
moment error in $\R^{N_E}$. The bounds $|D^j\chi|\le Cr_1^{-j}$ give
\[
\begin{aligned}
 \left\|(\chi-1)\sum_jd_j\Psi_j\right\|_{0,\alphaH;\nu-2;\Om}
 &\le C(r_0)r_1^{2-\nu}
       (\|f\|_{0,\alphaH;\nu-2;\Om}+|\mathbf m|),\\
 \|[L_\bgpar,\chi]z_1\|_{0,\alphaH;\nu-2;\Om}
 &\le C(r_0)r_1^{-\nu}
       (\|f\|_{0,\alphaH;\nu-2;\Om}+|\mathbf m|),\\
 \left|\int_M(1-\chi)z_1\Psi\vol\right|
 &\le C(r_0)r_1^n
       (\|f\|_{0,\alphaH;\nu-2;\Om}+|\mathbf m|).
\end{aligned}
\]
All three exponents are positive. Write $\mathcal T_0$ for the
approximate solution map and $\mathcal E$ for the combined error
in the equation and moment conditions. Fix $r_1$ so that $\|\mathcal E\|<1/2$.
Then $\mathcal T_0(I+\mathcal E)^{-1}$ is a solution operator with
\eqref{eq:augmented-bound}. The cutoffs and principal coefficients are
fixed in $\bgpar$. Differentiating the formula for $G_\bgpar$ and applying
\eqref{eq:inverse-parameter-calculus} to the compact inverse and
$I+\mathcal E$ proves the first $\bgpar$-derivative bound in
\eqref{eq:augmented-bound}.

\end{proof}

The zero-moment projected compact inverse maps smooth functions to
smooth functions by elliptic regularity; extend it to distributions by
transposition. Define $G_{i,\bgpar}$ to be the distributional solution
with zero spectral moments of

\[
 L_\bgpar G_{i,\bgpar}=(n-2)|\Sph^{n-1}|\delta_{p_i}
                 +\sum_jd^{\mathrm G}_{i,j}(\bgpar)\Psi_j,
 \qquad \int_MG_{i,\bgpar}\Psi_j\vol=0.
\]
Multiply $r^{2-n}$ by a smooth function that equals one near
$p_i$ and is supported in $B_{r_0}(p_i)$. Applying $L_\bgpar$
to this product gives the displayed Dirac mass and an
$O_{\rm sc}(r^{2-n})$ remainder near $p_i$.
For $X=re$ with $|e|=1$, the substitution $Y=rz$ gives
\[
 \int_{|Y|\le2r}|X-Y|^{2-n}|Y|^{2-n}\dd Y
 =r^{4-n}\int_{|z|\le2}|e-z|^{2-n}|z|^{2-n}\dd z
 \le Cr^{4-n}.
\]
The kernel of the compact projected inverse is locally bounded
by $C|X-Y|^{2-n}$, up to a smooth term. Integrating over
$|Y|\ge2r$ gives
\[
 \int_{B_{r_0}}|X-Y|^{2-n}|Y|^{2-n}\dd Y
 \le C\left(r^{4-n}+\int_{2r}^{r_0}s^{3-n}\dd s\right)
 \le Cr^{3-n},\qquad r=|X|\le r_0/4,\quad r_0\le1.
\]
The smooth source terms and their images under the compact augmented
inverse are bounded. Scaled Schauder estimates
on enlarged annuli therefore give
\begin{equation}
 G_{i,\bgpar}=r^{2-n}+O_{\mathrm{sc}}(r^{3-n}),
 \label{eq:green}
\end{equation}
where the remainder has scaled $C^{2,\alphaH}$ norm at most $Cr^{3-n}$.
The pole coefficient is fixed in $\bgpar$, and
\[
 L_\bgpar D_\bgpar G_{i,\bgpar}[\dot \bgpar]=-(D_\bgpar L_\bgpar[\dot \bgpar])G_{i,\bgpar}
       +\sum_jD_\bgpar d^{\mathrm G}_{i,j}[\dot \bgpar]\Psi_j,
 \qquad \int_MD_\bgpar G_{i,\bgpar}[\dot \bgpar]\Psi_j\vol=0.
\]
The differentiated source
$-(D_\bgpar L_\bgpar[\dot \bgpar])G_{i,\bgpar}+\sum_jD_\bgpar d^{\mathrm G}_{i,j}[\dot \bgpar]\Psi_j$
is $O_{\rm sc}(r^{2-n})$. The convolution and scaled Schauder estimates
give $D_\bgpar G_{i,\bgpar}[\dot \bgpar]=O_{\rm sc}(r^{3-n})$ for unit $\dot \bgpar$.
Both $G_{i,\bgpar}$ and $D_\bgpar G_{i,\bgpar}[\dot \bgpar]$ are uniformly bounded in
$C^{2,\alphaH}$ on every fixed compact set away from $p_i$,
including neighborhoods of the other punctures.

For $\psi\in C^{2,\alphaH}(\Sph^{n-1})$ with $\Pi_0\psi=0$,
let $Q_\rho\psi$ be the Euclidean harmonic extension to
$\R^n\setminus B_\rho$ with trace $\psi(\theta)$ at $X=\rho\theta$
and limit zero at infinity, multiplied by a fixed radial cutoff
that is one for $r\le r_0/2$ and zero for $r\ge r_0$.
Extend it by zero outside the coordinate chart, and use
$\|\psi\|=\|\psi\|_{C^{2,\alphaH}(\Sph^{n-1})}$. The estimate in
\cite[Proposition 2.5]{santos}, after scaling, gives
\[
 \norm{Q_\rho\psi}_{C^{2,\alphaH}_{\rm sc}(A_r)}
 \le C\norm\psi\rho^{n-1}r^{1-n},\qquad r\ge\rho.
\]
For each $(\bgpar,\mathsf m,\psi)$ in the ranges below, seek
\begin{equation}
 B_\bgpar=U_\bgpar+\sum_{i=1}^m\mathsf m_iG_{i,\bgpar}+\sum_{i=1}^mQ_{\rho_i}\psi_i+Z,
 \qquad \psi_i=\omega_i\cdot\theta+\vartheta_i,\label{eq:exterior-family}
\end{equation}
where $\omega_i\in\R^n$ and $\Pi_{\ge2}\vartheta_i=\vartheta_i$.
Work in the ranges
\begin{equation}
 |\mathsf m_i|\le C\rho^{n-3/2},\qquad
 \norm{\psi_i}\le C\rho^{1-\delta_0},\qquad
 \norm Z_{2,\alphaH;\nu;\Om}\le C\rho^{1-\nu}.\label{eq:exterior-boxes}
\end{equation}
Use the exterior fixed-point method of
\cite[Theorem 4.3]{santos}, imposing the augmented equation and
exact moments
\begin{equation}
 \Ff(B_\bgpar)=\sum_{j=1}^{N_E}d_j(\bgpar,\mathsf m,\psi)\Psi_j\quad\hbox{on }\Om,
 \qquad \int_{\Om} B_\bgpar\Psi_j\vol=\eta_j.
 \label{eq:exact-moments}
\end{equation}
The coefficients $d_j$ are solved for together with the correction $Z$.
Differentiating the moment condition in $\eta_l$ gives
$\int_{\Om}\Psi_j\partial_{\eta_l}B_\bgpar\vol=\delta_{jl}$.
Set $H=\sum_i\mathsf m_iG_{i,\bgpar}+\sum_iQ_{\rho_i}\psi_i+Z$ and
$\Lambda_i=\mathsf m_i\rho_i^{2-n}$. To estimate all parameter
derivatives together, use directions satisfying
\begin{equation}
 |\dot \bgpar|+\rho^{-1/2}|\dot\Lambda|
       +\rho^{\delta_0-1}\|\dot\psi\|_{C^{2,\alphaH}}\le1.
 \label{eq:low-exterior-parameter-directions}
\end{equation}
At this stage, derivatives with respect to $(\bgpar,\Lambda,\psi)$ keep $Z$ fixed.
Linearity in $\mathsf m,\psi$ and the uniform $C^1_\bgpar$ Green bounds
show that $H$ and $\dot H$ satisfy the same estimates:
\[
 \|H\|_{0,\alphaH;\nu;\Om}\le C\rho^{n-3/2},\qquad
 \|H\|_{0,\alphaH;-2;\Om}\le C\rho^{3/2},\qquad
 \|H\|_{0,\alphaH;0;\Om}\le C\rho^{1/2}.
\]
Use $r^{2-n-\nu}=r^{1/4}$ and $r^{1-n-\nu}=r^{-3/4}$.
On the annuli $\rho_i\le r\le r_0$ around $p_i$, use
$\rho_i\asymp\rho$ to obtain
\[
\begin{aligned}
 r^2\bigl(|\mathsf m_iG_{i,\bgpar}|+|Q_{\rho_i}\psi_i|+|Z|\bigr)
 &\le C\bigl[\rho^{n-3/2}r^{4-n}
       +\rho^{n-\delta_0}r^{3-n}+\rho^{n-3/4}r^{15/4-n}\bigr]\\
 &\le C\bigl[\rho^{3/2}+\rho^{3-\delta_0}+\rho^{9/4}\bigr],\\
 |\mathsf m_iG_{i,\bgpar}|+|Q_{\rho_i}\psi_i|+|Z|
 &\le C\bigl[\rho^{1/2}+\rho^{1-\delta_0}+\rho\bigr]
 \le C\rho^{1/2}.
\end{aligned}
\]
The smooth contributions from the other punctures satisfy these
bounds as well. Thus $U_\bgpar+tH$ stays in a fixed positive interval
for $0\le t\le1$.
The Taylor remainder
\[
 N_\bgpar(H)=-\frac{n(n-2)}4[(U_\bgpar+H)^p-U_\bgpar^p-pU_\bgpar^{p-1}H]
       =-\frac{n(n+2)}4(p-1)H^2\int_0^1(1-t)(U_\bgpar+tH)^{p-2}\dd t
\]
has a uniformly bounded scaled H\"older coefficient, also after one
parameter derivative in \eqref{eq:low-exterior-parameter-directions}.
Apply \eqref{eq:weighted-product-calculus} with weights $-2$ and
$\nu$ to $H^2$ and $H\dot H$. We obtain
\[
\|N_\bgpar(H)\|_{0,\alphaH;\nu-2;\Om}
 +\|\dot N_\bgpar(H)\|_{0,\alphaH;\nu-2;\Om}\le C\rho^n.
\]
For $Z_1,Z_2$ satisfying
$\|Z_k\|_{2,\alphaH;\nu;\Om}\le C\rho^{1-\nu}$, $k=1,2$, keep
$(\bgpar,\Lambda,\psi)$ fixed. The corresponding functions $H_1,H_2$
then satisfy $H_1-H_2=Z_1-Z_2$. Equation~\eqref{eq:weighted-product-calculus} gives
\[
 \|N_\bgpar(H_1)-N_\bgpar(H_2)\|_{0,\alphaH;\nu-2;\Om}
 \le C\rho^{3/2}\|Z_1-Z_2\|_{2,\alphaH;\nu;\Om}.
\]

Since $L_\bgpar G_{i,\bgpar}=\sum_jd^{\mathrm G}_{i,j}\Psi_j$ away from the pole,
write $\mathbf d_i^{\mathrm G}=(d^{\mathrm G}_{i,j})_{j=1}^{N_E}$ and put
$d_{\rm err}=d-d_0-\sum_i\mathsf m_i\mathbf d_i^{\mathrm G}$.
The equation for $Z$ and $d_{\rm err}$ is
\[
 L_\bgpar Z-\sum_jd_{{\rm err},j}\Psi_j=-L_\bgpar\sum_iQ_{\rho_i}\psi_i-N_\bgpar(H).
\]
The source $L_\bgpar\sum_iQ_{\rho_i}\psi_i$ has
$C^{0,\alphaH}_{\nu-2}(\overline\Om)$ norm $O(\rho^{n-\delta_0})$:
its geometric and potential terms are bounded by
$C\|\psi_i\|\rho_i^{n-1}r^{1-n}$, and the source weight gives
$r^{3-n-\nu}=r^{5/4}$. Cutoff errors lie on fixed annuli.
Substituting \eqref{eq:exterior-family} into the moment condition gives
the prescribed moments $\mathbf m_Z=\int_{\Om}Z\Psi\vol$:
\begin{equation}
\begin{aligned}
 \mathbf m_Z
 &=\bgpar-\int_{\Om}\left(U_\bgpar+\sum_i\mathsf m_iG_{i,\bgpar}
                         +\sum_iQ_{\rho_i}\psi_i\right)\Psi\vol\\
 &=\int_{M\setminus\Om}U_\bgpar\Psi\vol
   +\sum_i\mathsf m_i\int_{M\setminus\Om}G_{i,\bgpar}\Psi\vol
   -\sum_i\int_{\Om}Q_{\rho_i}\psi_i\Psi\vol.
\end{aligned}
\label{eq:low-moment-target}
\end{equation}
The second equality uses $\int_MU_\bgpar\Psi\vol=\bgpar$ and the zero
global Green moments. Radial integration bounds $\mathbf m_Z$ by
\[
O(\rho^n)+O(|\mathsf m|\rho^2)
             +O(\rho^{n-1}\|\psi\|)=O(\rho^{n-\delta_0}).
\]
Here $\int_{B_{\rho_i}(p_i)}|G_{i,\bgpar}\Psi_j|\vol\le C\rho_i^2$ and
$\int_{\rho}^{r_0}\|\psi\|\rho^{n-1}r^{1-n}r^{n-1}\dd r
\le C\rho^{n-1}\|\psi\|$. These source and moment bounds also
hold for parameter directions in
\eqref{eq:low-exterior-parameter-directions}: the lifts are linear,
the Green moments vanish identically in $\bgpar$, and the compact
coefficients have uniformly bounded first derivatives.

For each $Z$ in the stated ball, apply
Lemma~\ref{lem:exterior-inverse} to the source
$-L_\bgpar\sum_iQ_{\rho_i}\psi_i-N_\bgpar(H)$ and the moments
$\mathbf m_Z$. The returned function defines a map of $Z$ whose
value at zero is $O(\rho^{n-\delta_0})=o(\rho^{1-\nu})$. Its derivative with respect
to $Z$ is $O(\rho^{3/2})$. It is therefore a contraction on
$\|Z\|_{2,\alphaH;\nu;\Om}\le C\rho^{1-\nu}$.
Differentiate the fixed point equation for $Z$.
Lemma~\ref{lem:exterior-inverse} and
\eqref{eq:fixed-point-parameter-calculus} give
\[
 \|Z\|_{2,\alphaH;\nu;\Om}+|d_{\rm err}|
 +\|\dot Z\|_{2,\alphaH;\nu;\Om}+|\dot d_{\rm err}|
 \le C\rho^{n-\delta_0}.
\]
Here dots differentiate the solved family in the directions of
\eqref{eq:low-exterior-parameter-directions}. Let $\|Z\|_{\rm cap}$ be the maximum scaled
$C^{2,\alphaH}$ trace norm plus the scaled $C^{1,\alphaH}$ norm of
$\rho_i\partial_{\nu_g}Z$ over the spheres $S_{\rho_i}$. Multiplication by
$\rho^\nu$ and rescaling to unit directions give
\begin{equation}
\begin{aligned}
 \|Z\|_{\rm cap}+\|\dot Z\|_{\rm cap}
 &\le C\rho^{7/4-\delta_0}=o(\rho),\\
 \|D_\Lambda Z\|_{\rm cap}&\le C\rho^{5/4-\delta_0},\qquad
 \|D_\psi Z\|_{\rm cap}\le C\rho^{3/4}.
\end{aligned}
 \label{eq:low-exterior-boundary-estimates}
\end{equation}
In particular, $\|D_\bgpar Z\|_{\rm cap}=o(\rho)$, and the bounds for
$D_\Lambda Z$ and $D_\psi Z$ are $o(1)$. The estimates in
\eqref{eq:low-exterior-boundary-estimates} remain valid after
multiplication by $\Theta_i^{-1}$ and passage to $h_i$. Indeed,
$\Theta_i=1+O(\rho_i^2)$ and
$\rho_i\partial_{\nu_g}\Theta_i=O(\rho_i^2)$, and both are fixed
in the parameters.
Restoring $d=d_0+\sum_i\mathsf m_i\mathbf d_i^{\mathrm G}+d_{\rm err}$ gives the $C^1$ family
with exact moments.

\subsection{Joining the local and exterior solutions}
\label{subsec:low-boundary-matching}

We now determine the remaining local parameters and boundary traces.
On $S_{\rho_i}$ compare the local solution $V_i$ with
$\Theta_i^{-1}B_\bgpar$, so both are expressed relative to the same
local metric $h_i=\Theta_i^{4/(n-2)}g$.
Their values and increasing-radius normal derivatives must agree.
Equivalently, the boundary terms computed with the outward normals of
the two domains must cancel. We follow the matching construction of
\cite[Section 5]{santos}, using one contraction for the constant,
degree-one and higher-harmonic equations. The normal derivatives are
those of the local and exterior solutions with the chosen boundary values.

At the puncture $p_i$, suppress the index $i$ on the prescribed neck
size $\eps_i$ and the sphere radius $\rho_i$ in the following
local formulas. Represent $U_\bgpar$ in the fixed local metric by
\[
 U_\bgpar^{\mathrm{loc}}(\rho\theta)
  =c_i(\bgpar)+\rho e_i(\bgpar)\cdot\theta
        +\text{even }O(\rho^2)+\text{odd }O(\rho^3).
\]
The first $\bgpar$ derivatives of $U_\bgpar^{\mathrm{loc}}$ preserve the
even $O(\rho^2)$ and odd $O(\rho^3)$ remainders.
Also, $c_i(\bgpar)$ is bounded below by a positive constant. Choose a fixed
bounded box of Kelvin parameters containing
$\beta_i^0=e_i/[(n-2)c_i]$.

Let $Y$ collect the tuples $(c,\Lambda,B,\omega)$ over all
punctures. In the local formulas at $p_i$, the index is suppressed,
with $\Lambda=\mathsf m_i\rho_i^{2-n}$ and $B=\rho_i\beta$.
Write $\vartheta=(\vartheta_i)_i$ for the exterior high traces.
The pair $(c,\Lambda)$ controls the mean values and mean fluxes;
$(B,\omega)$ controls the degree-one data; $\vartheta$ contains the
remaining boundary functions.
For $B_{\bgpar,i}^{\rm loc}=\Theta_i^{-1}B_\bgpar$, equality of the
boundary values in harmonics of degree at least two gives
\[
 \phi_i=\Pi_{\ge2}(B_{\bgpar,i}^{\rm loc}-\uF_i)|_{S_{\rho_i}},
\]
where $\uF_i$ is the model $\uF_{\eps_i,c,\beta}$ at the $i$th puncture.
Here $\Pi_{\ge2}v_i|_{S_{\rho_i}}=0$.
Use the same local metric $h_i$ and increasing-radius direction to set
\[
 D_i=(V_i-B_{\bgpar,i}^{\rm loc})|_{S_{\rho_i}},\qquad
 N_i=\rho_i\partial_r(V_i-B_{\bgpar,i}^{\rm loc})|_{S_{\rho_i}}.
\]
Pull all traces back to the unit sphere. Collect the constant and
degree-one discrepancies in $\mathcal M$, and write the remaining
normal-derivative discrepancy as $\mathcal Z\vartheta+\mathcal R$:
\begin{equation}
\begin{aligned}
 \mathcal M(\vartheta,Y,\bgpar)&=(\Pi_0D_i,\Pi_0N_i,[\Pi_1D_i],[\Pi_1N_i])_i,\\
 \mathcal R(\vartheta,Y,\bgpar)&=(\Pi_{\ge2}N_i)_i-\mathcal Z\vartheta,\\
 (\mathcal Z\vartheta)_i
 &=\left.\rho_i\partial_r(P_{\rho_i}\vartheta_i)\right|_{r=\rho_i}
   -\left.\rho_i\partial_r(Q_{\rho_i}\vartheta_i)\right|_{r=\rho_i}.
\end{aligned}
\label{eq:low-matching-maps}
\end{equation}
Thus the high equation is $\mathcal Z\vartheta+\mathcal R=0$.
For every $n\ge3$, the interior and exterior harmonics of degree
$l\ge2$ have exponents $l$ and $-(l+n-2)$. Hence
\[
 \mathcal Z|_{l}=l+(l+n-2)=2l+n-2.
\]
The $C^{2,\alphaH}_{\ge2}\to C^{1,\alphaH}_{\ge2}$
isomorphism for $\mathcal Z$ is given in
\cite[Section 5.1, before Lemma 5.2]{santos}.
The rotation vector fields span the tangent spaces of $\Sph^{n-1}$.
They commute with the harmonic extension operators and $\mathcal Z^{-1}$.
Apply the $C^{2,\alphaH}$ estimate to four iterated rotational
derivatives to obtain
\begin{equation}
 \|\mathcal Z^{-1}f\|_{C^{j+1,\alphaH}}
 \le C_j\|f\|_{C^{j,\alphaH}},\qquad j=1,5.
 \label{eq:harmonic-matching-inverse}
\end{equation}
Applying four iterated rotational derivatives also upgrades the scaled
harmonic extension estimates in \cite[Propositions 2.4--2.5]{santos}
on $B_\rho\setminus B_{\rho/2}$ for the interior extension and
$B_{2\rho}\setminus B_\rho$ for the exterior extension.
The harmonic equation and interior estimates recover the remaining
radial derivatives, giving the $C^{6,\alphaH}$ extensions and
$C^{5,\alphaH}$ conormals.
We solve the high and low equations together. 
The radial and first-harmonic model contributions give the following reference system for the low modes.
\begin{equation}
 \begin{aligned}
 c+b_\eps(c)-\Lambda-c_i&=0,&
 (2-n)(b_\eps(c)-\Lambda)&=0,\\
 (n-2)cB-\omega-\rho e_i&=0,&
 (n-2)cB+(n-1)\omega-\rho e_i&=0.
 \end{aligned}\label{eq:low-leading}
\end{equation}
At the reference point
$Y_*=(c_i,b_\eps(c_i),\rho\beta_i^0,0)$, the diagonal Jacobian blocks
in $(c,\Lambda)$ and $(B_j,\omega_j)$ are
\[
 M_0(c)=\begin{pmatrix}1+b_\eps'(c)&-1\\
 (2-n)b_\eps'(c)&n-2\end{pmatrix},\qquad
 M_1(c)=\begin{pmatrix}(n-2)c&-1\\(n-2)c&n-1\end{pmatrix}.
\]
We have $\det M_0(c)=n-2$ and $\det M_1(c)=n(n-2)c$.
Since $c$ is bounded below and $b_\eps'(c)=O(\rho^{3/2})$,
the norms of $M_0(c)^{-1}$ and $M_1(c)^{-1}$ are uniformly bounded. The remaining amplitude coupling is
lower triangular and $O(|B|)=O(\rho)$.
For different punctures, the smooth Green contribution is
$O(|\mathsf m_j|)$ in the mean and $O(\rho_i|\mathsf m_j|)$
in the first harmonic. A $\Lambda_j$ derivative gains
$\rho_j^{n-2}$, so the off-diagonal blocks are $o(1)$.
Let $D_0$ be the Jacobian, with respect to $Y$, of the four
leading expressions in \eqref{eq:low-leading}, assembled over all
punctures and evaluated at $Y_*(\bgpar)$. It is uniformly invertible. The contributions from the other punctures
form an $o(1)$ perturbation.

For the maps in \eqref{eq:low-matching-maps}, on
$\|\vartheta\|_{C^{2,\alphaH}}+|Y-Y_*|\le c_0\rho$, the scaled Dirichlet and normal
trace estimates on $S_{\rho_i}$ give
\[
 \begin{aligned}
 \|D_{(\vartheta,Y)}(\mathcal Z^{-1}\mathcal R)\|&=o(1),\\
 \|D_0^{-1}(D_Y\mathcal M-D_0)\|
       +\|D_0^{-1}D_\vartheta\mathcal M\|&=o(1).
 \end{aligned}
\]
The explicit contributions to the high traces come from the Taylor
remainder of $U_\bgpar^{\mathrm{loc}}$, the quadratic Kelvin term and the
pole remainder. Each is $o(\rho)$, with $Y$ derivatives of size $o(1)$.
Proposition~\ref{prop:local} and
\eqref{eq:low-exterior-boundary-estimates} give values $o(\rho)$
and $Y$ derivatives $o(1)$ for the boundary traces of
$v_i$ and $\Theta_i^{-1}Z$.
The harmonic extensions $P_{\rho_i}\phi_i$ and
$Q_{\rho_i}\vartheta_i$ have zero constant and degree-one projections. At $(\vartheta,Y)=(0,Y_*)$, both
$\mathcal Z^{-1}\mathcal R$ and $D_0^{-1}\mathcal M$, including
their total $\bgpar$ derivatives along $\bgpar\mapsto(0,Y_*(\bgpar),\bgpar)$, are $o(\rho)$.

Fix $c_0>0$ sufficiently small. Then $c$ stays positive, and
$B/\rho=\beta_i^0+O(c_0)$ stays in the chosen box of Kelvin
parameters. The bounds $\Lambda,\omega,\vartheta=O(\rho)$ give
\[
 |\mathsf m|=O(\rho^{n-1})=o(\rho^{n-3/2}),\qquad
 \|\psi\|=O(\rho)=o(\rho^{1-\delta_0}).
\]
The boundary functions $\phi_i$ determined above satisfy
$\|\phi_i\|=O(\rho_i)=o(\rho_i^{1-\delta_0})$, so Proposition~\ref{prop:local} applies.
The choice of $c_0$ may depend on the fixed ratios $t_1,\ldots,t_m$. The map
\[
 (\vartheta,Y)\longmapsto
 \bigl(-\mathcal Z^{-1}\mathcal R(\vartheta,Y,\bgpar),\,
       Y-D_0^{-1}\mathcal M(\vartheta,Y,\bgpar)\bigr)
\]
is a contraction of this ball into itself. Its fixed point satisfies
\[
 \|\vartheta\|_{C^{2,\alphaH}}+|Y-Y_*|=o(\rho).
\]
To differentiate the matched parameters, write
$V=(\vartheta,Y)$, $V_*=(0,Y_*)$, and
\[
 \mathcal H(\bgpar,V)=
 \bigl(\vartheta+\mathcal Z^{-1}\mathcal R(\vartheta,Y,\bgpar),
       \mathcal M(\vartheta,Y,\bgpar)\bigr).
\]
Let $\mathcal M_0$ consist of the four leading expressions in
\eqref{eq:low-leading}, for all punctures, and put
$\mathcal H_0(\bgpar,V)=(\vartheta,\mathcal M_0(Y,\bgpar))$.
Set $\mathcal R_{\rm match}=\mathcal H-\mathcal H_0$.
The local and exterior estimates give
\[
 D_\bgpar\mathcal R_{\rm match}=o(\rho),\qquad
 D_c\mathcal R_{\rm match}=o(\rho),\qquad
 D_{(\Lambda,B,\omega,\vartheta)}\mathcal R_{\rm match}=o(1),
\]
with the other variables fixed. The explicit reference point satisfies
\[
 D_\bgpar Y_*=(O(1),O(\rho^{3/2}),O(\rho),0).
\]
Since $\mathcal H_0(\bgpar,V_*(\bgpar))=0$, its total $\bgpar$ derivative vanishes.
Changing $V_*$ to $V$ in this identity changes the explicit leading
coefficients by $O(\|\vartheta\|+|Y-Y_*|)=o(\rho)$.
The displayed bounds for $D_\bgpar\mathcal R_{\rm match}$ and
$D_V\mathcal R_{\rm match}$ therefore imply
\[
 D_\bgpar\mathcal H(\bgpar,V)+D_V\mathcal H(\bgpar,V)D_\bgpar V_*=o(\rho).
\]
Differentiating $\mathcal H(\bgpar,V(\bgpar))=0$ gives
\[
 D_V\mathcal H(\bgpar,V)[D_\bgpar V-D_\bgpar V_*]
 =-\{D_\bgpar\mathcal H(\bgpar,V)+D_V\mathcal H(\bgpar,V)D_\bgpar V_*\}.
\]
The uniform inverse bound for $D_V\mathcal H$ yields
$D_\bgpar(V-V_*)=o(\rho)$. Thus $D_\bgpar c=O(1)$, while
$D_\bgpar\Lambda$, $D_\bgpar B$, $D_\bgpar\omega$ and $D_\bgpar\vartheta$ are $O(\rho)$.

Let $U_\eps(\bgpar)$ equal $\Theta_iV_i$ on each punctured ball and
$B_\bgpar$ on $\Om$. The matched values and normal derivatives give a
$C^1$ function across each sphere. It satisfies
\begin{equation}
 \begin{gathered}
 \Ff(U_\eps(\bgpar))=\sum_{j=1}^{N_E}d_{\eps,j}(\bgpar)\Psi_j\quad\hbox{on }\Om,\qquad
 \Ff(U_\eps(\bgpar))=0\quad\hbox{on the ends},\\
 \int_{\Om}U_\eps(\bgpar)\Psi_j\vol=\eta_j,\qquad
 d_\eps\to d_0\quad\hbox{in }C^1.
 \label{eq:projected-global}
 \end{gathered}
\end{equation}
Along the matched family, $D_\bgpar\Lambda,D_\bgpar\psi=O(\rho)$ are bounded
in the directions of \eqref{eq:low-exterior-parameter-directions}.
The exterior estimates give
$\|d_{\rm err}\|_{C^1_\bgpar}\le C\rho^{n-\delta_0}$ and
$\|\mathsf m\|_{C^1_\bgpar}=O(\rho^{n-1})$. Hence
\[
 \|d_\eps-d_0\|_{C^1_\bgpar}
 =\left\|\sum_i\mathsf m_i\mathbf d_i^{\mathrm G}+d_{\rm err}\right\|_{C^1_\bgpar}
 \le C(\rho^{n-\delta_0}+\rho^{n-1})\longrightarrow0.
\]
For fixed $\eps$, the maps $\bgpar\mapsto\Theta_iV_i$ are $C^1$
with values in
$C^{2,\alphaH}_{-(n-2)/2}(\overline{B_{\rho_i}(p_i)}
\setminus\{p_i\})$, and $\bgpar\mapsto B_\bgpar$ is $C^1$
with values in $C^{2,\alphaH}_{\nu}(\overline\Om)$.
Only the background coefficients $\bgpar$ remain to be chosen.

\subsection{The action of the joined low-dimensional family}
\label{subsec:low-action}

For each fixed matching radius (and hence fixed neck), use the real phase
$\tau=\log R$, where $R^{-(n-2)/2}=2c/\eps$ and $c>0$.
Use this phase in the subtraction defining $S_{\eps,\rho}$ in
\eqref{eq:end-action}.
\begin{lemma}\label{lem:action}
Let $u_h$ be a local factor from Proposition~\ref{prop:local},
and put $u_g=\Theta u_h$.
The limit \eqref{eq:end-action} exists in $C^1$ with respect to
$(c,\beta,\phi)$.
Use \eqref{eq:action-metric-conversion} to express the action in $g$. Then
\begin{equation}
 DS_{\eps,\rho}(u_g)[\dot u_g]
    =\int_{S_\rho}\partial_{\nu_g}u_g\,\dot u_g\dd S_g,
 \label{eq:action-differential}
\end{equation}
where $\nu_g$ is the $g$-unit normal pointing from the end into the exterior.
\end{lemma}
\begin{proof}
Keep \(\varepsilon\) fixed while taking \(T\to\infty\) in the definition of the action.
 For $v(s,\theta)=r^{(n-2)/2}u_h(r,\theta)$, the Kelvin formula and the
weighted estimates in Proposition~\ref{prop:local} give
\begin{equation}
 \left\|D_{(c,\beta,\phi)}^j
      \bigl(v-v_\eps(\,\cdot+\tau)\bigr)
       \right\|_{C^1([s,s+1]\times\Sph^{n-1})}
 \le C(\eps)(e^{-s}+e^{-(n-2+2\mu)s/2}),\qquad j=0,1.
 \label{eq:low-terminal-mixed-estimate}
\end{equation}
For $j=1$, use the Fr\'echet operator norm from
$\R\times\R^n\times C^{2,\alphaH}_{\ge2}(\Sph^{n-1})$ to the
displayed spatial space. The functions are differentiated at fixed
physical coordinates, with $\eps$, $h$ and $S_\rho$ fixed.
The cylindrical metric $\gamma=r^{-2}h$ satisfies
$\gamma-(\dd s^2+g_{\Sph^{n-1}})=O(e^{-2s})$,
with the spatial derivatives used in the action density. Hence
\[
\begin{aligned}
 &\left\|D_{(c,\beta,\phi)}^j
  \{\mathcal L(s,c,\beta,\phi)
       -|\Sph^{n-1}|[(v_\eps')^2-h_\eps](s+\tau)\}\right\|\\
 &\hspace{2em}\le C(\eps)(e^{-s}+e^{-(n-2+2\mu)s/2}+e^{-2s})
       \in L^1(s_0,\infty),\qquad j=0,1.
\end{aligned}
\]
Estimate~\eqref{eq:low-terminal-mixed-estimate} proves the $C^1$
convergence in \eqref{eq:terminal-quadratic-error} and the convergence
in \eqref{eq:terminal-conormal-error}, uniformly in unit
$(c,\beta,\phi)$ directions. Since $\Ff_h(u_h)=0$ on the end,
\eqref{eq:renormalized-first-variation} gives the outer pairing
$\int_{S_\rho}\partial_{\nu_h}u_h\,\dot u_h\,\dd S_h$.
The inner boundary term in
\eqref{eq:smooth-conformal-boundary} satisfies
\[
 \left\|\int_{S_r}u_g^2\partial_{\nu_g}\log\Theta\dd S_g
       \right\|_{C^1(c,\beta,\phi)}
 \le C(\eps) r^{2-n}r\,r^{n-1}
 =C(\eps) r^2\longrightarrow0.
\]
The boundary term on $S_\rho$ in \eqref{eq:smooth-conformal-boundary}
converts this pairing into the $g$-conormal pairing
in \eqref{eq:action-differential}.
\end{proof}
For the actual matched family \eqref{eq:projected-global} set
\begin{equation}
 j_\eps(\bgpar)=J_g(U_\eps(\bgpar);\Om)+
             \sum_{i=1}^mS_{\eps_i,\rho_i}(U_\eps(\bgpar)|_{\text{end }i}).
 \label{eq:global-action}
\end{equation}
In the first variation of \eqref{eq:global-action}, the local
and exterior boundary integrals cancel on each $S_{\rho_i}$. The bulk residual vanishes on each end. The exact exterior moments therefore give
\begin{equation}
 \partial_{\eta_l}j_\eps(\bgpar)
 =\sum_{j=1}^{N_E}d_{\eps,j}(\bgpar)
       \int_{\Om}\Psi_j\partial_{\eta_l}U_\eps(\bgpar)\vol
 =d_{\eps,l}(\bgpar).\label{eq:exact-gradient}
\end{equation}
Since $Dj_\eps=d_\eps\in C^1$, we have $j_\eps\in C^2$.
Equation~\eqref{eq:projected-global} gives $Dj_\eps\to Dj_0$ in $C^1$.

There are constants $C_\eps$, independent of $\bgpar$, such that
\begin{equation}
 \sup_\bgpar|j_\eps(\bgpar)-C_\eps-j_0(\bgpar)|\to0.\label{eq:action-small}
\end{equation}
The difference $B_\bgpar-U_\bgpar$ tends to zero uniformly and in $H^1(\Om)$.
Indeed, radial integration and the harmonic extension estimate give
\[
 \begin{aligned}
 \left\|\sum_i\mathsf m_i\nabla G_{i,\bgpar}\right\|_2^2&\le C|\mathsf m|^2\rho^{2-n}
                              \le C\rho^{n-1},\\
 \left\|\nabla\sum_iQ_{\rho_i}\psi_i\right\|_2^2&\le C\rho^{n-2}\|\psi\|_{H^{1/2}}^2
                              \le C\rho^{n-2\delta_0},\\
 \|\nabla Z\|_2^2&\le o(\rho^{2-2\nu})\rho^{2\nu+n-2}
                              =o(\rho^n).
 \end{aligned}
\]
The last estimate uses $2\nu+n-2=3/2-n<0$.
The potential and nonlinear energy differences also tend to zero,
and the energy of $U_\bgpar$ on the removed balls is $O(\rho^n)$.

For each end, choose a reference tuple $\mathbf p_*$ in the connected
set of admissible parameters $\mathbf p=(c,\beta,\phi)$: $c$ lies in
a fixed compact positive interval, $\beta$ in a bounded convex box,
and $\phi$ in the high harmonic trace ball of Proposition~\ref{prop:local}.
Choose these sets to contain the matched parameters from every
component of $\K$. Join $\mathbf p_*$ to $\mathbf p$ by the straight
path $\mathbf p(t)$ in this set. Write $u_{g,\mathbf p(t)}$
for the corresponding factor in the global metric.
The boundary estimates give
$|\partial_{\nu_g}u_{g,\mathbf p(t)}|\le C/\rho$ and
$|\frac{d}{dt}u_{g,\mathbf p(t)}|\le C$ on $S_\rho$.
Equation~\eqref{eq:action-differential} yields
\[
\begin{aligned}
 S_{\eps,\rho}(u_{g,\mathbf p})-S_{\eps,\rho}(u_{g,\mathbf p_*})
 &=\int_0^1\int_{S_\rho}
       \partial_{\nu_g}u_{g,\mathbf p(t)}
       \frac{d}{dt}u_{g,\mathbf p(t)}\,\dd S_g\dd t,\\
 \bigl|S_{\eps,\rho}(u_{g,\mathbf p})-S_{\eps,\rho}(u_{g,\mathbf p_*})\bigr|
 &\le C\rho^{-1}\rho^{n-1}
       =C\rho^{n-2}\longrightarrow0.
\end{aligned}
\]
The real phase lift makes each action single valued on its chosen
parameter set. The sum of the reference actions therefore defines
one $\bgpar$-independent $C_\eps$ satisfying \eqref{eq:action-small}.

\begin{proposition}\label{prop:low-global-action}
Fix $t_i>0$ and take $\rho$ sufficiently small. Reindex the construction
by $\eps=\rho^{(2n-1)/4}$. The positive family
$U_\eps(\bgpar)$ is $C^1$ on the fixed compact parameter set chosen
at the start of this section.
Its limiting Fowler necks are $t_i\eps$. It satisfies
\eqref{eq:projected-global}, with matching values and physical
conormals. The action $j_\eps$ is $C^2$, with
$Dj_\eps=d_\eps$ and $j_\eps-C_\eps\to j_0$ in $C_\bgpar^2$ as $\eps\to0$
on the fixed compact parameter set. The constants $C_\eps$ are independent of $\bgpar$.
\end{proposition}
\begin{proof}
Combine \eqref{eq:projected-global}, \eqref{eq:exact-gradient} and
\eqref{eq:action-small}.
\end{proof}

\section{Metric expansions and model estimates}\label{sec:high-local}
\subsection{Metric jets, bubble parameters, and the three local regions}
\label{sec:jets}\label{subsec:bridge-geometry}

For $n\ge8$, we construct the local solution by joining an inner
Fowler end to a corrected smooth bubble, and then joining the outer
decay of the bubble to the regular background. We first describe the
metric and the bubble. The estimates from \cite{part1} supply the correction
of the bubble; the rest of this section concerns the equations on the
inner Fowler end.

Let $X$ be physical coordinates centered at the prescribed puncture.
As in Section~\ref{sec:low-construction}, fix a conformal representative
$h=\Theta^{4/(n-2)}g$ and its geodesic normal coordinates, with
$\Theta(0)=1$, $d\Theta(0)=0$, and $\det h\equiv1$ \cite{gunther}.
These choices are independent of the small scales and of the compact
background parameter $\bgpar$. The matrix logarithm of $h$ satisfies
\[
 h=\exp H,\qquad H(X)X=0,\qquad \operatorname{tr}H=0.
\]
Write $H_q$ for the homogeneous term of degree $q$ in the Taylor
expansion of $H$ at the puncture.

\begin{lemma}\label{lem:finite-jets}
Let $q\ge3$. If the Weyl derivatives vanish through order $q-3$,
then $h_{ij}=\delta_{ij}+O(|X|^q)$. If the first nonzero Weyl
derivative has order $q-2$, then $H_q$ has nonzero linearized Weyl
curvature. If all Weyl jets vanish, any fixed finite $q$ is available.
\end{lemma}
\begin{proof}
The exact normalization gives $\log\sqrt{\det h}=0$.
Suppose inductively that $\nabla^r\operatorname{Rm}(0)=0$ for $r<j$,
starting with $j=0$. The radial Jacobi equation gives
\[
 \log\sqrt{\det h(r\theta)}
 =-\frac{r^{j+2}}{(j+2)(j+3)j!}
 (\nabla_\theta^j\operatorname{Ric})(\theta,\theta)+O(r^{j+3}).
\]
Thus $\operatorname{Sym}\nabla^j\operatorname{Ric}(0)=0$, with
averaged symmetrization. Assume that $\nabla^rW(0)=0$ for
$0\le r\le j$. Covariant derivatives here use $h$. Set
\[
 P=\frac1{n-2}\left(\operatorname{Ric}
       -\frac{\operatorname{Scal}_h}{2(n-1)}h\right),\quad
 C_{ijk}=\nabla_kP_{ij}-\nabla_jP_{ik},\quad T=\nabla^jP(0).
\]
For $j\ge1$, the divergence identity for $W$ gives $\nabla^{j-1}C(0)=0$. The induction
hypothesis removes the curvature commutators. Hence $T$ is fully
symmetric. For $j=0$, use the symmetry of $P$.
Since $\operatorname{Ric}=(n-2)P+(\operatorname{tr}P)h$,
\[
 (n-2)T+\operatorname{Sym}(h\otimes\operatorname{tr}T)=0,
 \qquad (n-2)|T|^2+|\operatorname{tr}T|^2=0.
\]
Hence $T=0$ and $\nabla^j\operatorname{Rm}(0)=0$.
Induction through $j=q-3$ eliminates the metric terms of degrees
$2,\ldots,q-1$. The degree-$(q-2)$ Weyl term is then the linearized
Weyl curvature of $H_q$, and is nonzero when
$\nabla^{q-2}W(0)\ne0$.
\end{proof}

At each puncture, choose the least degree $q$ with nonzero linearized
Weyl curvature if $2\le q\le\lfloor(n-2)/4\rfloor$. At every other
puncture, choose $q=\lfloor(n-2)/4\rfloor+1$. For $q=2$, normal
coordinates give $h_{ij}=\delta_{ij}+O(|X|^2)$, and $W_h(0)$ is the
linearized Weyl curvature of $H_2$. For $q\ge3$, the preceding lemma
gives $h_{ij}=\delta_{ij}+O(|X|^q)$. In all cases,
\[
 H(X)=H_q(X)+O(|X|^{q+1}),\qquad 2\le q<\frac{n-2}{2}.
\]
The first choice has $n\ge4q+2$ and a nonzero $H_q$; the second has
$n<4q+2$ and allows $H_q=0$.

\medskip\noindent\textit{The bubble and its two parameters.}
The standard smooth bubble of width $\varrho$ and center $\xi$ is
\begin{equation}
 U_{\varrho,\xi}(X)
 =\left(\frac{2\varrho}{\varrho^2+|X-\xi|^2}\right)^{(n-2)/2},
 \qquad
 \varrho=\frac b{1+|\zeta|^2},\quad \xi=\frac{b\zeta}{1+|\zeta|^2}.
 \label{eq:finite-bubble}
\end{equation}
The scale $b$ is determined by the value at the prescribed puncture:
$U_{\varrho,\xi}(0)=(2/b)^{(n-2)/2}$. To explain the other parameter,
rescale by $y=X/b$ and write
\begin{equation}
 \begin{aligned}
 V_\zeta(y)
 &=b^{(n-2)/2}U_{\varrho,\xi}(by)
   =\left(\frac{2}{1-2\zeta\cdot y+(1+|\zeta|^2)|y|^2}\right)^{(n-2)/2}\\
 &=|y|^{2-n}U_{1,\zeta}\!\left(\frac y{|y|^2}\right).
 \end{aligned}
 \label{eq:normalized-bubble-parameters}
\end{equation}
Thus $\zeta$ is the center of the unit-width bubble before Kelvin
inversion. Its variation leaves the prescribed puncture at $y=0$ and
keeps the value there fixed. The Euclidean width and center of
$V_\zeta$ are $(1+|\zeta|^2)^{-1}$ and
$\zeta/(1+|\zeta|^2)$, respectively. At its two ends,
\begin{equation}
 \begin{gathered}
 V_\zeta(0)=2^{(n-2)/2},\qquad
 V_\zeta(y)=C_\infty(\zeta)|y|^{2-n}+O(|y|^{1-n}),\\
 C_\infty(\zeta)=2^{(n-2)/2}(1+|\zeta|^2)^{-(n-2)/2}.
 \end{gathered}
 \label{eq:outgoing}
\end{equation}
All local estimates are uniform when $\zeta$ belongs to a fixed compact
set. The size of this set will be chosen in the parameter-selection
argument.

\medskip\noindent\textit{The regions used to join the solutions.}
For a prescribed small neck $\eps$, choose
\begin{equation}
 \delta=\eps^{2/n},\qquad
 b^{\sigma_b}/\eps\in[c,C],\qquad
 \rho=\eps^{(n-2)/[\sigma_b(2n-1)]},\qquad
 \sigma_b=\min\{q,(n-2)/4\}.
 \label{eq:bridge-scales}
\end{equation}
The three regions are
\[
\begin{array}{c|l}
 0<|X|<b\delta & \text{the inner Fowler end},\\
 b\delta<|X|<b/\delta & \text{the middle region containing the corrected bubble},\\
 b/\delta<|X|<\rho & \text{the outer annulus joining to the regular background}.
\end{array}
\]
The two internal interfaces are fixed at $|y|=\delta,\delta^{-1}$
in the rescaled coordinates. Their physical radii are
$r_-=b\delta$ and $r_+=b/\delta$. The sphere $|X|=\rho$ is the
boundary of the entire local construction. For small $\eps$,
$b\delta\ll b/\delta\ll\rho$.

The coefficient of $|X|^{2-n}$ in the outgoing bubble is
$B_0=b^{(n-2)/2}C_\infty(\zeta)$. Substitution in
\eqref{eq:bridge-scales} gives
\begin{equation}
 \begin{gathered}
 b\asymp\rho^{(2n-1)/(n-2)},\qquad
 B_0\asymp\rho^{n-1/2},\qquad
 B_0\rho^{2-n}\asymp\rho^{3/2},\\
 B_0=O(\eps^2),\qquad \delta^n=\eps^2,\qquad
 \rho^n/\eps^2\le\rho^{1/2}\longrightarrow0.
 \end{gathered}
 \label{eq:bridge-scale-ledger}
\end{equation}
Parameter derivatives keep the physical chart, the reference metrics
$g$ and $h$, the neck $\eps$, and the outer radius $\rho$ fixed. Consequently $\delta$ is
fixed, whereas the two physical interfaces move with $b$. With
$\ell=\log b$, differentiation of a rescaled function gives
\[
 \partial_\ell\big|_X
 \bigl[b^{-(n-2)/2}W_b(X/b)\bigr]
 =b^{-(n-2)/2}
   \bigl(\partial_\ell W_b-\tfrac{n-2}{2}W_b-y\cdot\nabla_yW_b\bigr).
\]
The remaining parameters are $\zeta$ and the regular boundary data on
$|X|=\rho$. Boundary integrals are transformed with both the unit
normal derivative and the induced surface measure.

\subsection{The correction from \texorpdfstring{\cite{part1}}{Part I} and the middle-region operator}
\label{subsec:corrected-bubble}

\medskip\noindent\textit{The result from \cite{part1}.}
We first state the correction in the ordinary width and center
parameters $(\varrho,\xi)$ of a bubble. In this imported statement,
parameter derivatives keep $y$ and the tensor fixed. The metric
\[
 \gamma_{\varrho,\xi}=U_{\varrho,\xi}^{4/(n-2)}|dy|^2
\]
extends to a round metric of sectional curvature one on
$\R^n\cup\{\infty\}$.
Writing the corrected Euclidean factor as
$U_{\varrho,\xi}(1+s\psi)$ makes its linearized equation
$(-\Delta_{\gamma_{\varrho,\xi}}-n)\psi$ on this sphere.

For $R\ge2$, choose a smooth radial cutoff $\chi_R$, equal to one on $|y|\le R$ and
zero on $|y|\ge2R$, obtained by rescaling one fixed function, and set
$H_{q,R}=\chi_RH_q$. For $T=H_q$ or $H_{q,R}$, \cite{part1} supplies the
unique correction $\psi_{\varrho,\xi}[T]$ satisfying
\begin{equation}
 \begin{gathered}
 (-\Delta_{\gamma_{\varrho,\xi}}-n)\psi_{\varrho,\xi}[T]
 =-c(n)D\operatorname{Scal}_{\gamma_{\varrho,\xi}}
              [U_{\varrho,\xi}^{4/(n-2)}T],\\
 \psi_{\varrho,\xi}[T]\perp
       \ker(-\Delta_{\gamma_{\varrho,\xi}}-n)
       \quad\text{in }L^2(\gamma_{\varrho,\xi}).
 \end{gathered}
 \label{eq:part-i-first-response}
\end{equation}
Here $D\operatorname{Scal}_{\gamma}[k]$ is the first variation of
scalar curvature in the metric direction $k$. The kernel consists of
the first spherical harmonics. The source is orthogonal to this
kernel, and the inverse on its orthogonal complement gives
\[
 \|\psi_{\varrho,\xi}[T]\|_{H^1(\gamma_{\varrho,\xi})}
 \le C\|U_{\varrho,\xi}^{4/(n-2)}T\|_{H^1(\gamma_{\varrho,\xi})}.
\]
These statements are \cite[equations (6.12)--(6.14)]{part1}.
For $2\le q<(n-2)/2$, the uncut tensor has finite $H^1$ norm and
\cite[Lemma 6.2 and Proposition 6.4]{part1} give
\[
 \left|\partial_{(\varrho,\xi)}^\beta\partial_y^\alpha
     \bigl(U_{\varrho,\xi}\psi_{\varrho,\xi}[H_q]\bigr)\right|
 \le C(1+|y|)^{q+2-n-|\alpha|},
 \qquad |\alpha|\le2,\quad |\beta|\le1.
\]
The cutoff differences
\[
 U_{\varrho,\xi}^{4/(n-2)}(H_q-H_{q,R}),\qquad
 \psi_{\varrho,\xi}[H_q]-\psi_{\varrho,\xi}[H_{q,R}]
\]
have $H^1(\gamma_{\varrho,\xi})$ norm at most
$CR^{q-(n-2)/2}$, including one parameter derivative. For the
correction difference, this follows from the tensor estimate in
\cite[equation (6.7)]{part1} by applying and differentiating
\eqref{eq:part-i-first-response}; the undifferentiated estimate is
\cite[Lemma 6.8]{part1}. Constants are uniform on compact sets of
positive widths and finite centers. We apply this statement to
$V_\zeta$, whose width and center belong to such a set when
$\zeta$ varies in a fixed compact set.

The same construction gives the quadratic form
$Q^{\rm I}_{\varrho,\xi}(T,T)$ of
\cite[Proposition 6.6 and equation (6.34)]{part1}. It is continuous
in the $H^1$ norm of the compactified tensor
$U_{\varrho,\xi}^{4/(n-2)}T$. Thus it also applies to trace-free
polynomial tensors of degree at most $q$, even when they do not
satisfy the radial identity $T(y)y=0$. For the cutoff path, the
correction and metric use the same tensor, and
\begin{equation}
 \begin{aligned}
 &J_{\exp(sH_{q,R})}
   \bigl(U_{\varrho,\xi}(1+s\psi_{\varrho,\xi}[H_{q,R}]);\R^n\bigr)\\
 &\qquad=J_{g_{\rm E}}(U_{\varrho,\xi};\R^n)
       -s^2Q^{\rm I}_{\varrho,\xi}(H_{q,R},H_{q,R})+O_R(s^3).
 \end{aligned}
 \label{eq:part-i-cutoff-action}
\end{equation}
This is the fixed-cutoff expansion of \cite[equation (6.40)]{part1}.
The bubble in \cite{part1} is $2^{-(n-2)/2}U_{\varrho,\xi}$, and the
action used here is $-2^{n-2}$ times the action in \cite{part1} after this
rescaling. Its coefficient $2^{2-n}s^2Q^{\rm I}$ therefore becomes
$-s^2Q^{\rm I}$. The quotient of the scalar correction by its
bubble is unchanged by the rescaling; this quotient is
$\psi_{\varrho,\xi}$.
In \eqref{eq:part-i-cutoff-action}, $R$ is fixed as $s\to0$,
and the remainder constant may depend on $R$.
For parameter derivatives, regard these round metrics and tensors
as fields on the same sphere in fixed coordinates. The two-derivative
tensor estimate in \cite[equation (6.7)]{part1}, the differentiated
correction equation, and the smoothly varying metric coefficients
and measure give the same cutoff bound for the first derivative
of the quadratic form. Thus
\begin{equation}
 \left\|Q^{\rm I}_{\varrho,\xi}(H_{q,R},H_{q,R})
           -Q^{\rm I}_{\varrho,\xi}(H_q,H_q)\right\|_{C^1_{\varrho,\xi}}
 \le CR^{q-(n-2)/2}.
 \label{eq:tt-cutoff-parameter-bound}
\end{equation}
All these parameter bounds are on fixed compact sets. Finally,
\cite[Proposition 6.7]{part1} gives
$Q^{\rm I}_{\varrho,\xi}(H_q,H_q)>0$ when $H_q$ has nonzero
linearized Weyl curvature, uniformly on such parameter sets.

\medskip\noindent\textit{One identification with the fixed sphere used for gluing.}
The preceding result uses the round metric of each bubble. For the
boundary equations we use one fixed sphere, so that its measure and
the spheres $|y|=\delta,\delta^{-1}$ are fixed when $\zeta$ varies.
Set
\[
 V_0(y)=\left(\frac2{1+|y|^2}\right)^{(n-2)/2},\qquad
 g_{\rm rd}=V_0^{4/(n-2)}|dy|^2,\qquad u_\zeta=V_\zeta/V_0.
\]
A multiplicative correction $\psi$ of $V_\zeta$ is represented by
the additive correction $u_\zeta\psi$ on this fixed sphere:
\[
 V_\zeta(1+\psi)=V_0(u_\zeta+u_\zeta\psi).
\]
The later gluing equations all use $g_{\rm rd}$ and $\dd V_{\rm rd}$.

Let $h_b=b^{-2}\mathcal D_b^*h$, where $\mathcal D_b(y)=by$.
Its matrix logarithm on the local chart is $H(by)$. On
$\delta\le|y|\le\delta^{-1}$ put
$\gamma_b=V_0^{4/(n-2)}h_b$. Extend $\gamma_b$ to the sphere by
multiplying $H(by)$ by radial cutoffs before taking its exponential.
The cutoffs equal one on this annulus and zero near $y=0$ and
$y=\infty$; their transition annuli have round radius comparable
to $\delta$, and their rescaled profiles are fixed. The notation
$h_b$ still denotes the unmodified metric on the local chart. Put
\[
 K_b=\gamma_b-g_{\rm rd},\qquad
 K_q=V_0^{4/(n-2)}H_q,\qquad K_{q,R}=\chi_RK_q.
\]
The action gradient for a conformal factor $z$ on this fixed sphere is
\begin{equation}
 \Ff_b(z)=\frac{\dd V_{\gamma_b}}{\dd V_{\rm rd}}
   \left(-\Delta_{\gamma_b}z+c(n)\operatorname{Scal}_{\gamma_b}z
                        -\frac{n(n-2)}4z^p\right).
 \label{eq:compactified-action-source}
\end{equation}
In the estimates below, $d$ is round distance to the nearer of
$y=0$ and $y=\infty$. The notation $A_d$ denotes either annulus
with radii $d,2d$ about these points, and a supremum over $A_d$
includes both. The scaled annular norm is the ordinary H\"older norm
after rescaling the annulus to unit size; in the fixed region away
from these points it is the ordinary norm.

\begin{lemma}\label{lem:compactified-metric-bounds}
For $2\le q<(n-2)/2$, the compactified tensors satisfy
\[
 \begin{gathered}
 K_q\in H^1(\Sph^n),\qquad
 |D^jK_q|\le Cd^{-q-j}\quad(0\le j\le6),\qquad
 \|K_q-K_{q,R}\|_{H^1}\le CR^{q-(n-2)/2}.
 \end{gathered}
\]
\begin{equation}
 \begin{aligned}
 &\|K_b\|_{H^1}+\|\Ff_b(u_\zeta)\|_{H^{-1}}\\
 &\qquad+\sup_d d^{(n+2)/2}
           \|\Ff_b(u_\zeta)\|_{C^{4,\alphaH}_{\rm sc}(A_d)}
 \le Cb^q.
 \end{aligned}
 \label{eq:weak-metric-size}
\end{equation}
All bounds include one derivative in $(\ell,\zeta)$ and are uniform
on the fixed parameter sets. If $n<4q+2$, then $b^q=o(\eps)$.
\end{lemma}
\begin{proof}
The $H^1$ and cutoff estimates are the tensor estimates in \cite{part1} just
recalled, applied at the centered unit bubble. Near $y=\infty$,
Kelvin coordinates give $|D^jK_q|\le Cd^{-q-j}$; near $y=0$ the
tensor is smooth. Taylor expansion of the actual metric, including
the cutoff annuli, gives
\[
 \|K_b\|_{C^{6,\alphaH}_{\rm sc}(A_d)}\le Cb^qd^{-q},\qquad
 \|\Ff_b(u_\zeta)\|_{C^{4,\alphaH}_{\rm sc}(A_d)}
                                      \le Cb^qd^{-q-2}.
\]
Each cutoff derivative costs $d^{-1}$, whereas nonlinear metric
terms gain a factor $O((b/\delta)^q)=o(1)$. The first metric
variation of the action gradient is bounded from $H^1$ tensors to
$H^{-1}$: integrate its second metric derivatives once against a
test function. The nonlinear terms obey the same estimate using
\[
 \|\phi/d\|_{L^2}\le C\|\phi\|_{H^1},\qquad
 d|\nabla K_b|\le C(b/\delta)^q.
\]
Integration of the annular bounds proves the asserted estimate
because $q<(n-2)/2$. A derivative in $\ell$ replaces $H(by)$ by
$by\cdot\nabla H(by)$; derivatives in $\zeta$ act on the smooth
bubble coefficients. They satisfy the same bounds, with $\delta$
fixed. Finally, if $n<4q+2$, then
$b^q/\eps\asymp b^{q-(n-2)/4}\to0$.
\end{proof}

\medskip\noindent\textit{The kernel and the choice of correction.}
On the fixed sphere, the bubble linearization is
\[
 L_\zeta=-\Delta_{g_{\rm rd}}+\frac{n(n-2)}4
                         -\frac{n(n+2)}4u_\zeta^{p-1}.
\]
It is the operator $-\Delta-n$ on the bubble's round metric,
expressed in the fixed conformal representation. Its kernel is
spanned by the derivatives of the bubble with respect to its scale
and Kelvin center:
\[
 Z_0=-\frac{(n-2)V_\zeta+2y\cdot\nabla V_\zeta}{2V_0},\qquad
 Z_j=\frac{\partial_{\zeta_j}V_\zeta}{V_0}\quad(1\le j\le n).
\]
Fix a nonnegative smooth function $\chi$ supported in an annulus
away from $y=0,\infty$ and positive on an open subset. Prescribing
$\int\chi Z_jw\dd V_{\rm rd}$ determines the kernel component of
$w$. These localized integrals will be the same ones used in the
middle-region equation.

\begin{lemma}\label{lem:sphere-augmented-inverse}
The map
\[
 \begin{aligned}
 \mathscr A_\zeta:H^1(\Sph^n)\times\R^{n+1}
       &\longrightarrow H^{-1}(\Sph^n)\times\R^{n+1},\\
 (w,c)&\longmapsto
 \left(L_\zeta w-\sum_{j=0}^nc_j\chi Z_j,
       \left(\int\chi Z_iw\dd V_{\rm rd}\right)_{i=0}^n\right)
 \end{aligned}
\]
is an isomorphism. Uniformly on a fixed compact set of $\zeta$,
\[
 \|\mathscr A_\zeta^{-1}(f,\mathbf m)\|_{H^1\times\R^{n+1}}
 \le C\bigl(\|f\|_{H^{-1}}+|\mathbf m|\bigr).
\]
Its inverse is smooth in $\zeta$, with bounded first operator
derivatives on that set.
\end{lemma}
\begin{proof}
The matrix
$(G_{\rm J})_{ij}=\int\chi Z_iZ_j\dd V_{\rm rd}$ is positive
definite: a kernel function vanishing where $\chi>0$ vanishes
identically. Pairing
$L_\zeta w-\sum c_j\chi Z_j=f$ with the $Z_i$ first determines
$c$ from
\[
 G_{\rm J}c=-\bigl(\langle f,Z_i\rangle\bigr)_{i=0}^n.
\]
The resulting source $f+\sum c_j\chi Z_j$ is orthogonal to
$\ker L_\zeta$. The round-sphere inverse from \cite{part1}, transferred
by conformal covariance, gives a solution. Subtracting its
$L^2$ projection onto $\ker L_\zeta$ gives the unique solution
orthogonal to that kernel. Adding a linear combination of the
$Z_j$ and solving with $G_{\rm J}$ imposes the required localized
moments. This proves existence and uniqueness. The inverse on the
orthogonal complement and the Gram matrices have uniform bounds
on compact parameter sets. Their coefficients are smooth, and
the inverse differentiation formula
\eqref{eq:inverse-parameter-calculus} gives the asserted parameter
bounds.
\end{proof}

Let $\mathcal S_\zeta(K)$ denote the first metric variation of the
action gradient at the bubble, now in the fixed representation:
\[
 \mathcal S_\zeta(K)
 =\left.\frac{d}{ds}\right|_{s=0}
       \Ff_{g_{\rm rd}+sK}(u_\zeta).
\]
Multiplying by the volume density in
\eqref{eq:compactified-action-source} gives the same first variation,
since $\Ff_{g_{\rm rd}}(u_\zeta)=0$. Conformal covariance transfers
the correction equation in \cite{part1} and its compatibility to
\begin{equation}
 L_\zeta(u_\zeta\psi)=-\mathcal S_\zeta(K_q),\qquad
 \langle\mathcal S_\zeta(K_q),Z_j\rangle=0
                      \quad(0\le j\le n),
 \label{eq:metric-source-compatibility}
\end{equation}
where $\psi$ is the imported multiplicative correction for the
bubble $V_\zeta$ and tensor $H_q$. The same identities hold for
$H_{q,R},K_{q,R}$. The integration-by-parts estimate in the proof of
Lemma~\ref{lem:compactified-metric-bounds} also gives
$\|\mathcal S_\zeta(K)\|_{H^{-1}}\le C\|K\|_{H^1}$, including one
parameter derivative.

\begin{lemma}\label{lem:first-metric-correction}
Suppose $n\ge4q+2$ and write $\lambda=b^q$. There is a unique
$W_1\in H^1(\Sph^n)$ satisfying
\[
 L_\zeta W_1=-\lambda\mathcal S_\zeta(K_q),\qquad
 \int\chi Z_jW_1\dd V_{\rm rd}=0\quad(0\le j\le n).
\]
The function $W_1$ is the imported correction
$\lambda u_\zeta\psi$, with its kernel component chosen by these
localized moments. For $j=0,1$, it satisfies
\[
 \|D_{\rm par}^jW_1\|_{H^1}\le C\lambda,\qquad
 \|D_{\rm par}^jW_1\|_{C^{6,\alphaH}_{\rm sc}(A_d)}
                                      \le C\lambda d^{-q},
\]
where $D_{\rm par}$ is a derivative in a unit direction of
$(\ell,\zeta)$ and $D_{\rm par}^0$ is the identity. It is smooth
at $y=0$ and independent of the regular data on $|X|=\rho$.
For $n<4q+2$, we use $W_1=0$.
\end{lemma}
\begin{proof}
Equation~\eqref{eq:metric-source-compatibility} identifies the
imported correction as a solution. Subtracting its kernel component
using $G_{\rm J}$ imposes zero localized moments without changing
the equation. Lemma~\ref{lem:sphere-augmented-inverse} gives
uniqueness and the $H^1$ bound. The imported decay estimate, divided
by $V_0$, gives $|D_{\rm par}^{\le1}W_1|\le C\lambda d^{-q}$;
the subtracted kernel functions are smooth and satisfy the same
bound. Here $\partial_\ell\lambda=q\lambda$, and the width and
center of $V_\zeta$ range over a compact set.
Since the cutoff and measure are fixed, differentiation of a
localized moment condition gives
\begin{equation}
 \int\chi Z_jw_s\dd V_{\rm rd}
   =-\int\chi(Z_j)_sw\dd V_{\rm rd}
 \quad\text{if }\int\chi Z_jw\dd V_{\rm rd}=0.
 \label{eq:localized-moment-parameter-calculus}
\end{equation}
In particular, the differentiated moments of $W_1$ have the same
$O(\lambda)$ bound.

On an enlarged annulus, the source has scaled $C^{4,\alphaH}$ norm
at most $C\lambda d^{-q-2}$. Scaled Schauder estimates applied to
the equation, and then to its first parameter derivatives, give
the stated $C^{6,\alphaH}$ estimate. In a parameter direction $s$,
the differentiated equation is
\[
 L_\zeta(W_1)_s
 =-\partial_s(\lambda\mathcal S_\zeta(K_q))-(L_\zeta)_sW_1,
\]
where $(L_\zeta)_s$ is a bounded smooth potential.
The source is smooth at $y=0$, so the correction is smooth there.
Its Taylor expansion on $|y|=\delta$ has a constant and a linear
part, and the remaining boundary trace is
$O_{C^1}(\lambda\delta^2)$.

Changing the correction by a kernel function does not change the
quadratic action coefficient in \eqref{eq:part-i-cutoff-action}.
On the bubble's round sphere, with metric $\gamma$, this change is a first harmonic $Y$,
and its conformal metric variation is a Lie derivative:
\[
 \frac4{n-2}Y\gamma
   =\mathcal L_{-2\nabla_\gamma Y/(n-2)}\gamma.
\]
The quadratic form of \cite[Proposition 6.6]{part1} annihilates
these directions.
\end{proof}

\subsection{The Fowler profile and its Jacobi equations}
\label{sec:inner-tail}\label{subsec:tail-models}\label{def:tail-spaces}

On the inner region $0<|y|\le\delta$, put
\[
 r=|y|=\delta e^{-t},\qquad
 \mathcal C=[0,\infty)\times\Sph^{n-1},\qquad
 \delta=\eps^{2/n},\quad \lambda=b^q\le C\eps.
\]
Thus $t=0$ is the sphere joined to the middle region and
$t\to\infty$ approaches the puncture. Here $n\ge8$ and $q\ge2$.
Recall that $v_\eps(0)=\eps$ and $v_\eps'(0)=0$, so its
minimum occurs at zero. Let $T_\eps$ be its period. The profile
increases on $(0,T_\eps/2)$ and decreases on
$(T_\eps/2,T_\eps)$, from its maximum to its next minimum.
For small $\eps$, there is a unique point in this decreasing
half-period at which its value is $2^{(n-2)/2}\delta^{(n-2)/2}$.
Translate that point to $t=0$ and denote the resulting profile
by $v$. Thus
\[
 v(0)=2^{(n-2)/2}\delta^{(n-2)/2},\qquad v'(0)<0.
\]
As $t$ increases from zero, $r=\delta e^{-t}$ decreases toward
the puncture, and $v$ initially decreases to its next minimum;
it then increases again and continues periodically on $\R$.
The initial value is chosen so that
$\delta^{-(n-2)/2}v(0)=V_\zeta(0)$.
The profiles used on the inner end are restrictions and translates
of this periodic function on $t\ge0$. Its equation and Hamiltonian are
\begin{equation}
\begin{aligned}
 v''-\frac{(n-2)^2}{4}v+\frac{n(n-2)}4v^p&=0,\\
 h_\eps
 &=\tfrac12(v')^2-\frac{(n-2)^2}{8}v^2
                      +\frac{(n-2)^2}{8}v^{p+1}\\
 &=-\frac{(n-2)^2}{8}\eps^2(1-\eps^{4/(n-2)}).
\end{aligned}
 \label{eq:tail-energy}
\end{equation}
In particular, $\eps/v(0)\asymp\delta$ and
$v'(0)/v(0)=-(n-2)/2+O(\delta^2)$.

Since $|v'|/v\le(n-2)/2$, the times from the preceding maximum
to $t=0$ and from $t=0$ to the next minimum are at least
\[
 \frac2{n-2}\log\frac{\max v_\eps}{v(0)},\qquad
 \frac2{n-2}\log\frac{v(0)}\eps,
\]
respectively. Both tend to infinity as $\eps\to0$, since
$\eps\ll v(0)\ll\max v_\eps$.
For a phase $\tau$ in a fixed bounded interval, write
$v_\tau(t)=v(t+\tau)$. For small $\eps$, all these profiles have
$v_\tau'(0)<0$, $v_\tau(0)\asymp v(0)$ and
$|v_\tau'(0)|\asymp v(0)$. Indeed the energy identity and
$|v'|\le(n-2)v/2$ give these estimates on every fixed bounded
interval about zero. The norm weights below always use $v$;
$\eps,\delta$ and these weights are fixed under differentiation
in $\tau$ or the boundary data.

The cylindrical linearization at $v_\tau$ is
\[
 A_\tau=-\partial_t^2-\Delta_\theta+\frac{(n-2)^2}{4}
                     -\frac{n(n+2)}4v_\tau^{p-1}.
\]
On spherical harmonics of degree $l$, it is the scalar operator
\[
 A_{\tau,l}
 =-\frac{d^2}{dt^2}
   +\left(l+\frac{n-2}{2}\right)^2
   -\frac{n(n+2)}4v_\tau^{p-1}.
\]
The radial component will be solved with a prescribed limiting
Hamiltonian. For degree one there are two explicit positive
Jacobi functions,
\[
 j_-=r\left(\frac{n-2}{2}v-v'\right),\qquad
 j_+=r^{-1}\left(\frac{n-2}{2}v+v'\right).
\]
They satisfy $A_{0,1}j_\pm=0$ and
\[
 j_-j_+'-j_-'j_+=-2nh_\eps=:W_\eps\asymp\eps^2.
\]
The corresponding functions for $v_\tau$ are obtained by replacing
$v,v'$ in these formulas by $v_\tau,v_\tau'$ while keeping
$r=\delta e^{-t}$.

We record some factors that occur in the first-harmonic equation
and in the nonlinear products:
\[
 F=\frac{j_-}{r}=\frac{n-2}{2}v-v',\qquad
 G=rj_+=\frac{n-2}{2}v+v',\qquad
 \omega=\frac GF,\qquad \mathsf p=j_-^2.
\]
Thus $\mathsf p$ is the coefficient of the divergence-form equation
obtained by writing a degree-one correction as $j_-z$.
The energy identity gives
\begin{equation}
\begin{gathered}
 FG=\frac{(n-2)^2}{4}
       \bigl[v^{p+1}+\eps^2(1-\eps^{4/(n-2)})\bigr],\\
 \frac{\eps^2}{vF}\le C,\qquad
 \frac vF\le C\delta^{-2},\qquad
 \frac\eps F\le C\delta^{-1},\qquad
 \frac{v^2}F\le C.
\end{gathered}
 \label{eq:weight-bounds}
\end{equation}
For example, $G\le(n-2)v$ implies
$F\ge c(v^{p+1}+\eps^2)/v$. Maximizing the resulting bounds for
$v/F$ and $\eps/F$ gives the middle two estimates at
$v\asymp\eps^{2/(p+1)}=\delta^{(n-2)/2}$. The last follows from
$v^2/F\le Cv^{2-p}$ and $p<2$.
The Fowler equation gives bounded logarithmic derivatives of
$v,F,G$ of every fixed order.

Put $F_\tau=(n-2)v_\tau/2-v_\tau'$ and
$\mathsf p_\tau=r^2F_\tau^2$. The preceding bounds imply
\begin{equation}
\begin{gathered}
 c\le v_\tau/v\le C,\qquad c\le F_\tau/F\le C,\\
 \sum_{i=0}^1\sum_{k=0}^6
 \left\|\partial_\tau^i\partial_t^k
 \left(\frac{v_\tau}v,\frac v{v_\tau},
       \frac{F_\tau}F,\frac F{F_\tau}\right)
 \right\|_{C_{\rm b}^{0,\alphaH}}\le C.
\end{gathered}
 \label{eq:fowler-reference-multipliers}
\end{equation}
Here and below, $C_{\rm b}^{k,\alphaH}$ is the ordinary bounded
H\"older space on the half-cylinder, or on the half-line for
functions of $t$ alone. Its norm is equivalent to the supremum
of the same local norms on overlapping unit cylinders. All
Sobolev and Lebesgue norms use $\dd t\dd\theta$, or $\dd t$ in
the latter case. The symbol $D^j$ denotes all product covariant
derivatives of order $j$.

\subsection{The first and higher spherical harmonics}
\label{subsec:tail-inverses}

\medskip\noindent\textit{Degree one: variation of parameters.}
For the reference profile,
\[
 A_{0,1}(j_-z)=-j_-^{-1}(j_-^2z')'.
\]
The nonlinear center equation also contains a differentiated
source. We therefore consider given functions
$f,\mathsf j:[0,\infty)\to\R^n$ and the equation
\[
 -(\mathsf p_\tau z')'=f+(\mathsf p_\tau\mathsf j)',
 \qquad z(0)=0,\qquad \lim_{t\to\infty}z(t)\ \hbox{exists}.
\]
Integrating from infinity determines the flux:
\begin{equation}
 z'(t)=\mathsf p_\tau(t)^{-1}\int_t^\infty f(s)\dd s
                  -\mathsf j(t),\qquad z(0)=0.
 \label{eq:kelvin-tail-inverse}
\end{equation}
The Wronskian identity gives
\begin{equation}
 \int_0^u\mathsf p^{-1}\dd t
 =\frac{j_+(u)/j_-(u)-j_+(0)/j_-(0)}{W_\eps}.
 \label{eq:exact-capacity}
\end{equation}
Fubini's theorem and \eqref{eq:fowler-reference-multipliers}
therefore give
\begin{equation}
 \|z\|_\infty
 \le\frac C{W_\eps}\int_0^\infty\frac{j_+}{j_-}|f|\dd t
                         +\int_0^\infty|\mathsf j|\dd t.
 \label{eq:capacity-bound}
\end{equation}
If the right side is finite, \eqref{eq:kelvin-tail-inverse}
defines a locally absolutely continuous solution with $z'\in L^1$,
so $z$ has a finite limit. A nonzero homogeneous flux would add
$c\int_0^t\mathsf p_\tau^{-1}$, which diverges; hence these two
conditions determine the solution uniquely. The ratio
$j_+/j_-=r^{-2}\omega$ in the estimate is the moment against the
growing Jacobi solution.

We now specify the derivative and energy estimates used for this
equation. Fix
\[
 M=6,\qquad B_M=(1+t)^3,\qquad q_K=\eps^{-1}B_Mrv.
\]
The Kelvin coefficient space consists of functions $z$ with
$z(0)=0$ and finite norm
\begin{equation}
\begin{aligned}
 K(z):=\|z\|_{\mathcal K_{\rm tail}}
 ={}&\|z'\|_{W^{5,1}}+\|z'\|_{C_{\rm b}^{5,\alphaH}}\\
 &+\|q_Kz'\|_{H^5}+\|q_Kz'\|_{C_{\rm b}^{5,\alphaH}}.
\end{aligned}
 \label{eq:kelvin-six-norm}
\end{equation}
In particular $\|z\|_\infty\le\|z'\|_1\le K(z)$.
For stating source estimates, use the scalar or vector norm
\[
 \|h\|_{\mathcal K_{\rm tail,0}}
 =\|h\|_1+\|h\|_{C_{\rm b}^{0,\alphaH}}
       +\|q_Kh\|_2+\|q_Kh\|_{C_{\rm b}^{0,\alphaH}},
 \qquad
 e_f(t)=\mathsf p(t)^{-1}\int_t^\infty f(s)\dd s.
\]
The source space $\mathcal Y_K$ has norm
\[
\begin{aligned}
 \|(f,\mathsf j)\|_{\mathcal Y_K}
 ={}&\eps^{-2}\int_0^\infty\frac{j_+}{j_-}|f|\dd t
       +\|e_f\|_{\mathcal K_{\rm tail,0}}\\
 &+\sum_{k=0}^4
       \|\partial_t^k(f/\mathsf p)\|_{\mathcal K_{\rm tail,0}}
  +\sum_{k=0}^5
       \|\partial_t^k\mathsf j\|_{\mathcal K_{\rm tail,0}}.
\end{aligned}
\]
The functions $f,\mathsf j$ are locally $C^{4,\alphaH}$ and
$C^{5,\alphaH}$, respectively. Since the logarithmic derivatives
of $q_K$ are bounded, \eqref{eq:kelvin-six-norm} is uniformly equivalent
to $\sum_{k=1}^6\|\partial_t^kz\|_{\mathcal K_{\rm tail,0}}$.

\begin{lemma}\label{lem:tail-low}
The inverse \eqref{eq:kelvin-tail-inverse} satisfies
$K(z)\le C\|(f,\mathsf j)\|_{\mathcal Y_K}$, uniformly for
$\tau$ in a fixed bounded interval and small $\eps$.
Put
\[
 C_f=\eps^{-2}\int_0^\infty\frac{j_+}{j_-}|f|\dd t,
 \qquad I_f=\int_0^\infty B_M^2|f|\dd t.
\]
A sufficient bound for the data norm is
\begin{equation}
\begin{aligned}
 \|(f,\mathsf j)\|_{\mathcal Y_K}\le C\bigg[
 &C_f+\eps^{-1}\delta^{-2}(C_fI_f)^{1/2}\\
 &+\sum_{k=0}^4
    \|\partial_t^k(f/\mathsf p)\|_{\mathcal K_{\rm tail,0}}
  +\sum_{k=0}^5
    \|\partial_t^k\mathsf j\|_{\mathcal K_{\rm tail,0}}
 \bigg].
\end{aligned}
 \label{eq:kelvin-data-bound}
\end{equation}
The inverse is continuously differentiable in $\tau$ on the fixed
spaces, with the same operator bound for its derivative.
Moreover $\mathsf p_\tau(0)\asymp\mathsf p(0)\asymp\eps^2$.
\end{lemma}
\begin{proof}
It remains to estimate the derivatives of the explicit inverse.
For $\bar f=|f|$, put
$e_{\bar f}=\mathsf p^{-1}\int_t^\infty\bar f$ and
$z_{\bar f}(t)=\int_0^te_{\bar f}$.
Equation~\eqref{eq:exact-capacity} gives
$\|e_{\bar f}\|_1=\|z_{\bar f}\|_\infty\le CC_f$.
Integrating $(\mathsf p e_{\bar f})'=-\bar f$ against
$B_M^2z_{\bar f}$ yields
\[
\begin{aligned}
 \int_0^T B_M^2\mathsf p e_{\bar f}^2
 ={}&[B_M^2\mathsf p e_{\bar f}z_{\bar f}]_0^T
       +\int_0^T B_M^2\bar fz_{\bar f}\\
 &-\int_0^T(B_M^2)'\mathsf p e_{\bar f}z_{\bar f}.
\end{aligned}
\]
The last term is nonpositive. The term at zero vanishes, and
$B_M(T)^2\mathsf p(T)e_{\bar f}(T)
\le\int_T^\infty B_M^2\bar f\to0$ when $I_f<\infty$.
Since $|e_f|\le e_{\bar f}$ and
$q_K^2/\mathsf p=\eps^{-2}B_M^2(v/F)^2$, it follows that
\begin{equation}
 \int B_M^2\mathsf p|e_f|^2\le CC_fI_f,\qquad
 \|q_Ke_f\|_2
       \le C\eps^{-1}\delta^{-2}(C_fI_f)^{1/2}.
 \label{eq:kelvin-work}
\end{equation}
The identities
\[
 e_f'=-(\log\mathsf p)'e_f-f/\mathsf p,\qquad
 (q_Ke_f)'=((\log q_K)'-(\log\mathsf p)')q_Ke_f
                                      -q_Kf/\mathsf p
\]
first give uniform $C^0$ bounds for $e_f$ and $q_Ke_f$ by the
embeddings $W^{1,1},W^{1,2}\hookrightarrow C^0$ on enlarged
unit intervals. Substituting those bounds and the bounded source
terms $f/\mathsf p$, $q_Kf/\mathsf p$ into the same identities
bounds $e_f'$ and $(q_Ke_f)'$. The resulting Lipschitz bounds,
together with the supremum bounds, control their global
$C_{\rm b}^{0,\alphaH}$ norms. This proves
\eqref{eq:kelvin-data-bound}.
Four further differentiations control
$\partial_t^ke_f$ in $\mathcal K_{\rm tail,0}$ for $k\le5$.
The fixed-phase inverse then follows from $z'=e_f-\mathsf j$.

For general $\tau$, replace $e_f$ by $(F/F_\tau)^2e_f$.
The multipliers in \eqref{eq:fowler-reference-multipliers}
preserve all these estimates. For a differentiable family of
data and phases, differentiation gives
\begin{equation}
 (Dz)'=(F/F_\tau)^2e_{Df}
      +D\bigl((F/F_\tau)^2\bigr)e_f-D\mathsf j,\qquad Dz(0)=0.
 \label{eq:kelvin-parameter-identity}
\end{equation}
This proves the asserted derivative bound, and also the bound
with differentiated data. Finally
$F_\tau(0)\asymp v(0)\asymp\eps/\delta$, so
$\mathsf p_\tau(0)=\delta^2F_\tau(0)^2\asymp\eps^2$.
\end{proof}

\medskip\noindent\textit{Degrees at least two: the Green kernel.}
For $l\ge2$, the potential of $A_{\tau,l}$ is at least
$(n+2)/2$. Let $\psi_l$ be a positive solution of
$A_{\tau,l}\psi_l=0$ that decays at infinity. Reduction of order
gives the Dirichlet Green kernel
\[
 G_l(t,s)
 =\psi_l(t)\psi_l(s)
       \int_0^{\min(t,s)}\psi_l(u)^{-2}\dd u.
\]
Thus a mode with forcing $f_{l,m}$ and boundary value $h_{l,m}$
is given by
\[
 u_{l,m}(t)
 =\frac{\psi_l(t)}{\psi_l(0)}h_{l,m}
       +\int_0^\infty G_l(t,s)f_{l,m}(s)\dd s.
\]
We require a decay rate uniform as the neck shrinks:
\[
 0\le G_l(t,s)\le C_\kappa e^{-\kappa|t-s|},
 \qquad 0<\kappa<\frac{n+2}{2}.
\]
Here the constants are uniform in $l\ge2$ and the phase.
The mode decomposition is the one used in
\cite[Proposition 3]{kmps} and \cite[Proposition 2.3]{santos};
we give the uniform kernel estimate needed here.

Put $\eta=v_\tau^{4/(n-2)}$ and denote its minimum and maximum
by $\eta_-$ and $\eta_+$. The energy identity gives
\[
 (\eta')^2=4\eta^2Q_\eps(\eta),\qquad
 Q_\eps(z)=1-z-\eta_-^{(n-2)/2}(1-\eta_-)z^{-(n-2)/2}.
\]
For small necks, $\eta_-<1/4<1/2<\eta_+<1$ and
$Q_\eps(1/2)\ge3/8$. Concavity and the two chords from $1/2$
give $Q_\eps(z)\ge\frac34(z-\eta_-)(\eta_+-z)$.
Changing variables on the two half-periods therefore yields
\[
 \int_0^{T_\eps}\eta\dd t\le\frac{2\pi}{\sqrt3},
 \qquad
 T_\eps\ge\log\frac1{4\eta_-}\longrightarrow\infty.
\]
Every interval contains complete periods and at most two
incomplete ones, so
$\int_s^t\eta\le C_n(2+|t-s|/T_\eps)$ for $s<t$.
For the potential $q_l$ of $A_{\tau,l}$, set
$k(t)=(n+2)/2-n\eta(t)/2$.
The inequality $|\eta'|\le2\eta\sqrt{1-\eta}$ gives
$q_l-k^2-|k'|\ge(n-2)\eta/2\ge0$.
Consequently
\[
 E_s(t)=\tfrac12
 \exp\left(-\int_{\min(s,t)}^{\max(s,t)}k(u)\dd u\right)
\]
is a supersolution for the Green equation: its derivative jump
at $s$ is $-k(s)\le-1$.
The maximum principle on finite intervals, followed by exhaustion,
gives
\[
 0\le G_l(t,s)\le E_s(t)
 \le C\exp\left[-\left(\frac{n+2}{2}
                 -\frac{nC_n}{2T_\eps}\right)|t-s|\right].
\]
This proves the stated uniform estimate after decreasing the
neck threshold for any fixed $\kappa<(n+2)/2$.

The norms for the higher harmonics express two requirements:
decay faster than $j_-$, and control of their quadratic
contribution to the radial and first-harmonic equations.
Fix
\begin{equation}
 \gamma_*=\frac32.
 \label{eq:tail-buffers}
\end{equation}
The second requirement uses the single positive weight
\begin{equation}
 \beta_\eps
 =\eps^{-1}B_M
       \left(1+\delta^2\frac{j_+}{j_-}\right)^{1/2}.
 \label{eq:high-work-norm}
\end{equation}
For $j=4,6$, define
\begin{equation}
\begin{aligned}
 \mathcal W_j(u)
 &=\left\|\frac{e^{t/2}u}{j_-}\right\|_{C_{\rm b}^{j,\alphaH}}
       +\|\beta_\eps u\|_{C_{\rm b}^{j,\alphaH}}
       +\|\beta_\eps u\|_{H^j},
 \qquad \mathcal W=\mathcal W_6.
\end{aligned}
 \label{eq:relative-flux-norm}
\end{equation}
The first term is
$\delta^{-1}\|e^{\gamma_*t}u/F\|_{C_{\rm b}^{j,\alphaH}}$.
The square of the second weight is
$\eps^{-2}B_M^2(1+\delta^2r^{-2}\omega)$.
All derivatives of the weights divided by the weights are
uniformly bounded. Hence the squared Sobolev term controls, for
$0\le k\le j$,
\[
 \eps^{-2}\int B_M^2|D^ku|^2,\qquad
 \eps^{-2}\delta^2\int B_M^2\frac{j_+}{j_-}|D^ku|^2.
\]
The second integral contains the same Jacobi ratio as
\eqref{eq:capacity-bound}; it is used to bound the quadratic
sources in that equation.

Let $\mathcal X_{\rm hi}$ consist of fields with
$\Pi_0w=\Pi_1w=0$ and finite $\mathcal W(w)$, and let
$\mathcal X_{{\rm hi},0}$ impose $w(0,\theta)=0$.
The source space uses $\mathcal W_4$ and the same projections.
These are Banach spaces with the indicated norms. In particular,
\[
 \frac{|D^{\le6}w|}{v}\le C\delta\mathcal W(w),\qquad
 \frac{|D^{\le6}w|}{F}
       \le C\delta e^{-\gamma_*t}\mathcal W(w).
\]
Boundary data in the coupled construction will be measured by
$w(0,\theta)=\delta v(0)\phi(\theta)$, with
$\phi\in C^{6,\alphaH}_{\ge2}(\Sph^{n-1})$.

\begin{lemma}\label{lem:high-tail}
For sources $\Pi_{\ge2}f=f$ with $\mathcal W_4(f)<\infty$,
$A_\tau$ has a unique zero-Dirichlet inverse $R_{H,\tau}$
with values in $\mathcal X_{{\rm hi},0}$, and
\[
 \mathcal W_6(R_{H,\tau}f)\le C\mathcal W_4(f).
\]
For $h\in C^{6,\alphaH}_{\ge2}(\Sph^{n-1})$ with
$\|h\|_{C^{6,\alphaH}}\le\Gamma v(0)$, the homogeneous solution
$W_c$ with $W_c(0,\theta)=h(\theta)$ satisfies
\begin{equation}
 |D^iW_c|\le C\Gamma v(0)e^{-(2n+3)t/4}\quad(0\le i\le6),
 \qquad
 \mathcal W(W_c)\le C\Gamma/\delta.
 \label{eq:high-cap-bound}
\end{equation}
For fixed data, both bounds hold after one $\tau$ derivative
with the reference weights fixed. The solution depends linearly
and continuously on its boundary value. In particular,
\[
 \mathcal W(D_hW_c[\dot h])
       \le \frac{C\|\dot h\|_{C^{6,\alphaH}}}{\delta v(0)}.
\]
\end{lemma}
\begin{proof}
The Green kernel above defines the inverse mode by mode.
Fix the rate
\begin{equation}
 \sigma_0=\frac{4n+7}{8},\qquad
 \max\left\{\frac{2n+3}{4},\frac{n+4}{2}-\gamma_*\right\}
       <\sigma_0<\frac{n+2}{2}.
 \label{eq:high-decay-choice}
\end{equation}
Its kernel has operator norm at most
$Ce^{-\sigma_0|t-s|}$ on $L^2_{\ge2}(\Sph^{n-1})$.
Localizing the source to unit cylinders and applying interior
or zero-Dirichlet boundary Schauder estimates gives the same
decay for the localized map from $C^{4,\alphaH}$ to
$C^{6,\alphaH}$. Local Sobolev estimates give the corresponding
map from $H^4$ to $H^6$.

The Fowler identities imply
\begin{equation}
 F\ge cv(0)e^{-(n-2)t/2},\qquad
 (\log\omega)'=(n-2)
 -\frac{nv^{p+1}}{v^{p+1}+\eps^2(1-\eps^{4/(n-2)})},
 \quad -2\le(\log\omega)'\le n-2.
 \label{eq:omega-growth}
\end{equation}
Since $-(n-2)/2\le F'/F\le(n+2)/2$, the two weights in
\eqref{eq:relative-flux-norm} satisfy
\[
 \frac{e^{t/2}/j_-(t)}{e^{s/2}/j_-(s)}
       \le e^{(n+1)|t-s|/2},\qquad
 \frac{\beta_\eps(t)}{\beta_\eps(s)}
       \le(1+|t-s|)^3e^{n|t-s|/2}.
\]
Conjugating the Green kernel by either weight thus leaves
uniformly integrable rows and columns. Convolution and Schur's
test, followed by the local Schauder and Sobolev estimates,
prove $\mathcal W_6(R_{H,\tau}f)\le C\mathcal W_4(f)$.
Truncation gives the same inverse for general sources in this
space. For a homogeneous solution with zero trace,
$w,\partial_tw\in L^2$; integration by parts along cuts where
their boundary pairing tends to zero and the positive potential
give $w=0$.

For boundary data $h$, extend it by a fixed axial cutoff supported
in $[0,2]$, denoting this extension by $\mathcal E h$. Then
\[
 W_c=\mathcal E h-R_{H,\tau}(A_\tau\mathcal E h).
\]
On this fixed interval,
$v,F\asymp v(0)\asymp\eps/\delta$ and
$\omega\asymp\delta^2$. The two weights therefore give
\[
 \mathcal W_6(\mathcal E h)
       +\mathcal W_4(A_\tau\mathcal E h)\le C\Gamma/\delta.
\]
The Green estimate gives
$|D^{\le6}W_c|\le C\Gamma v(0)e^{-\sigma_0t}$.

Finally $\partial_\tau A_\tau$ is a uniformly bounded smooth
potential, including four derivatives. On the fixed spaces,
\[
 \partial_\tau R_{H,\tau}
       =-R_{H,\tau}(\partial_\tau A_\tau)R_{H,\tau}.
\]
This proves the derivative estimate for the inverse.
Differentiating the boundary solution gives source
$-(\partial_\tau A_\tau)W_c$ and zero boundary value when $h$
is fixed. Convolution bounds this derivative by
$C\Gamma v(0)(1+t)e^{-\sigma_0t}$, which implies
\eqref{eq:high-cap-bound}. Linearity gives the boundary-data
derivative.
\end{proof}

\subsection{The forced radial equation}

We seek a radial solution with two boundary conditions: value $m_0$
at $t=0$ and limiting Hamiltonian $h_\eps$. Its limiting phase
is determined by the equation. We keep $\eps$ and the reference
profile $v$ fixed, including when differentiating in $m_0$ or in
the forcing; thus all norm weights remain fixed.
The Hamiltonian is
\[
 \mathsf h(V,W)=\frac12W^2-\frac{(n-2)^2}{8}V^2
                         +\frac{(n-2)^2}{8}V^{p+1}.
\]
For a solution of
$V''-\frac{(n-2)^2}{4}V+\frac{n(n-2)}4V^p=\radsrc$,
\[
 \frac{d}{dt}\mathsf h(V,V')=V'\radsrc.
\]
Thus prescribing the limiting Hamiltonian amounts to requiring
\[
 \mathsf h(V(t),V'(t))-h_\eps
       =-\int_t^\infty V'(s)\radsrc(s)\dd s.
\]
To obtain uniform estimates as $\eps\to0$, we measure the forcing
relative to $|h_\eps|\asymp\eps^2$. For functions on $[0,\infty)$,
define
\begin{equation}
 \|\radsrc\|_{\mathcal Y_{\rm rad}}
 =\eps^{-2}\left[\|(1+t)^6v\radsrc\|_{W^{4,1}}
  +\|(1+t)^6v\radsrc\|_{C_{\rm b}^{4,\alphaH}}\right].
 \label{eq:master}
\end{equation}
The space $\mathcal Y_{\rm rad}$ consists of the locally
$C^{4,\alphaH}$ functions with finite displayed norm.  Its weights
depend on  $v$. Bounded logarithmic derivatives
of $v$ show that this norm is equivalent, uniformly in $\eps$, to
the sum of the $L^1$ and $C_{\rm b}^{0,\alphaH}$ norms of
$\eps^{-2}(1+t)^6v\partial_t^j\radsrc$, $0\le j\le4$.

\begin{lemma}\label{lem:radial-response}
Let $m_0/v(0)$ range over a fixed compact interval in $(0,\infty)$.
For sufficiently small $\eps$, let $\tau$ be the unique point in
the decreasing half-period of $v$ containing zero for which
$v(\tau)=m_0$. Thus $v_\tau(0)=m_0$ and $v_\tau'(0)<0$;
$\tau$ remains bounded uniformly in $\eps$ on the stated interval
of $m_0/v(0)$. There is a constant \(c>0\), independent of \(\varepsilon\), 
such that, if $\|\mathfrak f\|_{\mathcal Y_{\rm rad}}<c$, namely,
\[
 \|(1+t)^6v\radsrc\|_{W^{4,1}}
 +\|(1+t)^6v\radsrc\|_{C_{\rm b}^{4,\alphaH}}
 <c\eps^2,
\]
the construction below gives a positive radial solution
$V=R_\eps(m_0,\radsrc)$ of
\[
 V''-\frac{(n-2)^2}{4}V+\frac{n(n-2)}4V^p=\radsrc,
 \qquad V(0)=m_0,
 \qquad \lim_{t\to\infty}\mathsf h(V,V')=h_\eps.
\]
There is a real phase $\alpha_\infty$ such that
\[
 \sum_{j=0}^6
 \frac{|\partial_t^j[V(t)-v(t+\alpha_\infty)]|}{v(t)}
       \longrightarrow0.
\]
The solution and $\alpha_\infty$ depend jointly $C^1$ on
$\radsrc$ and $m_0/v(0)$.
Writing
$F_V=\frac{n-2}{2}V-V'$ and
$F_\tau=\frac{n-2}{2}v_\tau-v_\tau'$, we have
\begin{equation}
 \sum_{j=0}^1\left[
 \left\|[v(0)\partial_{m_0}]^j
                 \left(\frac{V}{v_\tau}-1\right)\right\|_{C_{\rm b}^{6,\alphaH}}
 +\left\|[v(0)\partial_{m_0}]^j
                 \left(\frac{F_V}{F_\tau}-1\right)\right\|_{C_{\rm b}^{5,\alphaH}}
 \right]
 \le C\|\radsrc\|_{\mathcal Y_{\rm rad}}.
 \label{eq:radial-background-first-difference}
\end{equation}
Here the boundary derivative holds $\radsrc$ fixed. The Fr\'echet
derivative in a source direction $\dot{\radsrc}$ satisfies
\[
 \left\|\frac{D_\radsrc V[\dot{\radsrc}]}v\right\|_{C_{\rm b}^{6,\alphaH}}
 +\left\|\frac{D_\radsrc F_V[\dot{\radsrc}]}F\right\|_{C_{\rm b}^{5,\alphaH}}
 \le C\|\dot{\radsrc}\|_{\mathcal Y_{\rm rad}}.
\]
All constants are uniform on the stated boundary-value interval.
\end{lemma}

\begin{proof}
\noindent\textit{Phase and Hamiltonian near the Fowler orbit.}
We use coordinates that remain regular at the maxima and minima
of $v$. The equation and its Hamiltonian identity give
\begin{equation}
 c_nv^2\le (v')^2+(v'')^2\le C_nv^2.
 \label{eq:radial-orbit-nondegeneracy}
\end{equation}
Indeed, if the lower bound failed as $\eps\to0$, the equation
would give $v^{4/(n-2)}\to(n-2)/n$, whereas the Hamiltonian
identity would give $(v'/v)^2\to(n-2)^2/(2n)>0$.

At the point $(v(s),v'(s))$, the vector
$(-v''(s),v'(s))=\nabla\mathsf h(v(s),v'(s))$ is normal to the
Fowler orbit. Define $Q(s,\energyvar)=(Q^V(s,\energyvar),Q^W(s,\energyvar))$ to be the
nearby point on this normal line with Hamiltonian
$h_\eps+|h_\eps|\energyvar$. Explicitly,
\[
 \begin{gathered}
 Q(s,\energyvar)=(v(s),v'(s))
                  +\beta(s,\energyvar)(-v''(s),v'(s)),\\
 \mathsf h(Q(s,\energyvar))=h_\eps+|h_\eps|\energyvar.
 \end{gathered}
\]
The coefficient $\beta$ is determined by these equations.
The derivative of the second equation with
respect to $\beta$ at zero is $(v')^2+(v'')^2$. Hence the implicit
function theorem, with distances divided by $v(s)$, gives a fixed
interval $|\energyvar|\le \energyvar_*$ on which
\[
 Q(s,0)=(v(s),v'(s)),\qquad
 |Q(s,\energyvar)-(v(s),v'(s))|\le C|h_\eps||\energyvar|/v(s).
\]
The constants are uniform in $s$ and small $\eps$. To check the
differentiated bounds, divide the scalar implicit equation by
$v(s)^2$. Its derivative in $\beta$ is bounded below, and all
required derivatives are bounded: the kinetic and quadratic terms
have this property directly, while the power term has the factor
$v(s)^{p-1}$ and a positive argument bounded away from zero after
this rescaling. Differentiation therefore gives, for every fixed
order used below,
\[
 \left|\partial_s^j\partial_\energyvar^kQ(s,\energyvar)\right|
 \le C_{j,k}|h_\eps|^k v(s)^{1-2k}\quad(k\ge1),\qquad
 \left|\partial_s^j\bigl[Q(s,\energyvar)-Q(s,0)\bigr]\right|
 \le C_j|h_\eps||\energyvar|/v(s).
\]
In particular, $Q^V\asymp v$ and
$\frac{n-2}{2}Q^V-Q^W\asymp F(s)$. The latter assertion and its
relative derivative bounds follow from
$|h_\eps|/(vF)\le C$ in \eqref{eq:weight-bounds}.

The determinant $\det(Q_s,Q_\energyvar)$ equals $|h_\eps|$ at $\energyvar=0$ and
is $|h_\eps|[1+O(|h_\eps||\energyvar|/v(s)^2)]$ on this interval.
Define the two coefficients
\[
 A(s,\energyvar)=\frac{|h_\eps|}{\det(Q_s,Q_\energyvar)}-1,
 \qquad B(s,\energyvar)=-\frac{Q_\energyvar^V}{\det(Q_s,Q_\energyvar)}.
\]
Writing $(V,V')=Q(t+\alpha(t),\energyvar(t))$, the forced equation is
equivalent to
\begin{equation}
 \alpha'=A(t+\alpha,\energyvar)+B(t+\alpha,\energyvar)\radsrc,
 \qquad \energyvar'=Q^W(t+\alpha,\energyvar)\radsrc/|h_\eps|.
 \label{eq:radial-chart}
\end{equation}
For this identity, $Q_s$ is tangent to a level curve of
$\mathsf h$, and
$\det(X_{\mathsf h}(Q),Q_\energyvar)=|h_\eps|$ for
$X_{\mathsf h}(V,W)=(W,\frac{(n-2)^2}{4}V-\frac{n(n-2)}4V^p)$.
Taking determinants of
$(1+\alpha')Q_s+\energyvar'Q_\energyvar=X_{\mathsf h}(Q)+(0,\radsrc)$
gives the first equation; applying $D\mathsf h$ gives the second.
The implicit estimates imply
\begin{equation}
 \begin{aligned}
 |\partial_s^j A(s,\energyvar)|&\le C_j\eps^2v(s)^{-2}|\energyvar|,\\
 |\partial_s^j\partial_\energyvar^k A(s,\energyvar)|&\le C_{j,k}\eps^2v(s)^{-2}
                       &&(k\ge1),\\
 |\partial_s^j\partial_\energyvar^k B(s,\energyvar)|&\le C_{j,k}/v(s),\\
 |\partial_s^j\partial_\energyvar^k Q^W(s,\energyvar)|&\le C_{j,k}v(s).
 \end{aligned}
 \label{eq:radial-drift-derivatives}
\end{equation}
Only finitely many orders, at most six, are needed. The factor
$|\energyvar|$ in the first line will control the dependence on the phase.

\medskip\noindent\textit{The two boundary conditions and a contraction.}
The equation $Q^V(\alpha(0),\energyvar(0))=m_0$ determines
$\alpha(0)=\alpha_{\rm in}(\energyvar(0),m_0)$. Indeed,
$Q_s^V(\tau,0)=v'(\tau)\le-cv(0)$, so the implicit function
theorem applies uniformly and gives
\[
 \alpha_{\rm in}(0,m_0)=\tau,\qquad
 |\partial_\energyvar\alpha_{\rm in}|\le C\eps^2/v(0)^2\le C\delta^2,
 \qquad |v(0)\partial_{m_0}\alpha_{\rm in}|\le C.
\]
These bounds hold with one further derivative in the scaled
boundary value.

For continuous $\alpha,\energyvar$, with $\alpha-\tau$ bounded and
$\energyvar\in C_{\rm b}\cap L^1$, define
\[
 \begin{aligned}
 \energyvar_{\rm new}(t)
 &=-|h_\eps|^{-1}\int_t^\infty
          Q^W(s+\alpha(s),\energyvar(s))\radsrc(s)\dd s,\\
 \alpha_{\rm new}(t)
 &=\alpha_{\rm in}(\energyvar_{\rm new}(0),m_0)\\
 &\quad+\int_0^t\bigl[
       A(s+\alpha(s),\energyvar_{\rm new}(s))
       +B(s+\alpha(s),\energyvar_{\rm new}(s))\radsrc(s)\bigr]\dd s.
 \end{aligned}
\]
The new energy is used in the phase equation. We use the norm
$\|\alpha\|_\infty+\|\energyvar\|_\infty+\|\energyvar\|_{L^1}$ on this product
space. Bounded phase shifts preserve the reference weights, and
Fubini's theorem gives
\[
 \|\energyvar_{\rm new}\|_\infty+\|\energyvar_{\rm new}\|_{L^1}
 \le C\eps^{-2}\int_0^\infty(1+s)v|\radsrc|\dd s
 \le C \|\radsrc\|_{\mathcal Y_{\rm rad}}.
\]
Also $\eps^2/v^2\le C$ and
$\int |\radsrc|/v\le\eps^{-2}\int v|\radsrc|$.
The coefficient bounds therefore give
$\|\alpha_{\rm new}-\tau\|_\infty
\le C \|\radsrc\|_{\mathcal Y_{\rm rad}}$.

For two input pairs, the difference of the new energies in
$C_{\rm b}\cap L^1$ is at most
$C\|\radsrc\|_{\mathcal Y_{\rm rad}}$ times their distance.
In the phase equation, the change of its first argument is
bounded by $C|\energyvar_{\rm new}||\Delta\alpha|$, using the first line
of \eqref{eq:radial-drift-derivatives}; its integral has the same
small factor. The change of the energy argument is controlled by
$\|\Delta \energyvar_{\rm new}\|_{L^1}$. The terms containing $B$ retain
$|\radsrc|/v$, and the initial phase contributes
$C\delta^2|\Delta \energyvar_{\rm new}(0)|$. Thus the displayed map is a
contraction on a ball of radius
$C\|\radsrc\|_{\mathcal Y_{\rm rad}}$ about $(\tau,0)$.
Its fixed point satisfies \eqref{eq:radial-chart}, $V(0)=m_0$,
and $\energyvar(t)\to0$. The phase derivative is integrable, so
$\alpha(t)$ has a finite limit.

\medskip\noindent\textit{Regularity and dependence on the data.}
The preceding integral map is $C^1$ on the indicated product
space. This follows from \eqref{eq:radial-drift-derivatives} and
the same integrable bounds used for its differences. Differentiating
the contraction gives bounded derivatives in $\radsrc$ and
$m_0/v(0)$. At fixed $\radsrc$, subtracting the derivative of
$(\tau,0)$ leaves a factor $\radsrc$ or $\energyvar$ in every integral.
Thus the derivative of $(\alpha-\tau,\energyvar)$ in the scaled boundary
value is $O(\|\radsrc\|_{\mathcal Y_{\rm rad}})$ in this product
norm.

The stronger weights and derivatives follow from the equations,
after the fixed point has been obtained. In fact, its energy
integral and the definition of the source norm give
\begin{equation}
 \begin{aligned}
 &\sum_{j=0}^5\|(1+t)^6\partial_t^j\energyvar\|_{C_{\rm b}^{0,\alphaH}}
       +\int_0^\infty(1+t)^5|\energyvar(t)|\dd t\\
 &\quad+\sum_{j=1}^5
          \|(1+t)^6\partial_t^j\alpha\|_{C_{\rm b}^{0,\alphaH}}
       +\|\alpha-\tau\|_\infty
 +\|\eps^{-1}(1+t)^3v\alpha'\|_{L^2}
 \le C \|\radsrc\|_{\mathcal Y_{\rm rad}}.
 \end{aligned}
 \label{eq:radial-terminal-estimate}
\end{equation}
For $\energyvar$ without derivatives, use Fubini on its integral formula.
For the remaining terms, differentiate \eqref{eq:radial-chart}
successively through order four. Each derivative of the energy
equation retains a source derivative multiplied by $O(v/\eps^2)$.
Each derivative of the phase equation retains either $\energyvar$ or one
of its derivatives, or a source derivative multiplied by $O(1/v)$.
The bounded logarithmic derivatives of $v$ and the inequalities
$\eps^2/v^2\le C$ close this induction, including the H\"older
seminorms. The same proof gives \eqref{eq:radial-terminal-estimate}
after a scaled boundary derivative, with the same right side, and
after a source derivative, with
$C\|\dot{\radsrc}\|_{\mathcal Y_{\rm rad}}$ on the right.
It also proves $C^1$ dependence of the limiting phase, by integrating
the differentiated phase equation.

The chart bounds transfer these estimates to $V/v$ and $F_V/F$.
The identity $V'=Q^W(t+\alpha,\energyvar)$ supplies the sixth derivative
of $V$ from the five derivatives of $\alpha,\energyvar$; the chart gives
five derivatives of $F_V$. At zero source, the constructed pair
is $(\alpha,\energyvar)=(\tau,0)$ and $V=v_\tau$. Taking differences in
the chart, including a scaled boundary derivative at fixed source,
proves \eqref{eq:radial-background-first-difference}. The
linearized contraction and the same chart estimates prove the
stated source derivative bound.

To prove the $L^2$ term in \eqref{eq:radial-terminal-estimate},
use $v^2|\alpha'|\le C\eps^2|\energyvar|+Cv|\radsrc|$. Its squared
norm is bounded by
\[
 C\int_0^\infty(1+t)^6\energyvar^2\dd t
  +C\eps^{-2}\int_0^\infty(1+t)^6\radsrc^2\dd t
 \le C\|\radsrc\|_{\mathcal Y_{\rm rad}}^2.
\]
The first term uses the energy estimates just proved; for the
second, write
$E=\eps^{-2}(1+t)^6v\radsrc$ and use
$\eps^2/v^2\le C$ to bound it by $C\|E\|_\infty\|E\|_{L^1}$.
The same argument applies after one scaled boundary derivative.
For a source variation it gives
\[
 \left\|\eps^{-1}(1+t)^3v\,
 D_{\radsrc}\alpha'[\dot{\radsrc}]\right\|_{L^2}
 \le C\|\dot{\radsrc}\|_{\mathcal Y_{\rm rad}}.
\]
\end{proof}

\section{Local families in higher dimensions}\label{sec:bridge}
\subsection{The coupled inner-end equations}\label{subsec:tail-equations}

We first solve the Yamabe equation on the inner Fowler end. We then
choose the boundary values at $|X|=b\delta$ and $|X|=b/\delta$
so that the solutions in the three regions have matching normal
derivatives.

On the inner end we use the cylinder
$\mathcal C=[0,\infty)\times\Sph^{n-1}$, with
$r=|y|=\delta e^{-t}$. Keep the neck $\eps$ fixed. Prescribe a phase $\tau$ in the
bounded interval of Section~\ref{sec:inner-tail}, a Kelvin parameter
$\kelpar_0$, and a higher-harmonic boundary value
$w(0,\theta)=\delta v(0)\phi(\theta)$.
Set $m_0=v_\tau(0)$. In this subsection the finite parameters are
$\vartheta=(\ell,\tau,\kelpar_0)$, where $\ell=\log b$ is the
logarithmic scale, $\tau$ is the Fowler phase just specified, and
$\kelpar_0$ is the initial Kelvin parameter. The trace $\phi$ is
a separate function-valued datum. The unknowns are a radial forcing
$\radsrc:[0,\infty)\to\R$, a function $w$ on $\mathcal C$
with $\Pi_{\ge2}w=w$, and a function
$z:[0,\infty)\to\R^n$ with $z(0)=0$. They determine
\[
 V=R_\eps(m_0,\radsrc),\qquad U=V+w,\qquad
 \kelpar(t)=\kelpar_0+z(t).
\]
Here $R_\eps(m_0,\radsrc)$ is the radial solution operator of
Lemma~\ref{lem:radial-response}; its unforced value is
$R_\eps(m_0,0)=v_\tau$. The varying Kelvin parameter
$\kelpar(t)$ accounts for the first spherical harmonics;
$\zeta$ remains the fixed parameter of the middle bubble.
We reconstruct from these functions a conformal factor relative to
$h_b$ and choose $(\radsrc,w,z)$ so that it solves the Yamabe equation.
The radial, first-harmonic and higher-harmonic solution operators
act on the fixed spaces of Section~\ref{sec:high-local}. We retain the divergence form of the
center equation in order to use the capacity estimate
\eqref{eq:capacity-bound}.

\medskip\noindent\textit{Coordinates and the equation to be solved.}
To separate the first variation in the metric from its quadratic
remainder, use $h_\upsilon=\exp(\upsilon H_b)$,
$0\le\upsilon\le1$; thus $h_0=g_{\rm E}$ and $h_1=h_b$.
For a constant Kelvin parameter $\kelpar$, set
\[
 E_\kelpar=1+2\kelpar\cdot y+|\kelpar|^2r^2,\qquad
 F_\kelpar(y)=\frac{y+\kelpar r^2}{E_\kelpar},\qquad
 \overline H_\kelpar
 =E_\kelpar^2(DF_\kelpar)^TH_b(F_\kelpar)DF_\kelpar.
\]
The map $F_\kelpar$ fixes the puncture. After this change of
coordinates, a further coordinate change and conformal factor remove
the metric terms that would force the radial and first-harmonic
equations at first order. Use the vector field $X_\kelpar$ and scalar
$\sigma_\kelpar$ constructed in Appendix~\ref{app:coordinates}, and put
\[
 T_\kelpar=\overline H_\kelpar-\mathcal L_{X_\kelpar}g_{\rm E}
                         -\sigma_\kelpar g_{\rm E}.
\]
This tensor is trace-free and satisfies, at every radius,
\[
 \fint (T_\kelpar)_{rr}=0,\qquad
 \fint (T_\kelpar)_{rr}\theta_i=0,\qquad
 \fint r^{-1}(T_\kelpar)_{rA}\nabla^A\theta_i=0
 \quad(1\le i\le n).
\]
Here the averages, angular indices and derivatives use the unit-sphere
metric. Moreover, $X_\kelpar=O(\lambda r^{q+1})$ and
$\sigma_\kelpar=O(\lambda r^q)$.
Let $Q_{\upsilon,\kelpar}$ be the flow at time one of
$-\upsilon X_\kelpar$. For the varying parameter $\kelpar(t)$, define
\begin{equation}
\begin{aligned}
 \Phi_\upsilon(t,\theta)
 &=(F_{\kelpar(t)}\circ Q_{\upsilon,\kelpar(t)})(r\theta),\\
 C_\upsilon(t,\theta)
 &=e^{-\upsilon\sigma_{\kelpar(t)}(r\theta)/2}
   E_{\kelpar(t)}\bigl(Q_{\upsilon,\kelpar(t)}(r\theta)\bigr).
\end{aligned}
\label{eq:finite-gauge-maps}
\end{equation}
Define the conformal factor $u_\upsilon$ relative to $h_\upsilon$
and the metric $\gamma_\upsilon$ on $\mathcal C$ by
\[
 u_\upsilon\circ\Phi_\upsilon
       =\left(\frac{C_\upsilon}{r}\right)^{(n-2)/2}U,
 \qquad
 \gamma_\upsilon
       =\left(\frac{C_\upsilon}{r}\right)^2\Phi_\upsilon^*h_\upsilon.
\]
Conformal covariance gives the transformed equation
\[
 \mathcal E_\upsilon(U)
 =-\Delta_{\gamma_\upsilon}U
  +c(n)\operatorname{Scal}_{\gamma_\upsilon}U
  -\frac{n(n-2)}4U^p.
\]
Consequently, the reconstructed factor $u_1$ solves the Yamabe equation
for $h_b$ precisely when $\mathcal E_1(V+w)=0$.

\medskip\noindent\textit{The range of the unknowns.}
Recall $F_V=\frac{n-2}{2}V-V'$, the unforced functions
$v_\tau,F_\tau$, and the fixed norm weights $v,F$ from
Section~\ref{sec:inner-tail}.
Fix a small constant $\eta_0>0$, depending only on the dimension
and the fixed parameter ranges. We work where
\begin{equation}
\begin{gathered}
 \frac12\le\frac V{v_\tau}\le2,\qquad
 \frac12\le\frac{F_V}{F_\tau}\le2,\qquad
 |w|\le\frac V2,\\
 r(|\kelpar|+|\kelpar'|)+\frac{|D^{\le1}w|}{F_V}
                  +\lambda r^q\frac V{F_V}\le\eta_0.
\end{gathered}
 \label{eq:inner-positive-neighborhood}
\end{equation}
Lemma~\ref{lem:radial-response} gives the first two bounds on a
small radial source ball, uniformly for the stated phase interval.
Equation~\eqref{eq:fowler-reference-multipliers} then gives
$V\asymp v$ and $F_V\asymp F$ for the fixed norm weights.
The other bounds follow from
\[
\begin{gathered}
 \frac{|w|}{V}+\frac{|D^{\le1}w|}{F_V}
 \le C\delta\mathcal W(w),\qquad
 r|\kelpar'|\le C\delta K(z),\qquad
 r|\kelpar|\le\delta(|\kelpar_0|+K(z)),\\
 \lambda r^qV/F_V\le C\lambda\delta^{q-2}\longrightarrow0.
\end{gathered}
\]
The parameters $\kelpar_0$ range over a fixed compact set.
By \eqref{eq:weight-bounds}, first choose small bounds for
$\|\radsrc\|_{\mathcal Y_{\rm rad}}$, $\mathcal W(w)$ and $K(z)$,
and then decrease $\eps$. These choices give
$U\ge V/2$, $F_V>0$ and $1+r\theta\cdot\kelpar'\ge1-\eta_0$.
When expressed in the $y$ coordinates, the differential of
$\Phi_\upsilon$ is
$I+O(r|\kelpar|+r|\kelpar'|+\lambda r^q)$, uniformly in $\upsilon$.
For small $\eta_0$, this map extends across zero to a bi-Lipschitz
map of the ball onto its image. Thus each such triple defines
$u_\upsilon$ on $\Phi_\upsilon(\mathcal C)$.

\medskip\noindent\textit{Separating the radial, first and higher harmonics.}
With Euclidean derivatives in $y$ and $\kelpar$ held fixed, define
\begin{equation}
 \mathcal M_{\rm stat}U=r^{(n+2)/2}\left[
 \operatorname{div}\bigl(T_\kelpar\nabla(r^{-(n-2)/2}U)\bigr)
 +c(n)\operatorname{div}\operatorname{div}T_\kelpar\,r^{-(n-2)/2}U\right],
 \label{eq:repair-static-metric}
\end{equation}
Put $\eta=r\theta\cdot\kelpar'$, $L=1+\eta$,
$\beta=\nabla_\theta\eta$, and $\mathscr D_t=\partial_t-\beta\cdot\nabla_\theta$.
The flat moving operator and its part linear in the center derivatives $\kelpar'$ and $\kelpar''$ are
\begin{align}
 \mathcal E_0(U)&=-L^{-1}\left\{\partial_t(\mathscr D_t U/L)
       +\operatorname{div}_\theta(L\nabla U-\beta\mathscr D_t U/L)\right\}
       +c(n)\operatorname{Scal}_{\gamma_0}U-\frac{n(n-2)}4U^p,\label{eq:repair-flat-moving}\\
 \mathcal TU&=2\eta U_{tt}+2\beta\cdot\nabla U_t
       +(\eta_t-(n-1)\eta)U_t+(\beta_t-\beta)\cdot\nabla U
       +\frac{n-2}{2}(\eta-\eta_t)U,\label{eq:repair-linear-motion}
\end{align}
where
\[
 \gamma_0=L^2\dd t^2+|\dd\theta+\beta\dd t|^2,\qquad
 c(n)\operatorname{Scal}_{\gamma_0}=\frac{(n-2)[n-2+2(n-1)\eta]}{4L^2}
                   -\frac{n-2}{2}\frac{\eta_t-|\beta|^2}{L^3}.
\]
The quadratic motion remainder is
\begin{equation}
 \mathcal R_{\rm mov}U
 =\mathcal E_0(U)-\mathcal E_{\rm flat}(U)-\mathcal TU,
 \label{eq:repair-motion-remainder}
\end{equation}
where $\mathcal E_{\rm flat}(U)=-U_{tt}-\Delta_\theta U+\frac{(n-2)^2}{4}U-\frac{n(n-2)}4U^p$.
The identical power terms cancel in the difference, making
$\mathcal R_{\rm mov}$ linear in $U$.

The remaining first metric variation and the second metric remainder
are, respectively,
\begin{equation}
\begin{aligned}
 \mathcal M_{\rm mov}U&=\left.\partial_\upsilon
                         \mathcal E_\upsilon(U)\right|_0
                                  -\mathcal M_{\rm stat}U,\\
 \mathcal R_{\rm met}U&=\int_0^1(1-\upsilon)
                       \partial_\upsilon^2\mathcal E_\upsilon(U)\dd\upsilon .
\end{aligned}
 \label{eq:repair-metric-remainders}
\end{equation}
Both are linear in $U$, since the nonlinear term $-\frac{n(n-2)}4U^p$ is
independent of the metric homotopy.
These homotopy derivatives use the actual center $\kelpar(t)$.
Spatial differentiation includes its derivatives; parameter
differentiation holds the reference weights fixed.
\begin{lemma}\label{lem:tail-gauge}
For every radial $V$,
\begin{equation}
 \Pi_0\mathcal M_{\rm stat}V=\Pi_1\mathcal M_{\rm stat}V=0.
 \label{eq:repair-low-selection}
\end{equation}
The metric operators satisfy, through four cylinder derivatives,
\begin{equation}
\begin{aligned}
\mathcal M_{\rm stat}U&=\sum_{m=0}^2O(\lambda r^q)D^mU,\\
\mathcal M_{\rm mov}U&=\sum_{j=1}^2\sum_{m=0}^2rz^{(j)}O(\lambda r^q)D^mU,\\
\mathcal R_{\rm met}U&=\sum_{m=0}^2O(\lambda^2r^{2q})D^mU.
\end{aligned}
\label{eq:coefficient-inventory}
\end{equation}
Each $O(\lambda r^q)$ or $O(\lambda^2r^{2q})$ denotes a
coefficient tensor contracting $D^mU$. After the displayed powers
of $\lambda$ and $r$, and the factors $rz^{(j)}$, are extracted,
the remaining coefficient tensors and their first parameter
derivatives are uniformly bounded in $C_{\rm b}^{4,\alphaH}$.

\end{lemma}
\begin{proof}
The scalar sector calculation and coefficient estimates are proved
in Appendix~\ref{app:coordinates}; see
\eqref{eq:scalar-gauge-euler}, \eqref{eq:metric-tensor-green-identity}
and \eqref{eq:finite-gauge-coefficient-bounds}.
\end{proof}

For the power remainder, Taylor's formula gives
\begin{equation}
 N_V(w)=(V+w)^p-V^p-pV^{p-1}w
 =p(p-1)w^2\int_0^1(1-s)(V+sw)^{p-2}\dd s.
 \label{eq:power-remainder-integral}
\end{equation}
By \eqref{eq:inner-positive-neighborhood}, $V+sw\ge V/2>0$ for
$0\le s\le1$. Thus the formula applies for every $p>1$.

Put
\[
 \radsrc=V''-\frac{(n-2)^2}{4}V+\frac{n(n-2)}4V^p,\qquad
 A_V=-\partial_t^2-\Delta_\theta+\frac{(n-2)^2}{4}-\frac{n(n+2)}4V^{p-1},
\]
and define
\begin{equation}
 \mathcal R_{\rm end}=\mathcal Tw+\mathcal R_{\rm mov}(V+w)
 +\mathcal M_{\rm stat}w+\mathcal M_{\rm mov}(V+w)
 +\mathcal R_{\rm met}(V+w)-\frac{n(n-2)}4N_V(w).
 \label{eq:repair-common-residual}
\end{equation}
The exact residual decomposition is
\[
 \mathcal E_1(V+w)=-\radsrc+A_Vw+\mathcal TV+\mathcal M_{\rm stat}V+\mathcal R_{\rm end}.
\]
The radial and high projections use
$\Pi_0\mathcal TV=\Pi_{\ge2}\mathcal TV=0$ and
\eqref{eq:repair-low-selection}; the center equation is kept in the
flux form below.

\medskip
Write $F_U=\frac{n-2}{2}U-U_t+\beta\cdot\nabla U$. The cylindrical Lagrangian depends on $U$ and its first derivatives:
\begin{equation}
 \mathcal L_{\rm can}=\frac{F_U^2}{2L}+\frac L2|\nabla U|^2
                  -\frac{(n-2)^2L}{8}U^{p+1}.
 \label{eq:cylindrical-lagrangian}
\end{equation}
Appendix~\ref{app:coordinates} differentiates $\mathcal L_{\rm can}$
with respect to $\kelpar'$. The normalized center momentum is
$\mathcal J_i=n\fint (\partial_{\kelpar_i'}\mathcal L_{\rm can})$, namely
\begin{equation}
 \mathcal J_i(U,\kelpar')=nr\fint \left[-\frac{\theta_iF_U^2}{2L^2}
 +\frac{F_U}{L}\nabla\theta_i\cdot\nabla U
 +\theta_i\left(\frac{|\nabla U|^2}{2}-\frac{(n-2)^2}{8}U^{p+1}\right)\right].
 \label{eq:repair-momentum}
\end{equation}
For a radial field, put
\[
 j_V=rF_V,\qquad \mathsf p_V=j_V^2.
\]
The center momentum then takes the form
$\mathcal J(V,\kelpar')=\mathsf p_VC(\kelpar')\kelpar'$, where
\begin{equation}
 C(\kelpar')=n\fint \left[\theta\otimes\theta
          \int_0^1(1+sr\theta\cdot\kelpar')^{-3}\dd s\right]
      =I+O(r^2|\kelpar'|^2).
 \label{eq:center-momentum-matrix}
\end{equation}
Here $cI\le C(\kelpar')\le CI$ on the prescribed neighborhood.
To isolate $\mathcal J(V+w,0)$, define
\[
 J_w=\mathcal J(V+w,0),\qquad
 \mathsf p_VA_w(\kelpar')=\int_0^1D_{\mathrm{vel}}
            [\mathcal J(V+w,s\kelpar')-\mathcal J(V,s\kelpar')]\dd s.
\]
Here $D_{\mathrm{vel}}$ differentiates the second argument of $\mathcal J$ at $s\kelpar'$.
The fundamental theorem of calculus gives
\begin{equation}
 \begin{aligned}
 \mathcal J(V+w,\kelpar')&=\mathsf p_V(C(\kelpar')+A_w(\kelpar'))\kelpar'+J_w,\\
 |J_w|&\le Cr|D^{\le1}w|^2,\qquad
 |A_w|\le C\left(\frac{|D^{\le1}w|}{F_V}
                       +\frac{|D^{\le1}w|^2}{F_V^2}\right).
 \end{aligned}
 \label{eq:repair-momentum-split}
\end{equation}
These bounds make $C(\kelpar')+A_w(\kelpar')$ uniformly invertible
on the chosen neighborhood. Appendix~\ref{app:coordinates} proves the
bounds, including the absence of terms linear in $w$ in $J_w$ and
linear in $\kelpar'$ in $C(\kelpar')-I$.

The flat center tangent is
\begin{equation}
 Z_i(U)=r\left(\frac{\theta_iF_U}{L}-\nabla\theta_i\cdot\nabla U\right).
 \label{eq:physical-center-variation}
\end{equation}
The identity obtained by varying the center for arbitrary states is
\begin{equation}
 -\mathcal J_i'=n\fint [L\mathcal E_0(U)Z_i(U)].
 \label{eq:repair-flat-ward}
\end{equation}
For high $w$, Appendix~\ref{app:coordinates} gives
$\fint (LZ_i(V+w))=0$. The first harmonic pairing matrix
\begin{equation}
 (G_{\rm lin})_{ij}=\frac n{j_V}\fint [LZ_i(V+w)\theta_j]
       =\delta_{ij}+O(|D^{\le1}w|/F_V),\qquad
 \|G_{\rm lin}-I\|\le C\delta\mathcal W(w).
 \label{eq:linear-harmonic-pairing}
\end{equation}
Choose the product neighborhood so that $\|G_{\rm lin}-I\|\le1/2$;
then $\|G_{\rm lin}^{-1}\|\le2$.

\medskip
The metric contribution to the center equation has components
\[
 (\mathcal B_H)_i=n\fint \bigl[L(\mathcal E_1-\mathcal E_0)(V+w)Z_i(V+w)\bigr],
 \qquad 1\le i\le n.
\]
The difference $\mathcal E_1-\mathcal E_0$ includes the finite conformal
normalization, since $h_0=g_{\rm E}$ and $Q_0=\operatorname{Id}$.
We solve the system
\begin{equation}
 \Pi_0\mathcal E_1=0,\qquad \Pi_{\ge2}\mathcal E_1=0,\qquad
 -\mathcal J'+\mathcal B_H=0.
 \label{eq:repair-final-system}
\end{equation}
For $1\le i\le n$, the $i$th component of $-\mathcal J'+\mathcal B_H$
is $n\fint [L\mathcal E_1(V+w)Z_i(V+w)]$ by \eqref{eq:repair-flat-ward}.
To display every contribution, set
\[
 Q_i(w)=\theta_i(\frac{n-2}{2}w-w_t+\beta\cdot\nabla w)-L\nabla\theta_i\cdot\nabla w.
\]
Then
\begin{align}
 \mathcal B_H={}&j_V[\Pi_1\{\mathcal M_{\rm stat}w+
            \mathcal M_{\rm mov}(V+w)+\mathcal R_{\rm met}(V+w)\}]
 \nonumber\\
 &+nr\fint \left[Q(w)\{\mathcal M_{\rm stat}V+\mathcal M_{\rm stat}w+
            \mathcal M_{\rm mov}(V+w)+\mathcal R_{\rm met}(V+w)\}\right].
 \label{eq:repair-force-ledger}
\end{align}
Apply the Kelvin inverse with the current coefficient
$\mathsf p_\tau$ to
\begin{equation}
 \begin{aligned}
 -(\mathsf p_\tau z')'&=f+(\mathsf p_\tau\mathsf j)',
                       \qquad f=-\mathcal B_H,\\
 \mathsf p_\tau\mathsf j
 &=\{\mathsf p_V(C+A_w)-\mathsf p_\tau\}z'+J_w.
 \end{aligned}
 \label{eq:repair-final-flux}
\end{equation}
When $\Pi_0\mathcal E_1=\Pi_{\ge2}\mathcal E_1=0$,
\eqref{eq:repair-flat-ward} gives
\[
 -\mathcal J'+\mathcal B_H=j_VG_{\rm lin}[\Pi_1\mathcal E_1],\qquad
 j_V>0,\quad \|G_{\rm lin}^{-1}\|\le2.
\]
Thus \eqref{eq:repair-final-system} is equivalent to
$\mathcal E_1(V+w)=0$; \eqref{eq:repair-momentum-split} gives its
center equation \eqref{eq:repair-final-flux}.

We also retain the boundary terms needed to compare actions.
Let $\mathcal L_\upsilon(U)$ be the density of the pulled-back action,
with respect to $\dd t\,\dd\theta$:
\[
\begin{aligned}
 \mathcal L_\upsilon(U)\,\dd t\,\dd\theta
 =\Phi_\upsilon^*\Bigl[\Bigl\{
 &\tfrac12\bigl(|\nabla_{h_\upsilon}u_\upsilon|^2
       +c(n)\operatorname{Scal}_{h_\upsilon}u_\upsilon^2\bigr)\\
 &-\frac{(n-2)^2}{8}u_\upsilon^{p+1}
 \Bigr\}\dd V_{h_\upsilon}\Bigr].
\end{aligned}
\]
The identity \eqref{eq:kelvin-endpoint-identity} combines the
cylindrical and conformal boundary terms into
\begin{equation}
 B_{\rm phys}
 =\frac{n-2}{2}\int U^2\frac{1+r\kelpar\cdot\theta}{E_\kelpar}\dd\theta.
 \label{eq:repair-physical-primitive}
\end{equation}
$B_{\rm phys}$ is independent of $\kelpar'$.
Putting $\mathcal L_H=\mathcal L_1-\mathcal L_0$ gives
exactly
\[
 J_{h_b}(u_1;\Phi_1([0,T]\times\Sph^{n-1}))
 =\int_0^T\!\int(\mathcal L_{\rm can}+\mathcal L_H)\dd\theta \dd t
      +B_{\rm phys}(T)-B_{\rm phys}(0).
\]
The boundary terms from the coordinate and conformal changes are
included in $\mathcal L_1-\mathcal L_0$.

\subsubsection{The coupled nonlinear map and its fixed point}
\label{subsec:tail-contraction}
We combine the three solution operators into a map whose fixed points
solve \eqref{eq:repair-final-system}.
For fixed boundary data at $t=0$, let $V=R_\eps(m_0,\radsrc)$ and
evaluate every residual at $(\radsrc,w,z)$. Use the Kelvin data $(f,\mathsf j)$ of
\eqref{eq:repair-final-flux}. Let $w_{\rm cap}$ be the homogeneous
high lift of $v(0)\delta\phi$, where
$\|\phi\|_{C^{6,\alphaH}}\le\kappa$, and write
\[
 \mathcal L_{\rm pert}w
 =\frac{n(n+2)}4(V^{p-1}-v_\tau^{p-1})w-\mathcal Tw-\mathcal R_{\rm mov}w
       -\mathcal M_{\rm stat}w-\mathcal M_{\rm mov}w-\mathcal R_{\rm met}w.
\]
The coupled map is
\begin{equation}
\begin{aligned}
 T_\radsrc&=\Pi_0\mathcal R_{\rm end},\\
 T_w&=w_{\rm cap}+R_{H,\tau}\Pi_{\ge2}\bigl[
  -(\mathcal M_{\rm stat}+\mathcal M_{\rm mov}+\mathcal R_{\rm met})V
  -\mathcal R_{\rm mov}V+\mathcal L_{\rm pert}w+\frac{n(n-2)}4N_V(w)\bigr],\\
 (T_z)'&=\mathsf p_\tau^{-1}\int_t^\infty f(s)\dd s-\mathsf j,
 \qquad T_z(0)=0.
\end{aligned}
\label{eq:one-map}
\end{equation}
The residual decomposition gives the equation for $T_w$, using
\eqref{eq:repair-low-selection}. The map preserves
$w(0)=v(0)\delta\phi$, $z(0)=0$.

Use the product norm
\begin{equation}
 \|(\radsrc,w,z)\|_{\rm prod}
 =\max\{\|\radsrc\|_{\mathcal Y_{\rm rad}},\mathcal W(w),K(z)\}.
 \label{eq:inner-product-norm}
\end{equation}
\begin{lemma}\label{lem:tail-source-map}
On the product neighborhood \eqref{eq:inner-positive-neighborhood}, set
\[
 \mathcal G=\|\radsrc\|_{\mathcal Y_{\rm rad}},\qquad W=\mathcal W(w),\qquad Z=K(z).
\]
For bounded small $\mathcal G,W,Z$, the map \eqref{eq:one-map} is jointly
$C^1$ into $\mathcal Y_{\rm rad}\times\mathcal X_{\rm hi}
\times\mathcal K_{\rm tail}$, and $(f,\mathsf j)\in\mathcal Y_K$.

At fixed boundary data, \eqref{eq:radial-background-first-difference}
bounds the relative field and flux changes from $R_\eps(m_0,0)$ to $R_\eps(m_0,\radsrc)$
through orders six and five by $C\mathcal G$, including one derivative
in $(\ell,\tau,\kelpar_0)$ at fixed $\radsrc$. Put
\begin{equation}
 c_{\rm pr}=C\left(\lambda^2\delta^{2q}
          +\lambda\delta^{q+1}Z+\delta^2Z^2\right),\qquad
 c_{\rm f}=C\left(\lambda\delta^{q-1}W+
 \lambda\delta^{q-2}Z+\lambda^2\delta^{2q-3}\right).
 \label{eq:principal-small}
\end{equation}
Then
\begin{align*}
 \|T_\radsrc\|_{\mathcal Y_{\rm rad}}\le{}&C(\delta^{2q}+\delta^qW+W^2+Z^2)
                          +c_{\rm pr}\mathcal G,\\
 \mathcal W(T_w)\le{}&C(\kappa+\delta^q)
 +C(\mathcal G+\lambda\delta^q+\delta Z)W+C\delta W^2
 +C\delta^q\mathcal G+C\delta(1+\mathcal G)Z^2,\\
 K(T_z)\le{}&C[\delta^{2q}+\delta^{q-3/2}W+\delta^{q-1}Z+W^2]
 +C(\mathcal G+\delta W+\delta^2Z^2)Z+c_{\rm f}\mathcal G.
\end{align*}
With the prescribed boundary data fixed, the state derivative in
\eqref{eq:inner-product-norm} satisfies
\begin{equation}
 \|D_{(\radsrc,w,z)}T\|
 \le C(\mathcal G+c_{\rm pr}+c_{\rm f}
             +\delta^{q-3/2}+\lambda\delta^q+W+Z).
 \label{eq:source-contraction-bound}
\end{equation}
Four derivatives in the cylinder variables of these source terms use
at most six derivatives of the fields and four derivatives in $t$ of
$\radsrc$. This is distinct from the single Fr\'echet derivative with
respect to the source or boundary datum.
\end{lemma}
\begin{proof}
The radial estimate follows from \eqref{eq:radial-product-estimates},
the forced multiplier bound \eqref{eq:principal-small}, and
$\delta^qZ\le(\delta^{2q}+Z^2)/2$. Equation~\eqref{eq:repair-low-selection}
gives $\Pi_0\mathcal M_{\rm stat}V=\Pi_1\mathcal M_{\rm stat}V=0$.
To estimate $T_w$, apply Lemma~\ref{lem:high-tail} to
\eqref{eq:high-total-sources}, \eqref{eq:high-power-bound} and
\eqref{eq:high-motion-source}. In particular, the exact radial motion
remainder retains the factor $\delta$ in its full $\mathcal W_4$ norm.
To estimate $T_z$, combine the complete force estimates
\eqref{eq:kelvin-unforced-data}, \eqref{eq:forced-small} with the
flux bounds \eqref{eq:kelvin-quadratic-flux}--\eqref{eq:kelvin-linear-flux}
and Lemma~\ref{lem:tail-low}. The factor $\delta^{q-3/2}$ includes
both the capacity bound and the weighted $L^2$ terms in
\eqref{eq:kelvin-work} for a general high field.
The first product rule in Appendix~\ref{app:weighted} gives
\eqref{eq:source-contraction-bound} in
$\mathcal Y_{\rm rad}\times\mathcal X_{\rm hi}\times\mathcal K_{\rm tail}$,
including the relative-$F$ entries.
\end{proof}

The following construction uses the same coupled-map organization as
\cite[Proposition 4.5]{part2}.
\begin{proposition}\label{prop:inner-graph}
Let $n\ge8$ and $q\ge2$. Fix compact sets of Fowler phases and
Kelvin centers. Assume that the smooth metric family satisfies
\[
 |D^iH_b|\le C_i\lambda r^q,\qquad \lambda\le C\eps,
\]
with the same bounds through total parameter order one.
There are $\kappa_0,\eps_0>0$ with the following property. For every
$0\le\kappa\le\kappa_0$ and $0<\eps\le\eps_0$, each
$\phi\in C^{6,\alphaH}_{\ge2}(\Sph^{n-1})$ with
$\|\phi\|_{C^{6,\alphaH}}\le\kappa$ determines a unique fixed
point of \eqref{eq:one-map} in the closed set
\begin{equation}
 \mathbb B_\phi=\{(\radsrc,w,z):\|(\radsrc,w,z)\|_{\rm prod}\le C_*(\kappa+\delta^q),
       \quad w(0)=v(0)\delta\phi,\quad z(0)=0\}.
 \label{eq:ball}
\end{equation}
The associated conformal factor $u_1$ is positive and solves the
Yamabe equation with the prescribed puncture and limiting neck $\eps$.

The family is jointly Fr\'echet $C^1$ in all finite parameters and scaled boundary
variables, with uniformly bounded first derivatives. The limiting
phase and bounded Kelvin coefficient are determined by the construction.
\end{proposition}
\needspace{6\baselineskip}
\begin{proof}
\leavevmode\par
\noindent\textit{Step 1. The fixed point.}

Evaluate $T(0)$ in the ambient product space.
With the boundary data fixed, Lemma~\ref{lem:tail-source-map} gives
on the ambient norm ball containing $\mathbb B_\phi$,
\[
 \|T(0)\|_{\rm prod}\le C_0(\kappa+\delta^q),\qquad
 \sup\|DT\|_{{\rm prod}\to{\rm prod}}
 \le C_{C_*}(\kappa+\delta^{q-3/2}).
\]
For the second bound use $\mathcal G=\|\radsrc\|_{\mathcal Y_{\rm rad}}$,
$\lambda\le C\delta^{n/2}$ and
\eqref{eq:source-contraction-bound}. Choose $C_*>2C_0$ large enough that
$(0,w_{\rm cap},0)\in\mathbb B_\phi$, using \eqref{eq:high-cap-bound}.
Then choose $\kappa_0,\eps_0$ so that the derivative bound is at most $1/2$.
For every state in the ball,
\[
 \|T(\radsrc,w,z)\|_{\rm prod}
 \le C_0(\kappa+\delta^q)+\tfrac12\|(\radsrc,w,z)\|_{\rm prod}
 \le C_*(\kappa+\delta^q).
\]
The output of $T$ has $w(0)=v(0)\delta\phi$ and $z(0)=0$.
Thus $T$ preserves the nonempty closed set $\mathbb B_\phi$.
Banach's theorem gives its unique fixed point.

\needspace{4\baselineskip}
\noindent\textit{Step 2. Parameter derivatives and the physical field.}

On each parameter patch, replace the input $w$ of $T$ by
$w+w_{\rm cap}$ and subtract $w_{\rm cap}$ from its $w$ component.
Denote the resulting map by $\widehat T$; its high component has
zero trace at $t=0$.
Lemma~\ref{lem:tail-source-map} and the bound \eqref{eq:high-cap-bound} give
$\|D_e\widehat T\|_{\rm prod}\le C$ for every unit finite or scaled
trace direction. Since $\|D\widehat T\|\le1/2$,
\eqref{eq:fixed-point-parameter-calculus} gives joint Fr\'echet
$C^1$ dependence and uniform first derivative bounds on the stated
parameter and trace domains.

We also need small derivatives along shrinking families of boundary data.
Suppose the scaled high trace $\phi$ and its first derivatives have
size $O(\eta)$, while the finite background directions are bounded.
At $(\radsrc,w,z)=0$, the sources and their finite background derivatives
are $O(\delta^q)$ in $\mathcal W_4$ and $O(\delta^{2q})$ in
$\mathcal Y_{\rm rad}$ and $\mathcal Y_K$, by \eqref{eq:repair-low-selection} and
\eqref{eq:coefficient-inventory}.
The lift and its derivative are $O(\eta)$. For general states,
\eqref{eq:radial-background-first-difference} and the first product
rule applied to \eqref{eq:one-map} therefore give
\[
 \|D_e\widehat T\|_{\rm prod}
 \le C(\eta+\delta^q+\mathcal G+W+Z).
\]
The value bounds \eqref{eq:ball}, with $\kappa=C\eta$, and
$\|(I-D\widehat T)^{-1}\|\le2$ yield
\begin{equation}
 \max_{j=0,1}\left\{
 \mathcal W(D_e^jw)+K(D_e^jz)+\|D_e^j\radsrc\|_{\mathcal Y_{\rm rad}}\right\}
 \le C(\eta+\delta^{q-3/2})\longrightarrow0.
 \label{eq:small-cap-parameter-bounds}
\end{equation}
Here $D_e^0$ is the identity; $D_e$ follows the specified boundary
family. For arbitrary trace directions $h$,
$\|D_\phi(\radsrc,w,z)[h]\|_{\rm prod}\le C\|h\|_{C^{6,\alphaH}}$.
The smaller bound \eqref{eq:small-cap-parameter-bounds} uses the
small trace derivative of the specified family.

The equivalence following \eqref{eq:repair-final-flux} identifies
the fixed point with $\mathcal E_1(V+w)=0$. The bound $|w|/V\le C\delta W$ gives
positivity. As \(t\to\infty\), \(\mathsf h(V(t),V'(t))\to h_\varepsilon\) and the phase displacement \(\alpha(t)\) converges;
 the Kelvin inverse gives a bounded limiting
Kelvin coefficient. 
On overlapping parameter sets, the neck, the selected unforced
profile $v_\tau$, the physical data at the initial boundary and
the gauge agree.
The fixed norms are uniformly equivalent.
Enlarging the ball constants by this equivalence factor and then
decreasing $\kappa,\eps$ gives a common contraction neighborhood.
Uniqueness identifies $u_1$ and its first parameter derivatives.
Fix one real lift of the asymptotic phase on the outer parameter box.
\end{proof}

\subsection{From geometric end parameters to full boundary values}
\label{subsec:tail-boundary}

The inner family is initially described by phase, Kelvin, and
higher-harmonic data. To join it to the middle annulus, we instead need
to prescribe an arbitrary small perturbation of the full boundary value
on $|y|=\delta$ and determine the resulting normal derivative.
The next proposition makes this change of boundary coordinates by
extending the end across a collar.
All boundary norms below use unit-sphere coordinates. Joint $C^1$
bounds include the value and every unit finite, trace or field
direction, with reference scales and measures fixed on each parameter patch.

In the centered $y=X/b$ coordinates for $h_b$, put $r=|y|$.
Here $u(y)=b^{(n-2)/2}u_{\rm phys}(by)$ is the rescaled field.
Write $u=r^{-(n-2)/2}\mathsf V(t,\theta)$, $t=-\log(r/\delta)$.
The inner domain's outer normal is increasing $r$, so
\[
 \delta\partial_ru|_{r=\delta}
       =\delta^{-(n-2)/2}[-\frac{n-2}{2}\mathsf V-\mathsf V_t]_{t=0}.
\]
For $\mathsf V=v$, the expression is $-\delta^{-(n-2)/2}G(0)$.
Freeze $v(0)$ and the two reference scales
\[
 e_0=v(0)\delta\asymp\eps,\qquad
 s_\eps=\delta^{-(n-2)/2}e_0=\delta^{(4-n)/2}v(0)\asymp\delta .
\]
Use $U=v_\tau$ and the fixed center $\kelpar_0$ in
\eqref{eq:finite-gauge-maps} at $\upsilon=1$ to define the reference
physical field $u_\vartheta^{\rm ref}$.
The profile uses $t=-\log(r/\delta)$.
Set $f_{\rm bg}(\vartheta)=u_\vartheta^{\rm ref}|_{r=\delta}$ and
prescribe the physical trace by
\[
 u_{\vartheta,\phi}|_{r=\delta}
       =f_{\rm bg}(\vartheta)+s_\eps\phi.
\]
Initially, $\phi\in C^{6,\alphaH}_{\ge2}(\Sph^{n-1})$.
Here $\vartheta=(\ell,\tau,\kelpar_0)$ specifies the reference
scale, Fowler phase, and Kelvin parameter; the full trace increment
$\phi$ is varied independently.
With the same normalization for trace and scaled normal derivative,
\begin{equation}
 \widetilde N(\vartheta,\phi)
 =s_\eps^{-1}[\delta\partial_ru_{\vartheta,\phi}
                    -\delta\partial_ru_{\vartheta,0}]_{r=\delta}
 =e_0^{-1}(-\frac{n-2}{2}-\partial_t)
             (\mathsf V_{\vartheta,\phi}-\mathsf V_{\vartheta,0})|_0.
 \label{eq:dn-normalization}
\end{equation}
For the regular harmonic lift $e_0e^{-(2l+n-2)t/2}\phi_l$, the operator
$e_0^{-1}(-\frac{n-2}{2}-\partial_t)$ gives $D_+\phi_l=l\phi_l$.
The term $-\frac{n-2}{2}$ comes from differentiating $r^{-(n-2)/2}$.
In normal coordinates for $h_b$, the scalar radial derivative is
the unit normal derivative. Under $X=by$, trace and scaled normal
derivative both acquire $b^{-(n-2)/2}$, which cancels in
\eqref{eq:dn-normalization}. Conormal densities additionally carry
the same physical surface density on both sides of each cut.

\begin{proposition}\label{prop:inner-dn}
The family of Proposition~\ref{prop:inner-graph} extends across a
fixed collar to $|y|=\delta$. For each reference background $\vartheta$,
the full trace determines the phase and Kelvin increments and the
higher-harmonic boundary data. The resulting map \eqref{eq:dn-normalization}
is jointly Fr\'echet $C^1$ on the full
trace ball $B_\Xtrace(0,\kappa_0)$, where
$\Xtrace=C^{6,\alphaH}(\Sph^{n-1})$ has no projection restriction.
Uniformly for $0\le\kappa\le\kappa_0$ and $0<\eps\le\eps_0$,
\[
 \begin{gathered}
 \|\widetilde N(\vartheta,\phi)-D_+\phi\|_
       {C^1(\vartheta,\phi;C^{5,\alphaH})}
 \le C\delta,\qquad D_+\phi_l=l\phi_l.
 \end{gathered}
\]
Here $\|\phi\|_\Xtrace\le\kappa$, and $\widetilde N(\vartheta,0)=0$.
\end{proposition}
\needspace{6\baselineskip}
\begin{proof}
Put $d=\delta/2$, $\widehat y=y/\delta$, and
\[
 \Ytrace=C^{5,\alphaH}(\Sph^{n-1}),\qquad
 \mathcal A=\{1/2<|\widehat y|<1\},\qquad
 \Zcollar=C^{6,\alphaH}(\overline{\mathcal A}).
\]
Subscripts $-$ and $+$ denote the boundaries $|\widehat y|=1/2$ and
$|\widehat y|=1$, respectively.
Apply the tail construction at $r=d$, using the same frozen scale
$s_\eps$ for trace and scaled conormal increments at both boundaries.
Retain their physical surface densities. At this boundary use
$t_d=-\log(r/d)$ and the shifted profile
$v_\tau(t_d+\log2)$. Thus $v_d=v_\tau(\log2)$ and
$v_d'=v_\tau'(\log2)$ belong to the same reference field
$u_\vartheta^{\rm ref}$. We write
$F_d=(n-2)v_d/2-v_d'=F_\tau(\log2)$.

\medskip
\noindent\textit{Step 1. The full trace and its tail conormal.}
Write $\trlow=(b_{{\rm tr},0},b_{{\rm tr},1},\ldots,b_{{\rm tr},n})$ for
the finite phase and Kelvin coordinates, with physical phase increment
$\delta b_{{\rm tr},0}$. Since $z(0)=0$, the initial boundary embedding depends
only on these finite coordinates. Pull all inner traces back by this
same embedding. For $h\in\Pi_{\ge2}\Xtrace$ the identity
$u\circ\Phi=mU$ then gives
$f_-=t_{\vartheta,\trlow}+s_\eps M_{\vartheta,\trlow}h$.
Here $t_{\vartheta,\trlow}$ is the trace for $h=0$. The cylindrical
high trace is $d^{(n-2)/2}s_\eps h$, so
$M_{\vartheta,\trlow}=d^{(n-2)/2}m|_{r=d}$ in these boundary
coordinates. The angular coordinates are fixed.
Finite parameter derivatives of $Mh$ act on $M$ and leave
the derivative order of $h$ unchanged.
The Fowler, Kelvin and gauge factors give
\begin{equation}
\begin{gathered}
 \|M-1\|_\Xtrace+\|D_{(\vartheta,\trlow)}M\|_\Xtrace
       +\|D_{\trlow}^2M\|_\Xtrace
       +\|D_\vartheta D_{\trlow}M\|_\Xtrace\le C\delta,\\
 D_{\trlow}(s_\eps^{-1}t)|_0
 =\left(\frac{\delta d^{-(n-2)/2}v_d'}{s_\eps},
        \frac{d^{(4-n)/2}F_d}{s_\eps}\theta_1,\ldots,
        \frac{d^{(4-n)/2}F_d}{s_\eps}\theta_n\right)+O_\Xtrace(\delta),\\
 \|D_{\trlow}^2(s_\eps^{-1}t)\|_\Xtrace\le C\delta,
 \qquad \|D_\vartheta D_{\trlow}(s_\eps^{-1}t)\|_\Xtrace\le C.
\end{gathered}
 \label{eq:collar-trace-coordinate-bounds}
\end{equation}
Here $v_d\asymp d^{(n-2)/2}$, $v_d'/v_d=-\frac{n-2}{2}+O(\delta^2)$ and
$F_d/v_d=(n-2)+O(\delta^2)$. Thus
$\delta d^{-(n-2)/2}v_d'/s_\eps$ and $d^{(4-n)/2}F_d/s_\eps$ are bounded
and bounded away from zero.

Choose the affine background from this same reference trace:
\[
 B_\vartheta(\widehat y)
 =\fint t_{\vartheta,0}+2[\Pi_1t_{\vartheta,0}]\cdot\widehat y,
 \qquad b_\vartheta=2\delta^{-1}[\Pi_1t_{\vartheta,0}].
\]
On the fixed collar, put $G_\tau=\frac{n-2}{2}v_\tau+v_\tau'$.
In the original variable $t=-\log(r/\delta)$, this collar is
$0\le t\le\log2$, where $v_\tau'(t)<0$.
The Hamiltonian identity therefore gives
\[
 \frac{G_\tau}{v_\tau}=O_{C^1_\vartheta}(\delta^2),\qquad
 r\partial_r(r^{-(n-2)/2}v_\tau)
       =-r^{-(n-2)/2}G_\tau.
\]
The radial physical factor therefore changes by $O_{C^1_\vartheta}(\delta^2)$
between $d$ and $\delta$. The Kelvin formula and the finite gauge give
\[
 \begin{aligned}
 t_{\vartheta,0}
 &=d^{-(n-2)/2}v_d[1+(n-2)d\kelpar_0\cdot\theta]
       +O_{C^1_\vartheta(\Xtrace)}(\delta^2),\\
 f_{\rm bg}
 &=d^{-(n-2)/2}v_d[1+(n-2)\delta\kelpar_0\cdot\theta]
       +O_{C^1_\vartheta(\Xtrace)}(\delta^2).
 \end{aligned}
\]
The gauge error is $O_{C^1}(\lambda\delta^q)$.
In particular, $\Pi_{\ge2}t_{\vartheta,0}=O_{C^1_\vartheta(\Xtrace)}(\delta^2)$
and $[\Pi_1t_{\vartheta,0}]=O_{C^1_\vartheta}(\delta)$.
The factor $2$ in $B_\vartheta$ matches the degree-one trace at
$|\widehat y|=1/2$. Consequently
\begin{equation}
 \|t_{\vartheta,0}-B_\vartheta|_-\|_{C^1_\vartheta(\Xtrace)}
 +\|f_{\rm bg}-B_\vartheta|_+\|_{C^1_\vartheta(\Xtrace)}\le C\delta^2.
 \label{eq:collar-affine-reference}
\end{equation}
After the same pullback, the reference scaled conormals have this error in
$C^1_\vartheta(\Ytrace)$ relative to
$\frac12\partial_{\widehat r}B_\vartheta|_-$ and
$\partial_{\widehat r}B_\vartheta|_+$, respectively.
This follows from the radial identity above and differentiation
of the Kelvin and gauge factors.
The lower bound for $d^{-(n-2)/2}v_d$ gives
$B_\vartheta(0)\ge2c_B>0$ for a fixed $c_B$.
Also $\|b_\vartheta\|_{C^1_\vartheta}\le C$, so
$B_\vartheta\ge c_B$ for small $\delta$.

On fixed small input balls, the trace map
\begin{equation}
 (\trlow,h)\longmapsto
 s_\eps^{-1}(t_{\vartheta,\trlow}-B_\vartheta|_-)
       +M_{\vartheta,\trlow}h
 \label{eq:tail-trace-coordinates}
\end{equation}
has derivative equal to its invertible diagonal low block and the
identity high block, up to $O(\delta)$. By
\eqref{eq:collar-affine-reference} and $s_\eps\asymp\delta$, its value
at $(0,0)$ is $O_{C^1_\vartheta(\Xtrace)}(\delta)$.
The quantitative inverse function theorem gives unique coordinates for every
$\|\trpar\|_\Xtrace\le r_*$, with
\[
 f_-=B_\vartheta|_-+s_\eps\trpar,\qquad
 |\trlow|+\|h\|_\Xtrace\le C(\|\trpar\|_\Xtrace+\delta),\qquad
 |D\trlow|+\|Dh\|_\Xtrace\le C.
\]
At fixed $\trpar$,
$|D_\vartheta\trlow|+\|D_\vartheta h\|_\Xtrace
\le C(\|\trpar\|_\Xtrace+\delta)$: the trace error at zero is $O_{C^1}(\delta)$,
and the mixed low coordinate derivatives are bounded.
Differentiating the explicit trace equation once in $\trpar$ and once
in $\vartheta$, or twice in $\trpar$, uses
\eqref{eq:collar-trace-coordinate-bounds} and the uniformly invertible
trace derivative. It gives
\[
 \|D_\vartheta D_\trpar(\trlow,h)\|
       +\|D_\trpar^2(\trlow,h)\|\le C,
\]
with the operator norms from $\Xtrace$ and $\Xtrace\times \Xtrace$ to
$\R^{n+1}\times\Pi_{\ge2}\Xtrace$, respectively.
Choose $r_*$ first and then $\eps_0$ so these coordinates remain
in the domain of Proposition~\ref{prop:inner-graph}.

After subtracting the conormal of $B_\vartheta$ and dividing by
$s_\eps$, the inner solution has conormal
$D_+\trpar+E(\vartheta,\trpar)$, with
$\|E(\vartheta,0)\|_{C^1_\vartheta(\Ytrace)}\le C\delta$.
The leading low coordinate block commutes with
$\operatorname{diag}(0,I_n)$, while the multiplier and nonlinear
coordinate errors are $O(\delta)$ by
\eqref{eq:collar-trace-coordinate-bounds}. Apply the
source estimates to the affine normalized conormal minus $D_+$
times the trace in \eqref{eq:tail-trace-coordinates}; denote that
residual in the old coordinates by $E^\sharp$.
The high, radial and Kelvin estimates
\eqref{eq:dn-high-products}--\eqref{eq:dn-source-gains},
\eqref{eq:dn-radial-gain} and
\eqref{eq:dn-kelvin-product}--\eqref{eq:actual-cap-flux}
control $E^\sharp$ and its first derivatives.
The coordinate composition gives
\[
\begin{aligned}
 E(\vartheta,\trpar)&=E^\sharp(\vartheta,\trlow(\vartheta,\trpar),h(\vartheta,\trpar)),\\
 \|E\|_{C^1}&\le C(1+\|D\trlow\|+\|Dh\|)\|E^\sharp\|_{C^1}.
\end{aligned}
\]
The bounded coordinate derivatives and the value bound at $\trpar=0$
therefore give, on the entire fixed trace ball,
\begin{equation}
 \|E\|_{C^1(\vartheta,\trpar;\Ytrace)}
 \le C\left(\delta+\lambda\delta^{q-2}+\delta^{q-1}\right)
 \le C\delta.
 \label{eq:fixed-tail-dn-small}
\end{equation}
Here the tail and coordinate derivatives are bounded, and
$\lambda\le C\delta^{n/2}$, $q\ge2$.

\medskip
\noindent\textit{Step 2. The collar and its outer conormal.}
Let
$\psi_-(\vartheta,\trlow,\theta)=\delta^{-1}\Phi_1(0,\theta)$
be the rescaled boundary parametrization of the tail constructed
at $r=d$. Extend this parametrization to $\mathcal A$ by
\[
 \Psi_{\vartheta,\trpar}(\widehat r\theta)
 =\widehat r\theta+\chi(\widehat r)
       [\psi_-(\vartheta,\trlow(\vartheta,\trpar),\theta)
                                                       -\tfrac12\theta],
\]
where $\widehat r=|\widehat y|\in(1/2,1)$, $\chi=1$ near $1/2$
and $\chi=0$ near $1$.
The gauge and coordinate bounds give
$\|\Psi-\operatorname{Id}\|_{C^1(C^{7,\alphaH})}\le C\delta$;
thus $\Psi$ is a diffeomorphism fixing a neighborhood of the physical
outer sphere. Write the pulled-back field as $B_\vartheta+s_\eps W$.
Let $\vartheta$ and $\trpar$ range over the fixed parameter set and trace
ball. Let $W$ belong to a fixed bounded ball in $\Zcollar$, and assume
$B_\vartheta+s_\eps W\ge c_B/2$. Equations
\eqref{eq:collar-nonlinear-map}--\eqref{eq:collar-raw-conormals} give
\[
\begin{gathered}
 -\Delta W+\mathcal Q(\vartheta,\trpar,W)=0,\qquad
 \|\mathcal Q\|_{C^1(\vartheta,\trpar,W;C^{4,\alphaH})}\le C\delta,\\
 \mathcal C_-=r_0\partial_{\widehat r}W|_-+\mathcal A_-,\qquad
 \mathcal C_+=\partial_{\widehat r}W|_++\mathcal A_+,\qquad
 \|\mathcal A_\pm\|_{C^1(\vartheta,\trpar,W;\Ytrace)}\le C\delta,
 \quad r_0=1/2.
\end{gathered}
\]
Both conormals use increasing radius, the physical surface density,
and subtraction of the Euclidean conormal of $B_\vartheta$ before division
by $s_\eps$. The spatial gradient of $B_\vartheta$ is
$O_{C^1_\vartheta}(\delta)$. Its contributions to $\mathcal A_\pm$
are therefore $O(\delta)$, including one parameter derivative. The $C^{4,\alphaH}$ source estimates use six derivatives
of $W$. The trace estimate then gives $C^{5,\alphaH}$ conormals.

Put $\mathcal B_-W=r_0\partial_{\widehat r}W|_--D_+\operatorname{tr}_-W$.
Let $G_D$ be the zero Dirichlet inverse of $-\Delta$ on $\mathcal A$.
Let $\mathsf E(\phi,\trpar)$ be the harmonic function on $\mathcal A$
with outer trace $\phi$ and inner trace $\trpar$. The radial harmonic basis $\widehat r^l,\widehat r^{2-n-l}$ gives the operator
$\mathsf L:\Xtrace\to \Ytrace$ defined by
\begin{equation}
 \mathsf L\trpar_l:=(r_0\partial_{\widehat r}\mathsf E(0,\trpar)|_--D_+\trpar)_l
       =-\frac{2l+n-2}{1-2^{-(2l+n-2)}}\trpar_l.
 \label{eq:collar-multiplier}
\end{equation}
Thus $\mathsf L^{-1}=-[2D_++n-2]^{-1}(I-2^{2-n}P_{1/4})$, where
$P_{1/4}$ is harmonic extension at radius $1/4$.
Equation~\eqref{eq:harmonic-matching-inverse}, with the low degrees
handled by $\frac1{n-2}$ and $\frac1n$, gives
$\|\mathsf L^{-1}\|_{\Ytrace\to \Xtrace}\le C$. Consequently
\begin{equation}
 \mathsf K(f,h)=G_Df+\mathsf E\left(0,
       \mathsf L^{-1}[h-r_0\partial_{\widehat r}G_Df|_-]\right)
 \label{eq:collar-fixed-inverse}
\end{equation}
is the unique inverse for $-\Delta u=f$, $u|_+=0$,
$\mathcal B_-u=h$, with
$\|\mathsf K(f,h)\|_\Zcollar\le C(\|f\|_{C^{4,\alphaH}}+\|h\|_\Ytrace)$
by substitution and Dirichlet Schauder estimates.
Regular harmonic extension $\mathsf H\phi$ has degree-$l$ term
$\widehat r^l\phi_l$, so $\mathcal B_-\mathsf H\phi=0$ and
$\|\mathsf H\phi\|_\Zcollar\le C\|\phi\|_\Xtrace$.

Matching the same physical densities at the inner interface gives
$\mathcal B_-W=E-\mathcal A_-$, with $\trpar=\operatorname{tr}_-W$.
The entire collar problem is therefore
\begin{equation}
 W=\mathsf H\phi+\mathsf K\bigl(
 -\mathcal Q(\vartheta,\operatorname{tr}_-W,W),
 E(\vartheta,\operatorname{tr}_-W)
       -\mathcal A_-(\vartheta,\operatorname{tr}_-W,W)\bigr).
 \label{eq:collar-single-equation}
\end{equation}
The resulting map for $u=W-\mathsf H\phi$ has joint
$C^1(\Zcollar)$ norm at most $C\delta$: $W=\mathsf H\phi+u$ and
$\trpar=\operatorname{tr}_-(\mathsf H\phi+u)$ are fixed bounded
affine maps. Choose $\kappa_0$ so that
$\|\operatorname{tr}_-\mathsf H\phi\|_\Xtrace\le r_*/3$ for
$\|\phi\|_\Xtrace\le2\kappa_0$, then $C_0>2C$, and finally $\eps$ small. The map preserves the ball $u|_+=0$,
$\|u\|_\Zcollar\le C_0\delta$, and its Lipschitz constant there is at most
$1/2$. On this ball, $\|\operatorname{tr}_-u\|_\Xtrace\le r_*/3$, so the
traces stay in the fixed positive domains. Banach's theorem gives a
unique fixed point $W=W(\vartheta,\phi)$.
Equation~\eqref{eq:fixed-point-parameter-calculus} gives
\begin{equation}
 \|W-\mathsf H\phi\|_{C^1(\vartheta,\phi;\Zcollar)}\le C\delta.
 \label{eq:collar-solution-c1}
\end{equation}
The outer conormal is $D_+\phi+\mathcal R_+$, where
\[
 \mathcal R_+=\partial_{\widehat r}(W-\mathsf H\phi)|_+
       +\mathcal A_+(\vartheta,\operatorname{tr}_-W,W).
\]
The trace bound and the full chain rule in the independent inputs
$(\vartheta,\trpar,W)$ prove
$\|\mathcal R_+\|_{C^1(\vartheta,\phi;\Ytrace)}\le C\delta$, since
$D\trpar,DW$ are bounded.
At the fixed outer sphere, conversion from density to scalar radial
derivative divides by a $\phi$-independent factor
$1+O_{C^1(C^{5,\alphaH})}(\lambda\delta^q)$; include it in
$\mathcal R_+$ without changing the estimate.

Finally write $f_{\rm bg}-B_\vartheta|_+=s_\eps\phi_\vartheta^{\rm ref}$.
Equation~\eqref{eq:collar-affine-reference} gives
$\|\phi_\vartheta^{\rm ref}\|_{C^1(\Xtrace)}\le C\delta$.
Decreasing $\eps$ keeps $\phi+\phi_\vartheta^{\rm ref}$ in
$B_\Xtrace(0,2\kappa_0)$ for $\|\phi\|_\Xtrace\le\kappa_0$.
An outer trace $f_{\rm bg}+s_\eps\phi$ determines the inner coordinate
\begin{equation}
 \widehat\trpar(\vartheta,\phi)
 =\operatorname{tr}_-W(\vartheta,\phi_\vartheta^{\rm ref}+\phi).
 \label{eq:collar-trace-transfer}
\end{equation}
Its finite tail coordinates are
$\trlow(\vartheta,\widehat\trpar(\vartheta,\phi))$.
Subtracting the conormal at the background trace $f_{\rm bg}$ gives
\[
 \widetilde N(\vartheta,\phi)=D_+\phi
 +\mathcal R_+(\vartheta,\phi+\phi_\vartheta^{\rm ref})
 -\mathcal R_+(\vartheta,\phi_\vartheta^{\rm ref}).
\]
$\widetilde N(\vartheta,0)=0$. The bounds for $\mathcal R_+$ and
$\phi_\vartheta^{\rm ref}$ give the estimate in Proposition~\ref{prop:inner-dn}.
The extension is positive on the collar because
$B_\vartheta+s_\eps W\ge c_B/2$. It agrees with the constructed
inner solution near the puncture, so the puncture, limiting phase,
and limiting Kelvin parameter are unchanged.
\end{proof}

\par\medskip

\subsection{The regular filling and the outer annulus}
\label{subsec:outer-annular-problem}

We compare the local singular solution with a smooth solution on the
whole ball having the same outer boundary values. This regular filling
will be used to compute their difference in flux and action.

Let $\mathcal B_s$ be a geodesic ball, with its metric $\gamma$
in a fixed smooth bounded family; use unit-ball rescaled norms.
\begin{lemma}\label{lem:regular-filling}
For $\varphi\ge c_*>0$ in a fixed bounded $C^{6,\alphaH}$ set
and sufficiently small $s$, the problem
\[
 (-\Delta_\gamma+c(n)\operatorname{Scal}_\gamma)R-\frac{n(n-2)}4R^p=0
       \quad\hbox{in }\mathcal B_s,\qquad
 R=\varphi\quad\hbox{on }\partial\mathcal B_s
\]
has a unique solution in the $C^{6,\alphaH}$ ball of radius $Cs^2$
about the Euclidean harmonic extension $P\varphi$. It is positive,
jointly $C^1$ in the metric and trace, and
$R-P\varphi=O_{C^1}(s^2)$ in $C^{6,\alphaH}$ at fixed $s$.
The scalar conormal
$\mathcal N_s^\gamma\varphi=(s\partial_{\nu_\gamma}R)|_{\partial\mathcal B_s}$,
pulled back to the unit sphere, differs from $D_+\varphi$ by $O(s^2)$ in
\[
 C^1(C^{6,\alphaH},C^{5,\alphaH}).
\]
For $-\Delta_{g_{\rm rd}}+\frac{n(n-2)}4-\frac{n(n+2)}4 u_\zeta^{p-1}$ on a small round
ball, the $O(s^2)$ operator bound relative to the harmonic
Dirichlet-to-Neumann map holds uniformly on compact $\zeta$ sets,
including its first $\zeta$ derivative.
\end{lemma}
\begin{proof}
Rescale normal coordinates by $x=sy$. With
$\gamma_s(y)=\gamma(sy)$, the metric and curvature bounds are
\[
 \|\gamma_s-g_{\rm E}\|_{C^{6,\alphaH}}
 +\|s^2\operatorname{Scal}_\gamma(s\cdot)\|_{C^{4,\alphaH}}
 \le Cs^2,
\]
including one metric derivative. Let $G_1$ be the zero Dirichlet
inverse of $-\Delta$ on $B_1$ and write $R=P\varphi+v$. Then
\[
 v=\mathcal Q_s(v):=G_1\left[
 (\Delta_{\gamma_s}-\Delta)(P\varphi+v)
 -s^2c(n)\operatorname{Scal}_\gamma(s\cdot)(P\varphi+v)
 +s^2\frac{n(n-2)}4(P\varphi+v)^p\right].
\]
Since $P\varphi\ge c_*$, \eqref{eq:relative-power-calculus} and
Dirichlet Schauder estimates give, on a fixed positive neighborhood,
\[
 \|\mathcal Q_s(0)\|_{C^{6,\alphaH}}
 +\|D_v\mathcal Q_s\|
 +\sup_i\|\partial_i\mathcal Q_s\|_{C^{6,\alphaH}}\le Cs^2,
\]
where $i$ ranges over unit metric and trace directions at fixed $v$.
The map $\mathcal Q_s$ preserves the ball of radius $C_0s^2$ for
fixed sufficiently large $C_0$. Its Lipschitz constant there is at most $1/2$.
At its fixed point,
\[
\begin{aligned}
 D_iv&=(I-D_v\mathcal Q_s)^{-1}\partial_i\mathcal Q_s,\qquad
 \|(I-D_v\mathcal Q_s)^{-1}\|\le2,\\
 \|v\|_{C^1(C^{6,\alphaH})}&\le Cs^2,\qquad R\ge c_*/2.
\end{aligned}
\]
The normal coefficients differ from their Euclidean values by
$O_{C^1(C^{5,\alphaH})}(s^2)$, so the trace estimate proves the scalar
conormal bound. The density version multiplies $\mathcal N_s^\gamma$
by the normalized surface Jacobian $1+O_{C^1(C^{5,\alphaH})}(s^2)$.

For the linearized bubble operator, replace the curvature and power
terms by $-V_{\zeta,s}(P\varphi+v)$, where
$\|V_{\zeta,s}\|_{C^1_\zeta(C^{4,\alphaH})}\le Cs^2$.
The Dirichlet inverse $G_1$ and
\eqref{eq:fixed-point-parameter-calculus} bound the conormal error
and its $\zeta$ derivative by $Cs^2\|\varphi\|_{C^{6,\alphaH}}$.
Linearity in $\varphi$ gives the boundary operator estimates.
\end{proof}

Apply Lemma~\ref{lem:regular-filling} in $B_\rho$ with outer trace
\[
 \varphi=A+\rho\mathbf v\cdot\theta+\rho^2\psi,\qquad
 A\in[c_*,C_*],\qquad \Pi_{\ge2}\psi=\psi,
\]
where $\mathbf v\in\R^n$ and $\psi$ range over fixed bounded sets.
The filling $R$ has center value
$R(0)=A+O_{C^1}(\rho^2)$.

For the outer annulus, work in the original coordinates $X$ and
write $r=|X|$. We solve on
$A_{r_+,\rho}=B_\rho\setminus\overline{B_{r_+}}$
with prescribed inner trace and outer trace
$A+\rho\mathbf v\cdot\theta+\rho^2\psi$. Use
\[
 H_+=A+B_0r^{2-n},\qquad A\in[c,C],\qquad
 B_0=b^{(n-2)/2}C_\infty(\zeta)\asymp b^{(n-2)/2}.
\]
Assume that, after pullback to the unit sphere, each prescribed
trace divided by the value of $H_+$ on that boundary belongs to
a fixed bounded $C^{6,\alphaH}$ set with a uniform positive
lower bound. For parameter-dependent traces, assume the same uniform
$C^{6,\alphaH}$ bounds for the first derivatives of these normalized
pullbacks. The coefficients $A$ and $B_0$ satisfy
$|D^jA|\le C$, $|D^jB_0|\le CB_0$, $j=0,1$, for logarithmic scale
and center derivatives. 

Let $G_D=(-\Delta)^{-1}_{A_{r_+,\rho},D}$ be the positive Euclidean
Green operator with both Dirichlet traces zero. For radial $Q\ge0$,
\[
 P_Q(r)=\frac1{n-2}\left[
 r^{2-n}\int_{r_+}^{r}s^{n-1}Q(s)H_+(s)\dd s
 +\int_r^\rho sQ(s)H_+(s)\dd s\right]
\]
satisfies $-\Delta P_Q=QH_+$ and is nonnegative at both boundaries.
Moreover,
$r^{2-n}s^{n-2}H_+(s)\le H_+(r)$ for $s\le r$ and
$H_+(s)\le H_+(r)$ for $s\ge r$.
The maximum principle therefore gives, also for signed nonradial $f$,
\begin{equation}
 |f|\le QH_+\quad\Longrightarrow\quad
 \left\|\frac{G_Df}{H_+}\right\|_\infty
 \le\frac1{n-2}\int_{r_+}^\rho sQ(s)\dd s.
 \label{eq:outer-green-comparison}
\end{equation}
The constant is independent of $A,B_0$ and $r_+$.

On fixed relative annuli $A_r$, truncated at a boundary when needed,
put
\begin{align*}
 Q(r)&=r^{q-2}+1+b^2r^{-4},\\
 \|u\|_{j,H}&=\sup_{r_+\le r\le\rho}
       H_+(r)^{-1}\|u\|_{C^{j,\alphaH}_{\rm sc}(A_r)},\\
 \|f\|_{4,QH}&=\sup_{r_+\le r\le\rho}
       [Q(r)H_+(r)]^{-1}\|f\|_{C^{4,\alphaH}_{\rm sc}(A_r)}.
\end{align*}
Reference weights are fixed on each parameter patch; nearby weights
are uniformly comparable with bounded first relative derivatives.
The scaled interior and Dirichlet Schauder estimates, followed by the
Green comparison, give
\begin{equation}
\begin{split}
 \|G_Df\|_{6,H}
 &\le C\left(\int_{r_+}^\rho rQ(r)\dd r
                +\sup_{r_+\le r\le\rho}r^2Q(r)\right)\|f\|_{4,QH}\\
 &\le C\vartheta_{\rm out}\|f\|_{4,QH},
 \qquad \vartheta_{\rm out}=\rho^q+\rho^2+\delta^2.
\end{split}
 \label{eq:outer-green-schauder}
\end{equation}
Indeed $H_+,Q$ are comparable on each enlarged annulus, and
$b^2r_+^{-2}=\delta^2$.  At each boundary radius $r\in\{r_+,\rho\}$, the scaled trace
estimate gives
\[
 \left\|\frac{r\partial_\nu(G_Df)}{H_+(r)}\right\|_{C^{5,\alphaH}}
 \le C\vartheta_{\rm out}\|f\|_{4,QH},
\]
where $\partial_\nu$ is the Euclidean outward normal derivative
and the boundary norm is taken after pullback to the unit sphere.

The metric coefficients satisfy
\[
 \|h-g_{\rm E}\|_{C^{6,\alphaH}_{\rm sc}(A_r)}\le Cr^q,
 \qquad
 \|\operatorname{Scal}_h\|_{C^{4,\alphaH}_{\rm sc}(A_r)}\le Cr^{q-2}.
\]
Thus the metric terms acting on $u$ are bounded by
$Cr^{q-2}H_+(r)\|u\|_{6,H}$ in scaled $C^{4,\alphaH}$.
Also $H_+^{p-1}\le C(1+b^2r^{-4})$.
Let $U_0$ be harmonic with the prescribed traces. The maximum principle
and boundary Schauder estimate give
$cH_+\le U_0\le CH_+$ and $\|U_0\|_{6,H}\le C$.
For $u=U_0+v$, solve
\[
 v=\mathcal Q(v):=G_D\bigl[(\Delta_h-\Delta)(U_0+v)
 -c(n)\operatorname{Scal}_h(U_0+v)+\frac{n(n-2)}4(U_0+v)^p\bigr],
 \qquad v|_{\partial A_{r_+,\rho}}=0.
\]
On a fixed positive relative neighborhood of $U_0$,
\eqref{eq:relative-power-calculus} and \eqref{eq:outer-green-schauder} imply
\[
 \|\mathcal Q(0)\|_{6,H}
       +\sup_v\|D_v\mathcal Q\|\le C\vartheta_{\rm out}.
\]
The map $\mathcal Q$ preserves the ball
$\|v\|_{6,H}\le C_0\vartheta_{\rm out}$ for fixed sufficiently large
$C_0$. Its
Lipschitz constant on this ball is at most $1/2$. The resulting
solution satisfies $u\ge(c/2)H_+>0$.

To differentiate the moving inner boundary, fix $b_0$, set
$s_0=b_0/\delta$, and use
\[
 \Phi_\beta(X)=e^{\beta\chi_0(|X|)}X,\qquad b=e^\beta b_0,
\]
where $\chi_0=1$ for $r\le2s_0$ and $\chi_0=0$ for $r\ge4s_0$.
For small bounded $\beta$, this identifies $A_{s_0,\rho}$ with
$A_{b/\delta,\rho}$ and fixes the outer sphere.
In the differentiated equations, $G_D$ denotes the Green operator
transported by this map to the reference annulus.
The pulled-back Euclidean Laplacian satisfies
\[
 \|(\partial_\beta\Delta_\beta)w\|_{C^{4,\alphaH}_{\rm sc}(A_r)}
 \le Cr^{-2}\|w\|_{C^{6,\alphaH}_{\rm sc}(A_r)}
 \quad(s_0\le r\le4s_0),
\]
and its derivative vanishes elsewhere. Its support weight is
$Q_0(r)=r^{-2}\boldsymbol1_{[s_0,4s_0]}(r)$, for which
$\int rQ_0\dd r=\log4$ and $\sup r^2Q_0\le1$.
Use a fixed smooth majorant on a slightly larger annulus in
H\"older norms. The Green comparison and Schauder estimates give
\begin{equation}
 \|G_Df\|_{6,H}\le C\|f\|_{4,Q_0H}.
 \label{eq:outer-moving-boundary-inverse}
\end{equation}
The constant is independent of the moving inner radius.

Write $g_-$ and $g_+$ for the prescribed inner and outer traces of
$U_0$. For a unit finite or Fr\'echet trace direction $i$, differentiate
the harmonic function $U_0$:
\[
 -\Delta_\beta(U_0)_i=(\partial_i\Delta_\beta)U_0,
 \qquad (U_0)_i|_{\partial_\pm}=(g_\pm)_i.
\]
The traces are measured in the fixed $C^{6,\alphaH}$ spaces after
pullback. Harmonic extension and the $Q_0$ bound give
$\|(U_0)_i\|_{6,H}\le C$.
On the reference annulus, write the equation for $v=u-U_0$ as
$-\Delta_\beta v=F(\beta,v)$. Differentiation gives
\[
 (-\Delta_\beta-F_v)v_i=F_i+(\partial_i\Delta_\beta)v,
 \qquad v_i|_{\partial_\pm}=0.
\]
Here $F_i$ holds $v$ fixed and includes $(U_0)_i$. The two
source bounds are
\[
 \|F_i\|_{4,QH}\le C,\qquad
 \|(\partial_i\Delta_\beta)v\|_{4,Q_0H}
       \le C\|v\|_{6,H}\le C\vartheta_{\rm out}.
\]
The inverse estimate with weight $Q$ has the factor
$\vartheta_{\rm out}$. The estimate with weight $Q_0$ is uniform,
and $\|(I-G_DF_v)^{-1}\|\le2$. Hence
\[
 \|D_{\rm par}^jv\|_{6,H}\le C\vartheta_{\rm out},
 \qquad j=0,1.
\]
The coefficient and positive power maps are jointly $C^1$ on the
entire parameter and trace neighborhood. The estimate for
$D_{\rm par}^jv$ holds jointly in the finite parameters and boundary
traces. The spatial regularity is six derivatives for the field, four
for the equation source, and five for its conormal trace; the parameter
regularity in this assertion is separately $C^1$.

The harmonic inner-data and outer-data
solutions $E_Dg$ and $P_{\rm out}g$ have inward conormals
\[
 \begin{aligned}
 [-r_+\partial_rE_Dg]_l
 &=\left((n-2)+l+\frac{(n-2+2l)(r_+/\rho)^{2l+n-2}}
                         {1-(r_+/\rho)^{2l+n-2}}\right)g_l,\\
 [-r_+\partial_rP_{\rm out}g]_l
 &=-\frac{(n-2+2l)(r_+/\rho)^l}{1-(r_+/\rho)^{2l+n-2}}g_l,
 \end{aligned}
\]
where the opposite trace is zero in each case. Expanding each denominator as a geometric series gives
\[
\begin{aligned}
 -r_+\partial_rE_D-((n-2)+D_+)
 &=\sum_{m\ge1}(r_+/\rho)^{(n-2)m}(n-2+2D_+)(r_+/\rho)^{2mD_+},\\
 -r_+\partial_rP_{\rm out}
 &=-\sum_{m\ge0}(r_+/\rho)^{(n-2)m}(n-2+2D_+)(r_+/\rho)^{(2m+1)D_+}.
\end{aligned}
\]
With $\rho$ fixed, Poisson smoothing, also after $r_+\partial_{r_+}$, gives for $j=0,1$
\[
\begin{aligned}
 \| (r_+\partial_{r_+})^j[-r_+\partial_rE_D-((n-2)+D_+)]
       \|_{C^{6,\alphaH}\to C^{5,\alphaH}}&\le C(r_+/\rho)^{n-2},\\
 \| (r_+\partial_{r_+})^j[-r_+\partial_rP_{\rm out}]
       \|_{C^{6,\alphaH}\to C^{5,\alphaH}}&\le C.
\end{aligned}
\]
On the $m$th summand, $r_+\partial_{r_+}$ introduces $(n-2)m+2mD_+$
or $(n-2)m+(2m+1)D_+$. The Poisson semigroup estimates and the
geometric series are uniform for $(r_+/\rho)\le1/2$.

Under inversion $y^\#=y/|y|^2$, the Kelvin transform satisfies
\[
 f^\#(y^\#)=|y^\#|^{2-n}f(y^\#/|y^\#|^2),\qquad
 \delta\partial_{|y^\#|}f^\#
 =\delta^{2-n}[-|y|\partial_{|y|}f-(n-2)f]
 \quad(|y|=\delta^{-1}).
\]
The Kelvin transformation subtracts the $(n-2)$ term from the leading
conormal coefficient $(n-2)+l$, leaving $D_+$.
For compactified trace $\phi_c$, the physical, Kelvin and round
traces are related by
\begin{equation}
 g=b^{-(n-2)/2}\delta^{n-2}V_0(\delta)\phi_c,
 \qquad \delta^{2-n}b^{(n-2)/2}g=V_0(\delta)\phi_c,
 \qquad H_+(r_+)\asymp b^{-(n-2)/2}\delta^{n-2}.
 \label{eq:compactified-trace-factors}
\end{equation}
Use \eqref{eq:compactified-trace-factors} to convert the traces
and conormals to the compactified factor. Express both the actual
and reference conormal densities relative to the same measure
$\dd\theta$ after pullback to the unit sphere.
For the round metric,
\begin{equation}
 \delta^{-(n-2)}\partial_{\nu_{\rm rd}}z\,
       \dd S_{\rm rd}/\dd\theta=V_0(\delta)^2\delta\partial_rz.
 \label{eq:round-conormal-density}
\end{equation}
For the scalar comparison, divide both conormal densities by
$V_0(\delta)^2$.
The change $u^\#=V_0z$ contributes
$-\delta V_0'(\delta)V_0(\delta)^{-1}z=O(\delta^2)z$.
Also $s=2\arctan\delta$ gives
$\delta\partial_\delta=(1+O(\delta^2))s\partial_s$, and the
normalized round surface measure differs from its geodesic Euclidean
reference by $1+O(\delta^2)$.

The bounded regular outer data contribute
$O_{C^1}(b^{(n-2)/2}\delta^{2-n})$ after compactification.
The bounds for $v=u-U_0$ and its first parameter derivatives give a
normalized conormal error $O(\vartheta_{\rm out})$:
an inner trace direction $g$ has relative size
$H_+(r_+)^{-1}\|g\|_{C^{6,\alphaH}}$, and the scaled boundary
derivative cancels $H_+(r_+)$. The physical normal and density errors
are $O_{C^1(C^{5,\alphaH})}(r_+^q)=O(\vartheta_{\rm out})$.
Let $f$ denote the round compactified trace on the inner boundary
of the outer annulus, and let $\vartheta$ collect the remaining
inputs $(\ell,\zeta,A,\mathbf v,\psi)$. Define
$\mathcal N_+^{\rm sc}(\vartheta,f)$ to be its conormal with the outward normal of the outer annulus,
transformed by \eqref{eq:compactified-trace-factors} and divided
by $V_0(\delta)^2$ as above. Combining the harmonic multipliers
with \eqref{eq:round-conormal-density} gives, on the prescribed
positive trace set,
\begin{equation}
 \mathcal N_+^{\rm sc}(\vartheta,f)=D_+f+E_+(\vartheta,f),\qquad
 \|E_+\|_{C^1(\vartheta,f;C^{5,\alphaH})}
 \le C(\vartheta_{\rm out}+\tau^{n-2}+b^{(n-2)/2}\delta^{2-n}).
 \label{eq:outer-affine-c1}
\end{equation}
The error $E_+$ includes
the regular outer trace and all normal and density perturbations.
By \eqref{eq:bridge-scale-ledger},
\[
 b^{(n-2)/2}\le C\delta^{n},\qquad \rho\le C\delta,\qquad
 \tau^{n-2}=(b^{(n-2)/2}\delta^{2-n})(b^{(n-2)/2}\rho^{2-n})
               =O(\delta^2\rho^{3/2}).
\]
Since $q\ge2$, $\|E_+\|_{C^1}\le C\delta^2$.
Let $DN_{L_\zeta}^{\rm reg,sc}$ map boundary values on the missing
round ball at the outer compactification point to the outward
conormal of their regular $L_\zeta$-harmonic extension, with the
same normalization as $\mathcal N_+^{\rm sc}$.
Lemma~\ref{lem:regular-filling} and
\eqref{eq:round-conormal-density} give
\[
 \begin{aligned}
 \|DN_{L_\zeta}^{\rm reg,sc}-D_+\|_{C^{6,\alphaH}\to C^{5,\alphaH}}
       &\le C\delta^2,\\
 \|D_\zeta DN_{L_\zeta}^{\rm reg,sc}\|_{C^{6,\alphaH}\to C^{5,\alphaH}}
       &\le C\delta^2.
 \end{aligned}
\]
In particular, $D_f\mathcal N_+^{\rm sc}$ differs from
$DN_{L_\zeta}^{\rm reg,sc}$ by $O(\vartheta_{\rm out}+\tau^{n-2})$
as an operator from $C^{6,\alphaH}$ to $C^{5,\alphaH}$.

\subsection{Solving the two internal transmission conditions}
\label{subsec:bridge-transmission}

We join the inner and outer solutions to the middle-region solution
by requiring equality of boundary values and cancellation of outward
conormal densities at $|y|=\delta$ and $|y|=\delta^{-1}$.

On $\Omega_\delta=\{y:\delta<|y|<\delta^{-1}\}\subset\Sph^n$,
use the compactified metric $\gamma_b$ and action gradient $\Ff_b$
defined in \eqref{eq:compactified-action-source}. Write
$\vartheta_b=\dd V_{\gamma_b}/\dd V_{\rm rd}$.
The physical conformal factor is
$u(X)=b^{-(n-2)/2}V_0(y)z(y)$. With the first metric correction, write
$z=u_\zeta+W_1+w$ and impose
\begin{equation}
 \Ff_b(z)=\sum_{j=0}^nc_j\chi Z_j,\qquad
 \int_{\Omega_\delta}\chi Z_jw\dd V_{\rm rd}=0
 \quad(0\le j\le n).
 \label{eq:bridge-projection}
\end{equation}
For small $\delta$, $\operatorname{supp}\chi\subset\Omega_\delta$.

For the transmission equations, we use the conormal density
of the compactified factor in the metric $\gamma_b$.
In the $y$ coordinates, parametrize the inner and outer boundary
spheres by
\[
 \iota_-(\theta)=\delta\theta,\qquad
 \iota_+(\theta)=\delta^{-1}\theta,\qquad
 \theta\in\Sph^{n-1},
\]
and write $\operatorname{tr}_\pm z=z\circ\iota_\pm$.
Define the central conormal maps by
\begin{equation}
 (\mathcal N_c^\pm z)\,\dd\theta
 =\delta^{2-n}\iota_\pm^*
       \bigl(\partial_{\nu_{\gamma_b}}z\,\dd S_{\gamma_b}\bigr),
 \label{eq:transmission-density}
\end{equation}
where $\nu_{\gamma_b}$ points out of $\Omega_\delta$ and
$\dd\theta$ is the standard surface measure on $\Sph^{n-1}$.
For each adjoining region, use the same normalization and
the normal pointing out of that region.

These conormal equations give the physical matching conditions.
Indeed, with $\mathcal D_b(y)=by$ and
$u(by)=b^{-(n-2)/2}V_0(y)z(y)$,
\[
 \mathcal D_b^*
 \bigl(\partial_{\nu_h}u\,\dd S_h\bigr)
 =b^{(n-2)/2}V_0^{-1}
 \bigl[\partial_{\nu_{\gamma_b}}z
       +z\partial_{\nu_{\gamma_b}}\log V_0\bigr]
       \dd S_{\gamma_b}.
\]
At a shared interface, the two factors have the same trace.
The terms involving $\partial_{\nu_{\gamma_b}}\log V_0$
therefore cancel with the opposite outward normals.

For the inner family of Proposition~\ref{prop:inner-dn},
$u_{\vartheta,\trpar(z)}$ denotes the conformal factor relative
to $h_b$ in the $y$ coordinates, and $r=|y|$.
Its boundary parameter and scalar radial derivative satisfy
\begin{equation}
\begin{aligned}
 \trpar(z)&=s_\eps^{-1}
       [V_0(\delta)\operatorname{tr}_-z-f_{\rm bg}(\vartheta)],\\
 \left.\delta\partial_ru_{\vartheta,\trpar(z)}\right|_{r=\delta}
 &=\left.\delta\partial_ru_{\vartheta,0}\right|_{r=\delta}
       +s_\eps\widetilde N(\vartheta,\trpar(z)).
\end{aligned}
\label{eq:inner-complete-flux}
\end{equation}
Its compactified factor is $u_{\vartheta,\trpar(z)}/V_0$.
On this sphere, the normal-coordinate and volume identities give
$\nu_{h_b}=\partial_r$ and
$\dd S_{h_b}=\delta^{n-1}\dd\theta$.
Thus the corresponding conormal map is
\[
\begin{aligned}
 \mathcal N_-(\operatorname{tr}_-z)
 &=V_0(\delta)^2
   \left.\delta\partial_r
       \left(\frac{u_{\vartheta,\trpar(z)}}{V_0}\right)
   \right|_{r=\delta}\\
 &=V_0(\delta)\left[
       \left.\delta\partial_ru_{\vartheta,0}\right|_{r=\delta}
       +s_\eps\widetilde N(\vartheta,\trpar(z))
       \right]
   -\delta V_0'(\delta)V_0(\delta)\operatorname{tr}_-z .
\end{aligned}
\]
For the outer solution, restore the factor divided out in
\eqref{eq:outer-affine-c1}:
\[
 \mathcal N_+(\operatorname{tr}_+z)
 =V_0(\delta)^2
       \mathcal N_+^{\rm sc}(\vartheta,\operatorname{tr}_+z).
\]
The transmission map is
\begin{equation}
 \mathcal T_b(w,c):=
 \begin{pmatrix}
 \Ff_b(u_\zeta+W_1+w)-\sum_{j=0}^nc_j\chi Z_j\\
 \mathcal N_c^-z+\mathcal N_-(\operatorname{tr}_-z)\\
 \mathcal N_c^+z+\mathcal N_+(\operatorname{tr}_+z)\\
 (\int\chi Z_jw\dd V_{\rm rd})_{j=0}^n
 \end{pmatrix}.
 \label{eq:bridge-transmission-map}
\end{equation}
The components of $\mathcal T_b$ are the interior residual,
the two conormal jumps and
the localized Jacobi moments.

Take fixed relative annuli $A_d$ around either compactification
point, one-sided at a cut, and include the fixed central region.
Define
\[
 \|w\|_{\mathcal X_\delta}
 =\|w\|_{H^1(\Omega_\delta)}
 +\sup_{d\ge c\delta}d^{(n-2)/2}\|w\|_{C^{6,\alphaH}_{\rm sc}(A_d)},
 \qquad
 \|(w,c)\|_{\widehat{\mathcal X}_\delta}=\|w\|_{\mathcal X_\delta}+|c|.
\]
Let $\widetilde f$ be the zero extension of the volume source as a
distribution on $\Sph^n$. For conormal jump data $\mathfrak j_-,\mathfrak j_+$ and moment data $\mathbf m\in\R^{n+1}$, set
\begin{align*}
 \|(f,\mathfrak j_-,\mathfrak j_+,\mathbf m)\|_{\mathcal Y_\delta}
 ={}&\|\widetilde f\|_{H^{-1}(\Sph^n)}
 +\sup_{d\ge c\delta}d^{(n+2)/2}\|f\|_{C^{4,\alphaH}_{\rm sc}(A_d)}+|\mathbf m|\\
 &+\delta^{(n-2)/2}\sum_\pm
 \|\mathfrak j_\pm\|_{C^{5,\alphaH}(\Sph^{n-1})}.
\end{align*}
The norms $\mathcal X_\delta$ and $\mathcal Y_\delta$ include the
fixed central norms. The combined forcing is
\[
 \phi\longmapsto\int_{\Omega_\delta}f\phi\dd V_{\rm rd}
 +\delta^{n-2}\sum_\pm\int_{\Sph^{n-1}}
          \mathfrak j_\pm(\phi\circ\iota_\pm)\dd\theta.
\]
This volume-and-surface distribution has $H^{-1}(\Sph^n)$ norm at
most $C\|(f,\mathfrak j_-,\mathfrak j_+,\mathbf m)\|_{\mathcal Y_\delta}$. The volume term is the zero extension.
Using $\|\mathfrak j_\pm\|_{H^{-1/2}}\le C\|\mathfrak j_\pm\|_{C^{5,\alphaH}}$ on the
unit sphere and the scaled trace inequality, we obtain
\[
 \begin{aligned}
 \left|\delta^{n-2}\int \mathfrak j_\pm(\phi\circ\iota_\pm)\dd\theta\right|
 &\le\delta^{n-2}\|\mathfrak j_\pm\|_{H^{-1/2}}
                    \|\phi\circ\iota_\pm\|_{H^{1/2}}\\
 &\le C\delta^{(n-2)/2}\|\mathfrak j_\pm\|_{C^{5,\alphaH}}\|\phi\|_{H^1(\Sph^n)}.
 \end{aligned}
\]
Localized multiplier sources have bounded data norm per unit
coefficient. Reference weights, collars, measures and trace
parametrizations remain fixed on each parameter patch; differentiate the conormal density in
\eqref{eq:transmission-density} before taking these norms.
Finally, put
\[
 \kappa_b=
 \begin{cases}
 C\delta,&n\ge4q+2,\\
 C(\delta+b^q/\eps),&n<4q+2.
 \end{cases}
\]
By \eqref{eq:weak-metric-size}, $\kappa_b\to0$ both for
$n\ge4q+2$ and for $n<4q+2$.

For the reference operator $L_\zeta$, use
\eqref{eq:transmission-density} with $\gamma_b$ replaced by
$g_{\rm rd}$.
Let $\mathcal N_{\zeta,c}^\pm$ denote the central outward
conormal maps and $\mathcal N_{\zeta,{\rm fill}}^\pm$ those
of the regular $L_\zeta$-harmonic fillings in the missing balls.
The central region and each filling use their respective
outward normals.
Define
\begin{equation}
 \mathcal A_\zeta(w,c)=
 \begin{pmatrix}
 L_\zeta w-\sum_{j=0}^nc_j\chi Z_j\\
 \mathcal N_{\zeta,c}^-w+\mathcal N_{\zeta,{\rm fill}}^-(\operatorname{tr}_-w)\\
 \mathcal N_{\zeta,c}^+w+\mathcal N_{\zeta,{\rm fill}}^+(\operatorname{tr}_+w)\\
 (\int\chi Z_jw\dd V_{\rm rd})_{j=0}^n
 \end{pmatrix}.
 \label{eq:reference-transmission}
\end{equation}
\begin{lemma}\label{lem:bridge-transmission}
The operator \eqref{eq:reference-transmission} is uniformly invertible
from $\widehat{\mathcal X}_\delta$ to $\mathcal Y_\delta$.
For fixed sufficiently large $R_0$ and small $\eps$, the map
\eqref{eq:bridge-transmission-map} has a unique zero in
\[
 \mathbb B_b=\{(w,c):
 \|(w,c)\|_{\widehat{\mathcal X}_\delta}\le R_0\kappa_b\eps\}.
\]
The solution $(w,c)$ is jointly $C^1$ in $(\ell,\zeta,A,\mathbf v,\psi)$.
We have
\begin{equation}
 \|D_{\rm par}^jw\|_{\mathcal X_\delta}+|D_{\rm par}^jc|
 \le C\kappa_b\eps,\qquad j=0,1.
 \label{eq:bridge-shrinking}
\end{equation}
\end{lemma}
\needspace{6\baselineskip}
\begin{proof}
\noindent\textit{Step 1. The reference inverse.}
For the combined volume and surface distribution $F_\delta$,
Lemma~\ref{lem:sphere-augmented-inverse} solves
\[
 L_\zeta\widetilde w-\sum_{j=0}^nc_j\chi Z_j=F_\delta,
 \qquad \int\chi Z_i\widetilde w\dd V_{\rm rd}=(\mathbf m)_i,
 \qquad
 \|\widetilde w\|_{H^1(\Sph^n)}+|c|
 \le C\|(f,\mathfrak j_-,\mathfrak j_+,\mathbf m)\|_{\mathcal Y_\delta}.
\]
Here the surface-source estimate supplies the factor $\delta^{(n-2)/2}$.
In each missing ball $F_\delta$ and $\chi$ vanish. The filling is
unique, since a zero boundary solution obeys
$\int|\nabla v|^2\le C\int v^2\le C\delta^2\int|\nabla v|^2$.
Restriction therefore gives the reference transmission solution with
the moment conditions $\int\chi Z_i\widetilde w\dd V_{\rm rd}=(\mathbf m)_i$.
Conversely, a homogeneous solution
joined to its fillings is a homogeneous augmented sphere solution.

H\"older and the sphere Sobolev inequality give, on every complete
or truncated annulus,
\[
 d^{-n/2}\|\widetilde w\|_{L^2(A_{2d})}
 \le Cd^{-(n-2)/2}\|\widetilde w\|_{H^1(\Sph^n)}.
\]
At each cut use the local compactification chart with radius $\delta$.
The jump is the ball conormal plus the central conormal.
By \eqref{eq:round-conormal-density}, density data $\mathfrak j_l$ of degree $l$
correspond to scalar data $\mathfrak j_l/V_0(\delta)^2$. The harmonic layer has
trace and radial factors
\[
 a_l=\frac{\mathfrak j_l}{V_0(\delta)^2(2l+n-2)},\qquad
 \mathscr S \mathfrak j_l=
 \begin{cases}
 a_l(r/\delta)^l,&r\le\delta,\\
 a_l(r/\delta)^{2-n-l},&r\ge\delta.
 \end{cases}
\]
Consequently
\[
 V_0(\delta)^2\delta
 [\partial_r\mathscr S \mathfrak j_l|_{\rm ball}
       -\partial_r\mathscr S \mathfrak j_l|_{\rm central}]=\mathfrak j_l.
\]
Multiply $\mathscr S\mathfrak j$ by a fixed smooth cutoff in $r/\delta$
that equals one near $r/\delta=1$, and retain the notation
$\mathscr S\mathfrak j$ for this localized layer. Its cutoff derivatives
and $L_\zeta\mathscr S\mathfrak j$ away from the interface enter the
volume sources on the two sides. The density jump remains exactly $\mathfrak j$.
Since $V_0(\delta)^2\asymp1$, for $d\asymp\delta$,
\[
 \|\mathscr S\mathfrak j\|_{C^{6,\alphaH}_{\rm sc}(A_d)}
 +d^2\|L_\zeta\mathscr S\mathfrak j\|_{C^{4,\alphaH}_{\rm sc}(A_d)}
 \le C\|\mathfrak j\|_{C^{5,\alphaH}}.
\]
Apply the local Schauder estimate to $\widetilde w-\mathscr S\mathfrak j$,
whose conormal jump is zero, and restore the layer after the estimate.
The moment conditions remain those of the global sphere solution
$\widetilde w$. Thus
\begin{equation}
 \|w\|_{C^{6,\alphaH}_{\rm sc}(A_d)}
 \le C\left[d^{-(n-2)/2}\|\widetilde w\|_{H^1(\Sph^n)}
 +d^2\|f+\sum_{j=0}^nc_j\chi Z_j\|_{C^{4,\alphaH}_{\rm sc}(A_{2d})}
 \right],
 \label{eq:transmission-schauder}
\end{equation}
with $C\|\mathfrak j_\pm\|_{C^{5,\alphaH}}$ added at a cut collar and
source norms taken separately on its two sides. Multiplying by
$d^{(n-2)/2}$ proves the asserted uniform inverse bound.

\medskip
\noindent\textit{Step 2. Estimating $\mathcal T_b(0,0)$.}
Put $z_0=u_\zeta+W_1$ and recall
$\Xtrace=C^{6,\alphaH}(\Sph^{n-1})$.
When $n\ge4q+2$, $\lambda=b^q\asymp\eps$ and $W_1$ cancels the
first metric source. Near the outer compactification point, the next
metric term is bounded by $C\lambda b\,d^{-q-3}$. The terms quadratic
in the metric perturbation, and those bilinear in the metric
perturbation and $W_1$, are bounded by $C\lambda^2d^{-2q-2}$.
The power remainder is bounded by $C\lambda^2d^{-2q}$. Each bound
includes four scaled spatial derivatives and one parameter derivative;
further factors are $O(\lambda\delta^{-q})=o(1)$. Hardy's inequality gives
\[
\begin{aligned}
 \|\lambda b\,d^{-q-3}\|_{H^{-1}}^2
 &\le C\lambda^2b^2\int_\delta^1t^{n-2q-5}\dd t
 \le C\lambda^2b^2,\\
 \|\lambda^2d^{-2q-2}\|_{H^{-1}}^2
 &\le C\lambda^4\int_\delta^1t^{n-4q-3}\dd t
 \le C\lambda^4L_\delta^2,\qquad
 L_\delta=\begin{cases}1,&n>4q+2,\\
 (1+|\log\delta|)^{1/2},&n=4q+2.
 \end{cases}
\end{aligned}
\]
The integral with exponent $n-2q-5$ and the strong annular bounds
use $\frac{n-2}{2}-q-1>0$. The integral with exponent $n-4q-3$ retains the
logarithm at $n=4q+2$. Since $\frac{n-2}{2}-2q\ge0$,
the full volume norm is $O_{C^1}(\lambda b+\lambda^2L_\delta)$ and
\[
 \frac{\lambda b+\lambda^2L_\delta}{\eps\delta}
 \le C(b/\delta+\delta^{(n-2)/2}L_\delta)\longrightarrow0.
\]
When $n<4q+2$, \eqref{eq:weak-metric-size} instead gives
$O_{C^1}(b^q)$ for the full volume norm.

For this application, use the reference profile $v$ fixed in
Section~\ref{sec:inner-tail}; thus $\tau=0$ and
\[
 \delta^{-(n-2)/2}v(0)=V_\zeta(0)=2^{(n-2)/2},
 \qquad \kelpar_0=\zeta.
\]
Here $t=0$ is $r=\delta$; the tail at $r=d$ uses the shift $\log2$.
The equality $\kelpar_0=\zeta$ is the choice of reference Kelvin data for
this matching problem; it does not identify the varying path
$\kelpar(t)$ with the fixed bubble parameter throughout the end.
The bubble has $\nabla V_\zeta(0)=(n-2)V_\zeta(0)\zeta$.
The reference expansion in \eqref{eq:collar-affine-reference} gives
$V_\zeta(\delta\,\cdot)-f_{\rm bg}=O_{C^1(\Xtrace)}(\delta^2)$.
Thus $V_0(\delta)\operatorname{tr}_-z_0-f_{\rm bg}
=O_{C^1(\Xtrace)}(\delta^2)$ after including the $W_1$ trace.
For $n\ge4q+2$, the inner $W_1$ trace obeys
$\|W_1\|_{C^1(\Xtrace)}\le C\lambda\delta^{-q}\le C\delta^3$;
for $n<4q+2$, $W_1=0$.

Since $s_\eps\asymp\delta$, the preceding trace comparison and
\eqref{eq:inner-complete-flux} give
\[
 \|\trpar(z_0)\|_{C^1(\Xtrace)}\le C\delta.
\]
Thus the background trace belongs to the full trace ball of
Proposition~\ref{prop:inner-dn}.

At each cut, expand the bubble in the corresponding smooth
compactification chart. The constant and linear terms give
opposite leading conormal contributions from the adjoining
regions. The remaining Taylor terms have scaled value and
normal traces $O_{C^1}(\delta^2)$.

At the inner cut, the Fowler slope error is
$O(\eps^2\delta^{2-n})=O(\delta^2)$. The reference conormal
comparison in the proof of Proposition~\ref{prop:inner-dn},
together with its full trace estimate, bounds the remaining
inner conormal error by
$s_\eps O_{C^1}(\delta)=O_{C^1}(\delta^2)$.
At the outer cut, \eqref{eq:outer-affine-c1} gives the same
$O_{C^1}(\delta^2)$ bound, including the regular outer data.
The metric conormal coefficients differ from their round
counterparts by
$O_{C^1}((b/\delta)^q)=O_{C^1}(\delta^2)$,
by \eqref{eq:bridge-scales}.
Lemma~\ref{lem:first-metric-correction} also bounds the traces
and scaled normal derivatives of $W_1$ at both cuts by
$C\lambda\delta^{-q}\le C\delta^3$ when $n\ge4q+2$;
for $n<4q+2$, $W_1=0$.
Consequently, for $j=0,1$,
\[
 \sum_\pm
 \left\|D_{\rm par}^j\left[
 \mathcal N_c^\pm z_0+
 \mathcal N_\pm(\operatorname{tr}_\pm z_0)
 \right]\right\|_{C^{5,\alphaH}(\Sph^{n-1})}
 \le C\delta^2.
\]
The contribution of these jumps to the data norm is at most
\[
 C\delta^{(n-2)/2}\delta^2=C\eps\delta.
\]
The moment component of $\mathcal T_b(0,0)$ is zero, together
with its parameter derivatives, since $w=0$.
Combining this with the volume estimates above gives
\[
 \|D_{\rm par}^j\mathcal T_b(0,0)\|_{\mathcal Y_\delta}
 \le C(\eps\delta+\boldsymbol1_{\{n<4q+2\}}b^q)
 \le C\kappa_b\eps,\qquad j=0,1.
\]

\medskip
\noindent\textit{Step 3. Linearization on $\mathbb B_b$.}
The metric agrees across each cut. Write the difference between
the linearized middle operator at $u_\zeta+W_1$ and $L_\zeta$ as
$-\operatorname{div}(A_b\nabla\,\cdot\,)+V_b$. These coefficients satisfy
\[
 \|A_b\|_\infty+\|d^2V_b\|_\infty\le C\nu_b,
 \qquad \nu_b=(b/\delta)^q.
\]
Hardy's inequality also holds for the nonzero boundary traces here:
subtract the mean on each adjacent annulus, extend on the fixed
rescaled annulus, and restore the mean. Scaled Poincar\'e gives
$\|\widetilde h\|_{H^1(\Sph^n)}\le C\|h\|_{H^1(\Omega_\delta)}$,
hence $\|h/d\|_2\le C\|h\|_{H^1(\Omega_\delta)}$. Thus
\[
\begin{aligned}
 \left|\int A_b\nabla h\cdot\nabla\phi\right|
   +\left|\int V_bh\phi\right|
 &\le C\nu_b\|h\|_{H^1}\|\phi\|_{H^1},\\
 d^{(n+2)/2}\left(\|\operatorname{div}(A_b\nabla h)\|_{C^{4,\alphaH}_{\rm sc}}
                 +\|V_bh\|_{C^{4,\alphaH}_{\rm sc}}\right)
 &\le C\nu_b\|h\|_{\mathcal X_\delta}.
\end{aligned}
\]
Integration by parts produces the conormal density in
\eqref{eq:transmission-density} in the jump component.
 The radially normalized boundary error produced by the
principal coefficient perturbation $A_b$ is at most
$C\nu_b\delta^{-(n-2)/2}\|h\|_{\mathcal X_\delta}$ in $C^{5,\alphaH}$;
the factor $\delta^{(n-2)/2}$ in the surface-source estimate controls both
the jump and the boundary term in the zero-extended volume source.
Both sides use the normalization
\eqref{eq:transmission-density}.

Fix $R_0$ and put $\mathbf u=(w,c)\in\mathbb B_b$.
The state norm gives
\[
 \|w\|_\infty\le C_{R_0}\kappa_b\delta,\qquad
 \|\operatorname{tr}_\pm w\|_\Xtrace
 \le C\delta^{-(n-2)/2}\|w\|_{\mathcal X_\delta}.
\]
The collar's outer input at $|y|=\delta$ is
$\trpar_w=\trpar(z_0+w)$ from \eqref{eq:inner-complete-flux}. Set
\[
 f_w=\operatorname{tr}_+(z_0+w),\qquad
 (\trpar_0,f_0)=(\trpar_w,f_w)|_{w=0}.
\]
The collar fixes $|y|=\delta$, and $w\mapsto\trpar_w$ is affine on $\Xtrace$.
Using the background trace estimate from Step 2, we obtain,
at fixed $w$,
\[
 \|\trpar_0\|_{C^1_\vartheta(\Xtrace)}\le C\delta,\qquad
 \|\trpar_w\|_\Xtrace+\|D_\vartheta\trpar_w\|_\Xtrace\le C_{R_0}\kappa_b,
 \qquad \|f_0\|_{C^1_\vartheta(\Xtrace)}\le C.
\]
For every $(w,c)\in\mathbb B_b$, $\trpar_w$ belongs to the full trace
ball of Proposition~\ref{prop:inner-dn}, and $f_w$ belongs to the
positive trace set used in \eqref{eq:outer-affine-c1}.

Set
$\widetilde R(\vartheta,\trpar)
 =\widetilde N(\vartheta,\trpar)-D_+\trpar$.
Thus the scalar conormal remainder in
\eqref{eq:inner-complete-flux} is
$s_\eps\widetilde R(\vartheta,\trpar)$.
Proposition~\ref{prop:inner-dn} bounds the joint $C^1$ norm of
$\widetilde R$ by $C\delta$ on the fixed trace ball;
\eqref{eq:outer-affine-c1} bounds the outer remainder $E_+$ by
$C\delta^2$. In a state direction $h$, the inner $s_\eps$ cancels
$s_\eps^{-1}$ in $D_w\trpar_w$, and the trace loss $\delta^{-(n-2)/2}$
cancels the surface-source factor $\delta^{(n-2)/2}$ at both cuts. The
jump operator errors are therefore
$C\delta\|h\|_{\mathcal X_\delta}$ and
$C\delta^2\|h\|_{\mathcal X_\delta}$.
For parameter derivatives at fixed $w$, use the full first chain rule
\[
 D_\vartheta[s_\eps\widetilde R(\vartheta,\trpar)]
 =s_\eps(\partial_\vartheta\widetilde R
                  +D_\trpar\widetilde R\,D_\vartheta\trpar)
\]
separately at $\trpar_w,\trpar_0$. Apply the chain rule also to
$E_+(\vartheta,f_w)$ and $E_+(\vartheta,f_0)$.
Since $\delta\le C\kappa_b$, this gives
\[
\begin{aligned}
 \delta^{(n-2)/2}\|D_\vartheta\{s_\eps[
 \widetilde R(\vartheta,\trpar_w)-\widetilde R(\vartheta,\trpar_0)]\}\|_{C^{5,\alphaH}}
 &\le C_{R_0}\delta^{(n-2)/2}s_\eps\delta
 \le C_{R_0}\eps\delta,\\
 \delta^{(n-2)/2}\|D_\vartheta[
 E_+(\vartheta,f_w)-E_+(\vartheta,f_0)]\|_{C^{5,\alphaH}}
 &\le C\delta^{(n+2)/2}=C\eps\delta.
\end{aligned}
\]
The two displayed parameter bounds for $\widetilde R$ and $E_+$
also control the $H^{-1/2}$ jumps. 
The conversion from the scalar boundary maps to the conormal
densities in \eqref{eq:transmission-density} preserves these
bounds. Differentiating its coefficients on linear trace
increments contributes at most $C\|w\|_{\mathcal X_\delta}$.
Lemma~\ref{lem:regular-filling} adds $O(\delta^2)$ to the state
operator error and bounds the reference DN parameter term by
$C\|w\|_{\mathcal X_\delta}$.

For the volume component, the fundamental theorem of calculus gives
\[
 (z_0+w)^{p-1}-z_0^{p-1}
 =(p-1)w\int_0^1(z_0+tw)^{p-2}\dd t.
\]
The coefficient $(p-1)\int_0^1(z_0+tw)^{p-2}\dd t$ has uniformly
bounded scaled $C^{4,\alphaH}$ norm.
The weak and strong products obey
\[
\begin{aligned}
 |\langle wh,\phi\rangle|
 &\le C_{R_0}\kappa_b\delta
       \|h\|_{\mathcal X_\delta}\|\phi\|_{H^1},\\
 d^{(n+2)/2}\|wh\|_{C^{4,\alphaH}_{\rm sc}(A_d)}
 &\le C_{R_0}\kappa_b\eps d^{(6-n)/2}\|h\|_{\mathcal X_\delta}
 \le C_{R_0}\kappa_b\delta^3\|h\|_{\mathcal X_\delta},
\end{aligned}
\]
using $n\ge8$, $d\ge c\delta$ and $\eps=\delta^{n/2}$.
The derivatives with respect to $c$ and the derivatives of the
moment conditions agree exactly with those of $\mathcal A_\zeta$. Combining
the volume bounds for $A_b$, $V_b$ and $wh$ with the two
conormal-jump bounds gives
\begin{equation}
 \sup_{\mathbf u\in\mathbb B_b}
 \|D_{\mathbf u}\mathcal T_b-\mathcal A_\zeta\|
 \le\eta_b:=C_{R_0}[(b/\delta)^q+\kappa_b]\longrightarrow0.
 \label{eq:transmission-full-ball-c1}
\end{equation}
At fixed state, the differentiated reference and central coefficients
contribute $C\|\mathbf u\|$. The power remainder and its parameter
derivative retain $w^2$ and are bounded by
$C_{R_0}\kappa_b^2\eps\delta$ by the preceding product bounds.
The full moment and multiplier derivatives contribute
$C(\|w\|_{H^1}+|c|)$ because $\chi Z_j$ are jointly $C^1$.
Together with the background residual estimate in Step 2 and the
parameter bounds for $\widetilde R$ and $E_+$, this proves
\begin{equation}
 \sup_{\mathbf u\in\mathbb B_b}
 \|\partial_\vartheta\mathcal T_b(\vartheta,\mathbf u)\|_{\mathcal Y_\delta}
 \le C_{R_0}\kappa_b\eps.
 \label{eq:transmission-first-parameter-source}
\end{equation}
The map $\mathcal T_b:\widehat{\mathcal X}_\delta\to\mathcal Y_\delta$
is jointly Fr\'echet $C^1$ in $\vartheta$ and $(w,c)$.

\medskip
\noindent\textit{Step 4. The fixed point and positivity.}
Let $\|\mathcal A_\zeta^{-1}\|\le M$ and
$\|\mathcal T_b(0)\|\le C_0\kappa_b\eps$. Choose
$R_0=2MC_0+1$, then $\eps$ so that $M\eta_b\le1/2$.
The map $\mathcal Q_b=\operatorname{Id}-\mathcal A_\zeta^{-1}\mathcal T_b$
is a contraction of $\mathbb B_b$ into itself, since
$\|\mathcal Q_b(0)\|\le MC_0\kappa_b\eps$ and its derivative
norm is at most $1/2$. Its unique fixed point is the required zero.
The inverse bound $\|(D_{\mathbf u}\mathcal T_b)^{-1}\|\le2M$
and \eqref{eq:fixed-point-parameter-calculus}, with
\eqref{eq:transmission-first-parameter-source}, give
\eqref{eq:bridge-shrinking}. Uniqueness identifies the maps on overlaps.

The function $u_\zeta$ has a uniform positive lower bound, and
\[
 \left\|(W_1+w)/u_\zeta\right\|_{L^\infty(\Omega_\delta)}
 \le C\lambda\delta^{-q}+C\kappa_b\eps\delta^{-(n-2)/2}
 =C\lambda\delta^{-q}+C\kappa_b\delta\longrightarrow0.
\]
For $n<4q+2$, $W_1=0$, so the $\lambda\delta^{-q}$ contribution
is absent. Together with the positivity of the inner and outer solutions, this
proves positivity throughout $B_\rho\setminus\{0\}$. The three solutions
have matching traces and opposite physical conormals at the two cuts.
Their interior residual is \eqref{eq:bridge-projection}.
Matched values and physical conormal densities cancel the interface
distributions. The multiplier source is supported away from the cuts,
so elliptic regularity for the common smooth metric applies across them.
\end{proof}

\subsection{The assembled local family and its residual}
\label{subsec:bridge-output}

\begin{proposition}\label{prop:bridge}
Fix the neck $\eps$ and outer radius $\rho$, let the bubble parameters
$(b,\zeta)$ satisfy \eqref{eq:bridge-scales} in the specified compact
ranges, and prescribe the outer trace
\[
 u|_{S_\rho}=A+\rho\mathbf v\cdot\theta+\rho^2\psi,
 \qquad A\in[c_*,C_*],\quad \Pi_{\ge2}\psi=\psi,
\]
with $\mathbf v$ and $\psi$ in the fixed bounded sets of
Lemma~\ref{lem:regular-filling}. These data determine a family of positive functions $u$ on
$0<|X|\le\rho$ with limiting Fowler neck $\eps$. It is locally
$C^{6,\alphaH}$ on each piece and jointly Fr\'echet $C^1$ in the
finite parameters and $C^{6,\alphaH}$ traces. Values and physical
conormals agree at $|X|=b\delta$ and $|X|=b/\delta$.
The equation is exact on the inner end and outer annulus. Its only
residual is supported in the middle region and, in the compactified
variables, equals $\sum_{j=0}^nc_j\chi Z_j$ as in
\eqref{eq:bridge-projection}.

There is a positive bounded coefficient $\mathcal B$ such that
\[
 \rho^{n-1/2}\mathcal B
       =b^{(n-2)/2}C_\infty(\zeta)[1+o_{C^1}(1)].
\]
Set $B=\rho^{n-1/2}\mathcal B$.

With boundary norms on the unit sphere, the Cauchy data satisfy
\begin{align}
 \fint u&=A,&
 \fint (\rho\partial_ru)
       &=-(n-2)\mathcal B\rho^{3/2}+O_{C^1}(\rho^2),\nonumber\\
 \Pi_1u&=\rho\mathbf v\cdot\theta,&
 \Pi_1(\rho\partial_ru)
       &=\rho\mathbf v\cdot\theta+
                   O_{C^1}(\rho^3+Bb\rho^{1-n}),\label{eq:bridge-cap}\\
 \Pi_{\ge2}u&=\rho^2\psi,&
 \Pi_{\ge2}(\rho\partial_ru)
       &=D_+(\rho^2\psi)+o_{C^1}(\rho^2).\nonumber
\end{align}
The normalized high remainder tends to zero jointly in $C^1$ from
$C^{6,\alphaH}$ to $C^{5,\alphaH}$. The coefficient vector $c$ and
the middle-region function $w$ obey \eqref{eq:bridge-shrinking}.
\end{proposition}
\begin{proof}
Let $R$ be the regular filling with the same outer trace as $u$,
given by Lemma~\ref{lem:regular-filling}. On the outer annulus
$r_+<|X|<\rho$, where $r_+=b/\delta$, both $u$ and $R$ solve the
Yamabe equation exactly. Subtracting their equations there gives
\[
 \begin{gathered}
 \left[-\Delta_h+c(n)\operatorname{Scal}_h
       -\frac{n(n+2)}4
        \int_0^1((1-s)R+su)^{p-1}\dd s\right](u-R)=0
 \quad\text{for }r_+<|X|<\rho,\\
 (u-R)|_{|X|=\rho}=0.
 \end{gathered}
\]
At $r_+=b/\delta$, the central trace and
\eqref{eq:bridge-shrinking} bound the inner data and their first
parameter derivatives by $Cb^{(n-2)/2}r_+^{2-n}$ in scaled
$C^{6,\alphaH}$. The Green comparison also holds with weight $b^{(n-2)/2}r^{2-n}$, since
its two monotonicity inequalities remain valid when the constant
term is zero. The difference equation has potential bounded by
$C(1+b^2r^{-4})$ and metric terms with weight $r^{q-2}$, including
first parameter derivatives. Apply \eqref{eq:outer-green-comparison}
and \eqref{eq:outer-green-schauder}, followed by the Neumann inverse.
For the derivative of the moving inner sphere, use
\eqref{eq:outer-moving-boundary-inverse}. For all unit finite and
Fr\'echet trace directions with $\rho$ fixed, we obtain
\[
\begin{gathered}
 \max_{0\le j\le1}\sup_{r_+\le r\le\rho}
   \frac{r^{n-2}}{b^{(n-2)/2}}
   \|D_{\rm par}^j(u-R)\|_{C^{6,\alphaH}_{\rm sc}(A_r)}\le C,\\
\begin{aligned}
 &\left\|\int_{A_{r_+,\rho}}
       \bigl[\tfrac{n(n-2)}4(u^p-R^p)
                 -c(n)\operatorname{Scal}_h(u-R)\bigr]
                \dd V_h\right\|_{C^1_{\rm par}}\\
 &\qquad\le Cb^{(n-2)/2}\int_{r_+}^{\rho}rQ(r)\dd r
        \le Cb^{(n-2)/2}(\delta^2+\rho^2+\rho^q).
\end{aligned}
\end{gathered}
\]
The lower-limit contribution is bounded by
$Cb^{(n-2)/2}(\delta^2+\rho^2+\rho^q)$, since
$|(r_+)_i|\le Cr_+$ and
$r_+^2Q(r_+)\le C(\delta^2+\rho^2+\rho^q)$.
Converting the physical flux to the Euclidean radial flux adds
$O_{C^1}(b^{(n-2)/2}\rho^q)$ at the two spheres.
Define the actual coefficient by the outer flux,
\[
 B=-\frac{\rho^{n-2}}{(n-2)}\fint \bigl(\rho\partial_r(u-R)\bigr).
\]
Put $B_0=b^{(n-2)/2}C_\infty(\zeta)$ and
$\mathfrak B_s(v)=-(n-2)^{-1}s^{n-2}\fint (s\partial_rv)$.
Then
\[
 B-B_0=
 [\mathfrak B_{r_+}(u-R)-B_0]
 +[\mathfrak B_\rho(u-R)-\mathfrak B_{r_+}(u-R)].
\]
At the inner cut, the bubble mean flux has relative error $O(\delta^2)$.
The terms $W_1$ and $w$ contribute
$O_{C^1}(\lambda\delta^{-q})$ and $O_{C^1}(\kappa_b\delta)$,
respectively.
The regular filling contributes at most
$Cb^{(n-2)/2}\delta^{2-n}=O(\delta^2)$ relative to the bubble.
Thus the cut bounds and the annular source integral above give
\[
\begin{aligned}
 \|\mathfrak B_{r_+}(u-R)/B_0-1\|_{C^1}
 &\le C(\delta^2+\lambda\delta^{-q}+\kappa_b\delta)=o(1),\\
 \|[\mathfrak B_\rho(u-R)-\mathfrak B_{r_+}(u-R)]/B_0\|_{C^1}
 &\le C(\delta^2+\rho^2+\rho^q)=o(1).
\end{aligned}
\]
The second bound includes conversion from physical to Euclidean
radial flux. All derivatives keep $\rho$ fixed. Hence
$B=B_0[1+o_{C^1}(1)]$.
Since $A$ is the prescribed boundary mean and the regular filling has
$\fint (\rho\partial_rR)=O_{C^1}(\rho^2)$,
\[
 \fint u=A,\qquad
 \fint (\rho\partial_ru)=-(n-2)B\rho^{2-n}+O_{C^1}(\rho^2).
\]
The identities for $\fint u$ and $\fint (\rho\partial_ru)$ are the
first two components of \eqref{eq:bridge-cap}, with
$B=\rho^{n-1/2}\mathcal B$.

It remains to estimate the nonconstant conormals. In the unit-ball
equation for the regular filling, the source is
\[
 \rho^2[\frac{n(n-2)}4A^p-c(n)\operatorname{Scal}_h(0)A]
       +O_{C^1(C^{4,\alphaH})}(\rho^3),
\]
because $P\varphi=A+\rho\mathbf v\cdot y+\rho^2P\psi$ and the
metric coefficients differ from the Euclidean ones by $O(\rho^2)$.
The leading term of $R-P\varphi$ is radial, so the nonconstant conormal error
is $O_{C^1}(\rho^3)$ by the ball Schauder estimate. The changes in the unit normal and surface density have the same order:
$O(\rho^2)\,O(\rho)=O(\rho^3)$, since the constant field has zero
gradient. The bounds use the scaled boundary variables
$(A,\mathbf v,\psi)$ and hold after one total derivative with respect
to these variables.

For $u-R$, retain its harmonic extension and write its Euclidean
source as
\[
 f=(\Delta_h-\Delta)(u-R)-c(n)\operatorname{Scal}_h(u-R)
                                  +\frac{n(n-2)}4(u^p-R^p).
\]
The relative $C^1(C^{6,\alphaH}_{\rm sc})$ estimate for $u-R$ gives
\[
 \|f\|_{C^{4,\alphaH}_{\rm sc}(A_r)}
 \le CB\bigl(r^{q-n}+r^{2-n}+b^2r^{-n-2}\bigr).
\]
For degree $l$, Green's identity tests this source against
\begin{equation}
 \frac{(r/\rho)^l-\tau^{2l+n-2}(\rho/r)^{l+n-2}}
      {1-\tau^{2l+n-2}},\qquad \tau=r_+/\rho.
 \label{eq:outer-harmonic-test}
\end{equation}
The harmonic test function in \eqref{eq:outer-harmonic-test} has modulus at
most $C(r/\rho)^l$. Poisson smoothing for
$r\le\rho/2$ and the boundary Schauder estimate on the outer
relative collar give the $C^{5,\alphaH}$ conormal bounds
\[
 \begin{aligned}
 l=1:\quad&
 CB\rho^{1-n}\int_{r_+}^\rho
       (r^q+r^2+b^2r^{-2})\,\dd r
 \le CB\rho^{2-n}\bigl[\rho^q+\rho^2+b\delta/\rho\bigr],\\
 l\ge2:\quad&
 CB\rho^{-n}\int_{r_+}^\rho
       (r^{q+1}+r^3+b^2r^{-1})\,\dd r
 \le CB\rho^{2-n}\bigl[\rho^q+\rho^2
                    +(b/\rho)^2\log(\rho/r_+)\bigr].
 \end{aligned}
\]
For first derivatives, pull back $\Delta_h-\Delta$ together;
its coefficient difference and its derivative retain the factor
$r^q$. The relative $C^1(C^{6,\alphaH}_{\rm sc})$ bound for $u-R$
therefore gives, for every total parameter direction $i$,
\[
 \|D_i f\|_{C^{4,\alphaH}_{\rm sc}(A_r)}
 \le CB\bigl(r^{q-n}+r^{2-n}+b^2r^{-n-2}\bigr).
\]
In Green's integral the radial pullback preserves each projection.
Differentiate the kernel and volume factor.
On the moving inner support, Poisson smoothing absorbs the extra
factor $l+\frac{n-2}{2}$. On the outer collar, the reflected kernel
retains $(r_+/r)^{2l+n-2}$. The kernel vanishes at $r_+$, so the endpoint
term is zero. The conormal estimates for $l=1$ and $l\ge2$ therefore
hold after one total parameter derivative, in the same trace spaces. Conversion to the physical conormal adds
$O_{C^1}(B\rho^{q+2-n})$.

Normalized by $Br_+^{2-n}$, the inner trace of $u-R$ has first harmonic
$O_{C^1}(\delta)$ and high part
$O_{C^1}(\delta^2+\lambda\delta^{-q}+\kappa_b\delta)
 \le C\kappa_b\delta$.
Here the $\lambda$ term is present only when $n\ge4q+2$;
then $\lambda\delta^{-q}\le C\delta^3$, and
$\kappa_b\ge c\delta$. Subtracting $R$ is absorbed because
$b^{(n-2)/2}\delta^{2-n}\le C\delta^2$.
Harmonic transfer therefore contributes
$O_{C^1}(Bb\rho^{1-n})$ in degree one and at most
\[
 CB\rho^{2-n}\kappa_b\delta(r_+/\rho)^2
\]
in the high modes. The degree-one source contribution is
$O_{C^1}(\rho^3+Bb\rho^{1-n})$. The contribution of degrees
$l\ge2$ is $o_{C^1}(\rho^2)$ by \eqref{eq:bridge-scale-ledger}.
For the logarithmic factor, use $b/\rho\asymp\rho^{(n+1)/(n-2)}$.
After division by $\rho^2$, its conormal contribution is at most
$C\rho^{3/2+6/(n-2)}|\log\rho|$. This bound tends to zero.
For the high harmonic transfer, division by $\rho^2$ gives
\[
 C\kappa_b b^2\delta^{-1}\rho^{-5/2}
 \le C\kappa_b\rho^{1/2+6/(n-2)+1/(2n)}=o(1).
\]
All estimates include unit $C^{6,\alphaH}$ directions in $\psi$,
so they prove the joint $C^1$ high remainder and the three
Cauchy-data identities in \eqref{eq:bridge-cap}.
\end{proof}

The internal matching is complete. The outer trace is still a free
input, and the residual coefficients $c_0,\ldots,c_n$ need not yet vanish.
We next show that variations of the bubble scale and center pair
nondegenerately with these coefficients, whereas variations of the outer
trace have zero bulk pairing. 
We now specify the normalization of the central residual in physical
coordinates. If $\mathcal F_h(u)=(-\Delta_h+c(n)\operatorname{Scal}_h)u-\frac{n(n-2)}4u^p$, then
\begin{equation}
 \mathcal F_h(u)(by)
 =b^{-(n+2)/2}V_0(y)^p\vartheta_b(y)^{-1}
                   \sum_{j=0}^nc_j\chi(y)Z_j(y).
 \label{eq:bridge-physical-residual}
\end{equation}
For a scalar variation parameter $s$, let
$\mathcal T_s=b^{(n-2)/2}V_0^{-1}\partial_s u|_X$.
The scale and center tangents are
\begin{equation}
 \mathcal T_\ell
 =\frac{\partial_\ell(V_0z)-\frac{n-2}{2}V_0z-y\cdot\nabla_y(V_0z)}{V_0},
 \qquad
 \mathcal T_{\zeta_j}=\partial_{\zeta_j}z.
 \label{eq:bridge-physical-tangents}
\end{equation}
The density factors cancel in the physical bulk pairing:
\begin{equation}
 \int\mathcal F_h(u)\,\partial_s u\dd V_h
 =\sum_{j=0}^nc_j\int\chi Z_j\mathcal T_s\dd V_{\rm rd}.
 \label{eq:bridge-density-pairing}
\end{equation}
Indeed the residual, variation and volume factors multiply to
\[
 (b^{-(n+2)/2}V_0^p\vartheta_b^{-1})(b^{-(n-2)/2}V_0)b^n\dd V_{h_b}
 =\vartheta_b^{-1}V_0^{p+1}\dd V_{h_b}=\dd V_{\rm rd},
\]
using $p+1=2n/(n-2)$.

\begin{lemma}\label{lem:bridge-pairing}
Fix $\eps,h,\rho,\chi$ and $\dd V_{\rm rd}$. For the local family
of Proposition~\ref{prop:bridge}, define $G_b$ with rows indexed by
$(\ell,\zeta_1,\ldots,\zeta_n)$ and columns indexed by
$(c_0,\ldots,c_n)$:
\[
 (G_b)_{0,j}=\int\chi Z_j\mathcal T_\ell\dd V_{\rm rd},\qquad
 (G_b)_{r,j}=\int\chi Z_j\mathcal T_{\zeta_r}\dd V_{\rm rd},
 \quad 0\le j\le n,\quad 1\le r\le n.
\]
Uniformly on the compact parameter boxes,
\[
 G_b=G_{\rm J}+o(1),
 \qquad \|G_b^{-1}\|\le C.
\]
For every $C^1$ parameter map
$s\mapsto(\ell,\zeta,A,\mathbf v,\psi)(s)$, differentiation at fixed
physical $X$ gives
\begin{equation}
 \int_{b\Omega_\delta}\mathcal F_h(u)\,D_su[\dot s]\dd V_h
 =(G_bc)\cdot
       \bigl(D_s\ell[\dot s],D_s\zeta[\dot s]\bigr).
 \label{eq:bridge-pullback-pairing}
\end{equation}
In particular, every boundary-data tangent at fixed $\ell,\zeta$ has zero bulk pairing.
\end{lemma}
\begin{proof}
On $\operatorname{supp}\chi$, Lemma~\ref{lem:first-metric-correction}
and \eqref{eq:bridge-shrinking} bound $W_1+w$, its first parameter
derivatives and its first spatial derivatives by
$C(\lambda+\kappa_b\eps)=o(1)$. Hence
\eqref{eq:bridge-physical-tangents} gives
\[
 \|\mathcal T_\ell-Z_0\|_{C^0}
 +\sum_{r=1}^n\|\mathcal T_{\zeta_r}-Z_r\|_{C^0}=o(1).
\]
By Lemma~\ref{lem:sphere-augmented-inverse}, $G_{\rm J}$ is uniformly
positive definite. Hence $G_b=G_{\rm J}+o(1)$ is uniformly invertible.

Let $(\dot A,\dot{\mathbf v},\dot\psi)$ be a boundary-data direction.
At fixed $\ell,\zeta$, the fields $u_\zeta,W_1$, the measure,
the cutoff and $Z_j$ are independent of the boundary data.
Differentiating the full moment condition \eqref{eq:bridge-projection}
therefore gives
\begin{equation}
\begin{aligned}
 &\int\chi Z_jD_{(A,\mathbf v,\psi)}z
                 [\dot A,\dot{\mathbf v},\dot\psi]\dd V_{\rm rd}\\
 &\qquad=D_{(A,\mathbf v,\psi)}
       \left(\int\chi Z_jw\dd V_{\rm rd}\right)
                 [\dot A,\dot{\mathbf v},\dot\psi]=0.
\end{aligned}
 \label{eq:cap-tangent-zero-moments}
\end{equation}
For the prescribed parameter map, the chain rule at fixed physical
$X$ yields
\[
 b^{(n-2)/2}V_0^{-1}D_su[\dot s]
 =\mathcal T_\ell D_s\ell[\dot s]
 +\sum_{r=1}^n\mathcal T_{\zeta_r}D_s\zeta_r[\dot s]
 +D_{(A,\mathbf v,\psi)}z[D_s(A,\mathbf v,\psi)[\dot s]].
\]
Substitution into \eqref{eq:bridge-density-pairing}, followed by
the identities in \eqref{eq:cap-tangent-zero-moments}, gives
\eqref{eq:bridge-pullback-pairing}.
\end{proof}

\section{Renormalized action and global matching}\label{sec:high-exterior}
The family in Proposition~\ref{prop:bridge} depends on the bubble scale
and center and on its boundary value at $S_\rho$. We compute its action
and join the local solutions to an exterior solution. Solving the
boundary equations leaves an action of the compact background and
bubble parameters; its critical points give solutions of the Yamabe equation.

\subsection{The quadratic curvature coefficient from \texorpdfstring{\cite{part1}}{Part I}}
\label{subsec:curvature-action}

After the dilation $X=by$, the leading metric perturbation is
$b^qH_q$. We apply the corrected-action expansion
\eqref{eq:part-i-cutoff-action} to this tensor and the bubble $V_\zeta$.

\begin{lemma}\label{lem:tt-coefficient}
Assume $2\le q<\frac{n-2}{2}$. Let $H_q$ be trace-free and homogeneous
of degree $q$, with $H_qX=0$ and nonzero linearized Weyl curvature.
The quadratic metric term in the corrected action has cutoff limit
$b^{2q}\Phi(\zeta)$, where
\begin{equation}
 \Phi(\zeta)
 =-Q^{\rm I}_{(1+|\zeta|^2)^{-1},\,\zeta/(1+|\zeta|^2)}(H_q,H_q).
 \label{eq:tt-action}
\end{equation}
This coefficient is unchanged by imposing the localized Jacobi
moments of the present construction. It is smooth, strictly negative
for every finite $\zeta$, and tends to zero as $|\zeta|\to\infty$.
\end{lemma}
\begin{proof}
For fixed $R$, apply \eqref{eq:part-i-cutoff-action} with $s=b^q$,
tensor $H_{q,R}$, and bubble $V_\zeta$. The cutoff estimate
\eqref{eq:tt-cutoff-parameter-bound} gives the stated limit.
Lemma~\ref{lem:first-metric-correction} shows that the localized
Jacobi normalization leaves the coefficient unchanged. The strict
positivity of the quadratic form of \cite{part1} gives $\Phi(\zeta)<0$.

For smoothness and decay, moving the bubble width and center into
the tensor argument and using homogeneity gives
\[
 -\Phi(\zeta)
 =(1+|\zeta|^2)^{-2q}
 Q^{\rm I}_{1,0}\bigl(H_q(\zeta+\cdot),H_q(\zeta+\cdot)\bigr).
\]
The translated tensor is a trace-free polynomial of degree at most
$q$. The general quadratic form from
Section~\ref{subsec:corrected-bubble} applies although this polynomial
need not satisfy the radial identity. Its coefficients depend
polynomially on $\zeta$ and are bounded by $C(1+|\zeta|)^q$.
Compactification maps this finite-dimensional polynomial space
continuously into $H^1$. Therefore
\[
 |\Phi(\zeta)|
 \le C(1+|\zeta|^2)^{-2q}(1+|\zeta|)^{2q}
 \le C(1+|\zeta|)^{-2q}.
\]
The same polynomial dependence proves smoothness.
\end{proof}

\subsection{The local relative action and its boundary terms}
\label{subsec:high-action}

Subtracting the action of the regular filling with the same boundary
value isolates the contribution of the singular end. The boundary data
are still independent parameters.

For regular fillings $R_g=\Theta R_h$ with
$R_g|_{S_\rho}=u_g|_{S_\rho}$, \eqref{eq:smooth-conformal-boundary} gives
\begin{equation}
 S^g_{\eps,\rho}(u_g)-J_g(R_g;B_\rho)
 =S^h_{\eps,\rho}(u_h)-J_h(R_h;B_\rho).
 \label{eq:relative-action-conformal}
\end{equation}

Fix $\eps$ and the physical outer radius $\rho$. We use the parameters
\[
 (\ell,\zeta,A,\mathbf v,\psi),\qquad
 \ell=\log b,\quad \psi\in C^{6,\alphaH}_{\ge2}(\Sph^{n-1}),
\]
and the induced inner traces from Proposition~\ref{prop:bridge}.
Choose a common convex product domain in
$\bigl(\ell-\sigma_b^{-1}\log\eps,\zeta,A,\mathbf v,\psi\bigr)$.
Its finite-dimensional factors are fixed compact intervals or balls,
with $A\ge c_*>0$; its trace factor is a fixed bounded H\"older ball.
Choose the domain within the local construction range and containing
the boundary data induced by the fixed compact parameter set for $\bgpar$.

The $C^1$ norms use $(\ell,\zeta,A,\mathbf v,\psi)$ at fixed $\eps$ and $S_\rho$.
First derivatives are evaluated in unit directions of
$\R\times\R^n\times\R\times\R^n\times C^{6,\alphaH}_{\ge2}$
and include the induced variations of the inner traces.
Pulled-back fields are differentiated at fixed reference coordinates.
For the physical field $u$, the dilation formula following
\eqref{eq:bridge-scale-ledger} converts these derivatives to derivatives
at fixed $X$ in \eqref{eq:high-local-action-differential}.
All bounds are uniform on these finite parameter sets and the bounded trace ball.

We apply Lemma~\ref{lem:end-first-variation}, the fixed-outer-sphere
version of the subtraction from \cite[Lemma 3.3]{part2}.
The limit towards the puncture is taken at fixed $\eps$; the subsequent
estimates concern $\eps\to0$. The phase used in the subtraction is the limiting phase of the
constructed end, including its dependence on the parameters.

\begin{proposition}\label{prop:bridge-action}
For the family of Proposition~\ref{prop:bridge}, let $R$ be the regular
filling of its outer trace. Fix $\eps$ and $\rho$, and choose one
$C^1$ real lift of the asymptotic Fowler phase on the parameter
domain specified above. The subtraction \eqref{eq:end-action}
defines $S_{\eps,\rho}(u)$ as a $C^1$ function of
$(\ell,\zeta,A,\mathbf v,\psi)$, where $\ell=\log b$.
There is a constant $C_\eps$, independent of
$(\ell,\zeta,A,\mathbf v,\psi)$, such that, when $n\ge4q+2$,
\begin{equation}
 I_{\rm rel}:=S_{\eps,\rho}(u)-J_h(R;B_\rho)-C_\eps
 =-|\Sph^{n-1}|h_\eps\ell+b^{2q}\Phi(\zeta)
       -(n-2)|\Sph^{n-1}|Ab^{(n-2)/2}C_\infty(\zeta)+o_{C^1}(\eps^2).
 \label{eq:bridge-action}
\end{equation}
For $n<4q+2$, omit $b^{2q}\Phi(\zeta)$.
The remainder uses the preceding $C^1$ parameter norm and is uniform
as $\eps\to0$ subject to \eqref{eq:bridge-scales}.
The exact first variation of $S_{\eps,\rho}$ is
\begin{equation}
 DS_{\eps,\rho}(u)[\dot u]
  =\sum_{j=0}^nc_j\int_{\Sph^n}\chi Z_j\dot z_{\rm phys}\,\dd V_{g_{\rm rd}}
    +\int_{S_\rho}\partial_{\nu_h}u\,\dot u\,\dd S_h,
 \label{eq:high-local-action-differential}
\end{equation}
where $\nu_h$ points out of the end, $\dot u$ is taken at fixed
physical $X$, and $\dot z_{\rm phys}(y)=b^{(n-2)/2}V_0(y)^{-1}\dot u(by)$.
The coefficients $c_j$ are those of Lemma~\ref{lem:bridge-pairing}.
\end{proposition}

\begin{proof}
\needspace{4\baselineskip}
\noindent\textit{Step 1. Subtraction at the puncture.}

In this step, put $d=\delta/2$ and $t=-\log(r/d)$ on the tail.
Shift both Fowler profiles by $\log2$, as in
Proposition~\ref{prop:inner-dn}. Keep the notation $v,v_\tau$
for these shifted profiles and use $v$ for the fixed weights.
This fixed translation gives uniformly equivalent reference norms.
We express the bounds using $\delta=2d$.
We verify the density and boundary assumptions of Lemma~\ref{lem:end-first-variation}.
The weighted estimates \eqref{eq:radial-terminal-estimate} give
\[
 \alpha_\infty=\alpha(0)+\int_0^\infty\alpha'(t)\dd t,
 \qquad
 \int_0^\infty\|\alpha'(t)\|_{C^1_{\rm par}}\dd t<\infty.
\]
Write the shifted reference profile as
  $v(t)=v_\eps(t+\tau_{\rm ref})$.
  The inner coordinates in \eqref{eq:finite-gauge-maps} satisfy
  \[
   X=b\Phi_1(t,\theta),\qquad
   |\Phi_1(t,\theta)|
   =d e^{-t}\bigl[1+O(r|\kelpar|+\lambda r^q)\bigr],
   \qquad r=d e^{-t}.
  \]
  Consequently, on the physical cutoff $|X|=e^{-T}$,
  \[
   t=T+\ell+\log d+o_{C^1_{\rm par}}(1)
   \quad\text{as }T\to\infty,
  \]
  with $\eps$ fixed and uniformly in the angular variable.
  The physical asymptotic phase is therefore
  \[
   \tau_{\rm phys}
   =\tau_{\rm ref}+\alpha_\infty+\ell+\log d.
  \]
  For fixed $\eps$, the corresponding change of cutoff contributes
  $o_{C^1_{\rm par}}(1)$ to the truncated action and the endpoint
  terms, by the bounds on the cylindrical density and conormal.
The difference in \eqref{eq:terminal-quadratic-error}
tends to zero in $C^1$. The difference in
\eqref{eq:terminal-conormal-error} tends to zero uniformly in unit
parameter directions. 

Let $\mathcal D_b(y)=by$ and $h_b=b^{-2}\mathcal D_b^*h$.
For the physical conformal factor $u$ on
$B_\rho\setminus\{0\}$, define
\[
 \widehat u(y)=b^{(n-2)/2}u(by),
 \qquad 0<|y|\le\rho/b.
\]
For the physical Fowler phase $\tau_{\rm phys}$, use
$\tau_{\rm phys}-\ell$ for $\widehat u$, where $\ell=\log b$.
The physical cutoff $|X|=e^{-T}$ corresponds to
$|y|=e^{-(T+\ell)}$, and finite-domain scaling gives
\[
 J_h(u;e^{-T}<|X|<\rho)
 =
 J_{h_b}(\widehat u;e^{-(T+\ell)}<|y|<\rho/b).
\]
In the two subtractions prescribed by \eqref{eq:end-action},
the Fowler profile has the same argument
$T+\tau_{\rm phys}=(T+\ell)+(\tau_{\rm phys}-\ell)$.
Only the term $-h_\eps T$ changes. Passing to the limit gives
\begin{equation}
 S^h_{\eps,\rho}(u)
 =
 S^{h_b}_{\eps,\rho/b}(\widehat u)
 -|\Sph^{n-1}|h_\eps\ell.
 \label{eq:renormalized-dilation}
\end{equation}

We estimate the error after the subtraction in \eqref{eq:end-action}.
Let $\mathcal W,K,\mathcal G$ be the maxima, for $j=0,1$, of
$\mathcal W_6(D_{\rm par}^jw)$,
$\|D_{\rm par}^jz\|_{\mathcal K_{\rm tail}}$ and
$\|D_{\rm par}^j\radsrc\|_{\mathcal Y_{\rm rad}}$ along the induced family of boundary data at $t=0$.
Transfer the central trace through
\eqref{eq:collar-trace-transfer}, then apply the finite trace coordinates
of Proposition~\ref{prop:inner-dn}.
Equations~\eqref{eq:collar-solution-c1} and \eqref{eq:bridge-shrinking}
give the induced scaled tail boundary data and their total first
derivatives of size $O(\kappa_b+\delta)$. Equation~\eqref{eq:small-cap-parameter-bounds}, including
the finite background variations, and Lemma~\ref{lem:radial-response} give
\[
 \max\{\mathcal W,K,\mathcal G\}
 \le C(\kappa_b+\delta+\delta^{q-3/2})\longrightarrow0.
\]
The phase and energy estimates in \eqref{eq:radial-terminal-estimate}
hold with right side $C\mathcal G$, for both their values and total
first derivatives along this family.
The radial chart identifies the energy defect as
$\mathsf h(V,V')-h_\eps=|h_\eps|\energyvar$ and gives
\[
 \int|\mathsf h(V,V')-h_\eps|\le C\int(1+t)v|\radsrc|,
 \qquad \int v^2|\alpha'|\le C\eps^2\|\radsrc\|_{\mathcal Y_{\rm rad}}.
\]
The radial cylindrical action density, before multiplication by
$|\Sph^{n-1}|$, is $(V')^2-\mathsf h(V,V')$.
Recall that $v(t)=v_\eps(t+\tau_{\rm ref})$.
The phase and energy coordinates in
Lemma~\ref{lem:radial-response} give
\[
 (V,V')=Q(t+\alpha,\energyvar),\qquad
 Q(t+\alpha,0)=(v(t+\alpha),v'(t+\alpha)).
\]
Consequently,
\[
\begin{aligned}
 &(V')^2-\mathsf h(V,V')
 -\frac{d}{dt}\left[
   \int_0^{t+\tau_{\rm ref}+\alpha(t)}
       (v_\eps'(s))^2\,\dd s-h_\eps t
   \right]\\
 &\qquad=(V')^2-[v'(t+\alpha)]^2
       -|h_\eps|\energyvar-[v'(t+\alpha)]^2\alpha'.
\end{aligned}
\]
The chart estimates and \eqref{eq:radial-chart} imply
\[
 |V'-v'(t+\alpha)|
       \le C\eps^2|\energyvar|/v,\qquad
 |V'|+|v'(t+\alpha)|\le Cv,\qquad
 v^2|\alpha'|\le C\eps^2|\energyvar|+Cv|\radsrc|.
\]
Thus the displayed difference is bounded by
$C\eps^2|\energyvar|+Cv|\radsrc|$.
Its first parameter derivatives are bounded by
\[
 C\eps^2\bigl(|\energyvar|+|D_{\rm par}\energyvar|\bigr)
 +Cv\bigl(|\radsrc|+|D_{\rm par}\radsrc|\bigr).
\]
At the upper endpoint, replacing $\alpha(t)$ by $\alpha_\infty$
in this integral changes its value by a quantity tending to zero in
$C^1$ as $t\to\infty$, with $\eps$ fixed. The resulting upper endpoint
agrees with the limiting-phase subtraction in \eqref{eq:end-action}
after the coordinate change recorded above. The lower-endpoint term
$\int_0^{\tau_{\rm ref}+\alpha(0)}(v_\eps'(s))^2\dd s$
is retained in the boundary calculation below. The quadratic
endpoint term is handled by \eqref{eq:terminal-quadratic-error}.

The first radial metric action vanishes by
\eqref{eq:metric-tensor-green-identity} with $\psi=u/2$ and the
degree-zero identities for $T_\kelpar$ in Appendix~\ref{app:coordinates}:
\[
 \begin{gathered}
 uu'\fint (T_\kelpar)_{rr}=0,\qquad
 u^2\fint (\operatorname{div}T_\kelpar)_r
 =u^2(\partial_r+n/r)\fint (T_\kelpar)_{rr}=0,\\
 -\frac12\int T_\kelpar(\dd u,\dd u)
 +\frac{c(n)}2\int\operatorname{div}\operatorname{div}T_\kelpar\,u^2=0.
 \end{gathered}
\]
For the removed gauge tangent, \eqref{eq:gauge-identity} includes
both field and domain variations. The two values of
$\mathcal B_{X,\sigma_{\rm g}}(u)$ cancel on each shared physical
cut, with opposite orientations. On the shrinking inner sphere \(r=d e^{-t}\), 
as \(t\to\infty\),
$|\mathcal B_{X,\sigma_{\rm g}}(u)|\le
C\lambda r^q(v^2+|v'|^2)\to0$. At \(r=d\)
 its absolute value is
$O(\lambda\delta^{q+n-2})=O(\eps^2\lambda\delta^{q-2})$.
We integrate the variable-center action once in first derivatives
and retain $\mathcal B_{X,\sigma_{\rm g}}(u)$ at both endpoints.

Put $B_M=(1+t)^{M/2}$ and $r=d e^{-t}$. The weighted norms give,
for $0\le j\le1$,
\begin{equation}
 \begin{aligned}
 |D^{\le1}D_{\rm par}^jw|&\le C\eps B_M^{-1}\mathcal W,\\
 rF|D_{\rm par}^jz'|&\le C\eps B_M^{-1}K,\\
 v|D_{\rm par}^j\radsrc|+\eps^2|D_{\rm par}^j\energyvar|
       &\le C\eps^2B_M^{-2}\mathcal G,\\
 |D_{\rm par}^j\alpha'|&\le CB_M^{-2}\mathcal G.
 \end{aligned}
 \label{eq:tail-weighted-pointwise}
\end{equation}
The reference weights remain fixed under differentiation. The metric
terms in \eqref{eq:repair-raw-density} satisfy
\[
\begin{aligned}
 \lambda^2r^{2q}v^2
 &\le C\eps^2\delta^{2q}e^{-2qt},\\
 \lambda r^qv|D^{\le1}w|
 &\le C\eps^2\delta^qe^{-qt}B_M^{-1}\mathcal W,\\
 \lambda r^{q+1}v^2|z'|
 &\le C\eps^2\delta^qe^{-qt}B_M^{-1}K.
\end{aligned}
\]
The first two inequalities use $\lambda\le C\eps$ and boundedness
of $v$. For the third, \eqref{eq:tail-weighted-pointwise} gives
$rF|z'|\le C\eps B_M^{-1}K$. Moreover,
$F(\frac{n-2}{2}v+v')=\frac{(n-2)^2}{4}v^{p+1}-2h_\eps$
gives $F\ge cv^p$. Hence $v^2/F\le Cv^{2-p}\le C$.
The same inequalities hold after one parameter derivative by the
relative background bounds and \eqref{eq:tail-weighted-pointwise}.

After the exact phase, center and boundary cancellations, the radial
energy and forcing contribute $C\eps^2|\energyvar|+Cv|\radsrc|$.
The terms quadratic in $w$ and $z'$ contribute
$C|D^{\le1}w|^2$ and $Cr^2F^2|z'|^2$, respectively.
For their mixed product, use
$2|D^{\le1}w|\,rF|z'|\le |D^{\le1}w|^2+r^2F^2|z'|^2$.
Thus, apart from the phase square, the density and its first parameter
derivatives have the common majorant
\begin{equation}
 C\eps^2\bigl[
 \delta^{2q}e^{-2qt}
 +\delta^q(\mathcal W+K)e^{-qt}B_M^{-1}
 +(\mathcal G+\mathcal W^2+K^2)B_M^{-2}\bigr].
 \label{eq:tail-common-majorant}
\end{equation}

For the phase square, the weighted $L^2$ estimate in
\eqref{eq:radial-terminal-estimate} and Cauchy--Schwarz give
\begin{equation}
 \int_T^\infty v^2
       |D_{\rm par}^{j_1}\alpha' D_{\rm par}^{j_2}\alpha'|\dd t
 \le C\eps^2\mathcal G^2(1+T)^{-M},\qquad j_1+j_2\le1.
 \label{eq:tail-phase-square}
\end{equation}
Since $M\ge6$, \eqref{eq:tail-common-majorant} and
\eqref{eq:tail-phase-square} give $C^1$ convergence in
\eqref{eq:end-action} as the inner cutoff approaches the puncture,
with $\eps$ fixed. After division by $\eps^2$, the integrals
over $t\ge T$, together with their first parameter derivatives,
tend to zero uniformly in the stated parameters as $T\to\infty$.
 Integration and
Young's inequality, using $\mathcal G\le1$, give the tail error
\begin{equation}
 C\eps^2\bigl[\mathcal G+\mathcal W^2+K^2+\delta^{2q}\bigr]
       =o(\eps^2)
 \label{eq:tail-action-error}
\end{equation}
in $C^1$. The induced variations of the boundary data at $t=0$ are included:
\eqref{eq:bridge-shrinking} gives
$\delta^{(n-2)/2}\|D\phi\|_{C^{6,\alphaH}}\le C\kappa_b\eps$, and
the inner conormal component of
\eqref{eq:transmission-first-parameter-source} obeys the same
scaled boundary estimate.
It remains to estimate the action contribution at the tail boundary
$t=0$ from the phase and center parameters.
Use the physical trace
$t_{\vartheta,\trlow}$ in the $y=X/b$ coordinates of
\eqref{eq:collar-trace-coordinate-bounds}. Fix the outer boundary data.
The middle-region corrections $w$ and $0$ induce collar boundary data
$\trpar_w$ and $\trpar_0$ at $|y|=\delta$, with
\[
 \trpar_w-\trpar_0=s_\eps^{-1}V_0(\delta)\operatorname{tr}_-w,
 \qquad
 \|\trpar_w-\trpar_0\|_{C^1_{\rm par}(\Xtrace)}\le C\kappa_b.
\]
The bound follows from \eqref{eq:bridge-shrinking}.
The collar transfers these data to its inner boundary:
\[
 \begin{aligned}
 \widehat\trpar_w&=\operatorname{tr}_-
       W(\vartheta,\phi_\vartheta^{\rm ref}+\trpar_w),&
 \trlow^{(w)}&=\trlow(\vartheta,\widehat\trpar_w),\\
 \widehat\trpar_0&=\operatorname{tr}_-
       W(\vartheta,\phi_\vartheta^{\rm ref}+\trpar_0),&
 \trlow^{(0)}&=\trlow(\vartheta,\widehat\trpar_0).
 \end{aligned}
\]
Write $W=\mathsf H\phi+\mathcal R_W$.
By \eqref{eq:collar-solution-c1}, each composition
$\mathcal R_W(\vartheta,\phi_\vartheta^{\rm ref}+\trpar_w)$ and
$\mathcal R_W(\vartheta,\phi_\vartheta^{\rm ref}+\trpar_0)$ is
$O_{C^1_{\rm par}(\Zcollar)}(\delta)$.
The bounded inner trace operator therefore gives
\[
 \|\widehat\trpar_w-\widehat\trpar_0\|_{C^1_{\rm par}(\Xtrace)}
 \le C(\kappa_b+\delta)\le C\kappa_b,
 \qquad
 \|\widehat\trpar_0\|_{C^1_{\rm par}(\Xtrace)}\le C\delta.
\]
The last bound uses $\|\trpar_0\|_{C^1(\Xtrace)}+
\|\phi_\vartheta^{\rm ref}\|_{C^1(\Xtrace)}\le C\delta$.

The explicit map
$(\vartheta,\trpar)\mapsto s_\eps^{-1}
t_{\vartheta,\trlow(\vartheta,\trpar)}$
has bounded $D_\trpar$, $D_\vartheta D_\trpar$ and $D_\trpar^2$ derivatives.
This follows from \eqref{eq:collar-trace-coordinate-bounds} and the
finite coordinate inverse. Integrate its $\trpar$ derivative along the
segment from $\widehat\trpar_0$ to $\widehat\trpar_w$ and differentiate once
in the parameters. The preceding estimates give
\[
 \|t_{\vartheta,\trlow^{(w)}}
          -t_{\vartheta,\trlow^{(0)}}\|_{C^1_{\rm par}(\Xtrace)}
 \le Cs_\eps\kappa_b\le C\delta\kappa_b.
\]
The physical phase increment is $\delta b_{{\rm tr},0}$.
The derivative bounds for the displayed trace map and
\eqref{eq:collar-solution-c1} include its dependence on the collar solution.

Pull the traces back to the fixed inner interface $\widehat r=1/2$,
corresponding to the tail boundary $r=d=\delta/2$.
Let $\mathcal N_{{\rm bg},-}$ be the physical scaled conormal there,
before division by $s_\eps$. The tail conormal formula gives
\[
 \begin{aligned}
 \mathcal N_{{\rm bg},-}
 &=\tfrac12\partial_{\widehat r}B_\vartheta|_-
       +s_\eps[D_+\widehat\trpar_0+E(\vartheta,\widehat\trpar_0)]\\
 &=\tfrac12\partial_{\widehat r}B_\vartheta|_-
       +O_{C^1_{\rm par}(\Ytrace)}(\delta^2)
 =O_{C^1_{\rm par}(\Ytrace)}(\delta).
 \end{aligned}
\]
The first equality uses the physical density on this same interface.
Along $\trpar_s=(1-s)\widehat\trpar_0+s\widehat\trpar_w$, $0\le s\le1$,
the full scaled conormal is bounded in $C^1_{\rm par}(\Ytrace)$ by
$C\delta(1+\kappa_b)$.
The finite trace derivative
$\partial_s t_{\vartheta,\trlow(\vartheta,\trpar_s)}$
has $C^1_{\rm par}(\Xtrace)$ norm at most $C\delta\kappa_b$.
Integrate the boundary first variation along this segment.
The area and radial factor is $d^{n-2}\asymp\delta^{n-2}$.
Thus the linear and quadratic boundary contributions are bounded in
$C^1$ by
\[
 C\delta^{n-2}\bigl[\delta(\delta\kappa_b)
                       +(\delta\kappa_b)^2\bigr]
 =C\eps^2(\kappa_b+\kappa_b^2)=o(\eps^2).
\]

\needspace{4\baselineskip}
\noindent\textit{Step 2. Fowler and homoclinic end comparison.}

We retain the full uncut bubble action, including its inner-ball
term of order $\delta^n=\eps^2$.
Use $t=-\log r$ and $t_0=-\log\delta$. Let $v_{\rm hom}$
be the cylindrical profile of a centered smooth bubble; its
Hamiltonian is zero. In this step, $v_\eps$ and $v_{\rm hom}$
denote translates chosen so that
\[
 \begin{gathered}
 v_\eps(t_0)=v_{\rm hom}(t_0)\asymp\delta^{(n-2)/2},\\
 v_\eps'(t_0)<0,\qquad v_{\rm hom}'(t_0)<0.
 \end{gathered}
\]
For $\Delta v=v_\eps-v_{\rm hom}$, the energies give
\[
 \Delta v(t_0)=0,\qquad
 |\Delta v'(t_0)|
 =\frac{2|h_\eps|}{|v_\eps'(t_0)+v_{\rm hom}'(t_0)|}
 \le C\delta^{(n+2)/2}.
\]
On $|t-t_0|\le C\delta$, $|v'|\le \frac{n-2}{2}v$ keeps both profiles
comparable to $\delta^{(n-2)/2}$. Their difference satisfies
\[
 (\Delta v)''=
 \left[\frac{(n-2)^2}{4}
 -\frac{n(n+2)}4\int_0^1(v_{\rm hom}+s\Delta v)^{p-1}\dd s\right]\Delta v.
\]
The coefficient is uniformly bounded. Integrating twice gives
$\|\Delta v\|_\infty\le C\delta^{(n+4)/2}
+C\delta^2\|\Delta v\|_\infty$. Absorbing the last term and
then using the ODE yields
\[
 \|\Delta v\|_\infty+\|\Delta v''\|_\infty
 \le C\delta^{(n+4)/2},\qquad
 \|\Delta v'\|_\infty\le C\delta^{(n+2)/2}.
\]

Let $\kelpar$ range over the fixed compact center set; constants in
this step may depend on that set. Put
\[
 D=1-2r\kelpar\cdot\theta+r^2|\kelpar|^2,\qquad F=\frac{n-2}{2}v-v',\quad G=\frac{n-2}{2}v+v'.
\]
At $r=\delta$,
$|\tfrac12\log D|+|D_t/(2D)|+|D_i/(2D)|\le C\delta$
for $1\le i\le n$, where $D_t=\partial_tD$ and
$D_i=\partial_{\kelpar_i}D$. In the following formulas, $v$
denotes either $v_\eps$ or $v_{\rm hom}$. The Kelvin-transformed
cylindrical factor is $\mathsf V$, its parameter derivative is
$\mathsf V_i=\partial_{\kelpar_i}\mathsf V$, and
$\mathsf G=\frac{n-2}{2}\mathsf V+\partial_t\mathsf V$:
\[
\begin{aligned}
 \mathsf V&=D^{-(n-2)/4}v(t+\tfrac12\log D),\\
 \mathsf G&=D^{-(n-2)/4}\left[G(t+\tfrac12\log D)
              -\frac{D_t}{2D}F(t+\tfrac12\log D)\right],\\
 \mathsf V_i&=-D^{-(n-2)/4}F(t+\tfrac12\log D)\frac{D_i}{2D}.
\end{aligned}
\]
Write $\Delta$ for the Fowler quantity minus its smooth-bubble
counterpart. All evaluation points satisfy $|t-t_0|\le C\delta$. Thus
\[
 |\Delta\mathsf V_i|\le C\delta^{(n+4)/2},\qquad
 |\Delta\mathsf G|\le C\delta^{(n+2)/2}.
\]
Also $F\asymp\delta^{(n-2)/2}$. The Hamiltonian identity
\[
 FG=\frac{(n-2)^2}{4}v_\eps^{p+1}-2h_\eps
\]
for the Fowler profile, and the same identity with $v_\eps$ replaced by $v_{\rm hom}$
and $h_\eps$ replaced by zero, give $G=O(\delta^{(n+2)/2})$
for either profile.
Consequently $|\mathsf V_i|+|\mathsf G|\le C\delta^{n/2}$.
Let $J_\eps^{\rm ren}(\kelpar)$ be the renormalized Euclidean
action, defined by \eqref{eq:end-action}, of the Kelvin transform
of the Fowler solution on $0<r<\delta$.
Let $J_0(\kelpar)$ be the ordinary Euclidean action on $B_\delta$
of the Kelvin transform of the corresponding smooth bubble.
Set
\[
 E(\kelpar)=J_\eps^{\rm ren}(\kelpar)-J_0(\kelpar).
\]
The boundary first-variation formulas give
\[
\begin{aligned}
 \partial_{\kelpar_i}J_\eps^{\rm ren}
 &=-\int_{\Sph^{n-1}}
       \mathsf G_\eps\mathsf V_{\eps,i}\,\dd\theta,\\
 \partial_{\kelpar_i}J_0
 &=-\int_{\Sph^{n-1}}
       \mathsf G_0\mathsf V_{0,i}\,\dd\theta.
\end{aligned}
\]
Here the subscripts $\eps$ and $0$ refer to the Fowler and
smooth-bubble profiles, respectively, and all boundary quantities
are evaluated at $r=\delta$. Therefore,
\[
 \partial_{\kelpar_i}E=-\int_{\Sph^{n-1}}
   (\Delta\mathsf G\,\mathsf V_{\eps,i}
                   +\mathsf G_0\,\Delta\mathsf V_i)\dd\theta,
 \qquad |\partial_{\kelpar_i}E|\le C\delta^{n+1}=C\eps^2\delta.
\]
Integration along the segment from $0$ to $\kelpar$ proves
$\|E-E(0)\|_{C^1}\le C\eps^2\delta=o(\eps^2)$.
After extracting the logarithmic scale term in
\eqref{eq:renormalized-dilation}, the reference value $E(0)$ depends
only on $(\eps,\delta)$ and is included in $C_\eps$.
The difference $E-E(0)$ enters the $o_{C^1}(\eps^2)$ remainder.
For the local family of Proposition~\ref{prop:inner-graph},
\eqref{eq:tail-action-error} controls the tail error and the derivatives
induced by its boundary traces.

\needspace{4\baselineskip}
\noindent\textit{Step 3. The outer interaction.}

The leading harmonic term of $V_\zeta$ at infinity is
\[
 H_\zeta(y)=C_\infty(\zeta)
       \left|y-\frac{\zeta}{1+|\zeta|^2}\right|^{2-n}.
\]
Put $w_{\rm harm}=(1-\chi_{\delta^{-1}})(H_\zeta-V_\zeta)$.
Since
\[
 \|\nabla D_{\rm par}^jw_{\rm harm}\|_2\le C\delta^{(n+2)/2}
 \quad(j\le1),\qquad DJ_{g_{\rm E}}(V_\zeta;\R^n)[w_{\rm harm}]=0,
\]
Taylor's formula gives
$\|J_{g_{\rm E}}(V_\zeta+w_{\rm harm};\R^n)-J_{g_{\rm E}}(V_\zeta;\R^n)\|_{C^1}
\le C\delta^{n+2}$.
Let $U=U_{\varrho,\xi}$ be the bubble in physical coordinates and
let $R$ be the regular filling in Proposition~\ref{prop:bridge-action}.
On the outer annulus $A_{r_+,\rho}$, put $E=u-U-R$. The annular construction gives
$u=U_0+v$, with $U_0$ harmonic and
$\|D_{\rm par}^jv\|_{6,H}\le C\vartheta_{\rm out}$, $j=0,1$.
Since $\rho\le C\delta$ and $q\ge2$, integration of its relative
value and gradient bounds gives
\[
 \left(\int_{A_{r_+,\rho}}
       (|\nabla D_{\rm par}^jv|^2+r^{-2}|D_{\rm par}^jv|^2)\right)^{1/2}
 \le C\vartheta_{\rm out}(\rho^{(n-2)/2}+Br_+^{-(n-2)/2})\le C\eps\delta.
\]
Here the harmonic extension of $U$ or $R$ means the Euclidean
harmonic function on $A_{r_+,\rho}$ with the same values on both
boundary spheres. The differences from $U$ and $R$ have norm at
most $C\eps\delta$ in the energy norm displayed above, including
one parameter derivative. These bounds follow from the bubble
estimate at large radius and Lemma~\ref{lem:regular-filling},
respectively, using the Dirichlet principle and Hardy inequality
to pass to the annular harmonic extensions.
The harmonic lift of the trace mismatch has weighted energy bounded by
$C[\kappa_b\eps+\lambda\delta^{(n-2)/2-q}+\eps\delta+B\rho^{-(n-2)/2}]$.
The traces contributed by $w$ and $W_1$ on the inner sphere
$S_{r_+}$ of this annulus have harmonic-lift norms
$C\kappa_b\eps$ and $C\lambda\delta^{(n-2)/2-q}$, respectively.
The trace $R|_{S_{r_+}}$ has lift norm
$Cr_+^{(n-2)/2}\le C\eps\delta$.
On the outer sphere $S_\rho$, the mismatch is $-U|_{S_\rho}$
and has lift norm
$CB\rho^{-(n-2)/2}$.
Thus, for $j=0,1$,
\begin{equation}
 \left(\int_{A_{r_+,\rho}}
       (|\nabla D_{\rm par}^jE|^2+r^{-2}|D_{\rm par}^jE|^2)\right)^{1/2}
 \le C[\eps(\kappa_b+\delta)+\lambda\delta^{(n-2)/2-q}+B\rho^{-(n-2)/2}]
 =o(\eps).
 \label{eq:outer-action-comparison}
\end{equation}
For moving cuts, \eqref{eq:outer-action-comparison} uses the
fixed-annulus pullback in Lemma~\ref{lem:bridge-transmission}.
For an annular integrand $f$, differentiation at fixed $\rho$ gives
\[
 D_{\rm par}\int_{A_{r_+,\rho}}f\dd X
 =\int_{A_{r_+,\rho}}D_{\rm par}f\dd X
 -r_+^{n-1}D_{\rm par}r_+
       \int_{\Sph^{n-1}}f(r_+,\theta)\dd\theta.
\]
The trace and scaled conormal estimates of
Lemma~\ref{lem:bridge-transmission}, with $|D_{\rm par}r_+|\le Cr_+$,
bound this boundary contribution by the products of the factors in
\eqref{eq:outer-action-comparison}.

Taylor expansion about the actual solution, with its boundary
variation retained, now gives
\[
 J_h(u;A_{r_+,\rho})-J_h(U+R;A_{r_+,\rho})
 =\int_{\partial A_{r_+,\rho}}\partial_{\nu_h}u\,E\,\dd S_h
       +o_{C^1}(\eps^2).
\]
The quadratic bound follows from \eqref{eq:outer-action-comparison}
and $r^2(U+R)^{p-1}\le C(\rho^2+\delta^2)$; relative first
parameter derivatives preserve this inequality.
On the middle region, expansion about $u_\zeta+W_1+w$ has
bulk pairing $\sum_jc_j\int\chi Z_jw\,\dd V_{\rm rd}=0$ by the
exact full moments. Differentiating the zero moment identity also
makes its total parameter derivative zero.

Combine the linear boundary integrals on $S_{r_+}$ before estimating
their sum. Use the outward normal $\nu_h$ of the central domain and
the common physical measure $\dd S_h$. For $X=by$, set
$W_{1,\rm phys}(X)=b^{-(n-2)/2}V_0(y)W_1(y)$ and
$w_{\rm phys}(X)=b^{-(n-2)/2}V_0(y)w(y)$.
Formula~\eqref{eq:smooth-conformal-boundary} gives their action
pairings in the physical metric.
The four boundary terms arise from $w_{\rm phys}$, the outer difference
$E$, the interaction of $U$ with $R$, and the linear variation in
$W_{1,\rm phys}$, respectively:
\begin{equation}
\begin{aligned}
 &\int_{S_{r_+}}\partial_{\nu_h}u\,w_{\rm phys}\dd S_h
 -\int_{S_{r_+}}\partial_{\nu_h}u
       (w_{\rm phys}+W_{1,\rm phys}-R)\dd S_h\\
 &\quad-\int_{S_{r_+}}\partial_{\nu_h}U\,R\dd S_h
 +\int_{S_{r_+}}\partial_{\nu_h}U\,W_{1,\rm phys}\dd S_h\\
 &=\int_{S_{r_+}}(\partial_{\nu_h}u-\partial_{\nu_h}U)
                         (R-W_{1,\rm phys})\dd S_h.
\end{aligned}
\label{eq:four-boundary-pairings}
\end{equation}
The trace and conormal estimates give
\[
\begin{aligned}
 r_+^{n/2}\|\partial_{\nu_h}u-\partial_{\nu_h}U
             \|_{C^1_{\rm par}(C^0)}
 &\le C(\eps\kappa_b+\lambda\delta^{(n-2)/2-q}+\eps\delta),\\
 r_+^{(n-2)/2}\|R-W_{1,\rm phys}\|_{C^1_{\rm par}(C^0)}
 &\le C(\lambda\delta^{(n-2)/2-q}+\eps\delta).
\end{aligned}
\]
Since $n-1-\frac n2-\frac{n-2}{2}=0$, the $C^1$ norm of the final integral in
\eqref{eq:four-boundary-pairings} is at most the product of these
two right-hand sides, which is $o(\eps^2)$.
When $n<4q+2$, set $W_{1,\rm phys}=0$ in this calculation.

On $S_\rho$, $E=-U$. The leading regular conormal is a first
harmonic. Its pairing with the constant leading pole is zero. The remaining regular, bubble dipole and pole contributions
are bounded by
$C(B\rho^2+Bb+B^2\rho^{2-n})=o(\eps^2)$.
This bound includes one parameter derivative and uses
$b\le C\rho^2$ from \eqref{eq:bridge-scale-ledger}.

For the nonlinear remainder, as in
\cite[the proof of Lemma 5.3, Step 1, and Lemma 4.3]{part2},
set, for $U,A>0$,
\[
 N(U,A)=(U+A)^{p+1}-U^{p+1}-A^{p+1}-(p+1)AU^p.
\]
Concavity of $s^{p-1}$ gives
\[
 0\le N_U=p(p+1)\int_0^A[(U+s)^{p-1}-U^{p-1}]\dd s
       \le(p+1)A^p.
\]
Integrating from $U=0$ and using homogeneity gives
\[
 AN_A=(p+1)N-UN_U,\qquad
 |N|+|UN_U|+|AN_A|\le CA^pU.
\]
By \eqref{eq:finite-bubble} and Lemma~\ref{lem:regular-filling},
$|D_{\rm par}U|\le CU$, $c_*\le R\le C_*$ and
$|D_{\rm par}R|\le C$, including unit trace directions. Hence
\[
\begin{aligned}
 D_{\rm par}N(U,R)&=N_U D_{\rm par}U+N_A D_{\rm par}R,\\
 |N(U,R)|+|D_{\rm par}N(U,R)|&\le CU\le CBr^{2-n}.
\end{aligned}
\]
At fixed $\rho$, differentiation of the moving annulus gives
\[
\begin{aligned}
 D_{\rm par}\int_{A_{r_+,\rho}}N(U,R)\,\dd X
 ={}&\int_{A_{r_+,\rho}}D_{\rm par}N(U,R)\,\dd X\\
 &-r_+^{n-1}D_{\rm par}r_+
       \int_{\Sph^{n-1}}N(U,R)|_{r=r_+}\,\dd\theta.
\end{aligned}
\]
Since $|D_{\rm par}r_+|\le Cr_+$ and $n-1-(n-2)=1$,
\[
 \left\|\int_{A_{r_+,\rho}}N(U,R)\,\dd X\right\|_{C^1}
 \le CB\left(\int_{r_+}^\rho r\,\dd r+r_+^2\right)
 \le CB\rho^2=o(\eps^2).
\]
Replacing $R(X)$ by $R(0)$ in the linear interaction gives
an error $O_{C^1}(Bb)$, using
$\int|y|V_\zeta^p\dd y<\infty$ and
$b\asymp\rho^{(2n-1)/(n-2)}$.

The equation $-\Delta U=\frac{n(n-2)}4U^p$ identifies the linear interaction:
\[
 J_{g_{\rm E}}(U+R;D)-J_{g_{\rm E}}(U;D)-J_{g_{\rm E}}(R;D)
 =\int_{\partial D}R\partial_{\nu_{g_{\rm E}}}U\dd S_{g_{\rm E}}
                    -\frac{(n-2)^2}{8}\int_DN(U,R).
\]
With the outward normal of the central domain, the inner boundary
term is $-\int_{S_{r_+}}R\partial_{\nu_h}U\dd S_h$ after the
physical conversion. It is the third term of
\eqref{eq:four-boundary-pairings}. The full bubble flux is
\[
 \frac{n(n-2)}4\int_{\R^n}U^p\dd X
 =(n-2)|\Sph^{n-1}|b^{(n-2)/2}C_\infty(\zeta).
\]
Hence the interaction is
$-(n-2)|\Sph^{n-1}|R(0)b^{(n-2)/2}C_\infty(\zeta)$.
After subtracting the regular filling, the metric and scalar curvature
errors, of orders $r^q$ and $r^{q-2}$, give
\[
 C\left(B\rho^q+B^2\int_{r_+}^{\rho}r^{q+1-n}\dd r\right)
       \le C(B\rho^q+B^2r_+^{q+2-n}).
\]
The nonlinear bubble density satisfies
$U^{p+1}\le C(b/r)^2r^{-2}U^2$ on $A_{r_+,\rho}$, so it gains the
factor $b^2/r_+^2=\delta^2$ relative to the quadratic density
$r^{-2}U^2$. Relative parameter bounds preserve
$C(B\rho^q+B^2r_+^{q+2-n})$. The moving-boundary contribution is
at most $C(Br_+^{q-1}+B^2r_+^{q+1-n})|D_{\rm par}r_+|$.
Since $|D_{\rm par}r_+|\le Cr_+$, this contribution is also bounded
by $C(B\rho^q+B^2r_+^{q+2-n})$.
Finally,
\[
 \frac{B^2r_+^{q+2-n}}{\eps^2}
 \asymp b^q\delta^{-q-2}
 \asymp
 \begin{cases}
 \delta^{(n-4)/2-q},&n\ge4q+2,\\
 \delta^{q+4q/(n-2)-2},&n<4q+2.
 \end{cases}
\]
Both exponents are positive: the first is at least $q-1$ and the
second at least $4q/(n-2)$. Since $B=O(\eps^2)$, the geometric error is
$o_{C^1}(\eps^2)$.

\needspace{4\baselineskip}
\noindent\textit{Step 4. The quadratic metric term.}

For $n\ge4q+2$, put $\lambda=b^q$ and use the cutoff from
Lemma~\ref{lem:compactified-metric-bounds} to define
\[
 W_{1,\delta}=\chi_{2/\delta}W_1.
\]
This field is smooth on $\Sph^n$. It equals $W_1$ near
$|y|=\delta,\delta^{-1}$ and on $\operatorname{supp}\chi$.
Thus it has the same traces, derivatives and zero localized moments.
All parameter derivatives keep $\delta$ fixed.
Lemma~\ref{lem:first-metric-correction} gives
\[
 \|W_{1,\delta}\|_{C^1_{\rm par}(L^\infty)}
       \le C\lambda\delta^{-q},\qquad
 \|W_{1,\delta}-W_1\|_{C^1_{\rm par}(H^1)}
       \le C\lambda\delta^{(n-2)/2-q}.
\]
For the second bound, integrate the squared gradient estimate
near the point of $\Sph^n$ corresponding to $|y|=\infty$:
\[
 C\lambda^2\int_0^{C\delta}s^{n-3-2q}\dd s
 \le C\lambda^2\delta^{n-2-2q}.
\]
The cutoff-gradient term has the same bound, since
$|\nabla_{g_{\rm rd}}\chi_{2/\delta}|\le C\delta^{-1}$
on the transition annulus.
The bubbles $u_\zeta$ have a uniform positive minimum on the fixed
center set. Since $\lambda\delta^{-q}\le C\delta^{n/2-q}\to0$,
\[
 u_\zeta+sW_{1,\delta}
 \ge\tfrac12\min_{\Sph^n}u_\zeta>0,\qquad 0\le s\le1.
\]
For $K_b=\gamma_b-g_{\rm rd}$,
Lemma~\ref{lem:compactified-metric-bounds} gives
$\|K_b\|_{C^1_{\rm par}(L^\infty)}\le C\eps\delta^{-q}$.
Both $K_b$ and $W_{1,\delta}$ have
$C^1_{\rm par}(H^1)$ norm at most $C\eps$.
Integrate the scalar curvature term by parts once and keep its
physical boundary term. The cubic contractions are controlled by products such as
\[
\begin{aligned}
 \int|K_b||\nabla W_{1,\delta}|^2
 &\le\|K_b\|_\infty\|W_{1,\delta}\|_{H^1}^2,\\
 \int|W_{1,\delta}||\nabla K_b|^2
 &\le\|W_{1,\delta}\|_\infty\|K_b\|_{H^1}^2.
\end{aligned}
\]
Mixed products use Cauchy--Schwarz. For example, the contraction
$T=\int K_b(\nabla W_{1,\delta},\nabla W_{1,\delta})\dd V_{\rm rd}$ obeys
\[
 |D_iT|\le C\|D_iK_b\|_\infty\|W_{1,\delta}\|_{H^1}^2
 +C\|K_b\|_\infty\|D_iW_{1,\delta}\|_{H^1}\|W_{1,\delta}\|_{H^1}.
\]
The same product rule applies to each cubic contraction.
The stated $C^1(L^\infty)$ and $C^1(H^1)$ bounds therefore give
$C\eps^3\delta^{-q}=C\eps^2\delta^{n/2-q}=o(\eps^2)$.
After applying the fixed cutoffs, the compactified Taylor
remainder beyond $b^qH_q(y)$ has $C^1_{\rm par}(H^1)$ norm
$O(b^{q+1})=O(\eps b)$, since $q+1<(n-2)/2$ in the present
range. Its quadratic pairing with the leading metric term is
therefore $O_{C^1}(\eps^2b)$.

The expansion about $u_\zeta+W_1+w$ in Step~3 has zero
bulk pairing with $w$. The extension and Hardy estimate in
Lemma~\ref{lem:bridge-transmission}, followed by
\eqref{eq:bridge-shrinking}, bound the quadratic remainder in $w$ by
\[
 C\|w\|_{C^1_{\rm par}(H^1)}^2\le C\kappa_b^2\eps^2.
\]
For the corrected bubble, retain the linear boundary term at $u_\zeta$ in the round metric:
\[
 DJ_{g_{\rm rd}}(u_\zeta;\Omega_\delta)[W_1]
 =\int_{\partial\Omega_\delta}\partial_{\nu_{\rm rd}}u_\zeta\,W_1
                      \,\dd S_{\rm rd}.
\]
We record the boundary terms when restricting the complete-space
calculation to the central domain $D_c=\{r_-<|X|<r_+\}$.
For a smooth metric $g_*$, a symmetric tensor $K$ and a hypersurface
$\Sigma$ with chosen unit normal $\nu_{g_*}$, define the oriented metric
boundary integral
\begin{equation}
\begin{aligned}
 \mathfrak M_\Sigma(g_*,K,u)
 =\frac{c(n)}2\int_\Sigma \bigl[&
 u^2(\operatorname{div}_{g_*}K-d\operatorname{tr}_{g_*}K)(\nu_{g_*})\\
 &-K(\nu_{g_*},\nabla_{g_*}u^2)
 +(\operatorname{tr}_{g_*}K)\partial_{\nu_{g_*}}u^2\bigr]\dd S_{g_*}.
\end{aligned}
\label{eq:metric-action-boundary-current}
\end{equation}
Indeed,
$D\operatorname{Scal}_{g_*}[K]=
-\langle\operatorname{Ric}_{g_*},K\rangle
+\operatorname{div}_{g_*}\operatorname{div}_{g_*}K
-\Delta_{g_*}\operatorname{tr}_{g_*}K$.
Twice integrating the last two terms against $c(n)u^2/2$ gives
\eqref{eq:metric-action-boundary-current}.
The gradient and power terms of the action have no derivatives of
$K$, so they add no metric boundary integral.

Use the smooth metric extension defining $K_b$ and set
\[
 g_s=(1-s)g_{\rm E}+sh,\qquad
 u_s(X)=U(X)+s b^{-(n-2)/2}V_0(X/b)W_{1,\delta}(X/b),
 \quad 0\le s\le1.
\]
The compactified field is $u_\zeta+sW_{1,\delta}$.
Near both physical cuts, $u_s=U+sW_{1,\rm phys}$.
The restrictions from the two sides therefore have the same metric,
field and first derivatives on each cut.
Apply \eqref{eq:metric-action-boundary-current} along this path.
The boundary contribution to the derivative of $J_{g_s}(u_s;D_c)$
on each cut is
\[
 \int_{S_{r_\pm}}\partial_{\nu_{g_s}}u_s\,\partial_su_s\dd S_{g_s}
 +\mathfrak M_{S_{r_\pm}}(g_s,\partial_sg_s,u_s).
\]
Integrate in $s$ and retain the terms of degrees one and two in the
metric perturbation and $W_{1,\rm phys}$. This gives precisely the
finite-domain boundary terms in the Taylor expansion.
At the trace-free flat linearization, the boundary expression in
\eqref{eq:metric-tensor-green-identity}, evaluated at $\psi=u/2$, satisfies
\[
 \mathfrak M_{\partial\mathcal A}(g_{\rm E},T_\kelpar,u)
 =\left.\mathcal B_{\mathcal A}\right|_{\psi=u/2}
 -\frac12\int_{\partial\mathcal A}
       uT_\kelpar(\nu_{g_{\rm E}},\nabla u)\dd S_{g_{\rm E}}.
\]
This follows by substituting $\psi=u/2$ into the three terms defining
$\mathcal B_{\mathcal A}$. In gauge coordinates the coordinate and
conformal tangent contribute the explicit expression
$\mathcal B_{X,\sigma_{\rm g}}$ in \eqref{eq:gauge-identity}.

The conformal conversion has its own exact boundary difference.
For $h_* =\Theta^{4/(n-2)}g_*$ and $u_{g_*}=\Theta u_{h_*}$, write
\[
 \mathfrak C_\Sigma(g_*,\Theta,u_{g_*})
 =-\frac12\int_\Sigma u_{g_*}^2\partial_{\nu_{g_*}}\log\Theta\dd S_{g_*}.
\]
Formula~\eqref{eq:smooth-conformal-boundary} adds
$\mathfrak C_\Sigma$ when converting the action from $g_*$ to $h_*$.
The conformal boundary terms for the corrected and uncorrected
fields differ by
\[
 -\frac12\int_\Sigma 
 \bigl(2U W_{1,\rm phys}+W_{1,\rm phys}^2\bigr)
             \partial_{\nu_{g_*}}\log\Theta\dd S_{g_*}.
\]
The metric and conformal integrals use the same physical surface and
measure as the field integral. Reversing the normal reverses all
three integrals. For the central and adjacent domains, denoted by
superscripts $c$ and $\mathrm{adj}$, the cancellation on either cut is therefore
\[
 \mathfrak M_{S_{r_\pm}}^c+\mathfrak M_{S_{r_\pm}}^{\mathrm{adj}}=0,
 \qquad
 \mathfrak C_{S_{r_\pm}}^c+\mathfrak C_{S_{r_\pm}}^{\mathrm{adj}}=0,
 \qquad
 \mathcal B_{X,\sigma_{\rm g}}^c+
             \mathcal B_{X,\sigma_{\rm g}}^{\mathrm{adj}}=0.
\]
These equalities hold for the common extended path at every $s$,
and hence for each of its Taylor coefficients.
On $S_{r_-}$ the field pairing also cancels its inner-end counterpart:
\[
 \int_{S_{r_-}}\partial_{\nu_h^c}U\,W_{1,\rm phys}\dd S_h
 +\int_{S_{r_-}}\partial_{\nu_h^{\mathrm{adj}}}U\,W_{1,\rm phys}\dd S_h=0.
\]
On $S_{r_+}$ we retain
$\int_{S_{r_+}}\partial_{\nu_h^c}U\,W_{1,\rm phys}\dd S_h$;
it is the fourth term of \eqref{eq:four-boundary-pairings}.
After adding the adjacent extended model actions, the other
oriented boundary integrals cancel. The finite-domain expansion
therefore contains the complete-space quadratic coefficient, the
field pairing on $S_{r_+}$, and the cutoff and Taylor remainders.
Steps~1--3 estimate the differences between the actual fields and
these extended model fields. At $S_{r_+}$ the outer
trace difference is $w_{\rm phys}+W_{1,\rm phys}-R$, so its linear
field pairing is exactly the second term of
\eqref{eq:four-boundary-pairings}. The metric differences involving
$R$ are bounded, together with their boundary contributions, by
$C(B\rho^q+B^2r_+^{q+2-n})$ as computed in Step~3.
Equation~\eqref{eq:outer-action-comparison} bounds the quadratic
field differences. The inner-end differences are bounded by
\eqref{eq:tail-action-error} and $C\eps^2\delta$ from Step~2.
Comparison with the actual fields therefore preserves the boundary
pairing in \eqref{eq:four-boundary-pairings} and contributes an error
$o_{C^1}(\eps^2)$.

To identify the degree-$q$ quadratic term, take $R_\delta=2/\delta$.
Let $W_{1,R_\delta}$ be the field component of
$\mathscr A_\zeta^{-1}(-\lambda\mathcal S_\zeta(K_{q,R_\delta}),0)$.
The compatibility in \eqref{eq:metric-source-compatibility} for
$K_{q,R_\delta}$ makes the augmented multiplier zero. Hence
\[
 L_\zeta W_{1,R_\delta}=-\lambda\mathcal S_\zeta(K_{q,R_\delta}),\qquad
 \int\chi Z_jW_{1,R_\delta}\dd V_{\rm rd}=0\quad(0\le j\le n).
\]
Lemmas~\ref{lem:compactified-metric-bounds}
and~\ref{lem:sphere-augmented-inverse}, including their first
parameter bounds, give
\[
 \|W_{1,R_\delta}-W_1\|_{C^1_{\rm par}(H^1)}
 +\|W_{1,\delta}-W_{1,R_\delta}\|_{C^1_{\rm par}(H^1)}
 \le C\lambda\delta^{(n-2)/2-q}.
\]
The degree-$q$ term of the metric extension agrees with
$\lambda K_{q,R_\delta}$ on $\Omega_\delta$. On the omitted caps, their difference has
$C^1_{\rm par}(H^1)$ norm at most
$C\lambda\delta^{(n-2)/2-q}$ by
Lemma~\ref{lem:compactified-metric-bounds}.
The quadratic Taylor form is continuous in the $H^1$ norms of the
metric and field. Replacing this metric term and $W_{1,\delta}$ by
$\lambda K_{q,R_\delta}$ and $W_{1,R_\delta}$ therefore changes the quadratic term
by $O_{C^1}(\lambda^2\delta^{(n-2)/2-q})$.

For fixed $R_\delta$, apply \eqref{eq:part-i-cutoff-action} to
$H_{q,R_\delta}=\chi_{R_\delta}H_q$ and the scalar first correction
computed from this same tensor. For the bubble $V_\zeta$ in
\eqref{eq:normalized-bubble-parameters}, the cutoff estimate
\eqref{eq:tt-cutoff-parameter-bound} and the definition
\eqref{eq:tt-action} give
\[
 \begin{aligned}
 &-b^{2q}Q^{\rm I}_{(1+|\zeta|^2)^{-1},\,\zeta/(1+|\zeta|^2)}
                 (H_{q,R_\delta},H_{q,R_\delta})\\
 &\qquad=b^{2q}\Phi(\zeta)
   +O_{C^1}\bigl(b^{2q}\delta^{(n-2)/2-q}\bigr)
  =b^{2q}\Phi(\zeta)+o_{C^1}(\eps^2).
 \end{aligned}
\]
Here $b^{2q}\asymp\eps^2$ in the range $n\ge4q+2$ under
consideration, and $(n-2)/2-q>0$.
The cutoff radius $R_\delta$ is fixed under parameter
differentiation, since $\eps$ and $\delta$ are fixed, and tends to
infinity as $\eps\to0$. The estimates above control the Taylor
remainder uniformly in the changing cutoff;
\eqref{eq:part-i-cutoff-action} identifies the quadratic term.

If $n<4q+2$, then $W_1=0$, $\|K_b\|_{H^1}\le Cb^q=o(\eps)$,
and $\|K_b\|_\infty=O((b/\delta)^q)=o(1)$.
The quadratic and cubic metric terms, the quadratic remainder in
$w$, and their mixed products are therefore bounded in $C^1$ by
$C(b^q+\kappa_b\eps)^2=o(\eps^2)$. In either range the remaining
profile, outer boundary and overlap errors satisfy
\[
 (\eps^4+b^{n-2})\delta^{2-n}=O(\eps^2\delta^2),\qquad
 B^2\rho^{2-n}=O(\eps^2\rho^{3/2}),\qquad
 B(\rho^2+b+\rho^q)=o(\eps^2).
\]
Each displayed remainder bound holds after one unit derivative in
$(\ell,\zeta,A,\mathbf v,\psi)$. The bounds apply after cancellation
of the gauge and moving-interface boundary terms.

\needspace{4\baselineskip}
\noindent\textit{Step 5. First variation and interface cancellation.}

Use finite-domain first variation on the assembled field in the
fixed physical annulus $e^{-T}<|X|<\rho$, then pass to the limit
by \eqref{eq:renormalized-first-variation}.
Fix a unit parameter direction. Let $\Sigma_i$ be an artificial
interface, parametrized by $X_i$ along this variation.
Let $u_+$ and $u_-$ be the fields on its two sides.
Value matching gives tangential gradient matching; physical conormal
matching gives the remaining normal component. Hence
$\nabla u_+=\nabla u_-$ on the interface. Differentiating
$u_+(X_i)=u_-(X_i)$ along this direction yields
\[
 (\dot u_+-\dot u_-)
 +(\nabla u_+-\nabla u_-)\cdot\dot X_i=0,
 \qquad \dot u_+=\dot u_-.
\]
Here the dots on the fields mean differentiation at fixed physical
$X$. For the opposite outward normals $\nu_+=-\nu_-$, the two
conormal integrals therefore satisfy
\[
 \int_{\Sigma_i}\bigl(\partial_{\nu_+}u_+\dot u_+
                   +\partial_{\nu_-}u_-\dot u_-\bigr)\dd S_h=0.
\]
Write $e_h(u)=\tfrac12(|\nabla u|_h^2+c(n)\operatorname{Scal}_hu^2)
-\tfrac{(n-2)^2}{8}u^{p+1}$ for the energy density.
The moving-domain contributions cancel separately:
\[
 \int_{\Sigma_i}\bigl[e_h(u_+)-e_h(u_-)\bigr]
                    (\dot X_i\cdot\nu_+)\dd S_h=0.
\]
The equality follows from the matching values and full gradients.
Consequently the weak Euler residual has no interface measure.
The gauge and conformal boundary terms also occur on these shared
physical surfaces with opposite orientations.
Step 1 verifies the density bound and the two convergence assumptions for the boundary integrals in Lemma~\ref{lem:end-first-variation}. Lemma~\ref{lem:bridge-pairing} identifies
the remaining central bulk term, proving
\eqref{eq:high-local-action-differential} with its full outer
conormal term. Formula~\eqref{eq:smooth-conformal-boundary}
expresses \eqref{eq:high-local-action-differential} in the global metric $g$.

Since $R(0)=A+O_{C^1}(\rho^2)$, the interaction in Step~3 is
\[
 -(n-2)|\Sph^{n-1}|R(0)b^{(n-2)/2}C_\infty
 =-(n-2)|\Sph^{n-1}|Ab^{(n-2)/2}C_\infty+O_{C^1}(B\rho^2).
\]
The error is $o_{C^1}(\eps^2)$.

Fix the reference phase and action at one point of the common convex
product domain chosen before Proposition~\ref{prop:bridge-action}. The action of the untruncated bubble, including its
$\delta^n=\eps^2$ inner-ball term, and the regular filling subtraction
fix $C_\eps$ there. Subtract the differentials of the leading terms
computed in Steps~1--4 from the exact first variation. Integrating
this difference along straight segments in the translated scale and
outer parameters gives one $C_\eps$ on the product domain.
The segment lengths are uniformly bounded in the specified product norm. Steps~1--4 and
\eqref{eq:bridge-scales}--\eqref{eq:bridge-scale-ledger} therefore give
\eqref{eq:bridge-action} with a parameter-independent $C_\eps$ and
uniform $o_{C^1}(\eps^2)$ remainder.
\end{proof}

\subsection{The exterior problem with prescribed spectral moments}
\label{subsec:high-exterior-problem}

On the complement of the puncture balls, we prescribe the pole and
boundary data independently of the compact background. We solve the
Yamabe equation up to a residual in $E$, with prescribed integrals
against the basis functions $\Psi_j$.

We use the family $U_\bgpar$ from Lemma~\ref{lem:compact-graph} in
dimension $n\ge8$. The local exponents $\sigma_i$ and fixed neck
ratios $t_i$ determine the radii:
\[
 \eps_i=t_i\eps,\quad \sigma_i=\min(q_i,(n-2)/4),\quad
 \rho_i=\eps_i^{(n-2)/[\sigma_i(2n-1)]},\qquad
 \Omega=M\setminus\bigcup_iB_{\rho_i}(p_i).
\]
Let $L_i\in\R$, $\Omega_i\in\R^n$ and
$\psi_i^{\mathrm{ext}}\in C^{6,\alphaH}_{\ge2}(\Sph^{n-1})$ range over
fixed bounded sets. Their pole coefficient and boundary data are
\[
 \mu_i=\rho_i^{n-1/2}L_i,\qquad
 \rho_i\Omega_i\cdot\theta,\qquad \rho_i^2\psi_i^{\mathrm{ext}},
\]
respectively.
Let $T$ be the cutoff sum of the poles $\mu_i r^{2-n}$ and the
decaying harmonic extensions of $\rho_i\Omega_i\cdot\theta$ and
$\rho_i^2\psi_i^{\mathrm{ext}}$,
multiplied by $\Theta_i$. We use the extensions and their higher
regularity established after \eqref{eq:harmonic-matching-inverse}.
We seek $Z$ with zero boundary value and $d\in\R^{N_E}$ such that

\begin{equation}
 U_{\rm ext}=U_\bgpar+T+Z,\quad Z|_{\partial\Omega}=0,\qquad
 \Ff(U_{\rm ext})=\sum_{j=1}^{N_E}d_j\Psi_j,\qquad
 \int_{\Omega} U_{\rm ext}\Psi_j\vol=\eta_j.
 \label{eq:high-exterior-family}
\end{equation}
Thus $U_{\rm ext}$ has boundary value $U_\bgpar+T$.
The unperturbed multiplier is $d_0$, defined by
$\Ff(U_\bgpar)=\sum_jd_{0,j}\Psi_j$.

Define the nonlinear remainder by
\[
\begin{aligned}
 N_\bgpar(W)
 &=\Ff(U_\bgpar+W)-\Ff(U_\bgpar)-L_\bgpar W\\
 &=-\frac{n(n-2)}4\bigl[(U_\bgpar+W)^p-U_\bgpar^p
                         -pU_\bgpar^{p-1}W\bigr].
\end{aligned}
\]
The equation for $(Z,d-d_0)$ is
\begin{equation}
 L_\bgpar Z-\sum_{j=1}^{N_E}(d_j-d_{0,j})\Psi_j=-L_\bgpar T-N_\bgpar(T+Z).
 \label{eq:actual-exterior-correction}
\end{equation}

Fix disjoint coordinate balls $B_{4r_0}(p_i)$.
Write $K_0=M\setminus\bigcup_iB_{r_0}(p_i)$.
For $\rho_i\le r\le r_0$, set
$\mathcal A_{i,r}=\{r\le|X|\le2r\}\cap\overline\Omega$.
Pull its norm back by $X=ry$; at $r=\rho_i$ this includes the
Dirichlet boundary. Introduce the parameter-independent sizes
\begin{align}
 \widehat\mu_i&=\rho_i^{n-1/2},&
 \sigma_{\rm ext}&=\sum_{i=1}^m(\widehat\mu_i+\rho_i^n),\nonumber\\
 e_{\rm ext}&=\sigma_{\rm ext}+
 \left(\sum_{i=1}^m[\widehat\mu_i^2\rho_i^{2-n}+\rho_i^n+\rho_i^{n+2}]\right)^{1/2}.
 \label{eq:exterior-norm-sizes}
\end{align}

The scale relations \eqref{eq:bridge-scale-ledger} give
\[
 \sigma_{\rm ext}=O(\eps^2),\qquad e_{\rm ext}=o(\eps).
\]
For $\rho_i\le r\le r_0$, put
\[
 \omega_i(r)=\sigma_{\rm ext}+\widehat\mu_i r^{3-n},\qquad
 \omega_i^\sharp(r)=\sigma_{\rm ext}+\widehat\mu_i r^{1-n}.
\]
The weights $\omega_i$ and $\omega_i^\sharp$ control $Z$ and
$L_\bgpar Z$, respectively. They satisfy
$r^2\omega_i^\sharp\le C\omega_i$ for scaled Schauder estimates.
For $Z|_{\partial\Omega}=0$, $e\in\R^{N_E}$ and a source $F$, define
\begin{align}
 \|(Z,e)\|_{\mathcal X_{\rm ext}}
  ={}&e_{\rm ext}^{-1}\|Z\|_{H^1}
   +\sigma_{\rm ext}^{-1}
      (\|Z\|_{L^1}+|e|+\|Z\|_{C^{6,\alphaH}(K_0)})\nonumber\\
  &+\max_i\sup_{\rho_i\le r\le r_0}
      \omega_i(r)^{-1}\|Z\|_{C^{6,\alphaH}_{\rm sc}(\mathcal A_{i,r})},
 \label{eq:exterior-response-norm}\\
 \|F\|_{\mathcal Y_{\rm ext}}
  ={}&e_{\rm ext}^{-1}\|F\|_{H^{-1}}
   +\sigma_{\rm ext}^{-1}
      (\|F\|_{L^1}+\|F\|_{C^{4,\alphaH}(K_0)})\nonumber\\
  &+\max_i\sup_{\rho_i\le r\le r_0}
      \omega_i^\sharp(r)^{-1}
      \|F\|_{C^{4,\alphaH}_{\rm sc}(\mathcal A_{i,r})}.
 \label{eq:exterior-source-norm}
\end{align}
Here $H^{-1}=(H^1_0(\Omega))^*$, using $g$.
A fixed finite cover of $K_0$ gives equivalent H\"older norms.
The spaces $\mathcal X_{\rm ext}$ and $\mathcal Y_{\rm ext}$ consist of
pairs $(Z,e)$ and sources $F$ with finite respective norms; they are
Banach spaces for each fixed $\eps$.
The weights remain fixed under exterior parameter variations.
Before matching, these parameters are $(\bgpar,(L_i,\Omega_i,\psi_i^{\mathrm{ext}})_i)$;
first derivatives include unit directions in every
$C^{6,\alphaH}_{\ge2}$ trace factor.

\begin{lemma}\label{lem:exterior-weighted-inverse}
The augmented Dirichlet problem
\[
 L_\bgpar Z-\Psi e=F,\qquad Z|_{\partial\Omega}=0,\qquad
 \int_{\Omega} Z\Psi\vol=\mathbf m
\]
has a unique solution with the uniform bound
\begin{equation}
 \|(Z,e)\|_{\mathcal X_{\rm ext}}
 \le C\bigl(\|F\|_{\mathcal Y_{\rm ext}}
                   +\sigma_{\rm ext}^{-1}|\mathbf m|\bigr).
 \label{eq:exterior-intersection-inverse}
\end{equation}
The inverse is $C^1$ in $\bgpar$. Its $\bgpar$ derivative satisfies
\eqref{eq:exterior-intersection-inverse}, with $(F,\mathbf m)$ fixed.
\end{lemma}
\begin{proof}
Extend zero trace functions by zero across the holes. The zero
moment subspace lies in the fixed coercive space
\eqref{eq:high-coercive}.
Choose a fixed smooth cutoff $\chi_0$, zero on small puncture
neighborhoods and one outside slightly larger neighborhoods.
Since $(\Psi_j)$ is orthonormal, these neighborhoods can be chosen so that
\[
 \mathsf M_{jk}=\int_M\chi_0\Psi_k\Psi_j\vol,\qquad
 \|\mathsf M-I\|<\tfrac12.
\]
For all sufficiently small holes,
\[
 Z_{\mathbf m}=\sum_k(\mathsf M^{-1}\mathbf m)_k\chi_0\Psi_k,\qquad
 \int_{\Omega} Z_{\mathbf m}\Psi\vol=\mathbf m,\qquad
 \|Z_{\mathbf m}\|_{H^1}\le C|\mathbf m|.
\]
Subtract $Z_{\mathbf m}$ and apply Lax--Milgram on the zero
moment subspace. Testing the resulting equation against
$\chi_0\Psi_j$, $1\le j\le N_E$, determines $e$ and gives uniqueness and
\[
 \|Z\|_{H^1}+|e|\le C(\|F\|_{H^{-1}}+|\mathbf m|).
\]
Since $\sigma_{\rm ext}\le e_{\rm ext}$, this proves the energy part
of \eqref{eq:exterior-intersection-inverse}.

We first prove a uniform supremum bound for bounded data.
Choose $C_*>0$ and let $G_\Omega$ be the Dirichlet inverse of
$P=-\Delta_g+C_*$. Let $G_M$ be the Green kernel of $P$ on $M$.
The positive Dirichlet kernel satisfies
$G_\Omega(X,Y)\le G_M(X,Y)\le C(d(X,Y)^{2-n}+1)$.
Writing $L_\bgpar=-\Delta_g+V_\bgpar$, split the potential term in
$Z=G_\Omega[F+\Psi e+(C_*-V_\bgpar)Z]$ at distance $r_*$ to obtain
\[
 \|Z\|_\infty\le C(\|F\|_\infty+|e|)
       +Cr_*^2\|Z\|_\infty+C_{r_*}\|Z\|_{L^1}.
\]
Elliptic regularity makes $Z$ bounded for each fixed hole domain.
Choose $r_*$ small to absorb $Cr_*^2\|Z\|_\infty$.
The preceding $H^1$ estimate controls the remaining $L^1$ term.
Thus $\|Z\|_\infty+|e|\le C(\|F\|_\infty+|\mathbf m|)$,
uniformly in the holes.

For adjoint data $(H,\widetilde{\mathbf m})$, let $(V,f)$ solve
the augmented Dirichlet problem with those data. Self-adjointness gives
\[
 \int_{\Omega} ZH=\int_{\Omega} FV+e\cdot\widetilde{\mathbf m}-f\cdot \mathbf m.
\]
Taking the supremum over $\|H\|_\infty+|\widetilde{\mathbf m}|\le1$
and using the supremum estimate for $(V,f)$ yields
\[
 \|Z\|_{L^1}+|e|\le C(\|F\|_{L^1}+|\mathbf m|).
\]
Put $A=\|F\|_{\mathcal Y_{\rm ext}}+|\mathbf m|/\sigma_{\rm ext}$.
Interior estimates starting from
$\|Z\|_{L^1}+|e|\le C(\|F\|_{L^1}+|\mathbf m|)$ give
$\|Z\|_{C^{6,\alphaH}}\le C\sigma_{\rm ext}A$ on fixed shells
and $K_0$, since both local source weights there are
$O(\sigma_{\rm ext})$.

Choose $r_0$ small enough that the first Dirichlet eigenvalue on
$B_{2r_0}(p_i)$ dominates the negative part of the potential.
The normal coordinate bounds give
$L_\bgpar(r^{3-n})\ge cr^{1-n}$; its metric error is $O(r^{3-n})$.
For a fixed large $C_0$ and then small $r_0$,
$b_0(r)=1+C_0(r_0^2-r^2)$ is positive and $L_\bgpar b_0\ge1$.
A sufficiently large fixed multiple
\[
 B_i(r)=CA[\sigma_{\rm ext}b_0(r)+\widehat\mu_i r^{3-n}]
\]
therefore satisfies $L_\bgpar B_i\ge|F|+|\Psi e|$ and
$B_i\ge|Z|$ on both annular boundaries. The maximum principle
applied to $B_i\pm Z$ gives $|Z|\le CA\omega_i(r)$.
Scaled interior and Dirichlet Schauder estimates give
$\|Z\|_{C^{6,\alphaH}_{\rm sc}(\mathcal A_{i,r})}
\le CA\omega_i(r)$, because
$r^2\omega_i^\sharp\le C\omega_i$. This proves
\eqref{eq:exterior-intersection-inverse}.

A derivative of $L_\bgpar$ with respect to $\bgpar$ is multiplication by a
smooth bounded function. Such multiplication maps
$\mathcal X_{\rm ext}$ to $\mathcal Y_{\rm ext}$ uniformly, by
$\omega_i\le C\omega_i^\sharp$ and the $H^1$, $L^1$ and scaled
H\"older product estimates. Formula
\eqref{eq:inverse-parameter-calculus} now gives the first derivative
bound from $\mathcal Y_{\rm ext}$ to $\mathcal X_{\rm ext}$.
\end{proof}

For the trace and scaled conormal of $Z$ on $S_{\rho_i}$, use
\[
 \|\mathcal C_iZ\|
 :=\|\Theta_i^{-1}Z\|_{C^{6,\alphaH}_{\rm sc}(\mathcal A_{i,\rho_i})}
   +\|\rho_i\partial_{\nu_{h_i}}(\Theta_i^{-1}Z)
                       \|_{C^{5,\alphaH}(S_{\rho_i})}.
\]
Pull the trace norm back to the unit sphere. Multiplication by the
fixed surface density converts the scalar conormal to the conormal
density used in the matching equations.

We now solve \eqref{eq:high-exterior-family} for all bounded inputs
$(\bgpar,L_i,\Omega_i,\psi_i^{\mathrm{ext}})$, including the full
$C^{6,\alphaH}$ trace balls.
The prescribed moment of $Z$ is
\[
 \mathbf m=\bgpar-\int_{\Omega}(U_\bgpar+T)\Psi\vol
   =\int_{M\setminus\Omega}U_\bgpar\Psi\vol-\int_{\Omega} T\Psi\vol.
\]
The pole lift, degree-one harmonic lift, and harmonic lift in degrees
at least two have, respectively, the annular envelopes
$|\mu_i|r^{2-n}$, $\rho_i^nr^{1-n}$ and $\rho_i^{n+2}r^{-n}$,
through six scaled derivatives. Their integrals are bounded by
$C|\mu_i|$, $C\rho_i^n$ and
$C\rho_i^{n+2}\log(r_0/\rho_i)$, respectively. Their squared
$H^1$ norms are bounded by
\[
 C\widehat\mu_i^2\rho_i^{2-n},\qquad
 C\rho_i^n,\qquad C\rho_i^{n+2}.
\]
The three lift bounds hold after differentiation in any unit
direction of $(\bgpar,(L_i,\Omega_i,\psi_i^{\mathrm{ext}})_i)$.
The integrals over $B_{\rho_i}(p_i)$ in the formula for $\mathbf m$
are $O_{C^1}(\rho_i^n)$, and hence
\[
 \|\mathbf m\|_{C^1}\le C\sigma_{\rm ext},\qquad
 \max_{j\le1}\|D_{\rm par}^jT\|_{H^1}\le Ce_{\rm ext}.
\]

Through four scaled derivatives, the geometric and potential
residuals of the pole, degree-one, and high harmonic lifts have the
respective envelopes
\[
 |\mu_i|r^{2-n},\qquad \rho_i^nr^{1-n},\qquad
 \rho_i^{n+2}r^{-n}.
\]
Their comparison with the source weight is
\[
 \frac{|\mu_i|r^{2-n}}{\omega_i^\sharp(r)}\le Cr,
 \qquad
 \frac{\rho_i^nr^{1-n}}{\omega_i^\sharp(r)}\le C\rho_i^{1/2},
 \qquad
 \frac{\rho_i^{n+2}r^{-n}}{\omega_i^\sharp(r)}
       \le C\rho_i^{5/2}/r.
\]
All three ratios are uniformly bounded for $\rho_i\le r\le r_0$.
The cutoff errors are supported on fixed shells and have
$C^1$ size $O(\sigma_{\rm ext})$.
Integrating the three source envelopes gives the $L^1$ bound for
$L_\bgpar T$. Integration by parts against
$\phi\in H^1_0(\Omega)$ gives
$|\langle L_\bgpar T,\phi\rangle|\le C\|T\|_{H^1}\|\phi\|_{H^1}$,
which supplies the $H^{-1}$ component.
Consequently,
\[
 \max_{j\le1}\|D_{\rm par}^j(L_\bgpar T)\|_{\mathcal Y_{\rm ext}}\le C,
 \qquad
 \max_{j\le1}\left(
   \|D_{\rm par}^jT\|_{L^1}
   +\|D_{\rm par}^jT\|_{C^{6,\alphaH}(K_0)}\right)
 \le C\sigma_{\rm ext}.
\]

Put $\kappa_{\rm ext}=\max_i\rho_i+\max_i\omega_i(\rho_i)\to0$.
Let $F$ be $T$ or one of its first parameter derivatives, or a function
satisfying $\|(F,0)\|_{\mathcal X_{\rm ext}}\le C$. Then
\[
 \begin{gathered}
 \|F\|_{H^1}\le Ce_{\rm ext},\qquad
 \|F\|_{L^1}+\|F\|_{C^{4,\alphaH}(K_0)}\le C\sigma_{\rm ext},\\
 \|F\|_\infty\le C\kappa_{\rm ext},\qquad
 \|F\|_{C^{4,\alphaH}_{\rm sc}(\mathcal A_{i,r})}
 \le C\min\{\kappa_{\rm ext},\omega_i^\sharp(r)\}.
 \end{gathered}
\]
For the three lifts, their ratios to $\omega_i^\sharp$ are at most
$r$, $\rho_i^{1/2}$ and $\rho_i^{5/2}/r$, and their maximum amplitudes
are $\rho_i^{3/2}$, $\rho_i$ and $\rho_i^2$. For $\|(F,0)\|_{\mathcal X_{\rm ext}}\le C$, use
$\omega_i(r)\le\omega_i(\rho_i)$ and $\omega_i\le C\omega_i^\sharp$.
For two such fields $F,G$, H\"older's integral inequality gives
\[
 \|FG\|_{H^{-1}}\le C\|F\|_\infty\|G\|_{L^2}
       \le C\kappa_{\rm ext}e_{\rm ext},\qquad
 \|FG\|_{L^1}\le\|F\|_\infty\|G\|_{L^1}
       \le C\kappa_{\rm ext}\sigma_{\rm ext}.
\]
The product estimate in scaled H\"older spaces controls the other
two components. Hence
\begin{equation}
 \|FG\|_{\mathcal Y_{\rm ext}}\le C\kappa_{\rm ext}.
 \label{eq:exterior-product-bound}
\end{equation}
For arbitrary field directions $F,G$, the right-hand side is
$C\kappa_{\rm ext}\|(F,0)\|_{\mathcal X_{\rm ext}}
\|(G,0)\|_{\mathcal X_{\rm ext}}$.
Multiplication by bounded scaled $C^{4,\alphaH}$ coefficients
preserves \eqref{eq:exterior-product-bound}.

For $\|(Z,d-d_0)\|_{\mathcal X_{\rm ext}}\le C$,
$W=T+Z\to0$ uniformly and $U_\bgpar\ge c>0$.
By \eqref{eq:relative-power-calculus}, $N_\bgpar(W)$ is $W^2$
times a coefficient with bounded scaled H\"older norm and first
derivatives. At fixed $W$, differentiation in $\bgpar$ keeps the
factor $W^2$. At fixed $Z$, differentiation in a lift parameter
gives terms containing $W D_{\rm par}T$.
Equation~\eqref{eq:exterior-product-bound} gives joint Fr\'echet
$C^1$ dependence, including every $C^{6,\alphaH}$ trace direction.
It also gives $\|N_\bgpar(T)\|_{\mathcal Y_{\rm ext}}\le C\kappa_{\rm ext}$ and
\begin{equation*}
 \|N_\bgpar(T+Z)-N_\bgpar(T+\widetilde Z)\|_{\mathcal Y_{\rm ext}}
 \le C\kappa_{\rm ext}
        \|(Z-\widetilde Z,0)\|_{\mathcal X_{\rm ext}}.
\end{equation*}
Let $\mathscr G_\bgpar$ denote the inverse in
Lemma~\ref{lem:exterior-weighted-inverse}. Applying it to the equation
for $Z$ and the prescribed moment vector $\mathbf m$ gives
\[
 (Z,d-d_0)=\mathscr G_\bgpar\bigl(-L_\bgpar T-N_\bgpar(T+Z),\mathbf m\bigr).
\]
Choose a ball centered at zero in $\mathcal X_{\rm ext}$ whose fixed
radius is at least twice the uniform norm bound for the image of zero. The source and inverse estimates
show that the displayed map sends this ball into itself and is a
jointly $C^1$ contraction.
Equation~\eqref{eq:fixed-point-parameter-calculus} gives
\[
 \max_{0\le j\le1}
 \|D_{\rm par}^j(Z,d-d_0)\|_{\mathcal X_{\rm ext}}\le C.
\]
The solution $U_{\rm ext}=U_\bgpar+T+Z$ is positive. Since
$\omega_i(\rho_i)=\sigma_{\rm ext}+\rho_i^{5/2}$,
$e_{\rm ext}=o(\eps)$ and $\eps^2=o(\rho_i^2)$, it satisfies
\begin{align}
 \|\mathcal C_iZ\|&\le C(\rho_i^{5/2}+\eps^2)=o(\rho_i^2),
 \label{eq:localized-cap}\\
 \|T\|_{H^1}+\|Z\|_{H^1}&\le Ce_{\rm ext}=o(\eps).
 \label{eq:sharp-exterior-energy}
\end{align}
Both estimates hold after one unit derivative in
$(\bgpar,(L_i,\Omega_i,\psi_i^{\mathrm{ext}})_i)$, including all
$C^{6,\alphaH}$ trace directions.

\subsection{Joining each local family to the global exterior}
\label{subsec:high-exterior-matching}

We join the local solution to $U_{\rm ext}$ on $S_{\rho_i}$, expressing
both conformal factors in the metric $h_i$. Their values and
increasing-radius normal derivatives must agree.
The free variables are $(\bgpar,\ell_i,\zeta_i)$, with $\ell_i=\log b_i$;
the boundary variables are
$(A_i,L_i,\mathbf v_i,\Omega_i,\psi_i,\psi_i^{\mathrm{ext}})$.
Let $c_i(\bgpar)>0$ and $e_i(\bgpar)$ be the constant and linear Taylor
coefficients of $\Theta_i^{-1}U_\bgpar$ at $p_i$, and write its high
boundary trace as $\rho_i^2C_i(\bgpar)$. Value matching in the high
modes gives exactly $\psi_i=C_i(\bgpar)+\psi_i^{\mathrm{ext}}$.
Use the affine coordinates
\[
 \alpha_i=A_i-c_i(\bgpar),\qquad
 \mathbf u_i=\mathbf v_i-e_i(\bgpar),\qquad
 \widehat L_i=L_i-\rho_i^{1/2-n}b_i^{(n-2)/2}C_\infty(\zeta_i).
\]
Proposition~\ref{prop:bridge} gives, in all independent inputs,
\[
 \mathcal B_i-\rho_i^{1/2-n}b_i^{(n-2)/2}C_\infty(\zeta_i)=o_{C^1}(1).
\]

Before solving the boundary variables, the mean equations in the
coordinates $(\alpha_i,\widehat L_i)$ are
\begin{equation}
\begin{aligned}
 \alpha_i-\rho_i^{3/2}L_i+E_i^{\rm val}&=0,\\
 (n-2)\widehat L_i
 +(n-2)[\rho_i^{1/2-n}b_i^{(n-2)/2}C_\infty(\zeta_i)-\mathcal B_i]
 +\rho_i^{-3/2}E_i^{\rm nor}&=0,\\
 E_i^{\rm val},E_i^{\rm nor}&=O_{C^1}(\rho_i^2+\sigma_{\rm ext}).
\end{aligned}
\label{eq:high-mean-matching}
\end{equation}
The Taylor errors of the regular filling and $\Theta_i^{-1}U_\bgpar$
are $O_{C^1}(\rho_i^2)$. Equation~\eqref{eq:localized-cap} bounds
the traces of $\Theta_i^{-1}Z$. The value equation is unscaled; only the
conormal equation has been divided by $\rho_i^{3/2}$.

The first harmonic equations, divided by $\rho_i$, are
\[
 \mathbf u_i-\Omega_i+R_i^{1,\rm val}=0,
 \qquad \mathbf u_i+(n-1)\Omega_i+R_i^{1,\rm nor}=0,
\]
with
\[
 \|R_i^{1,\rm val}\|_{C^1}+\|R_i^{1,\rm nor}\|_{C^1}
       \le C(\rho_i^{3/2}+\sigma_{\rm ext}\rho_i^{-1}).
\]
For each first-harmonic component, the $2\times2$ coefficient matrix has determinant $n$.
After division by $\rho_i^2$, the high conormal equation is
$K\psi_i^{\mathrm{ext}}+R_i^{\rm hi}=0$, with $K=D_+-D_-$.
Here $D_-$ is the increasing-radius conormal of the decaying
exterior harmonic extension. The bounded inverse
$C^{5,\alphaH}_{\ge2}\to C^{6,\alphaH}_{\ge2}$ is
\eqref{eq:harmonic-matching-inverse}.
The local remainder is jointly $C^1$ small by
Proposition~\ref{prop:bridge}; the exterior contribution is
$O_{C^1}(\rho_i^{1/2}+\sigma_{\rm ext}\rho_i^{-2})=o(1)$.
For $f_i=\Theta_i^{-1}U_\bgpar$ in the fixed $h_i$ coordinates, write
\[
 f_i=c_i+e_i\cdot X+Q_i+R_i,\qquad
 \|R_i\|_{C_\bgpar^1(C^{6,\alphaH}_{\rm sc}(A_r))}\le Cr^3,
\]
where $Q_i$ is homogeneous of degree two. Its high projection has
degree two, so $D_+\Pi_{\ge2}Q_i=2\Pi_{\ge2}Q_i$. Thus, on the
unit-sphere pullback of $S_{\rho_i}$,
\begin{equation*}
 \|\Pi_{\ge2}(\rho_i\partial_rf_i)-D_+\Pi_{\ge2}f_i\|_
       {C_\bgpar^1(C^{5,\alphaH})}\le C\rho_i^3.
\end{equation*}
The normal and density errors are $O(\rho_i^2)$ and multiply
$\rho_i\nabla f_i=O(\rho_i)$, so they have the same order.
After division by $\rho_i^2$, the background contribution is $O_{C^1}(\rho_i)$.
The matching equations use $h_i$. Once the value equations hold,
the local factor $u_{{\rm loc},i}$ has trace $\Theta_i^{-1}U_{\rm ext}$.
The extra conormal term on returning to $g$ is then
$(\partial_{\nu_g}\Theta_i)(u_{{\rm loc},i}-\Theta_i^{-1}U_{\rm ext})=0$.

For the boundary variables
\[
 Y=(\alpha_i,\widehat L_i,\mathbf u_i,\Omega_i,
       \psi_i^{\mathrm{ext}})_i,
\]
use the maximum over punctures of the finite-dimensional norms and
$\|\psi_i^{\mathrm{ext}}\|_{C^{6,\alphaH}}$. The high equation is
measured in $C^{5,\alphaH}$.
The fixed block operator is
\[
 \mathcal A_{\rm match}=\bigoplus_i\left[
 \operatorname{diag}(1,n-2)\oplus
 \begin{pmatrix}1&-1\\1&n-1\end{pmatrix}\otimes I_n\oplus K\right].
\]
Its inverse is uniformly bounded.
The mean remainders are $o_{C^1}(1)$ by the pole expansion and
$\rho_i^{1/2}+\sigma_{\rm ext}\rho_i^{-3/2}\to0$.
All other normalized errors also tend to zero in $C^1$, since
$\sigma_{\rm ext}=O(\eps^2)=o(\rho_i^2)$ at every puncture;
the constants are uniform for independently chosen radii
$\rho_i$ satisfying \eqref{eq:bridge-scales}.
Uniform positivity of $c_i$ and boundedness of the reference data
place one fixed small ball in the product space of $Y$ inside the
local and exterior parameter domains. Including the small
off-diagonal terms in $\mathcal R_\eps$, the matching system is
\begin{equation}
 \mathcal A_{\rm match}Y+\mathcal R_\eps(\bgpar,\ell,\zeta,Y)=0,
 \qquad Y=-\mathcal A_{\rm match}^{-1}\mathcal R_\eps(\bgpar,\ell,\zeta,Y).
 \label{eq:high-matching-fixed-point}
\end{equation}
The normalized error bounds make the last map a contraction on that ball.
Equation~\eqref{eq:fixed-point-parameter-calculus} gives its unique solution
$Y=o_{C^1(\bgpar,\ell,\zeta)}(1)$ in this product norm.
Uniqueness identifies overlapping parameter patches.
From this point, $C^1$ estimates for the matched family use
$(\bgpar,\ell,\zeta)$ as parameters.

The solved boundary parameters have bounded first derivatives. The unscaled
value equation in \eqref{eq:high-mean-matching}, the degree-one
equations and the high conormal equation give
\begin{equation}
 \begin{aligned}
 A_i&=c_i(\bgpar)+O_{C^1}(\rho_i^{3/2}),&
 \mathbf v_i&=e_i(\bgpar)+o_{C^1}(1),\\
 L_i&=\mathcal B_i+o_{C^1}(1),&
 \Omega_i,\psi_i^{\mathrm{ext}}&=o_{C^1}(1).
 \end{aligned}
 \label{eq:high-matched-data}
\end{equation}

The matching gives a positive family $U_\eps(\bgpar,\ell,\zeta)$ on the
prescribed compact, scale and center parameter sets. The Fowler necks
are fixed, and the family is $C^1$ in the stated local and exterior
norms. Its residual and exact moments are
\eqref{eq:high-exterior-family} and \eqref{eq:bridge-projection}.
Values and physical conormals match on $S_{\rho_i}$ and every
artificial interface. Thus $U_\eps$ is $C^1$ across these surfaces
in physical coordinates.

\subsection{The global action and elimination of the residuals}
\label{subsec:global-high-action}

The boundary variables are now functions of $(\bgpar,\ell,\zeta)$ at
fixed necks. Parameter derivatives include this dependence. The
exterior residual is $\sum_jd_j\Psi_j$. On the $i$th middle region,
the action gradient relative to $\dd V_{\rm rd}$ is
$\sum_jc_{i,j}\chi Z_j$.

Let $S_i$ be the full local action of Proposition~\ref{prop:bridge-action},
expressed in the global metric by \eqref{eq:action-metric-conversion}.
Use \eqref{eq:relative-action-conformal} to restore the regular filling action,
and set
\begin{equation}
 I_\eps=J_g(U_\eps;\Omega)+\sum_{i=1}^mS_i,\qquad
 S_i=J_{{\rm regular},i}+I_{{\rm rel},i}+C_{\eps_i}.
 \label{eq:full-global-action}
\end{equation}
Here $J_{{\rm regular},i}$ is the action of the regular filling
on $B_{\rho_i}(p_i)$, expressed in the metric $g$.

Let $d$ be the exterior multiplier vector and write
$c=(c^{(1)},\ldots,c^{(m)})$, with $c^{(i)}=(c_{i,0},\ldots,c_{i,n})$.
Order the local parameters by puncture as $(\ell_i,\zeta_i)$ and let
$G_{b_i}$ be the matrix associated with the $i$th puncture in
Lemma~\ref{lem:bridge-pairing}, evaluated at the matched boundary data.
\begin{proposition}\label{prop:global-high-action}
The action $I_\eps(\bgpar,\ell,\zeta)$ is single-valued and $C^1$.
Uniformly on the fixed compact parameter sets for $\bgpar$, $\zeta_i$
and $b_i^{\sigma_i}/\eps_i$,
\begin{equation}
 I_\eps=C_\eps+j_0(\bgpar)+\sum_{i=1}^m
       [\mu_i^H\ell_i+b_i^{2\sigma_i}\mathcal P_i(\bgpar,\zeta_i)]
       +o_{C^1}(\eps^2),\qquad \mu_i^H=-|\Sph^{n-1}|h_{\eps_i}>0,
 \label{eq:global-high-action}
\end{equation}
where $C_\eps$ is independent of $(\bgpar,\ell,\zeta)$ and
\begin{equation}
 \mathcal P_i(\bgpar,\zeta_i)=
 \begin{cases}
 \Phi_i(\zeta_i),&n>4q_i+2,\\
 \Phi_i(\zeta_i)-(n-2)|\Sph^{n-1}|c_i(\bgpar)C_\infty(\zeta_i),&n=4q_i+2,\\
 -(n-2)|\Sph^{n-1}|c_i(\bgpar)C_\infty(\zeta_i),&n<4q_i+2.
 \end{cases}\label{eq:leading-potentials}
\end{equation}
The exact differential of $I_\eps$ is
\begin{equation}
 \begin{aligned}
 D_\bgpar I_\eps&=d,\\
 D_{(\ell,\zeta)}I_\eps&=\operatorname{diag}(G_{b_1},\ldots,G_{b_m})c.
 \end{aligned}
 \label{eq:full-multiplier-pairing}
\end{equation}
Each block is uniformly invertible at the matched data, so a critical
point of $I_\eps$ has $c=0$ and $d=0$.
\end{proposition}
\begin{proof}
Substituting \eqref{eq:high-matched-data} in
Proposition~\ref{prop:bridge-action} preserves its $o_{C^1}(\eps_i^2)$
remainders, since all recovered outer data have bounded first derivatives.

Put $W=T+Z$ on $\Omega$. Taylor's formula, the exact moments and
\eqref{eq:sharp-exterior-energy} give
\begin{align*}
 J_g(U_\bgpar+W;\Omega)-J_g(U_\bgpar;\Omega)
 ={}&\sum_jd_{0,j}\int_{M\setminus\Omega}U_\bgpar\Psi_j\vol
   +\int_{\partial\Omega}\partial_{\nu_g}U_\bgpar\,W\,\dd S_g
   +O_{C^1}(e_{\rm ext}^2).
\end{align*}
Uniform positive bounds give the quadratic remainder bound
$C\|W\|_{H^1}^2$. Its first derivative is bounded by
$C\|W\|_{H^1}(\|W\|_{H^1}+\|D_{\rm par}W\|_{H^1})$.
Thus the remainder is $o_{C^1}(\eps^2)$, while the moment term is
$O_{C^1}(\sum_i\rho_i^n)=o_{C^1}(\eps^2)$.

On $S_{\rho_i}$,
$\Theta_i^{-1}U_\bgpar=c_i+e_i\cdot X+O_{C^1}(|X|^2)$,
with exterior conormal $-e_i\cdot\theta+O_{C^1}(\rho_i)$.
The constant leading pole trace has zero pairing with
$-e_i\cdot\theta$, and the high trace is orthogonal to the first
spherical harmonics. Pole, first and high traces thus contribute,
respectively,
\[
 C|\mu_i|\rho_i^2,\qquad C\rho_i^n,\qquad C\rho_i^{n+2}
\]
in $C^1$; fixed conformal factors and densities preserve these bounds.
Both $J_g(U_\bgpar;B_{\rho_i}(p_i))$ and $J_{{\rm regular},i}$
are $O_{C^1}(\rho_i^n)$; Lemma~\ref{lem:regular-filling} gives
the bound for $J_{{\rm regular},i}$. These contributions are
$o_{C^1}(\eps^2)$.
Moreover, \eqref{eq:high-matched-data} gives
\[
 \left\|\sum_iB_i[A_i-c_i(\bgpar)]\right\|_{C^1}
 \le C\sum_iB_i\rho_i^{3/2}=o(\eps^2).
\]
We may therefore sum the local expansions with $A_i=c_i(\bgpar)$.
For $n>4q_i+2$,
$b_i^{(n-2)/2}/\eps_i^2\asymp b_i^{(n-2)/2-2q_i}\to0$.
At $n=4q_i+2$, the interaction and curvature terms combine.
For $n<4q_i+2$, the curvature term belongs to the remainder. This proves
\eqref{eq:global-high-action} and \eqref{eq:leading-potentials}.
The local constants were fixed on connected boundary-data sets
containing all traces induced by the background parameter $\bgpar$. Their
sum is therefore independent of the connected component of the
minimizing set containing the limiting background.

For the exact differential, combine
\eqref{eq:high-local-action-differential} with
\eqref{eq:bridge-pullback-pairing} at the matched local parameters.
Decompose each local variation into its direct $(\ell_i,\zeta_i)$
variation and its boundary-data variation.
For the latter, differentiate the local family in
$(A_i,\mathbf v_i,\psi_i)$ while holding $(\ell_i,\zeta_i)$ fixed.
Equation~\eqref{eq:cap-tangent-zero-moments} makes its bulk pairing zero.
The direct $(\ell_i,\zeta_i)$ variation gives the $G_{b_i}$ pairing
in \eqref{eq:bridge-pullback-pairing}.
Conformal change preserves the bulk pairing, since
$\Theta_i^{-p}\Theta_i^{-1}\Theta_i^{p+1}=1$.
On each fixed sphere $S_{\rho_i}$, traces, trace variations and physical
conormals match with opposite orientations, so the boundary terms
cancel.

Fix a direction $(\dot \bgpar,\dot\ell,\dot\zeta)$ and let $\dot U_\eps$
be the corresponding variation at fixed physical coordinates.
The exterior moments give
\[
 \int_{\Omega}\Psi_j\dot U_\eps\vol=\dot\eta_j.
\]
The first variation is therefore
\[
 DI_\eps[(\dot \bgpar,\dot\ell,\dot\zeta)]
 =d\cdot\dot \bgpar+\sum_i(G_{b_i}c^{(i)})\cdot(\dot\ell_i,\dot\zeta_i).
\]
This proves \eqref{eq:full-multiplier-pairing}; uniform invertibility
from Lemma~\ref{lem:bridge-pairing} then gives $c=d=0$ at every
critical point.
\end{proof}

\section{Variational selection and the proof of the theorem}\label{sec:selection}
We select the remaining parameters by finding a critical point of the
action. In dimensions three through seven, its $C^2$ estimate permits
elimination of the amplitude by the implicit function theorem. In
higher dimensions, we use the variational argument below, which
requires only the $C^1$ estimate.

\subsection{A boundary gap and a finite-dimensional variational lemma}

Use the amplitude and normalized spectral coordinates of
Section~\ref{subsec:amplitude}, $\bgpar=\amp\bgshape$ with
$\bgshape\in\mathcal H_E$. 
Choose a compact neighborhood $D$ of $\mathcal S$ contained
in the domain of $\amp_0$, with boundary disjoint from $\mathcal S$.
Since $\mathcal S$ is exactly the minimum set and $\varphi_0$ is continuous,
\begin{equation}
 \min_{\partial D}\varphi_0>\min_D\varphi_0.
 \label{eq:core-boundary-gap}
\end{equation}
Boundaries are taken in $\mathcal H_E$; if this affine space has dimension
zero, $D$ is a point and the boundary condition is vacuous.
Lemma~\ref{lem:negative-scale} gives
$\partial_\amp^2j_0(\amp \bgshape)\le-\eta_0<0$ near
$\{(\amp_0(\bgshape),\bgshape):\bgshape\in D\}$.
The next lemma gives a critical point when the leading action is
strictly concave in some variables and has a positive boundary gap in
the remaining variables.

\begin{lemma}\label{lem:c1-minimax}
Let $\Omega_0$ be a compact neighborhood in a finite-dimensional
Euclidean space, $h,\gamma,R>0$ and $N\ge1$.
Suppose $F\in C^1$ near $\Omega_0$ satisfies
\begin{equation}
 m_*=\min_{\Omega_0}F,\qquad
 F\ge m_*+\gamma h\quad\hbox{on }\partial\Omega_0.
 \label{eq:c1-minimax-gap}
\end{equation}
On a neighborhood of $\Omega_0\times\overline B_{2R}^{\,N}$, let
\[
 J(s,z)=F(s)-hK(s,z)+h\rho(s,z),
\]
where $K,D_zK,\rho$ are $C^1$ and
\begin{equation}
 \begin{gathered}
 K(s,0)=0,\qquad D_zK(s,0)=0,\\
 c\,\operatorname{Id}\le D_z^2K(s,z)\le C\,\operatorname{Id}
 \end{gathered}
 \label{eq:c1-minimax-convexity}
\end{equation}
for fixed $c,C>0$. If
$\delta=\|\rho\|_{C^1(\Omega_0\times\overline B_{2R}^{\,N})}$
is sufficiently small in terms of $\gamma,c,C,R$, then $J$ has a critical point
$(s_h,z_h)\in\operatorname{int}\Omega_0\times B_R^N$ satisfying
\begin{equation}
 \begin{gathered}
 |J(s_h,z_h)-m_*|\le\delta h,\qquad |z_h|\le\delta/c,\\
 0\le F(s_h)-m_*
 \le\left(2\delta+\frac{C\delta^2}{2c^2}\right)h.
 \end{gathered}
 \label{eq:c1-minimax-output}
\end{equation}
Boundaries are relative to the ambient Euclidean space;
the gap condition is vacuous when the boundary is empty.
\end{lemma}
\begin{proof}
Integrating \eqref{eq:c1-minimax-convexity} gives
\begin{equation*}
 \frac c2|z|^2\le K(s,z)\le\frac C2|z|^2,
 \qquad D_zK(s,z)[z]\ge c|z|^2.
\end{equation*}
Choose $0<\beta<\min\{\gamma/4,cR^2/8\}$ and put
\[
 \Omega=\{s\in\operatorname{int}\Omega_0:F(s)<m_*+\beta h\},
 \qquad X=\overline\Omega\times\overline B_R^N.
\]
Then $\overline\Omega\Subset\operatorname{int}\Omega_0$ and
$F=m_*+\beta h$ on $\partial\Omega$. Choose $s_*\in\Omega$
with $F(s_*)=m_*$. For sufficiently small $\delta$,
\begin{equation}
 J\le m_*-bh\quad\hbox{on }\overline\Omega\times\partial B_R^N,
 \qquad b=cR^2/4.
 \label{eq:c1-minimax-ball-boundary}
\end{equation}
Let $\Gamma$ consist of continuous maps
$q:\overline B_R^N\to X$ with $q(z)=(s_*,z)$ on $\partial B_R^N$.
By Brouwer's no-retraction theorem, for each $q\in\Gamma$ there is
$z_q\in\overline B_R^N$ such that the second component of $q(z_q)$
is zero. Hence $J(q(z_q))\ge m_*-\delta h$.
Conversely $q_0(z)=(s_*,z)$ satisfies $\max_{z\in\overline B_R^N}J(q_0(z))\le m_*+\delta h$.
Thus
\begin{equation}
 c_*:=\inf_{q\in\Gamma}\max_{z\in\overline B_R^N}J(q(z))
 \quad\hbox{satisfies}\quad |c_*-m_*|\le\delta h.
 \label{eq:c1-minimax-level}
\end{equation}

If $F(s)>m_*+\beta h/2$ and $|J(s,z)-m_*|\le\beta h/4$,
then $K(s,z)\ge\beta/8$ for $\delta\le\beta/8$, hence
$|z|^2\ge\beta/(4C)$. For smaller $\delta$,
\begin{equation}
 DJ(s,z)[(0,z)]
 \le h(-c|z|^2+\delta|z|)
 \le-\frac{c\beta h}{8C}.
 \label{eq:c1-minimax-boundary-descent}
\end{equation}
This gives a descent direction with zero $s$ component.

Suppose there is no critical point at level $c_*$ in $X$.
By compactness, choose $\eta>0$ and set
$A=\{p\in X:|J(p)-c_*|\le3\eta\}$.
The set $A$ has no critical point.
It is disjoint from $\overline\Omega\times\partial B_R^N$.
It lies in the energy band $\{|J-m_*|<\beta h/4\}$. Choose a smooth cutoff $\chi:\R\to[0,1]$
equal to zero on $(-\infty,5/8]$ and one on $[3/4,\infty)$, and set
\[
 \chi_F(s)=\chi\bigl((F(s)-m_*)/(\beta h)\bigr),\qquad
 W=-(1-\chi_F)\nabla J+\chi_F(0,z).
\]
On $\{\chi_F>0\}\cap A$, \eqref{eq:c1-minimax-boundary-descent}
gives $DJ[(0,z)]<0$. On $\{\chi_F=0\}\cap A$,
$DJ[-\nabla J]=-|DJ|^2<0$. Compactness therefore yields
$DJ[W]\le-\nu<0$ on $A$.
Moreover, the $s$-component of $W$ vanishes near $\partial\Omega$.
Let $\psi:\R\to[0,1]$ be a smooth cutoff equal to one on
$[c_*-\eta,c_*+\eta]$ and zero outside
$(c_*-2\eta,c_*+2\eta)$. For a fixed sufficiently small $t>0$,
\[
 T(p)=p+t\psi(J(p))W(p)
\]
maps $X$ into itself and fixes a neighborhood of
$\overline\Omega\times\partial B_R^N$.
The support of $\psi\circ J$ is disjoint from this boundary.
The $s$ component of $T-\operatorname{Id}$ vanishes outside the
compact subset $\{F\le m_*+3\beta h/4\}$ of $\Omega$.
Uniform continuity of $DJ$ and Taylor's integral formula give
\[
 J(Tp)\le J(p)-\tfrac12\nu t\psi(J(p))
\]
for small $t$.
Reduce $t$ so that $\nu t<\eta$ and choose $q\in\Gamma$ with
$\max_{z\in\overline B_R^N}J(q(z))<c_*+\nu t/4$. For every $z\in\overline B_R^N$,
\[
 J(Tq(z))\le
 \begin{cases}
 J(q(z))-\nu t/2<c_*-\nu t/4,&J(q(z))\ge c_*-\eta,\\
 J(q(z))<c_*-\eta,&J(q(z))<c_*-\eta.
 \end{cases}
\]
Since $T\circ q\in\Gamma$, this gives the contradiction
\[
 c_*\le\max_{z\in\overline B_R^N}J(T(q(z)))
 \le\max\{c_*-\nu t/4,c_*-\eta\}<c_*.
\]
There is therefore a critical point at level $c_*$. Equations
\eqref{eq:c1-minimax-ball-boundary} and
\eqref{eq:c1-minimax-boundary-descent} place it in $\Omega\times B_R^N$.

At the critical point $(s_h,z_h)$, $D_zK=D_z\rho$, so $c|z_h|^2\le\delta|z_h|$.
Consequently $|z_h|\le\delta/c$, and
\[
 0\le F(s_h)-m_*
 =J(s_h,z_h)-m_*+hK(s_h,z_h)-h\rho(s_h,z_h)
 \le\left(2\delta+\frac{C\delta^2}{2c^2}\right)h.
\]
Together with \eqref{eq:c1-minimax-level}, this proves
\eqref{eq:c1-minimax-output}.
\end{proof}

\subsection{Selection of the parameters and completion of the construction}

\begin{proof}[Proof of Theorem~\ref{thm:main} and Corollary~\ref{cor:neck-ratios}]
For $3\le n\le7$, Proposition~\ref{prop:low-global-action} gives
$j_\eps-C_\eps\to j_0$ in $C^2$. By Lemma~\ref{lem:negative-scale},
the equation
$\partial_\amp j_\eps(\amp\bgshape)=0$ has a $C^1$ solution
$\amp=\amp_\eps(\bgshape)$, with
$\amp_\eps\to\amp_0$ uniformly on $D$.
The functions $j_\eps(\amp_\eps(\bgshape)\bgshape)-C_\eps$ converge uniformly
to $\varphi_0$. The boundary gap \eqref{eq:core-boundary-gap}
places their minima in $\operatorname{int}D$, where the chain rule
and amplitude equation give $Dj_\eps=0$. Equation~\eqref{eq:exact-gradient}
then removes all residual multipliers.

For $n\ge8$, use the action of Proposition~\ref{prop:global-high-action}.
Its coefficients satisfy
\[
 \eps_i=t_i\eps,\qquad \mu_i^H=-|\Sph^{n-1}|h_{\eps_i}\asymp\eps_i^2>0,
 \qquad \sigma_i=\min(q_i,(n-2)/4).
\]
Each $\mathcal P_i$ in \eqref{eq:leading-potentials} is negative.
At $\zeta_i=0$, $-\mathcal P_i$ has a uniform positive lower bound.
As $|\zeta_i|\to\infty$, $\mathcal P_i\to0$ uniformly for $\bgpar$ in
the fixed compact parameter set. These assertions follow from
Lemma~\ref{lem:tt-coefficient}, the formula for $C_\infty$ in
\eqref{eq:outgoing}, and the positive bounds for $c_i(\bgpar)$. Choose $K_i=\overline B_{R_i}$ so that
\begin{equation}
 \sup_{\bgpar,\,|\zeta_i|=R_i}[-\mathcal P_i(\bgpar,\zeta_i)]
 \le\tfrac12\inf_\bgpar[-\mathcal P_i(\bgpar,0)].
 \label{eq:center-face-gap}
\end{equation}
For $\bgpar$ in its fixed compact set and $\zeta_i\in K_i$, the equation
$\mu_i^H+2\sigma_i e^{2\sigma_i\ell_i}\mathcal P_i(\bgpar,\zeta_i)=0$
has roots
\[
 \ell_i^*(\bgpar,\zeta_i)=\frac1{2\sigma_i}
 \log\frac{\mu_i^H}{2\sigma_i[-\mathcal P_i(\bgpar,\zeta_i)]}
 \quad\hbox{with}\quad
 \ell_i^*-\sigma_i^{-1}\log\eps_i=O(1).
\]
Choose fixed intervals for
$\ell_i-\sigma_i^{-1}\log\eps_i$ so that
$\ell_i^*(\bgpar,\zeta_i)-\sigma_i^{-1}\log\eps_i$
has a fixed positive distance from their boundaries.
The scales $e^{\ell_i^*(\bgpar,\zeta_i)}$ satisfy
\eqref{eq:bridge-scales}, so the local and exterior constructions
apply for sufficiently small $\eps$.
Write $r_\eps=o_{C^1}(\eps^2)$ for the remainder in
\eqref{eq:global-high-action}.

To apply Lemma~\ref{lem:c1-minimax}, set
$s=(\bgshape,\zeta)$ (the normalized compact parameter and all bubble centers), $\bgpar_0(\bgshape)=\amp_0(\bgshape)\bgshape$ and
$\ell_i^0(s)=\ell_i^*(\bgpar_0(\bgshape),\zeta_i)$. For
$z=(z_0,\ldots,z_m)\in B_{2R}^{m+1}$, use
\[
 \amp=\amp_0(\bgshape)+\eps z_0,\qquad \ell_i=\ell_i^0(s)+z_i.
\]
Choose $2R$ smaller than the fixed distance of $\ell_i^*$ from
the boundary of the scale interval. For small $\eps$, the map
\[
 \Xi_\eps(s,z)=
 \bigl((\amp_0(\bgshape)+\eps z_0)\bgshape,
       (\ell_i^0(s)+z_i)_{i=1}^m,\zeta\bigr)
\]
takes values in the prescribed $(\bgpar,\ell,\zeta)$ parameter
sets and has uniformly bounded first derivatives and invertible differential.
Put
\[
\begin{aligned}
 C'_\eps&=C_\eps+
       \sum_i\frac{\mu_i^H}{2\sigma_i}(\log\mu_i^H-\log(2\sigma_i)-1),\\
 F_\eps(s)&=C'_\eps+\varphi_0(\bgshape)
       -\sum_i\frac{\mu_i^H}{2\sigma_i}
                     \log[-\mathcal P_i(\bgpar_0(\bgshape),\zeta_i)],\\
 K_\eps(s,z)&=
       \frac{\varphi_0(\bgshape)-j_0((\amp_0(\bgshape)+\eps z_0)\bgshape)}{\eps^2}
       +\sum_i\frac{\mu_i^H}{2\sigma_i\eps^2}(e^{2\sigma_iz_i}-1-2\sigma_iz_i).
\end{aligned}
\]
Then
\[
 I_\eps\circ\Xi_\eps=F_\eps-\eps^2K_\eps+\eps^2\rho_\eps,
 \qquad \|\rho_\eps\|_{C^1}
 \le C\eps^{-2}\|r_\eps\|_{C^1}+C\eps=o(1).
\]
The full remainder satisfies
\[
 \eps^2\rho_\eps
 =r_\eps\circ\Xi_\eps+
 \sum_i e^{2\sigma_i(\ell_i^0+z_i)}
 [\mathcal P_i(\bgpar_0+\eps z_0\bgshape,\zeta_i)
       -\mathcal P_i(\bgpar_0,\zeta_i)].
\]
The exponential factor is $O_{C^1}(\eps^2)$.
The difference of the two values of $\mathcal P_i$ is
$O_{C^1}(\eps)$, because $\mathcal P_i$ is $C^2$ on the fixed
compact parameter sets.
Stationarity gives $K_\eps(s,0)=D_zK_\eps(s,0)=0$.
The Hessian of $K_\eps$ with respect to $z$ is diagonal, with entries
\[
 -\partial_\amp^2j_0((\amp_0(\bgshape)+\eps z_0)\bgshape),\qquad
 \frac{2\sigma_i\mu_i^H}{\eps^2}e^{2\sigma_iz_i},\quad 1\le i\le m.
\]
Thus $K_\eps,D_zK_\eps$ are $C^1$ and
$c'\operatorname{Id}\le D_z^2K_\eps\le C'\operatorname{Id}$.

Let $\Omega_0=D\times\prod_iK_i$ and
$m_\eps=\min_{\Omega_0}F_\eps$. On $\{|\zeta_i|=R_i\}$, replacing
$\zeta_i$ by zero lowers $F_\eps$ by at least
$\mu_i^H\log2/(2\sigma_i)\ge c\eps^2$, by \eqref{eq:center-face-gap}.
The boundary $\partial D\times\prod_iK_i$ has a fixed positive gap by
\eqref{eq:core-boundary-gap} and
$F_\eps=C'_\eps+\varphi_0+O(\eps^2)$.
Hence $F_\eps\ge m_\eps+\gamma\eps^2$ on $\partial\Omega_0$.
Apply Lemma~\ref{lem:c1-minimax} with
$(J,F,K,\rho,h,N)=(I_\eps\circ\Xi_\eps,F_\eps,K_\eps,\rho_\eps,\eps^2,m+1)$.
It supplies an interior critical point with $z_\eps=o(1)$ and
$F_\eps(s_\eps)-m_\eps=o(\eps^2)$.
Invertibility of $D\Xi_\eps$ gives $DI_\eps=0$, so
\eqref{eq:full-multiplier-pairing} removes every residual multiplier.

In both dimension ranges, the selected parameters satisfy
\begin{equation}
 \amp_\eps-\amp_0(\bgshape_\eps)=o(1),\qquad
 \varphi_0(\bgshape_\eps)-\min_D\varphi_0=o(1).
 \label{eq:selection-parameter-location}
\end{equation}
The first estimate follows from $\amp_\eps\to\amp_0$ in low dimensions
and from $\amp=\amp_0+\eps z_0$ in high dimensions. The second follows
from the bounds on the selected action values.

Let $u$ be the assembled conformal factor at the selected
critical point. Its values and physical conormals agree on every
interface. Integration by parts therefore gives
\[
 \int_{M\setminus P}
 (\langle\nabla u,\nabla\phi\rangle
  +c(n)\operatorname{Scal}_gu\phi-\frac{n(n-2)}4u^p\phi)\vol=0,
 \qquad \phi\in C_c^\infty(M\setminus P).
\]
The local relative estimates and the positive exterior family
$U_{\rm ext}$ give $u>0$;
elliptic regularity \cite[Theorem 2.4]{lp} then gives
$u\in C^\infty(M\setminus P)$.
Near $p_i$, the bound $u\ge c\eps_i r^{-(n-2)/2}$ implies
\[
 L_{u^{4/(n-2)}g}(\gamma)
 =\int u^{2/(n-2)}|\dot\gamma|_g\dd s
 \ge c\eps_i^{2/(n-2)}\int\frac{|r'(s)|}{r(s)}\dd s=\infty
\]
for every locally rectifiable curve approaching $p_i$.
A continuous curve leaving every compact subset of $M\setminus P$
eventually stays in one of the disjoint puncture neighborhoods and
approaches its puncture. This proves completeness and, together with
Proposition~\ref{prop:necessity-round}, Theorem~\ref{thm:main}.

The Fowler necks are exactly $t_i\eps$ for every small $\eps$.
Take any sequence $\eps_j\downarrow0$.
Compactness of $D$ and \eqref{eq:selection-parameter-location} give
$\operatorname{dist}(\bgshape_\eps,\mathcal S)\to0$.
Choose a subsequence with $\bgshape_{\eps_j}\to \bgshape_*\in\mathcal S$. Then
\[
 \bgpar_{\eps_j}\longrightarrow\amp_0(\bgshape_*)\bgshape_*,
 \qquad U_{\amp_0(\bgshape_*)\bgshape_*}\in\K
\]
by \eqref{eq:reduced-minimum}. Lemma~\ref{lem:compact-graph} gives
smooth convergence of $U_{\bgpar_{\eps_j}}$ on $M$.
The exterior estimates in Sections~\ref{sec:exterior-inverse}
and~\ref{subsec:high-exterior-problem} give
$u_{\eps_j}-U_{\bgpar_{\eps_j}}\to0$ uniformly on compact subsets
of $M\setminus P$. Thus $u_{\eps_j}$ converges there to the same limit.
The limiting minimizer is positive, so uniform convergence gives a
common positive lower bound on each fixed compact subset away from
$P$. The equation and interior elliptic estimates first give uniform
local $C^{2,\alphaH}$ bounds and then, since the power is smooth on this
positive range, bounds of every order. Applying these estimates on
nested compact subsets upgrades the convergence to
$C^\infty_{\rm loc}(M\setminus P)$; compare
\cite[Lemma A.4(iii), proof Step 4]{part2}.
This proves Corollary~\ref{cor:neck-ratios}.
\end{proof}

\appendix
\refstepcounter{section}
\section*{Appendix A. Coordinate identities for the inner end}
\label{app:coordinates}
\addcontentsline{toc}{section}{Appendix A. Coordinate identities for the inner end}
We prove the coordinate identities in Lemma~\ref{lem:tail-gauge} and
\eqref{eq:repair-flat-ward} for arbitrary inputs in
\eqref{eq:inner-positive-neighborhood}, using the fixed spaces and
weights of Section~\ref{sec:high-local}. The scalar cancellation in
fixed conformal normal coordinates is recorded in
\cite[Lemma 2.8]{santos}. Here the Kelvin center may vary, and we
include the coordinate and conformal changes that remove the scalar
tensor components of degrees zero and one.

\subsection{Removal of the low scalar tensor components}

We first remove the scalar tensor components of degrees zero and one.
Use the notation of Section~\ref{subsec:tail-equations} and set
$\tau_\kelpar=\operatorname{tr}\overline H_\kelpar$,
$\overline H_\kelpar^{\rm tf}=\overline H_\kelpar-(\tau_\kelpar/n)g_{\rm E}$.
The orthogonal projection of $\overline H_\kelpar^{\rm tf}$ onto
the scalar tensor components of degrees zero and one is
\begin{equation}
\begin{aligned}
 (H_\kelpar^{\rm low})_{rr}&=A_0+\sum_{i=1}^nA_i\theta_i,&
 (H_\kelpar^{\rm low})_{rA}&=r\sum_{i=1}^nB_i\nabla_A\theta_i,\\
 (H_\kelpar^{\rm low})_{AB}
 &=-\frac{A_0+\sum_{i=1}^nA_i\theta_i}{n-1}r^2\sigma_{AB}.
\end{aligned}
\label{eq:scalar-tensor-projections}
\end{equation}
Here $\sigma$ is the unit-sphere metric, and angular indices are raised
with $\sigma$. The coefficients depend on $(r,\kelpar)$ and are
\begin{equation}
\begin{aligned}
 A_0&=\fint (\overline H_\kelpar^{\rm tf})_{rr},&
 A_i&=n\fint [(\overline H_\kelpar^{\rm tf})_{rr}\theta_i],\\
 B_i&=\frac n{n-1}\fint [r^{-1}(\overline H_\kelpar^{\rm tf})_{rA}
                                      \nabla^A\theta_i].
\end{aligned}
\label{eq:scalar-projection-coefficients}
\end{equation}

For degree zero, write $X_0=r\alpha_g\partial_r$. The $rr$ and
angular equations for $\mathcal L_{X_0}g_{\rm E}+\sigma_{{\rm g},0}g_{\rm E}$ are
\[
 2\alpha_g+2r\alpha_g'+\sigma_{{\rm g},0}=A_0,\qquad
 2\alpha_g+\sigma_{{\rm g},0}=-A_0/(n-1).
\]
Subtracting the equations gives $r\alpha_g'=c_gA_0$,
$c_g=n/[2(n-1)]$. The solution satisfying
$\alpha_g=O(\lambda r^q)$ as $r\to0$ is
\[
 \alpha_g=c_g\int_0^r A_0(s)\frac{\dd s}{s},\qquad
 \sigma_{{\rm g},0}=-A_0/(n-1)-2\alpha_g.
\]

For degree one, put
$X_1=\sum_{i=1}^n(\chi_i\theta_i\partial_r+\psi_i\nabla\theta_i)$,
where $\nabla\theta_i$ is written in angular coordinates.
The $rr$, $rA$, and angular equations are, respectively,
\[
 2\chi_i'+(\sigma_1)_i=A_i,\qquad
 \chi_i+r^2\psi_i'=rB_i,\qquad
 2\chi_i/r-2\psi_i+(\sigma_1)_i=-A_i/(n-1).
\]
Set $\chi_i=rz_i^g$. Eliminating $(\sigma_1)_i$ gives
$\psi_i=c_gA_i-r(z_i^g)'$, and the $rA$ equation becomes
\begin{equation}
 r^2(z_i^g)''+r(z_i^g)'-z_i^g=c_g rA_i'-B_i.
 \label{eq:scalar-gauge-euler}
\end{equation}
Thus, with $f_i=c_g rA_i'-B_i$,
\[
\begin{gathered}
 z_i^g=\frac r2\int_0^r f_i(s)s^{-2}\dd s
                 -\frac1{2r}\int_0^r f_i(s)\dd s,\\
 \chi_i=rz_i^g,\qquad \psi_i=c_gA_i-r(z_i^g)',\qquad
 \sigma_1=\sum_{i=1}^n(A_i-2\chi_i')\theta_i.
\end{gathered}
\]
Differentiating the displayed expression for $z_i^g$ verifies
\eqref{eq:scalar-gauge-euler}. The $rr$ and angular equations then
follow from $1-2c_g=-1/(n-1)$.

The homogeneous solutions $r^{\pm1}$ are excluded by
$z_i^g=O(\lambda r^q)$, $q\ge2$. A degree-$q$ term is divided
by $q^2-1$ in the formula for degree one and by $q$ in the formula
for degree zero. The integrals preserve $O(\lambda r^{q+1})$ Taylor
remainders and every fixed number of $r\partial_r$ derivatives.
Set $X_\kelpar=X_0+X_1$ and
$\sigma_\kelpar=\sigma_{{\rm g},0}+\sigma_1+\tau_\kelpar/n$ in
\eqref{eq:finite-gauge-maps}. The component equations make the
resulting $T_\kelpar$ trace-free with zero scalar projections in
degrees zero and one.

We now compute the action density defined in
Section~\ref{subsec:tail-equations}. For this calculation write
\[
 m_\upsilon=(C_\upsilon/r)^{(n-2)/2},\qquad
 K_\upsilon=\Phi_\upsilon^*h_\upsilon,\qquad
 \mu_\upsilon=\frac{\dd V_{K_\upsilon}}{\dd t\,\dd\theta}.
\]
Thus $u_\upsilon\circ\Phi_\upsilon=m_\upsilon U$.
After pullback by $\Phi_\upsilon$, the density of
$J_{h_\upsilon}(u_\upsilon)$ with respect to $\dd t\,\dd\theta$
has the coefficient fields
\begin{gather*}
 A_\upsilon=\mu_\upsilon m_\upsilon^2K_\upsilon^{-1},\qquad
 b_\upsilon=\mu_\upsilon m_\upsilon K_\upsilon^{-1}dm_\upsilon,\qquad
 d_\upsilon=\mu_\upsilon m_\upsilon^{p+1},\\
 c_\upsilon=\mu_\upsilon\bigl(|dm_\upsilon|_{K_\upsilon}^2
       +c(n)\operatorname{Scal}_{h_\upsilon}(\Phi_\upsilon)m_\upsilon^2\bigr).
\end{gather*}
Here $A_\upsilon$ is a symmetric matrix, $b_\upsilon$ is a vector,
and $c_\upsilon,d_\upsilon$ are scalars. The density is
\begin{equation}
 \mathcal L_\upsilon(U)
 =\tfrac12\langle A_\upsilon DU,DU\rangle
  +U\,b_\upsilon\cdot DU+\tfrac12c_\upsilon U^2
  -\frac{(n-2)^2}{8}d_\upsilon U^{p+1}.
 \label{eq:repair-raw-density}
\end{equation}
Indeed, the product rule gives
\[
\begin{aligned}
 \partial_\mu(m_\upsilon U)\partial_\nu(m_\upsilon U)
 ={}&m_\upsilon^2U_\mu U_\nu
 +m_\upsilon U\bigl((\partial_\mu m_\upsilon)U_\nu
                    +(\partial_\nu m_\upsilon)U_\mu\bigr)\\
 &+(\partial_\mu m_\upsilon)(\partial_\nu m_\upsilon)U^2.
\end{aligned}
\]
Substitution in the physical action yields the coefficients in
\eqref{eq:repair-raw-density}, with scalar curvature
$\operatorname{Scal}_{h_\upsilon}(\Phi_\upsilon)$.
On the finite cylinder $\mathcal C_T=[0,T]\times\Sph^{n-1}$,
\[
 \int_{\mathcal C_T}b_\upsilon\cdot DU\,U
 =-\frac12\int_{\mathcal C_T}(\operatorname{div}b_\upsilon)U^2
   +\frac12\left[\int_{\Sph^{n-1}}b_\upsilon^tU^2\dd\theta\right]_0^T.
\]
Varying $U$ and dividing by $d_\upsilon$ gives
\begin{equation}
 \mathcal E_\upsilon(U)
 =d_\upsilon^{-1}\Bigl\{-\operatorname{div}(A_\upsilon DU)
 +(c_\upsilon-\operatorname{div}b_\upsilon)U
 -\frac{n(n-2)}4d_\upsilon U^p\Bigr\}.
 \label{eq:repair-exact-operator}
\end{equation}
Here divergence is taken relative to $\dd t\,\dd\theta$.
This is the coefficient form of the operator for
$\gamma_\upsilon=(C_\upsilon/r)^2K_\upsilon$ used in
Section~\ref{subsec:tail-equations}.
The conormal density paired with a variation $\dot U$ is
\[
 \mathcal N_\upsilon(U)
 =\sum_{\nu=0}^{n-1}(A_\upsilon)^{t\nu}U_\nu
       +(b_\upsilon)^tU.
\]
Its first-variation term has positive sign at the upper $t$ endpoint
and negative sign at the lower endpoint. The retained action term
$b_\upsilon^tU^2/2$ has the same endpoint sign convention.

Fix a smooth bounded domain $\Omega$ and put $J_0(u)=J_{g_{\rm E}}(u;\Omega)$.
For the following variation identity, let $X$ be a vector field and let $\sigma_{\rm g}$ be a scalar function. Define
\[
 h^{\rm g}=\mathcal L_Xg_{\rm E}+\sigma_{\rm g}g_{\rm E},\qquad
 J_1(u)=\left.\frac{\mathrm d}{\mathrm d\upsilon}J_{g_{\rm E}+\upsilon h^{\rm g}}(u;\Omega)\right|_{\upsilon=0},\qquad
 \mathsf D u=X\cdot\nabla u-\frac{n-2}{4}\sigma_{\rm g}u.
\]

Let $\Phi_\upsilon^X$ be the flow of $X$. Differentiating
\[
 J_{(\Phi_\upsilon^X)^*g_{\rm E}}(u\circ\Phi_\upsilon^X;\Omega)
 =J_{g_{\rm E}}(u;\Phi_\upsilon^X(\Omega))
\]
shows that the coordinate contribution to
$J_1(u)+DJ_0(u)[\mathsf Du]$ is
\[
 \int_{\partial\Omega}(X\cdot\nu_{g_{\rm E}})
       \left(\tfrac12|\nabla u|^2-\tfrac{(n-2)^2}{8}u^{p+1}\right)\dd S_{g_{\rm E}}.
\]
For the conformal contribution, use
$\left.\partial_\upsilon\operatorname{Scal}_{e^{\upsilon\sigma_{\rm g}}g_{\rm E}}
\right|_0=-(n-1)\Delta\sigma_{\rm g}$,
$c(n)(n-1)=(n-2)/4$, and $n/(p+1)=\frac{n-2}{2}$.
The gradient and power terms cancel between the metric variation
and the field variation $-\frac{n-2}{4}\sigma_{\rm g}u$, leaving
\[
 -\frac{n-2}{8}\int_{\Omega}
       \bigl(u^2\Delta\sigma_{\rm g}+2u\nabla u\cdot\nabla\sigma_{\rm g}\bigr)\dd X
 =-\frac{n-2}{8}\int_{\partial\Omega}u^2\partial_{\nu_{g_{\rm E}}}\sigma_{\rm g}\dd S_{g_{\rm E}}.
\]
Adding the two contributions gives

\begin{align}
 J_1(u)+DJ_0(u)[\mathsf D u]
 &=\mathcal B_{X,\sigma_{\rm g}}(u),\label{eq:gauge-identity}\\
 \mathcal B_{X,\sigma_{\rm g}}(u)
 &=\int_{\partial\Omega}
 \left[(X\cdot\nu_{g_{\rm E}})\left(\frac{|\nabla u|^2}{2}
                         -\frac{(n-2)^2}{8}u^{p+1}\right)
             -\frac{n-2}{8}u^2\partial_{\nu_{g_{\rm E}}}\sigma_{\rm g}\right]\dd S_{g_{\rm E}}.\nonumber
\end{align}
Here $\nu_{g_{\rm E}}$ is outward on every boundary component.

\subsection{The moving Kelvin coordinates and the action identities}

For a varying Kelvin center $\kelpar(t)$, recall
\[
 \eta=r\theta\cdot\kelpar',\qquad L=1+\eta,\qquad
 \beta=\nabla_\theta\eta,\qquad
 \mathscr D_t=\partial_t-\beta\cdot\nabla_\theta.
\]
Since $F_{\kelpar+\xi}=F_\kelpar\circ F_\xi$, the differential of
the varying-center map $\Phi_0$ factors into $DF_\kelpar$ and the
motion of $r\theta$.
At $\xi=0$ the Kelvin generator is
$r^2\xi-2r^2(\theta\cdot\xi)\theta$. Therefore, after applying
$(DF_\kelpar)^{-1}$,
\[
 d\Phi_0=DF_\kelpar\bigl[-rL\theta\,\dd t+
                       r(\dd\theta+\beta\,\dd t)\bigr].
\]
The identity $(DF_\kelpar)^TDF_\kelpar=E_\kelpar^{-2}I$ gives the
metric $\gamma_0=L^2dt^2+|d\theta+\beta dt|^2$ and its inverse
\[
 \gamma_0^{tt}=L^{-2},\quad
 \gamma_0^{tA}=-L^{-2}\beta^A,\quad
 \gamma_0^{AB}=\sigma^{AB}+L^{-2}\beta^A\beta^B.
\]
Using $E_\kelpar=1+2r\kelpar\cdot\theta+r^2|\kelpar|^2$, we find
\[
 \nabla_\theta E_\kelpar=2r[\kelpar-(\kelpar\cdot\theta)\theta],\qquad
 \mathscr D_tE_\kelpar
 =-2r\kelpar\cdot\theta-2r^2|\kelpar|^2
                      +2\eta(1+r\kelpar\cdot\theta),
\]
hence
\begin{equation}
 2E_\kelpar+\mathscr D_tE_\kelpar
 =2+2r\kelpar\cdot\theta+2\eta(1+r\kelpar\cdot\theta)
 =2L(1+r\kelpar\cdot\theta).
 \label{eq:kelvin-endpoint-identity}
\end{equation}
At $\upsilon=0$, the physical metric is Euclidean,
$K_0=(r/E_\kelpar)^2\gamma_0$, $m_\upsilon|_{\upsilon=0}=(E_\kelpar/r)^{(n-2)/2}$ and $d_0=L$.
Substitution in \eqref{eq:repair-raw-density} expresses the difference
of the cylindrical action densities as a divergence:
\[
\begin{split}
 \mathcal L_0-\mathcal L_{\rm can}
 ={}&\frac{n-2}{4}\partial_t
       \left[\frac{U^2}{L}(2+\mathscr D_t\log E_\kelpar)\right]\\
 &+\frac{n-2}{4}\operatorname{div}_\theta\left\{
 U^2\left[L\nabla\log E_\kelpar
       -\frac\beta L(2+\mathscr D_t\log E_\kelpar)\right]\right\}.
\end{split}
\]
To verify it, expand the squared gradients and use
$\Delta_\theta E_\kelpar=-2(n-1)r\kelpar\cdot\theta$ and
\eqref{eq:kelvin-endpoint-identity}. The coefficients of $U\nabla U$
and $UU_t$ agree directly. For the $U^2$ coefficients, use
$E_\kelpar=1+2r\kelpar\cdot\theta+r^2|\kelpar|^2$,
$L=1+\eta$ and $\operatorname{div}_\theta\beta=-(n-1)\eta$ to obtain
\[
\begin{aligned}
 &\partial_t\frac{2+\mathscr D_t\log E_\kelpar}{L}
 +\operatorname{div}_\theta\!\left[
 L\nabla_\theta\log E_\kelpar
 -\frac\beta L(2+\mathscr D_t\log E_\kelpar)\right]\\
 &=\frac{2(n-2)(\eta-r\kelpar\cdot\theta)}{E_\kelpar}\\
 &=\frac{n-2}{2}\left[
 \frac{2\mathscr D_t\log E_\kelpar+(\mathscr D_t\log E_\kelpar)^2}{L}
 +L|\nabla_\theta\log E_\kelpar|^2\right].
\end{aligned}
\]
This proves equality of the $U^2$ coefficients for every $U$.
The angular divergence integrates to zero on each sphere.
For the interior variation, the density $\mathcal L_{\rm can}$ has
derivatives
\[
 \partial_{U_t}\mathcal L_{\rm can}=-F_U/L,\qquad
 \partial_{\nabla U}\mathcal L_{\rm can}=L\nabla U+\beta F_U/L,
 \qquad
 \partial_U\mathcal L_{\rm can}=\frac{n-2}{2}F_U/L-\frac{n(n-2)}4LU^p.
\]
The Euler--Lagrange expression of $\mathcal L_{\rm can}$, divided by $L$ and using
$\Delta_\theta L=-(n-1)\eta$, is \eqref{eq:repair-flat-moving}.
Expanding $L^{-j}=1-j\eta+O(\eta^2)$, $j=1,2,3$, gives
\eqref{eq:repair-linear-motion}. The identical power terms cancel
in \eqref{eq:repair-motion-remainder}.
If $\mathcal E_0^{[s]}$ denotes the expression with
$(\kelpar',\kelpar'')$ replaced by $(s\kelpar',s\kelpar'')$, then
\[
 \mathcal R_{\rm mov}U
 =\int_0^1(1-s)\partial_s^2\mathcal E_0^{[s]}(U)\dd s.
\]
Each coefficient of $\mathcal R_{\rm mov}$ contains at least two factors from $r\kelpar'$ and $r\kelpar''$.
For an arbitrary radial profile,
\[
 \mathcal TV=-[j_V\kelpar''+2j_V'\kelpar']\cdot\theta,
 \qquad
\mathcal E_0(V)=\frac{-V''+\frac{(n-2)^2}{4}V+(n-1)\eta F_V}{L^2}
                    -F_V\frac{\eta_t-|\beta|^2}{L^3}-\frac{n(n-2)}4V^p.
\]

At fixed $\kelpar$, the first variation of the conformal Laplacian
in the trace-free direction $T_\kelpar$ is
\eqref{eq:repair-static-metric}; Taylor's formula in $\upsilon$
gives \eqref{eq:repair-metric-remainders}.
To prove \eqref{eq:repair-low-selection}, set $u=r^{-(n-2)/2}V$ and use a test
function $\psi$. On a Euclidean annulus $\mathcal A$ integration by
parts gives
\begin{equation}
\begin{aligned}
 &\int_{\mathcal A}\psi\sum_{i,j=1}^n\left\{
       \partial_i(T_\kelpar^{ij}u_j)
                   +c(n)(\partial_i\partial_jT_\kelpar^{ij})u\right\}\dd X\\
 &\quad=\int_{\mathcal A}\sum_{i,j=1}^nT_\kelpar^{ij}
 \left[c(n)u\psi_{ij}+c(n)\psi u_{ij}
       +(2c(n)-1)u_{(i}\psi_{j)}\right]\dd X+\mathcal B_{\mathcal A},
\end{aligned}
\label{eq:metric-tensor-green-identity}
\end{equation}
where $u_{(i}\psi_{j)}=(u_i\psi_j+u_j\psi_i)/2$ and
\[
 \mathcal B_{\mathcal A}
 =\int_{\partial\mathcal A}\sum_{i,j=1}^n\left[
 \psi(\nu_{g_{\rm E}})_iT_\kelpar^{ij}u_j+
 c(n)u\psi(\nu_{g_{\rm E}})_i\partial_jT_\kelpar^{ij}
       -c(n)(\nu_{g_{\rm E}})_jT_\kelpar^{ij}\partial_i(u\psi)\right]\dd S_{g_{\rm E}}.
\]
For radial $u$ and $\psi=\psi_0(r)$, the bracketed tensor paired
with $T_\kelpar^{ij}$ in \eqref{eq:metric-tensor-green-identity} is a
radial combination of $\theta\otimes\theta$ and
$g_{\rm E}-\theta\otimes\theta$. For
$\psi=\psi_1(r)\theta_\ell$, its $rr$, $rA$ and $AB$ components
are radial multiples of $\theta_\ell$, $\nabla_A\theta_\ell$ and
$\theta_\ell\sigma_{AB}$, respectively.
We compute these components from
$d\psi=\psi_1'\theta_\ell\,dr+\psi_1d\theta_\ell$ and
$\nabla_\theta^2\theta_\ell=-\theta_\ell\sigma$.
The test tensors obtained from $\psi_0(r)$ and
$\psi_1(r)\theta_\ell$ pair to zero with $T_\kelpar$: its degree-zero $rr$
and angular trace averages vanish, as do the degree-one $rr$,
$rA$ and angular trace pairings. For compactly supported
$\psi_0,\psi_1$, the boundary term also vanishes. Since these radial
test functions are arbitrary, \eqref{eq:repair-low-selection} holds
at every radius. On finite annuli we retain $\mathcal B_{\mathcal A}$.
Identity \eqref{eq:repair-low-selection} holds for every smooth
radial $V$, including solutions of the forced radial equation.

For \eqref{eq:repair-flat-ward}, vary the center with compact support.
At fixed $U$, the physical center variation divided by the
conformal factor is $Z_i(U)$ in
\eqref{eq:physical-center-variation}. Indeed the Kelvin
generator in cylinder coordinates is
$r\theta_i\partial_t+r\nabla\theta_i$; the inverse differential of
the moving chart changes it to
\[
 (r\theta_i/L)\partial_t+
       (r\nabla\theta_i-r\theta_i\beta/L)\cdot\nabla_\theta.
\]
Combining the displayed vector field with the conformal weight gives
$Z_i(U)$ in \eqref{eq:physical-center-variation}.
For a compactly supported center variation, the momentum in
\eqref{eq:repair-momentum} is the normalized spherical average
\[
 \mathcal J_i=n\fint \bigl(\partial_{\kelpar_i'}\mathcal L_{\rm can}\bigr).
\]
Here $\mathcal L_{\rm can}$ is defined in
\eqref{eq:cylindrical-lagrangian}. Integration by parts gives
$-\mathcal J_i'$ for the Euler--Lagrange expression in the center
coordinate, with this normalization. The physical field variation gives
$n\fint [L\mathcal E_0(U)Z_i(U)]$, proving
\eqref{eq:repair-flat-ward}.
If $\Pi_{\ge2}w=w$, the remaining term in $\fint (LZ_i(V+w))$ is
$\fint [(\theta_i\nabla\eta-\eta\nabla\theta_i)\cdot\nabla w]$.
It vanishes by integration by parts, since
\[
 \operatorname{div}(\theta_i\nabla\eta-\eta\nabla\theta_i)
 =\theta_i\Delta\eta-\eta\Delta\theta_i=0.
\]
Using $LZ_i(V+w)=j_V\theta_i+rQ_i(w)$ and
$n\fint (\theta_i\theta_j)=\delta_{ij}$ gives
\[
\begin{aligned}
 (G_{\rm lin})_{ij}-\delta_{ij}
 &=\frac{nr}{j_V}\fint [Q_i(w)\theta_j],\\
 |(G_{\rm lin})_{ij}-\delta_{ij}|
 &\le C\frac{|D^{\le1}w|}{F_V}\le C\delta W.
\end{aligned}
\]

For radial $V$, \eqref{eq:repair-momentum} gives
\[
 \mathcal J(V,\kelpar')=-\frac{nrF_V^2}{2}
       \fint [\theta(1+r\theta\cdot\kelpar')^{-2}].
\]
Differentiate $\mathcal J(V,\kelpar')$ with respect to $\kelpar'$, then
integrate from $0$ to $\kelpar'$. This gives
$\mathcal J(V,\kelpar')=\mathsf p_VC(\kelpar')\kelpar'$, with the
positive matrix $C$ in \eqref{eq:center-momentum-matrix}. The linear term of $C-I$
is a third spherical moment, hence zero; thus
$C-I=O(r^2|\kelpar'|^2)$. Integrating the derivative of
$\mathcal J(V+w,\kelpar')-\mathcal J(V,\kelpar')$ from zero velocity
to $\kelpar'$ gives the term $\mathsf p_VA_w\kelpar'$.
At zero velocity, $\mathcal J(V,0)=0$ and
$D_U\mathcal J(V,0)[w]=0$: the linear terms are multiples of
$\fint (\theta_iw)$, $\fint (\theta_iw_t)$ and
$\fint (\nabla\theta_i\cdot\nabla w)=(n-1)\fint (\theta_iw)$.
Taylor's formula gives
\begin{equation}
 J_w=\int_0^1(1-s)D_U^2\mathcal J(V+sw,0)[w,w]\dd s.
 \label{eq:kelvin-momentum-taylor}
\end{equation}
The second derivative of \eqref{eq:repair-momentum} with respect to
$U$ contains products of the variations and their first cylinder
derivatives, with bounded power coefficient $(V+sw)^{p-1}$. Hence
$|J_w|\le Cr|D^{\le1}w|^2$.
For fixed $V,w$, replacing the two field directions in
\eqref{eq:kelvin-momentum-taylor} by $h_1,h_2$ defines the symmetric
bilinear form $\mathcal B_{V,w}$ used in Appendix~\ref{app:weighted}.

The divergence identity for the action densities gives the boundary
primitive
\[
 \frac{n-2}{4}\int_{\Sph^{n-1}}
       U^2(2+\mathscr D_t\log E_\kelpar)/L\,\dd\theta.
\]
Equation~\eqref{eq:kelvin-endpoint-identity} reduces this primitive to
\eqref{eq:repair-physical-primitive}, independent of $\kelpar'$.
On each finite annulus, retain both the term in
\eqref{eq:repair-physical-primitive} and the boundary integral in
\eqref{eq:gauge-identity} on both boundary components.

Four cylinder derivatives of the source equations use at most six
physical field derivatives. For $\radsrc\in C^{4,\alphaH}$, the radial
field is $C^{6,\alphaH}$ and its phase and energy coordinates are
$C^{5,\alphaH}$. The equivalence of the radial, high harmonic and
center equations with the complete equation follows from
\eqref{eq:repair-final-system} and the invertibility of $G_{\rm lin}$
in \eqref{eq:linear-harmonic-pairing}.

\subsection{Coefficient and parameter-derivative bounds}

The integral formulas for $X_\kelpar$ and $\sigma_\kelpar$ give
\begin{equation}
 X_\kelpar=O(\lambda r^{q+1}),\qquad
 \partial_\kelpar X_\kelpar=O(\lambda r^{q+2}),\qquad
 \partial_\kelpar\sigma_\kelpar=O(\lambda r^{q+1}).
\label{eq:finite-gauge-coefficient-bounds}
\end{equation}
More generally, $D_\kelpar^jX_\kelpar=O(\lambda r^{q+1+j})$ for
$j\le2$. Differentiate the flow equation
$\partial_sQ_s=-\upsilon X_\kelpar(Q_s)$, $Q_0=\operatorname{Id}$,
with respect to $\kelpar$.
The linearized equations have coefficients $O(\lambda r^q)$.
The first and second center sources have sizes
$O(\lambda r^{q+2})$ and $O(\lambda r^{q+3})$, respectively.
Products with an extra factor $\lambda r^q$ are smaller.
Gronwall's inequality yields
\[
 \sup_{0\le s\le1}|\partial_{\kelpar_i}Q_s|\le C\lambda r^{q+2},\qquad
 \sup_{0\le s\le1}|\partial_{\kelpar_i}\partial_{\kelpar_j}Q_s|
 \le C\lambda r^{q+3}.
\]

A center derivative of the Kelvin map is $O(r^2)$, adding a factor
$r$ to a normalized metric coefficient. Differentiation in the homotopy
parameter $\upsilon$ introduces a factor $O(\lambda r^q)$.
Since $\mathcal M_{\rm mov}$ vanishes at $z'=z''=0$, the fundamental
theorem of calculus in these variables gives its factors
$rz^{(j)}O(\lambda r^q)$, $j=1,2$.
The metric coefficients depend on the background and center only;
radial profile factors remain in $D^\alpha U$.

Expanding the divergence in \eqref{eq:repair-exact-operator},
the coefficient fields of its linear part are formed from
\[
 d_\upsilon^{-1}A_\upsilon,\qquad
 d_\upsilon^{-1}\operatorname{div}A_\upsilon,\qquad
 d_\upsilon^{-1}(c_\upsilon-\operatorname{div}b_\upsilon).
\]
The definitions give
\[
 d_\upsilon=\frac{\dd V_{\gamma_\upsilon}}{\dd t\,\dd\theta},
 \qquad d_\upsilon^{-1}A_\upsilon=\gamma_\upsilon^{-1},
 \qquad
 d_\upsilon^{-1}(c_\upsilon-\operatorname{div}b_\upsilon)
       =c(n)\operatorname{Scal}_{\gamma_\upsilon}.
\]
These coefficient fields are smooth functions of the evaluated metric
derivatives $\partial_y^{\le2}H_b(\Phi_\upsilon)$ and the derivatives
through order two of $\Phi_\upsilon$ and $C_\upsilon$.
Their dependence on $r\kelpar''$ is affine.
Relative to the product metric $\dd t^2+g_{\Sph^{n-1}}$,
\eqref{eq:inner-positive-neighborhood} gives
\[
 c(\dd t^2+g_{\Sph^{n-1}})\le\gamma_\upsilon
       \le C(\dd t^2+g_{\Sph^{n-1}}),\qquad
 c\le d_\upsilon\le C,\qquad L\ge1-\eta_0.
\]
The scaled coordinate derivative $r^{-1}D\Phi_\upsilon$, viewed
from the product tangent space to Euclidean space, and its inverse
are uniformly bounded as well.

For example, a first center variation of an evaluated metric factor
contains $\partial_yH_b(\Phi_\upsilon)D_\vartheta\Phi_\upsilon$.
After four cylinder derivatives, its terms have the form
\[
 \partial_y^{j+1}H_b(\Phi_\upsilon)
 \left[ D^{\alpha_1}\Phi_\upsilon,\ldots,
        D^{\alpha_j}\Phi_\upsilon,
        D^{\alpha_0}D_\vartheta\Phi_\upsilon\right],
 \quad
 \sum_{i=0}^j|\alpha_i|=4,\quad |\alpha_i|\ge1\ (i\ge1).
\]
Since $j\le4$, at most five metric derivatives occur in this
example. The operator $\partial_y$ differentiates in Euclidean
coordinates, whereas $D$ differentiates in the cylinder variables.
Since $\partial_y^{m}H_b=O(\lambda r^{q-m})$,
$D^\alpha\Phi_\upsilon=O(r)$ for $|\alpha|\ge1$, and
$D^\alpha D_\vartheta\Phi_\upsilon=O(r^2)$ for a center direction,
each term is $O(\lambda r^{q+1})$.
The second derivatives of the metric appearing in
\eqref{eq:repair-exact-operator} require at most six metric derivatives
after four cylinder derivatives, and at most seven after one parameter
variation. The center and field factors require at most six spatial
derivatives. The fixed $C^{q+8,\alphaH}$ bound controls these terms.
Applying the product, inverse and composition rules to each summand
preserves the powers of $r$ in \eqref{eq:coefficient-inventory} and
bounds its first Fr\'echet derivative in the fixed reference spaces.

\refstepcounter{section}
\section*{Appendix B. Weighted source and boundary estimates}
\label{app:weighted}
\addcontentsline{toc}{section}{Appendix B. Weighted source and boundary estimates}
All spaces and reference weights are fixed as in
Section~\ref{sec:high-local}. The unforced profile is $v_\tau$,
where $v_\tau(0)=m_0$; derivatives of the prescribed data keep
$\eps$ and the norm weights fixed. This appendix proves the estimates used
in Lemma~\ref{lem:tail-source-map} and Proposition~\ref{prop:inner-dn}.
The inputs lie in the small ball of the product norm
\eqref{eq:inner-product-norm} and satisfy
\eqref{eq:inner-positive-neighborhood}. Write
$\mathcal G=\|\radsrc\|_{\mathcal Y_{\rm rad}}$,
$W=\mathcal W_6(w)$ and $Z=K(z)$.

\subsection{Radial source and nonlinear product estimates}

The norm \eqref{eq:master} controls the following source quantities
for the nonlinear radial response in Lemma~\ref{lem:radial-response}:
\begin{equation}
 \sum_{j=0}^4\left[
 \eps^{-2}\int B_M^2v|D_t^j\radsrc|+\int B_M^2|D_t^j\radsrc|/v
 +\eps^{-1}\|B_MD_t^j\radsrc\|_2
 +\|B_M^2D_t^j\radsrc/v\|_{C_{\rm b}^{0,\alphaH}}
 +\|B_M^2D_t^j\radsrc/F\|_{C_{\rm b}^{0,\alphaH}}\right].
 \label{eq:radial-source-norm}
\end{equation}

Adjoining $\sum_{j=0}^4\|E_j\|_{C_{\rm b}^{0,\alphaH}}$, where
$E_j=B_M^2\eps^{-2}vD_t^j\radsrc$, gives a norm equivalent to $\|\cdot\|_{\mathcal Y_{\rm rad}}$.
Indeed, the bounded relative multipliers in \eqref{eq:weight-bounds}
give
\[
 \begin{aligned}
 B_M^2D_t^j\radsrc/v&=(\eps^2/v^2)E_j,&
 B_M^2D_t^j\radsrc/F&=(\eps^2/(vF))E_j,\\
 \eps^{-2}\int B_M^2|D_t^j\radsrc|^2
 &=\int\frac{\eps^2E_j^2}{v^2B_M^2}
 \le C\|E_j\|_\infty\|E_j\|_1.
 \end{aligned}
\]
The displayed identities bound \eqref{eq:radial-source-norm}
by $\|\radsrc\|_{\mathcal Y_{\rm rad}}$.
The $L^1$ and bounded H\"older norms of $E_j$ recover
the reverse bound.

We now estimate the products that enter the radial equation.
Let $H=D^\beta w$, where $|\beta|\le2$ and $w$ has only spherical
harmonics of degree at least two. For $i,j\in\{1,2\}$, write
$z^{(i)}=\partial_t^i z$. Then
$\mathcal W_4(H)\le C\mathcal W_6(w)$.
Set
\[
 h_H=B_MH/\eps,\qquad k_i=B_Mrvz^{(i)}/\eps,\qquad
 m_q=(\lambda/\eps)B_Mr^qv.
\]
The factors $h_H,k_i,m_q$ lie in
$L^2\cap C_{\rm b}^{0,\alphaH}$, with bounds
$\mathcal W_4(H)$, $K(z)$, and $C\delta^q$, respectively.
For $s=1,2$, write
$\|h\|_{L^s\cap C_{\rm b}^{0,\alphaH}}=\|h\|_{L^s}+\|h\|_{C_{\rm b}^{0,\alphaH}}$.
Let $H^{(1)}$ and $H^{(2)}$ be two functions of the same form as $H$.
The weighted products factor as
\[
\begin{aligned}
 B_M^2\eps^{-2}v(\lambda r^qH)&=m_qh_H,\\
 B_M^2\eps^{-2}v(H^{(1)}H^{(2)}/v)&=h_{H^{(1)}}h_{H^{(2)}},\\
 B_M^2\eps^{-2}v(rz^{(i)}H)&=k_i h_H,\\
 B_M^2\eps^{-2}v(\lambda^2r^{2q}v)&=m_q^2,\\
 B_M^2\eps^{-2}v(r^2vz^{(i)}z^{(j)})&=k_i k_j,\\
 B_M^2\eps^{-2}v(\lambda r^{q+1}vz^{(i)})&=m_qk_i.
\end{aligned}
\]
Cauchy--Schwarz and H\"older multiplication give
\[
 \|AB\|_1\le\|A\|_2\|B\|_2,\qquad
 \|AB\|_{C_{\rm b}^{0,\alphaH}}\le C\|A\|_{C_{\rm b}^{0,\alphaH}}\|B\|_{C_{\rm b}^{0,\alphaH}},
\]
and hence the weighted product bounds below. In these bounds only,
$\|\cdot\|_{\mathcal Y_{\rm rad}}$ denotes the expression in
\eqref{eq:master} also for functions of $(t,\theta)$, with $L^1$
integration over $\mathcal C$ and componentwise norms for vector or
tensor products. Angular averaging of a scalar contraction then gives
a radial source in $\mathcal Y_{\rm rad}$.
\begin{equation}
\begin{aligned}
 \|\lambda^2r^{2q}v\|_{\mathcal Y_{\rm rad}}&\le C\delta^{2q},&
 \|\lambda r^qH\|_{\mathcal Y_{\rm rad}}&\le C\delta^q\mathcal W_4(H),\\
 \|H^{(1)}H^{(2)}/v\|_{\mathcal Y_{\rm rad}}&\le C\mathcal W_4(H^{(1)})\mathcal W_4(H^{(2)}),&
 \|rz^{(i)}H\|_{\mathcal Y_{\rm rad}}&\le CK(z)\mathcal W_4(H),\\
 \|r^2vz^{(i)}z^{(j)}\|_{\mathcal Y_{\rm rad}}&\le CK(z)^2,&
 \|\lambda r^{q+1}vz^{(i)}\|_{\mathcal Y_{\rm rad}}&\le C\delta^qK(z).
\end{aligned}
\label{eq:radial-product-estimates}
\end{equation}
These bounds remain valid after contraction with a tensor coefficient
whose cylinder derivatives through order four are uniformly bounded.
Angular averaging is bounded in both the $L^1$ and bounded
H\"older components.
For $0\le\ell\le4$, differentiated products use
\begin{equation}
 D_t^\ell(ABC)=
 \sum_{\ell_1+\ell_2+\ell_3=\ell}
       \frac{\ell!}{\ell_1!\ell_2!\ell_3!}
       (D_t^{\ell_1}A)(D_t^{\ell_2}B)(D_t^{\ell_3}C).
 \label{eq:weighted-leibniz-calculus}
\end{equation}
For $0\le\ell\le4$, the fixed profile satisfies
$|D_t^\ell v|\le C_\ell v$ and
$|D_t^\ell(v^{-1})|\le C_\ell v^{-1}$;
differentiating a power of $r$ preserves that power.
Since $i,j\le2$, four derivatives of these products
use at most $D^6w,z^{(6)}$.
Each summand of \eqref{eq:weighted-leibniz-calculus} obeys the
weighted $L^1$ and bounded H\"older bound for the applicable
one of the six products in \eqref{eq:radial-product-estimates}.
For $|w|\le cV$ and $V\asymp v$, the power remainder satisfies
$|N_V(w)|\le Cw^2/v$, including $p<2$.
The terms containing $\radsrc$ have multipliers $\lambda^2r^{2q}$,
$\lambda r^{q+1}z^{(j)}$ and $r^2z^{(i)}z^{(j)}$.
Their local norms through four derivatives are bounded by
$c_{\rm pr}$ in \eqref{eq:principal-small}; multiplication in
\eqref{eq:master} therefore gives $c_{\rm pr}\|\radsrc\|_{\mathcal Y_{\rm rad}}$.

Before estimating the radial source, use
$\Pi_0\mathcal M_{\rm stat}V=0$ from
\eqref{eq:repair-low-selection}.

\subsection{Kelvin force, displacement, and momentum estimates}

For the fixed reference profile $v$, its chosen value and
negative derivative at $t=0$ give
$F(0)\asymp\delta^{(n-2)/2}$ and
$G(0)=O(\delta^{(n+2)/2})$. Duhamel's formula yields
\[
 v(t)\le C\min\{1,\delta^{(n-2)/2}e^{-(n-2)t/2}+\delta^{(n+2)/2}e^{(n-2)t/2}\}.
\]
For $m\ge2$, use $\min\{1,x\}\le x^{2/n}$ for the growing
term. Direct integration gives
\begin{equation}
\begin{aligned}
 \int_0^\infty B_M^2r^mv^2\dd t
 &\le C\delta^{m+(n-2)}\int_0^\infty
          B_M^2e^{-(m+(n-2))t}\dd t\\
 &\quad+C\delta^{m+2+4/n}\int_0^\infty
          B_M^2e^{-(m-2+4/n)t}\dd t\\
 &\le C_m\delta^{m+2}.
\end{aligned}
\label{eq:train-moment}
\end{equation}
Both exponential rates are positive for $m\ge2$.
The bound is uniform on the fixed boundary-data set.
Estimate~\eqref{eq:train-moment} also holds with $v^2$ replaced by
$vF,vG,F^2$, or $G^2$, since $0<F,G\le(n-2)v$.
We estimate the Kelvin forces in the fixed norm
$\mathcal Y_K$, so divided source terms use $f/\mathsf p$.
Lemma~\ref{lem:tail-low} applies these bounds to the inverse with
coefficient $\mathsf p_\tau$. We keep $\eps,\delta$ and all
reference weights fixed in parameter derivatives.

For the Kelvin equation, insert $LZ_i(V+w)=j_V\theta_i+rQ_i(w)$ into
$n\fint [L(\mathcal E_1-\mathcal E_0)(V+w)Z_i(V+w)]$, for $1\le i\le n$.
The term $j_V\theta\,\mathcal M_{\rm stat}V$ vanishes after angular
averaging because $\Pi_1\mathcal M_{\rm stat}V=0$.
The pairings of $j_V\theta$ with $\mathcal M_{\rm stat}w$ and of
$rQ(w)$ with $\mathcal M_{\rm stat}V$ give the first two bounds below.
The radial terms in $\mathcal M_{\rm mov}$ and
$\mathcal R_{\rm met}$ give the third and fourth bounds.
We retain the $\radsrc$ terms when substituting the forced equation
for $V''$.
The displayed expressions record bounds for these source terms;
their smooth angular coefficients have uniformly bounded cylinder
derivatives through order four and first parameter derivatives on
the fixed input ball.
\begin{equation}
 \begin{split}
 &\lambda r^{q+1}F D^{\le2}w,
 \qquad\lambda r^{q+1}(v+|\radsrc|)D^{\le1}w,\\
 &\lambda r^{q+2}F(v+|\radsrc|)z^{(j)},
 \qquad\lambda^2r^{2q+1}F(v+|\radsrc|).
 \end{split}
 \label{eq:four-forces}
\end{equation}

Let $f$ denote a source term represented by either of the first two
bounds in \eqref{eq:four-forces}, evaluated at $\radsrc=0$. Then
\begin{equation}
 |f|\le C\lambda r^{q+1}v|D^{\le2}w|,
 \label{eq:selector}
\end{equation}
since $F\le(n-2)v$. The bound for $w/F$ in $\mathcal W_6(w)$ and
\eqref{eq:weight-bounds} give, for $0\le\ell\le4$,
\[
\begin{aligned}
 |\partial_t^\ell(f/\mathsf p)|
 &\le C\lambda\delta^{q-2}e^{-(q-1+\gamma_*)t}W,\\
 q_K|\partial_t^\ell(f/\mathsf p)|
 &\le C(\lambda/\eps)\delta^{q+1}B_Me^{-(q+\gamma_*)t}W.
\end{aligned}
\]
The displayed bounds for $\partial_t^\ell(f/\mathsf p)$ and
$q_K\partial_t^\ell(f/\mathsf p)$ also hold in bounded
H\"older norms. Their positive exponential rates control the
unweighted $L^1$ and weighted $L^2$ components,
so the sum of the $\mathcal K_{\rm tail,0}$ norms over
$\ell=0,\ldots,4$ is bounded by
$C[\lambda\delta^{q-2}+(\lambda/\eps)\delta^{q+1}]W
\le C\delta^{q-1}W$.
For the two integrals $C_f,I_f$ in Lemma~\ref{lem:tail-low},
split the weight from \eqref{eq:high-work-norm} into
\[
 \beta_1=\eps^{-1}B_M,\qquad
 \beta_2=\eps^{-1}B_M\frac{\delta}{r}\sqrt{\omega},
 \qquad \beta_1^2+\beta_2^2=\beta_\eps^2.
\]
Cauchy--Schwarz with $\beta_2D^{\le2}w$ and
$\beta_1D^{\le2}w$ gives
\begin{equation}
 \begin{aligned}
 C_f&\le C(\lambda/\eps)\delta^{-1}
       \left(\int B_M^{-2}r^{2q}v^2\omega\right)^{1/2}W
       \le C(\lambda/\eps)\delta^{q-1}W,\\
 I_f&\le C\lambda\eps
       \left(\int B_M^2r^{2q+2}v^2\right)^{1/2}W
       \le C\lambda\eps\delta^{q+2}W.
 \end{aligned}
 \label{eq:selector-data-estimates}
\end{equation}
Here $v^2\omega\le C$ and \eqref{eq:train-moment} apply.
The term $\eps^{-1}\delta^{-2}(C_fI_f)^{1/2}$ in
\eqref{eq:kelvin-data-bound} is therefore bounded by
$C(\lambda/\eps)\delta^{q-3/2}W$. Thus
$\|(f,0)\|_{\mathcal Y_K}\le C\delta^{q-3/2}W$.
Let $z_f$ be the solution in Lemma~\ref{lem:tail-low} with data
$(f,0)$. The same bound holds for $K(z_f)$ and tends to zero
relative to $W$ as $\delta\to0$, since $q\ge2$.

With $\radsrc=0$, let $f_{mz}$ denote a term represented by the
third bound in \eqref{eq:four-forces}, and let $f_{mm}$ denote a term
represented by the fourth. They satisfy
$|f_{mz}|\le C\lambda r^{q+2}vF|z^{(j)}|$, $j=1,2$, and
$|f_{mm}|\le C\lambda^2r^{2q+1}vF$.
The quotient $f_{mz}/\mathsf p$ is bounded by
$C\lambda r^q(v/F)|z^{(j)}|$; its multiplier has four bounded
local derivatives of size $C\lambda\delta^{q-2}$.
For $f_{mm}/\mathsf p$, division gives the bounds
$C\lambda^2\delta^{2q-3}e^{-(2q-1)t}$ and,
after multiplication by $q_K$,
$C(\lambda^2/\eps)\delta^{2q}B_Me^{-2qt}$.
The two bounds for $f_{mm}/\mathsf p$ and
$q_Kf_{mm}/\mathsf p$ remain valid after four $t$ derivatives. Hence
\[
\begin{aligned}
 \sum_{\ell=0}^4\|\partial_t^\ell(f_{mz}/\mathsf p)\|_{\mathcal K_{\rm tail,0}}
 &\le C\lambda\delta^{q-2}Z,\\
 \sum_{\ell=0}^4\|\partial_t^\ell(f_{mm}/\mathsf p)\|_{\mathcal K_{\rm tail,0}}
 &\le C[\lambda^2\delta^{2q-3}+(\lambda^2/\eps)\delta^{2q}].
\end{aligned}
\]
The four quantities $C_{f_{mz}},I_{f_{mz}},C_{f_{mm}},I_{f_{mm}}$
are bounded using \eqref{eq:train-moment}:
\[
\begin{aligned}
 C_{f_{mz}}&\le C(\lambda/\eps)
    \|B_M^{-1}r^{q-1}G\|_2\|q_Kz^{(j)}\|_2
       \le C(\lambda/\eps)\delta^qZ,\\
 I_{f_{mz}}&\le C\lambda\eps
    \|B_Mr^{q+1}F\|_2\|q_Kz^{(j)}\|_2
       \le C\lambda\eps\delta^{q+2}Z,\\
 C_{f_{mm}}&\le C(\lambda/\eps)^2\int r^{2q-1}vG
       \le C(\lambda/\eps)^2\delta^{2q+1},\\
 I_{f_{mm}}&\le C\lambda^2\int B_M^2r^{2q+1}vF
       \le C\lambda^2\delta^{2q+3}.
\end{aligned}
\]
Inserting the four displayed bounds for
$C_{f_{mz}},I_{f_{mz}},C_{f_{mm}},I_{f_{mm}}$ in
\eqref{eq:kelvin-data-bound} gives the full data
bounds $C(\lambda/\eps)\delta^{q-1}Z$ and
$C(\lambda/\eps)^2\delta^{2q}$, respectively.
Four source derivatives require only $D^6w,z^{(6)}$.
Thus the terms with $\radsrc=0$ in \eqref{eq:four-forces} satisfy
\begin{equation}
 \|(f,0)\|_{\mathcal Y_K}
 \le C[\delta^{q-3/2}W+\delta^{q-1}Z+\delta^{2q}].
 \label{eq:kelvin-unforced-data}
\end{equation}

The three terms containing $\radsrc$ in \eqref{eq:four-forces} satisfy
one common estimate. Write $f=\mathsf p A(\radsrc/F)$.
The coefficient $A$ includes angular averaging and the bounded
tensor factors.
For $E_j=B_M^2\eps^{-2}vD_t^j\radsrc$, \eqref{eq:weight-bounds} gives
\[
 |D_t^j(\radsrc/F)|\le C B_M^{-2}\sum_{l\le j}|E_l|,
 \qquad |q_KD_t^j(\radsrc/F)|\le C B_M^{-1}\sum_{l\le j}|E_l|.
\]
Indeed, $\eps^2/(vF)\le C$ and $\eps r/F\le C$.
The inequality $\|E_l/B_M\|_2^2\le\|E_l\|_\infty\|E_l\|_1$
and the relative weight derivative bounds prove all four entries of
\[
 \sum_{j=0}^4\|D_t^j(f/\mathsf p)\|_{\mathcal K_{\rm tail,0}}
 \le C\|A\|_{C_{\rm b}^{4,\alphaH}}\mathcal G.
\]
Since $G/v\le C$, the definitions and \eqref{eq:kelvin-work} give
\[
\begin{aligned}
 C_f&=\eps^{-2}\int G|A\radsrc|\le C\|A\|_\infty\mathcal G,\\
 I_f&=\int B_M^2r^2F|A\radsrc|
       \le C\eps^2\|r^2(F/v)A\|_\infty\mathcal G,\\
 \|q_Ke_f\|_2&\le C\delta^{-2}
 \bigl(\|A\|_\infty\|r^2(F/v)A\|_\infty\bigr)^{1/2}\mathcal G.
\end{aligned}
\]
For these three source terms, the coefficient bounds are
\[
\begin{array}{ccc}\toprule
 A\text{, up to bounded coefficients}
 &\|A\|_{C_{\rm b}^{4,\alphaH}}
 &\|r^2(F/v)A\|_\infty\\ \midrule
 \lambda r^{q-1}D^{\le1}w/F
 &C\lambda\delta^qW&C\lambda\delta^{q+2}W\\
 \lambda r^q(v/F)z^{(j)}
 &C\lambda\delta^{q-2}K(z)&C\lambda\delta^{q+2}K(z)\\
 \lambda^2r^{2q-1}(v/F)
 &C\lambda^2\delta^{2q-3}&C\lambda^2\delta^{2q+1}\\
\bottomrule
\end{array}
\]
Here $j\le2$; four derivatives use at most six field or center
derivatives. The first row uses $\|D^{\le5}w/F\|_{C_{\rm b}^{0,\alphaH}}\le C\delta W$.
The other rows retain the loss $v/F\le C\delta^{-2}$, which
cancels in the moment multiplier $r^2(F/v)A$.
Substituting these bounds in \eqref{eq:kelvin-data-bound} gives, for
the response $z_f$ to the sum of these forces with $z_f(0)=0$,
\begin{equation}
 K(z_f)\le
 C\left[\lambda\delta^{q-1}W
        +\lambda\delta^{q-2}K(z)+\lambda^2\delta^{2q-3}\right]\|\radsrc\|_{\mathcal Y_{\rm rad}}.
 \label{eq:forced-small}
\end{equation}
Fix a unit parameter direction. A dot denotes differentiation in
this direction at fixed reference coordinates; reference weights
and neck size remain fixed. Thus
$\dot f=\mathsf p(\dot A \radsrc+A\dot \radsrc)/F$.
Apply \eqref{eq:weighted-leibniz-calculus} to
$\dot A \radsrc+A\dot \radsrc$ to estimate $e_{\dot f}$ in the fixed norms.
The parameter identity \eqref{eq:kelvin-parameter-identity} and the
multiplier bounds \eqref{eq:fowler-reference-multipliers} then
estimate the derivative of $z_f$.

For the quadratic momentum $J_w$ in
\eqref{eq:kelvin-momentum-taylor}, fix $V,w$, put $d_{h_j}=\frac{n-2}{2}h_j-\partial_t h_j$ and define
\[
\begin{aligned}
 (\mathcal B_{V,w}(h_1,h_2))_i
 =nr\int_0^1(1-s)\fint \bigl[
 &-\theta_i d_{h_1}d_{h_2}
 +d_{h_1}\nabla\theta_i\cdot\nabla h_2
 +d_{h_2}\nabla\theta_i\cdot\nabla h_1\\
 &+\theta_i\nabla h_1\cdot\nabla h_2
 -\frac{n(n+2)}4\theta_i(V+sw)^{p-1}h_1h_2\bigr]\dd s.
\end{aligned}
\]
Then $J_w=\mathcal B_{V,w}(w,w)$; the subscripts record the
dependence of the coefficients on $V,w$.
Under \eqref{eq:inner-positive-neighborhood}, the coefficient
$(V+sw)^{p-1}$ and its required normalized derivatives are bounded.
Hence
\[
 |\mathcal B_{V,w}(h_1,h_2)|/\mathsf p
 \le \frac C r\frac{|D^{\le1}h_1|}{F}
                   \frac{|D^{\le1}h_2|}{F}.
\]
The $\mathcal W_6$ bounds on $h_1/F,h_2/F$ and
$2\gamma_*-1>0$ give
\[
 \|\mathcal B_{V,w}(h_1,h_2)/\mathsf p\|_{L^1\cap C_{\rm b}^{0,\alphaH}}
 \le C\delta\mathcal W_6(h_1)\mathcal W_6(h_2).
\]
For the $L^2$ and bounded H\"older norms after
multiplication by $q_K$, we have
\[
\begin{aligned}
 q_K|\mathcal B_{V,w}(h_1,h_2)|/\mathsf p
 &\le C\delta^{-1}\|D^{\le1}h_1/F\|_\infty
       \frac{B_M}{\eps}
       \sqrt{1+\delta^2\omega/r^2}\,|D^{\le1}h_2|,\\
 \|q_K\mathcal B_{V,w}(h_1,h_2)/\mathsf p\|_{L^2\cap C_{\rm b}^{0,\alphaH}}
 &\le C\mathcal W_6(h_1)\mathcal W_6(h_2).
\end{aligned}
\]
Here $v/F=(1+\omega)/(n-2)$, $\omega\le C\delta^{-2}$ and
$r\le\delta$ justify the pointwise estimate for
$q_K|\mathcal B_{V,w}(h_1,h_2)|/\mathsf p$.
Exchanging $h_1,h_2$ gives the symmetric bound. For $\ell\le5$, the highest differentiated factors are
$D^{\ell+1}h_1$, $D^{\ell+1}h_2$, and the normalized derivatives of
$(V+sw)^{p-1}$ through order $\ell$. They require at most six
field derivatives. The product rule therefore gives
\begin{equation}
 \sum_{\ell=0}^5
 \|\partial_t^\ell(J_w/\mathsf p)\|_{\mathcal K_{\rm tail,0}}
 \le CW^2.
 \label{eq:kelvin-quadratic-flux}
\end{equation}
The bilinear version has right-hand side
$C\mathcal W_6(h_1)\mathcal W_6(h_2)$.

The formula for $C(\kelpar')$, the vanishing third spherical moment,
and the definition of $A_w$ give
\[
 \begin{aligned}
 \sum_{\ell=0}^5\|\partial_t^\ell(C-I)\|_{C_{\rm b}^{0,\alphaH}}
 &\le C\delta^2K(z)^2,\\
 \sum_{\ell=0}^5\|\partial_t^\ell A_w\|_{C_{\rm b}^{0,\alphaH}}
 &\le C(\delta W+\delta^2W^2).
 \end{aligned}
\]
Apply the product rule, placing derivatives of $C-I$ and $A_w$ in
the bounded H\"older norm and derivatives of $z'$ in
$\mathcal K_{\rm tail,0}$. This gives
\begin{equation}
 \begin{aligned}
 \sum_{\ell=0}^5\|\partial_t^\ell[(C-I)z']\|_{\mathcal K_{\rm tail,0}}
 &\le C\delta^2K(z)^3,\\
 \sum_{\ell=0}^5\|\partial_t^\ell[A_wz']\|_{\mathcal K_{\rm tail,0}}
 &\le C\delta WK(z).
 \end{aligned}
 \label{eq:kelvin-linear-flux}
\end{equation}
Only $z^{(6)}$ is needed. Since $R_\eps(m_0,0)=v_\tau$,
\eqref{eq:radial-background-first-difference} gives
\[
 \left\|\frac{\mathsf p_V}{\mathsf p_\tau}-1\right\|_
       {C^1_\vartheta(C_{\rm b}^{5,\alphaH})}
 =\left\|\left(\frac{F_V}{F_\tau}\right)^2-1\right\|_
       {C^1_\vartheta(C_{\rm b}^{5,\alphaH})}
 \le C\mathcal G.
\]
Here the background derivative keeps $\radsrc$ fixed.
For $\radsrc=0$, $\mathsf p_V=\mathsf p_\tau$.
The momentum term $\mathsf j$ in \eqref{eq:repair-final-flux} is
\[
 \mathsf j=
 \left(\frac{\mathsf p_V}{\mathsf p_\tau}-1\right)z'
 +\frac{\mathsf p_V}{\mathsf p_\tau}(C-I+A_w)z'
 +(F/F_\tau)^2\frac{J_w}{\mathsf p}.
\]
Use \eqref{eq:kelvin-quadratic-flux}, \eqref{eq:kelvin-linear-flux}
and the multiplier bounds to obtain
\[
 \sum_{\ell=0}^5\|\partial_t^\ell\mathsf j\|_{\mathcal K_{\rm tail,0}}
 \le C\bigl[W^2+(\mathcal G+\delta W+\delta^2Z^2)Z\bigr].
\]
The same product estimates include one background derivative at fixed
$(\radsrc,w,z)$. A source derivative of $\mathsf p_V/\mathsf p_\tau$
is bounded by $C\|\dot \radsrc\|_{\mathcal Y_{\rm rad}}$ and retains the
factor $z'$ in the momentum.

The products estimated in \eqref{eq:radial-product-estimates},
\eqref{eq:kelvin-quadratic-flux}, \eqref{eq:kelvin-linear-flux},
\eqref{eq:kelvin-unforced-data} and \eqref{eq:forced-small}
cover the exact residual \eqref{eq:repair-common-residual}.
To see this, substitute
$V''=\frac{(n-2)^2}{4}V-\frac{n(n-2)}4V^p+\radsrc$.
Then use \eqref{eq:repair-metric-remainders} and the motion Taylor
formula before applying $\Pi_0$ and $\Pi_{\ge2}$. The terms
$\mathcal R_{\rm mov}w$ have an additional bounded multiplier
$rz^{(i)}$. The low selection
identity \eqref{eq:repair-low-selection} also holds after one
variation.

\subsection{Differentiated source estimates and the higher harmonics}
\label{app:source-polarization}

We differentiate the source estimates as maps between the fixed
weighted spaces. The variation replaces a field factor by its direction
and differentiates its coefficients; it is not an additional derivative
in the cylinder variables.

Differentiate each product with the reference weights fixed and
$|\dot\lambda|\le C\lambda$. Then use \eqref{eq:kelvin-parameter-identity} to differentiate the
Kelvin response $z$.

For the metric coefficient tensors in
\eqref{eq:coefficient-inventory}, remove the displayed factors
$\lambda r^q$, $\lambda^2r^{2q}$ and $rz^{(j)}$ before differentiating
the remaining coefficients. These coefficients are independent of
$\radsrc,w$. By Appendix~\ref{app:coordinates}, their derivatives in
a $z$ direction are bounded by $C\delta K(\dot z)$, and their finite
parameter derivatives are bounded by $C$.

The coefficients involving the fields use relative factors such as
$V/v$, $F_V/F$ and $D^{\le1}w/F_V$.
With the boundary data fixed, Lemma~\ref{lem:radial-response}
bounds variations of the first two factors by
$C\|\dot\radsrc\|_{\mathcal Y_{\rm rad}}$ in the required local norms.
The bound
$\|D^{\le1}\dot w/F_V\|_{C_{\rm b}^{0,\alphaH}}\le C\delta\mathcal W_6(\dot w)$
controls the field variation of the third factor; varying its
denominator uses the same relative radial estimate.
Thus unit directions are measured in the component norms of
\eqref{eq:inner-product-norm}.
Replacing one factor by its variation and applying these bounds gives
the first derivative estimates in $\mathcal Y_{\rm rad}$,
$\mathcal W_4$, $\mathcal Y_K$ and $\mathcal K_{\rm tail}$.

In the power remainder, write
$(V+sw)^{p-2}=V^{p-2}(1+sw/V)^{p-2}$.
Since $|w|\le V/2$ and $0\le s\le1$, the last factor is smooth
with uniformly bounded derivatives. Together with the relative
estimates for $V$ and $w$, this gives continuity of the
Fr\'echet derivatives in the fixed spaces.

For $R_\eps(m_0,\radsrc)-R_\eps(m_0,0)$, a background variation at fixed $\radsrc$ retains
$\mathcal G$ by \eqref{eq:radial-background-first-difference};
a source variation uses $D_\radsrc R_\eps(m_0,\cdot)$.

The factorization in \eqref{eq:radial-product-estimates} includes
the highest motion term:
\[
 \|B_M^2\eps^{-2}v\,r z^{(6)}w_t\|_{L^1\cap C_{\rm b}^{0,\alphaH}}
 \le C\|q_Kz^{(6)}\|_{L^2\cap C_{\rm b}^{0,\alphaH}}
            \|\beta_1w_t\|_{L^2\cap C_{\rm b}^{0,\alphaH}}
 \le CK(z)W.
\]

To estimate $\mathcal M_{\rm stat}V$ in $\mathcal W_4$, first
consider the term $\lambda r^qv$. Its two weighted integral terms
and its bound relative to $F$ use
\[
\begin{aligned}
 \beta_1\lambda r^qv
   &=(\lambda/\eps)\delta^qB_Me^{-qt}v,\\
 \beta_2\lambda r^qv
   &=(\lambda/\eps)\delta^qB_Me^{-(q-1)t}v\sqrt\omega,\\
 \delta^{-1}e^{\gamma_*t}\lambda r^qv/F
   &\le C\lambda\delta^{q-3}e^{-(q-\gamma_*)t}.
\end{aligned}
\]
Since $q\ge2>\gamma_*$, $v^2\omega\le C$ and
$\lambda\le C\eps=C\delta^{n/2}$ with $n\ge8$, all local and
integral entries are bounded by $C\delta^q$. Relative weight
derivatives are bounded, so this holds through four derivatives.
The terms $\mathcal M_{\rm mov}V$ and $\mathcal R_{\rm met}V$
contribute the additional factors $rz^{(j)}$ and $\lambda r^q$,
respectively. The coefficient bounds in \eqref{eq:coefficient-inventory} give
\begin{equation}
\begin{aligned}
 \mathcal W_4\bigl(\Pi_{\ge2}(\mathcal M_{\rm stat}
       +\mathcal M_{\rm mov}+\mathcal R_{\rm met})V\bigr)
 &\le C\delta^q(1+\delta Z+\lambda\delta^q)(1+\mathcal G),\\
 \mathcal W_4(\Pi_{\ge2}\mathcal L_{\rm pert}w)
 &\le C(\mathcal G+\lambda\delta^q+\delta Z)W.
\end{aligned}
 \label{eq:high-total-sources}
\end{equation}
The linear coefficients act as bounded local multipliers on every
weighted integral and relative-$F$ entry.
For $N_V(w)$, use \eqref{eq:power-remainder-integral} and
$\|D^{\le6}w/v\|_{C_{\rm b}^{0,\alphaH}}\le C\delta W$. Put one relative
factor in its bounded local norm and the other in each weighted
source norm. The bound $1+sw/V\ge1/2$ controls
$(1+sw/V)^{p-2}$, also for $p<2$, and hence
\begin{equation}
 \mathcal W_4(\Pi_{\ge2}N_V(w))\le C\delta W^2.
 \label{eq:high-power-bound}
\end{equation}
Differentiating \eqref{eq:power-remainder-integral} and
applying the product bounds gives a continuous Fr\'echet derivative
with values in $\mathcal W_4$.

For the radial motion source, the exact identity is
\begin{equation}
\begin{split}
 \mathcal R_{\rm mov}V={}&(\frac{n(n-2)}4V^p-\radsrc)(L^{-2}-1+2\eta)
 +(n-1)F_V\eta(L^{-2}-1)\\
 &-F_V\eta_t(L^{-3}-1)+F_V|\beta|^2L^{-3}.
\end{split}
\label{eq:radial-motion-remainder}
\end{equation}
Taylor's formula gives $L^{-2}-1+2\eta=O(\eta^2)$ and
$L^{-j}-1=O(\eta)$, $j=2,3$, with first derivatives $O(\eta)$
and $O(1)$, respectively. Since $\eta=r\theta\cdot z'$,
$\eta_t=r\theta\cdot(z''-z')$ and $V^p/F_V\le C$,
$\mathcal R_{\rm mov}V$ is bounded by
$Cr^2(F+|\radsrc|)(|z'|^2+|z'z''|)$.
For each product $r^2Fz^{(i)}z^{(j)}$, $i,j\le2$, we have
\[
\begin{aligned}
 |\beta_1r^2Fz^{(i)}z^{(j)}|
 &\le Cr\,|q_Kz^{(i)}|\,|z^{(j)}|,\\
 |\beta_2r^2Fz^{(i)}z^{(j)}|
 &\le C\delta\,|q_Kz^{(i)}|\,|z^{(j)}|,\\
 \delta^{-1}e^{\gamma_*t}r^2|z^{(i)}z^{(j)}|
 &\le C\delta e^{-(2-\gamma_*)t}|z^{(i)}z^{(j)}|.
\end{aligned}
\]
The bound for $\beta_2r^2Fz^{(i)}z^{(j)}$ uses
$\sqrt{FG}\le Cv$. Put one center factor in its unweighted local
norm and the other in its Kelvin norm.
The estimate for $\delta^{-1}e^{\gamma_*t}r^2z^{(i)}z^{(j)}$
controls the relative-$F$ component because $\gamma_*<2$.
Four derivatives use at most $z^{(6)}$; a factor $\radsrc/F$ contributes
$\mathcal G$ by \eqref{eq:radial-source-norm}. Thus
\begin{equation}
 \mathcal W_4(\Pi_{\ge2}\mathcal R_{\rm mov}V)
 \le C\delta(1+\mathcal G)Z^2.
 \label{eq:high-motion-source}
\end{equation}
The bounds \eqref{eq:high-total-sources},
\eqref{eq:high-power-bound} and \eqref{eq:high-motion-source},
and their first Fr\'echet derivatives give the corresponding solution
bounds after applying the inverse on harmonics of degree at least two
from Lemma~\ref{lem:high-tail}.

The remaining factors $\lambda r^q$, $rz^{(j)}$ and
$D^{\le2}w/F$ in \eqref{eq:four-forces} have bounded local
derivatives. Differentiating each factor before applying the
weighted estimates preserves the bounds for $C_f$, $I_f$ and
$\|\partial_t^\ell(f/\mathsf p)\|_{\mathcal K_{\rm tail,0}}$,
and hence the complete
$\mathcal Y_K$ and $\mathcal K_{\rm tail}$ estimates.
The flat center terms are exactly \eqref{eq:repair-final-flux}.

For derivatives with respect to $(\radsrc,w,z)$, fix the boundary
data at $t=0$ and use the product norm
\eqref{eq:inner-product-norm}.
The metric terms and sources in \eqref{eq:selector} retain $\delta^q$,
$\delta^{2q}$ and $\delta^{q-3/2}$; quadratic terms and radial
coefficient variations retain $W$ or $Z$.

The coefficient $\mathsf p_\tau$ is fixed under these
derivatives. Its difference from $\mathsf p_V$ contributes
$C\mathcal G$ to the center derivative and $CZ$ to the radial
source derivative. The other linear coefficients give
$\lambda\delta^q+\delta Z$, while the terms proportional to $\radsrc$
give $c_{\rm pr}$ or $c_{\rm f}$.

Combining these bounds with $W^2\le W$, $Z^2\le Z$ on the small
ball and the bounds for the three solution operators proves
\eqref{eq:source-contraction-bound}.

Finite parameter derivatives and
full boundary-trace derivatives have uniform bounds on the input
ball. The small bound \eqref{eq:source-contraction-bound} applies
to derivatives in $(\radsrc,w,z)$ with the boundary data fixed.

\subsection{Boundary estimates for the differentiated local solutions}

Differentiate the equation for $U$ in a fixed unit parameter
direction. In the derivatives of its coefficients, keep $U$ fixed.
Use the coefficients in \eqref{eq:repair-exact-operator}, and put
$C_0=c-\operatorname{div}b$. The field linearization and coefficient
source are
\[
\begin{aligned}
 \mathcal L&=-\operatorname{div}(A D\cdot)+C_0-\frac{n(n+2)}4dU^{p-1},\\
 Q^{\rm c}&=-\operatorname{div}(\dot A\,DU)+\dot C_0U-\frac{n(n-2)}4\dot d\,U^p.
\end{aligned}
\]
Differentiating
$d\mathcal E=-\operatorname{div}(ADU)+C_0U-\frac{n(n-2)}4dU^p$ gives
\begin{equation}
 d\dot{\mathcal E}=\mathcal L\dot U+Q^{\rm c}-\dot d\,\mathcal E.
 \label{eq:local-jets}
\end{equation}
Taking four spatial derivatives of \eqref{eq:local-jets} uses
at most six derivatives of $U$ and $\dot U$.
The derivative of $U$ is obtained from the fixed zero-trace
parameter equation in Proposition~\ref{prop:inner-graph}.

We estimate the boundary values of the differentiated solutions
in the fixed trace spaces. In this subsection, for an unprojected
source $f$, the notation $R_{H,\tau}f$ means
$R_{H,\tau}(\Pi_{\ge2}f)$.
For the high equation we estimate
$e_0^{-1}\partial_tR_{H,\tau}f(0)$ and its first parameter derivative.
Use the fixed $\sigma_0$ from \eqref{eq:high-decay-choice}.
For the high harmonic inverse with zero trace at $t=0$, we have
\begin{equation}
 \|\partial_tR_{H,\tau}f(0)\|_{C^{5,\alphaH}}
 \le C\sum_{j\ge0}e^{-\sigma_0j}
                   \|f\|_{C^{4,\alphaH}([j,j+1])}.
 \label{eq:boundary-row}
\end{equation}
One-sided boundary Schauder estimates control $[0,1]$,
where $\frac{n(n+2)}4v_\tau^{p-1}$ and its first parameter derivative are
$O(\delta^2)$; the exponential Green bound and homogeneous
boundary regularity control the slabs $[j,j+1]$, $j\ge1$. Thus
$e_0^{-1}\|\partial_tR_{H,\tau}f(0)\|\le C\mathcal W_4(f)$.

For $h_1,h_2$ satisfying $\Pi_{\ge2}h_j=h_j$,
\begin{equation}
\begin{aligned}
 |D^{\le6}h_1|/v&\le C\delta e^{-\gamma_*t}\mathcal W_6(h_1),\\
 \mathcal W_4(h_1h_2/v)&\le
       C\delta\mathcal W_6(h_1)\mathcal W_6(h_2).
\end{aligned}
\label{eq:dn-high-products}
\end{equation}
The power remainder has derivative
\[
 D_wN_V(w)[h_2]
 =p(p-1)wh_2\int_0^1(V+sw)^{p-2}\dd s.
\]
Relative to $wh_2/v$, positivity bounds the coefficient by
$CV^{p-1}$, so \eqref{eq:dn-high-products} applies with $h_1=w$.
Applying \eqref{eq:boundary-row} to $D_wN_V(w)[h_2]$ gives its
contribution to the boundary normal derivative.

A potential coefficient variation has the form
\[
 \dot{(\frac{n(n+2)}4V^{p-1})}=\frac{n(n+2)}4(p-1)V^{p-2}\dot V=O(v^{p-1}),
\]
by the relative first derivative bound in
Lemma~\ref{lem:radial-response}.
Use
\[
 v^{p-1}\le C\min(\delta^2e^{2t},1)\le C\delta e^t,\qquad
 F(t)\le CF(0)e^{(n+2)t/2}.
\]
For a scalar coefficient $a$ satisfying
$\|a/v^{p-1}\|_{C_{\rm b}^{4,\alphaH}}\le C$,
\eqref{eq:boundary-row} and the lower bound in
\eqref{eq:high-decay-choice} give
\begin{equation}
\begin{aligned}
 e_0^{-1}\|\partial_tR_{H,\tau}(ah)(0)\|_{C^{5,\alphaH}}
 &\le C\delta\mathcal W_6(h)
       \int_0^\infty e^{-(\sigma_0-(n+4)/2+\gamma_*)t}\dd t\\
 &\le C\delta\mathcal W_6(h).
\end{aligned}
 \label{eq:dn-potential-gain}
\end{equation}

For the remaining estimates, $a$ denotes a scalar coefficient with
norm $\|a\|=\|a\|_{C_{\rm b}^{4,\alphaH}}$.
The other displayed products include the angular tensor coefficients
of the actual sources, whose $C_{\rm b}^{4,\alphaH}$ norms are uniformly
bounded. Multiplication and contraction with these coefficients are
suppressed in the formulas; in particular, the term written
$\lambda r^q\radsrc$ may have a nonconstant angular coefficient
before the projection $\Pi_{\ge2}$ is applied.
The boundary estimates are
\begin{align}
 \mathcal W_4(\lambda r^qD^{\le2}h)&\le C\lambda\delta^q\mathcal W_6(h),\nonumber\\
 e_0^{-1}\|\partial_tR_{H,\tau}(\lambda r^qva)(0)\|
       &\le C\lambda\delta^{q-1}\|a\|,\nonumber\\
 e_0^{-1}\|\partial_tR_{H,\tau}(\lambda r^q\radsrc)(0)\|
       &\le C\lambda\delta^{q-1}\|\radsrc\|_{\mathcal Y_{\rm rad}},\nonumber\\
 \mathcal W_4(rz^{(j)}D^{\le2}h)&\le C\delta K(z)\mathcal W_6(h),\nonumber\\
 e_0^{-1}\|\partial_tR_{H,\tau}(r^2Fz_1^{(i)}z_2^{(j)})(0)\|
       &\le C\delta K(z_1)K(z_2),\nonumber\\
 e_0^{-1}\|\partial_tR_{H,\tau}(\lambda r^{q+1}vz^{(j)})(0)\|
       &\le C\lambda\delta^qK(z).
 \label{eq:dn-source-gains}
\end{align}
All boundary normal-derivative norms in this display are
$C^{5,\alphaH}(\Sph^{n-1})$.
For the source $\lambda r^qva$, use $v\le Cv(0)e^{(n-2)t/2}$ in
\eqref{eq:boundary-row}. For $\lambda r^q\radsrc$, use the norm of $\radsrc/F$.
For $r^2Fz_1^{(i)}z_2^{(j)}$, use
$r^2F/e_0\le C\delta e^{(n-2)t/2}$.
After extracting the displayed prefactors, the remaining integrals are
\[
 \int_0^\infty e^{-(\sigma_0+q-(n-2)/2)t}\dd t,\qquad
 \int_0^\infty e^{-(\sigma_0+q-(n+2)/2)t}\dd t,\qquad
 \int_0^\infty e^{-(\sigma_0-(n-2)/2)t}\dd t.
\]
All exponents are positive by \eqref{eq:high-decay-choice} and $q\ge2$.
The sources $\lambda r^qD^{\le2}h$ and
$rz^{(j)}D^{\le2}h$ use bounded multipliers in $\mathcal W_4$.
The source $\lambda r^{q+1}vz^{(j)}$ uses
$\int_0^\infty e^{-(\sigma_0+q-(n-2)/2)t}\dd t$ with the additional
factor $rz^{(j)}$.
Replacing $F$ by $\radsrc$ in the source $r^2Fz_1^{(i)}z_2^{(j)}$, we obtain
\[
 e_0^{-1}\|\partial_tR_{H,\tau}(r^2\radsrc z_1^{(i)}z_2^{(j)})(0)\|_{C^{5,\alphaH}}
 \le C\delta\|\radsrc\|_{\mathcal Y_{\rm rad}} K(z_1)K(z_2).
\]
The prefactors in \eqref{eq:dn-source-gains} and in the displayed
bound for $R_{H,\tau}(r^2\radsrc z_1^{(i)}z_2^{(j)})$ are
$O(\delta)$ for
$\lambda\le C\eps=C\delta^{n/2}$ and $q\ge2$.

For the radial boundary estimate, use the unforced comparison
profile $\bar V=v_\tau=R_\eps(m_0,0)$. It has
$\bar V(0)=V(0)=m_0$, $\mathsf h(\bar V,\bar V')=h_\eps$,
and $\bar V'(0)<0$. Put
$D=V-\bar V$, $I_\radsrc=\int_0^\infty V'\radsrc$, and
$s_V=V'(0)+\bar V'(0)$. 
Integrating \(\frac{d}{dt}\mathsf h(V,V')=V'\mathfrak f\) from \(t\) to infinity and using \(\mathsf h(V,V')\to h_\varepsilon\) gives
\[
 \mathsf h(V,V')-h_\eps=-\int_t^\infty V'\radsrc,
 \qquad D'(0)s_V=-2I_\radsrc,
\]
because $D(0)=0$. Since $-\bar V'(0)\asymp v(0)$,
the relative comparison and parameter bounds in
Lemma~\ref{lem:radial-response} give
$|s_V|\asymp v(0)$ and $|\dot s_V|\le Cv(0)$.
The relative response bounds and the weighted source integral give
\[
 |I_\radsrc|\le C\eps^2\|\radsrc\|_{\mathcal Y_{\rm rad}},\qquad
 |\dot I_\radsrc|\le C\eps^2(\|\radsrc\|_{\mathcal Y_{\rm rad}}+\|\dot \radsrc\|_{\mathcal Y_{\rm rad}}).
\]
Differentiating $D'(0)=-2I_\radsrc/s_V$ and using
$\eps^2/(e_0v(0))\asymp\delta$, we obtain
\begin{equation}
 \begin{aligned}
 |D'(0)|&\le C\eps^2v(0)^{-1}\|\radsrc\|_{\mathcal Y_{\rm rad}},\\
 e_0^{-1}|D'(0)|&\le C\delta\|\radsrc\|_{\mathcal Y_{\rm rad}},\\
 e_0^{-1}|\dot D'(0)|&\le C\delta(\|\radsrc\|_{\mathcal Y_{\rm rad}}+\|\dot \radsrc\|_{\mathcal Y_{\rm rad}}).
 \end{aligned}
 \label{eq:dn-radial-gain}
\end{equation}
The bounds in \eqref{eq:dn-radial-gain} are conormal bounds
because $D(0)=\dot D(0)=0$.
All reference weights and the neck remain fixed.

For the Fowler profile with boundary value $\bar V(0)=m_0\asymp v(0)$,
\[
 \bar V'(0)=-m_0
   \left(\frac{(n-2)^2}{4}-\frac{(n-2)^2}{4}m_0^{4/(n-2)}+2h_\eps/m_0^2\right)^{1/2}.
\]
Both $m_0^{4/(n-2)}$ and $h_\eps/m_0^2$ are $O(\delta^2)$,
including their first derivatives in the normalized trace directions at $t=0$. Hence
$\partial_{m_0}\bar V'(0)=-\frac{n-2}{2}+O(\delta^2)$.
The radial conormal for this Fowler profile therefore differs from
the degree-zero harmonic Dirichlet-to-Neumann map by $O(\delta^2)$,
including one parameter derivative.

At $t=0$, the center coefficients satisfy
$\mathsf p_\tau(0)\asymp\mathsf p(0)\asymp\eps^2$.
We use the fixed denominator $\mathsf p(0)$ in the bilinear estimates.
Multiplication by $(F(0)/F_\tau(0))^2$ replaces this denominator by
$\mathsf p_\tau(0)$, including one parameter derivative.
The pointwise bilinear estimate preceding
\eqref{eq:kelvin-quadratic-flux} gives, for $h_1,h_2$ with
$\Pi_{\ge2}h_j=h_j$,
\begin{equation}
 |\mathcal B_{V,w}(h_1,h_2)(0)|/\mathsf p(0)
 \le C\delta\mathcal W_6(h_1)\mathcal W_6(h_2).
 \label{eq:dn-kelvin-product}
\end{equation}

Only the power coefficient of $\mathcal B_{V,w}$ depends on $V,w$.
For a joint direction $(V_1,h)$,
\[
 DJ_w[V_1,h]=2\mathcal B_{V,w}(w,h)
 +(D_{(V,w)}\mathcal B_{V,w}[V_1,h])(w,w).
\]
The coefficient derivative $D_{(V,w)}(V+sw)^{p-1}[V_1,h]$ is
$(p-1)(V+sw)^{p-2}(V_1+sh)$. Positivity and the relative response
bounds control this coefficient by $CV^{p-1}$ in a unit direction.
For a variation of $w$, use $|w|/V\le C$ in the coefficient term
and \eqref{eq:dn-kelvin-product} in the bilinear term to obtain
\[
 |D_wJ_w[h](0)|/\mathsf p(0)
 \le C\delta\mathcal W_6(w)\mathcal W_6(h).
\]
For a parameter derivative of the solved inner family,
$|D^{\le1}w(0)|+|D^{\le1}\dot w(0)|\le C\delta F(0)$.
The displayed formula for $DJ_w[V_1,h]$ gives
$|\dot J_w(0)|/\mathsf p(0)\le C\delta$, with all reference
scales fixed.

The term $e_{f,\vartheta}(0)$ in the center equation satisfies
\[
 |e_{f,\vartheta}(0)|
 =\mathsf p_\tau(0)^{-1}\left|\int_0^\infty f\right|
 \le C\eps^{-2}I_f,\qquad I_f=\int_0^\infty B_M^2|f|.
\]
Since $|D_i\log\mathsf p_\tau(0)|\le C$, its derivative is
bounded by $C\eps^{-2}(I_f+I_{D_i f})$.
For the source in \eqref{eq:selector},
\eqref{eq:selector-data-estimates}
gives $C(\lambda/\eps)\delta^{q+2}W$. For metric motion and square
terms the bounds are
$C[(\lambda/\eps)\delta^{q+2}Z+(\lambda/\eps)^2\delta^{2q+3}]$.
The forced terms contribute at most
$C[\lambda\delta^{q+2}W+\lambda\delta^{q+2}Z+
\lambda^2\delta^{2q+1}]\mathcal G$.
Here $\lambda\le C\eps$.
Every term and its first input derivative is $O(\delta)$ on the
fixed bounded ball. For the solved inner family,
$\mathsf p_\tau(0)[z'(0)+\mathsf j(0)]=\int_0^\infty f$.
Substituting \eqref{eq:repair-final-flux} gives the exact center flux
\begin{equation}
 \mathsf p_V(0)\left.[C(\kelpar')+A_w(\kelpar')]\right|_{t=0}z'(0)
       =\int_0^\infty f-J_w(0).
 \label{eq:actual-cap-flux}
\end{equation}
The matrix
$\mathsf p_V(0)\left.[C(\kelpar')+A_w(\kelpar')]\right|_{t=0}$
and its first derivative
are $O(\eps^2)$, with inverse
$O(\eps^{-2})$. The force and bilinear bounds give
\[
 \left|\int f\right|+\left|\int \dot f\right|
       +|J_w(0)|+|\dot J_w(0)|\le C\eps^2\delta.
\]
Solving \eqref{eq:actual-cap-flux} and then differentiating it,
with the neck and reference weights fixed, yields
\[
 |z'(0)|\le C\delta,\qquad
 |\dot z'(0)|\le C\eps^{-2}
       (\eps^2\delta+\eps^2|z'(0)|)\le C\delta.
\]

The contribution of $z'(0)$ to the degree-one coefficient of the
normalized conormal is $-(F_V(0)/v(0))z'(0)$ by
\eqref{eq:dn-normalization}.
The coefficient $F_V(0)/v(0)$ and its first parameter derivative are
bounded, so this term is $O_{C^1}(\delta)$.
For the degree-one Jacobi functions
$j_\tau(t)\theta_i$, where $j_\tau=rF_\tau$, we have
\[
 -\frac{\frac{n-2}{2}j_\tau+j_\tau'}{j_\tau}
 =1-\frac{n(n-2)}4\frac{v_\tau^p}{F_\tau}
 =1+O(\delta^2)\quad\text{at }t=0.
\]
The last bound includes one parameter derivative by
\eqref{eq:fowler-reference-multipliers}. The remaining metric and
forced contributions satisfy \eqref{eq:dn-source-gains}.

The physical conormal, with the fixed reference density and endpoint
orientation of \eqref{eq:repair-raw-density}, is
\[
 \mathcal N=\sum_{\nu=0}^{n-1}(A)^{t\nu}U_\nu+(b)^tU.
\]
Differentiate both the field and its coefficients. The bounds in
Appendix~\ref{app:coordinates},
\eqref{eq:dn-high-products}--\eqref{eq:dn-source-gains},
\eqref{eq:dn-radial-gain} and \eqref{eq:actual-cap-flux} give the
$C^{5,\alphaH}$ conormal bounds of
Proposition~\ref{prop:inner-dn}, including one parameter derivative.

We finish on the fixed annulus
$\mathcal A=\{1/2<|\widehat y|<1\}$, where $\widehat y=y/\delta$.
Use the chart $\Psi$ and the trace coordinates from the proof of
Proposition~\ref{prop:inner-dn}, and the affine function
$B=B_\vartheta$ in \eqref{eq:collar-affine-reference}. Put
\[
 g_{\vartheta,\trpar}=\Psi^*[h_b(\delta\,\cdot)],\qquad
 c=\delta^2c(n)\operatorname{Scal}_{h_b}(\delta\Psi).
\]
In this collar calculation, the dot fixes $W$ and differentiates in a
unit direction of $(\vartheta,\trpar)$.
The map $\trlow(\vartheta,\trpar)$ of finite boundary coordinates
has bounded
first derivative, and the chart estimates give
\[
 \begin{aligned}
 \|g_{\vartheta,\trpar}-I\|_{C^1(\vartheta,\trpar;C^{5,\alphaH})}&\le C\delta,
 &\|c\|_{C^1(\vartheta,\trpar;C^{4,\alphaH})}&\le C\lambda\delta^q,\\
 \|\nabla B\|_{C^1(\vartheta,\trpar;C^{5,\alphaH})}&\le C\delta,
 &\|B\|_{C^1(\vartheta,\trpar;C^{6,\alphaH})}&\le C.
 \end{aligned}
\]
For $U=B+s_\eps W$, division by the fixed $s_\eps\asymp\delta$ gives
\begin{equation}
 \begin{aligned}
 -\Delta W+\mathcal Q(\vartheta,\trpar,W)&=0,\\
 \mathcal Q(\vartheta,\trpar,W)
 &=(-\Delta_{g_{\vartheta,\trpar}}+\Delta+c)W\\
 &\quad+\frac{-\Delta_{g_{\vartheta,\trpar}}B+cB-\frac{n(n-2)}4\delta^2B^p}{s_\eps}\\
 &\quad-\frac{n(n-2)\delta^2}{4s_\eps}[(B+s_\eps W)^p-B^p].
 \end{aligned}
 \label{eq:collar-nonlinear-map}
\end{equation}
Choose $c_B>0$ with $B\ge c_B$, and consider every input in a fixed
bounded ball in $C^{6,\alphaH}(\overline{\mathcal A})$
with $B+s_\eps W\ge c_B/2$.
Since $D^2B=D^2\dot B=0$ and $DB,D\dot B=O(\delta)$,
\[
 \|\Delta_{g_{\vartheta,\trpar}}B\|_{C^1(\vartheta,\trpar;C^{4,\alphaH})}
 \le C\delta^2,\qquad
 c/s_\eps=O(\lambda\delta^{q-1}),\quad \delta^2/s_\eps=O(\delta).
\]
The first, second and third summands in
\eqref{eq:collar-nonlinear-map} have the respective
$C^{4,\alphaH}$ bounds
\[
 C(\delta+\lambda\delta^q)\|W\|_{C^{6,\alphaH}},
 \qquad C(\delta+\lambda\delta^{q-1}),
 \qquad C\delta^2\|W\|_{C^{6,\alphaH}}.
\]
The first and second summands obey their displayed bounds after one
parameter derivative. For the third summand,
the field derivative is $-\frac{n(n+2)}4\delta^2(B+s_\eps W)^{p-1}H=O(\delta^2)H$,
and at fixed $W$,
\[
 \dot{[(B+s_\eps W)^p-B^p]}
 =p[(B+s_\eps W)^{p-1}-B^{p-1}]\dot B
 =O(s_\eps|W||\dot B|).
\]
H\"older multiplication and positivity make
$D_W\mathcal Q,D_\vartheta\mathcal Q,D_\trpar\mathcal Q$ continuous. On the whole input ball,
\begin{equation}
 \|\mathcal Q\|_{C^1(\vartheta,\trpar,W;C^{4,\alphaH})}\le C\delta.
 \label{eq:collar-first-source-bound}
\end{equation}

Use Euclidean surface measure on the fixed reference spheres
$|\widehat y|=1/2,1$, with increasing-radius orientation.
Write $\widehat r=|\widehat y|$, $\widehat r_-=1/2$ and
$\widehat r_+=1$.
The physical conormal density and its Euclidean reference are
\[
 \mathcal D_\pm^{g_{\vartheta,\trpar}}V
 =\left.\widehat r_\pm\sqrt{\det g_{\vartheta,\trpar}}
       \sum_{i,j}g_{\vartheta,\trpar}^{ij}\theta_i\partial_jV\right|_\pm,
 \qquad \mathcal D_\pm^0V=\widehat r_\pm\partial_{\widehat r}V|_\pm.
\]
The tail uses this Euclidean reference surface density at
$|\widehat y|=1/2$.
Since
\[
 \|\sqrt{\det g_{\vartheta,\trpar}}\,g_{\vartheta,\trpar}^{-1}-I
       \|_{C^1(\vartheta,\trpar;C^{5,\alphaH})}\le C\delta,
\]
the operator $\mathcal D_\pm^{g_{\vartheta,\trpar}}-\mathcal D_\pm^0$
is a first-order differential operator with
$O_{C^1(C^{5,\alphaH})}(\delta)$ coefficients. For every input
in the stated ball, where $B+s_\eps W\ge c_B/2$,
\begin{equation}
 \begin{aligned}
 \mathcal C_\pm(\vartheta,\trpar,W)
 &=s_\eps^{-1}[\mathcal D_\pm^{g_{\vartheta,\trpar}}(B+s_\eps W)
                                    -\mathcal D_\pm^0B]\\
 &=\mathcal D_\pm^0W+\mathcal A_\pm(\vartheta,\trpar,W),\\
 \mathcal A_\pm(\vartheta,\trpar,W)
 &=(\mathcal D_\pm^{g_{\vartheta,\trpar}}-\mathcal D_\pm^0)W
 +s_\eps^{-1}(\mathcal D_\pm^{g_{\vartheta,\trpar}}-\mathcal D_\pm^0)B.
 \end{aligned}
 \label{eq:collar-raw-conormals}
\end{equation}
Using $DB,D\dot B=O(\delta)$ in the affine term gives
\[
 \|\mathcal A_\pm\|_{C^1(\vartheta,\trpar,W;C^{5,\alphaH})}
 \le C\delta(1+\|W\|_{C^{6,\alphaH}})+Cs_\eps^{-1}\delta^2\le C\delta.
\]

\refstepcounter{section}
\section*{Appendix C. Parameter differentiation and nonlinear estimates}
\label{app:calculus}
\addcontentsline{toc}{section}{Appendix C. Parameter differentiation and nonlinear estimates}
We collect the identities used to differentiate inverses and fixed
points and to estimate nonlinear powers.

Parameter derivatives are taken on fixed Banach spaces after pulling
back moving domains. For a differentiable family of bounded invertible
operators $A$, differentiating $AG=I$ gives
\begin{equation}
 G_i=-GA_iG,\qquad G=A^{-1}.
 \label{eq:inverse-parameter-calculus}
\end{equation}
For a jointly $C^1$ contraction $u=T(\vartheta,u)$ with
$\|T_u\|\le\theta_0<1$, its fixed point is $C^1$ and
\begin{equation}
 (I-T_u)u_i=T_i,\qquad
 \|u_i\|\le(1-\theta_0)^{-1}\|T_i\|.
 \label{eq:fixed-point-parameter-calculus}
\end{equation}
Here $T_i$ is the derivative with respect to $\vartheta_i$, with
$u$ fixed. Its norm is taken in the Banach space containing the solution $u$.

For $p>1$ and $|s|\le1/2$, Taylor's theorem applied to
$H_p(s)=(1+s)^p-1-ps$ gives
\begin{equation}
 |H_p(s)|\le Cs^2,\qquad |H_p'(s)|\le C|s|,\qquad
 |H_p^{(j)}(s)|\le C_j\quad(j\ge2).
 \label{eq:relative-power-calculus}
\end{equation}
For $V>0$ and $|w|\le V/2$, the exact remainder is
\[
 (V+w)^p-V^p-pV^{p-1}w
 =p(p-1)w^2\int_0^1(1-s)(V+sw)^{p-2}\dd s
 =V^pH_p(w/V).
\]
This is the Taylor formula used in
\cite[Lemma 3.2 and the proof of Lemma 5.3]{part2}.
For the scaled annular norms of Section~\ref{sec:low-construction},
H\"older multiplication gives
\begin{equation}
 \|fg\|_{j,\alphaH;\beta+\gamma;\rho}
 \le C\|f\|_{j,\alphaH;\beta;\rho}\|g\|_{j,\alphaH;\gamma;\rho}.
 \label{eq:weighted-product-calculus}
\end{equation}
The same estimate holds for the exterior norms with $;\rho$ replaced
by $;\Om$, with a constant independent of the small matching radii.

\section*{AI assistance and interactions}
\phantomsection
\addcontentsline{toc}{section}{AI assistance and interactions}

This work was a natural  continuation of our previous work as posted in  Part I \cite{part1}
and Part II \cite{part2}. In this  process, we used GPT, including 
 ChatGPT for Academic Researchers at a later stage, as interactive research-assistance tools for literature reading, exploratory proof searches, calculations, checking arguments, and manuscript preparation. The coauthors 
determined the research questions, supplied and verified the citations and attributions, and selected, checked, and substantially revised every adopted argument, and are responsible for all mathematical claims, citations, and final content. 

\section*{Added note on subsequent work related to Theorem~\texorpdfstring{\ref{thm:main}}{1.1}}\label{sec:addednote}
\phantomsection
\addcontentsline{toc}{section}{Added note on subsequent work related to Theorem~\texorpdfstring{\ref{thm:main}}{1.1}}

We obtained the finite-point existence theorem and its construction in early October 2026. The introduction to Part~II~\cite{part2}, whose first arXiv version was submitted on October~17, 2026,
announced existence with Fowler asymptotics at each point of any
prescribed nonempty finite set in a nonround positive Yamabe conformal class.
On October 22, 2026, UTC, Chen and Ma sent us a preliminary manuscript, dated October 19, concerning the same problem. In that version, their existence result was restricted to any prescribed finite set of two or more points either in dimensions \(3\le n\le 7\), or in dimensions \(n\ge8\) and the Weyl tensor is nonzero at each prescribed point. In our response, we informed them of the statement of our Theorem~\ref{thm:main},  and that we were in the process of revising the present paper. Later that day, they sent us a substantially revised manuscript containing results closer to those subsequently posted on arXiv  on
October 22, 2026, UTC, as  \cite{CM26}.

Chen and Ma~\cite[Theorem~1.1]{CM26} prove that, on a smooth closed connected manifold of dimension $n\ge3$ with $Y(M,[g])>0$, the complement of any prescribed finite set of at least two points admits a complete conformal metric of positive constant scalar curvature.
That paper acknowledges the announcement in \cite{part2} and uses the
geometric quadratic-form analysis developed in \cite{part1}.
Its construction uses Fowler-type ends and an exterior conformal factor
modeled on a small positive sum of Green functions, tending to zero
away from the punctures. In contrast, the present construction is made
near the entire compact set of normalized Yamabe minimizers, without
requiring an individual minimizer to be isolated or nondegenerate.
It prescribes arbitrary fixed positive ratios of the exact asymptotic
Fowler neck sizes and yields subsequential convergence away from the
punctures to a smooth positive Yamabe minimizing metric.
Theorem~\ref{thm:main} gives the full existence characterization,
including a single prescribed singularity in a nonround positive
conformal class.

\section*{Appendix. Announcement of the extension to closed singular sets}\label{sec:closed-set-announcement}
\addcontentsline{toc}{section}
{Appendix. Announcement of the extension to closed singular sets}

During the revision of this paper, we obtained an extension from finite singular sets to arbitrary nonempty closed sets with empty interior. We announce the principal results here; their proofs will appear in \cite{part4}. Although the forthcoming work includes the finite-set existence characterization, we publish the present paper first to give a separate account of the construction for prescribed isolated singularities. This construction is also used in  \cite{part4}, and presenting it separately allows the reader to follow its methods before encountering the additional arguments needed to treat accumulation points. The existence criterion is expressed in terms of the Bessel capacity defined below.

For $s>0$ and $1<q<\infty$, let $\mathcal G_s$ denote the Bessel kernel on
$\R^n$ with Fourier multiplier $(1+|\xi|^2)^{-s/2}$. For a compact set
$K\subset\R^n$, the Bessel capacity $C_{s,q}(K)$ is defined as the infimum of
$\|f\|_{L^q(\R^n)}^q$ over all nonnegative functions $f\in L^q(\R^n)$ satisfying
$\mathcal G_s*f\ge1$ on an open neighborhood of $K$. On a smooth manifold $M$,
we say that a compact set $\Lambda\subset M$ has zero $C_{s,q}$-capacity, and
write $C_{s,q}(\Lambda)=0$, if $\varphi(K)$ has zero $C_{s,q}$-capacity in
$\R^n$ for every smooth coordinate chart $(U,\varphi)$ on $M$ and every compact
subset $K\subset\Lambda\cap U$.

Write
  \[
   L_g=-\Delta_g+\frac{n-2}{4(n-1)}\operatorname{Scal}_g,
   \qquad \kappa_n=\frac{n(n-2)}4,
   \qquad p=\frac{n+2}{n-2}.
  \]
The theorem and proposition numbers below refer to~\cite{part4}. 
\begin{namedthm}{Theorem 1.1~\cite{part4}}
Let $(M^n,g)$ be a smooth, connected, closed Riemannian manifold of
dimension $n\ge3$, and let $\Lambda\subsetneq M$ be a nonempty closed set
with empty interior. Then the following assertions are equivalent:
\begin{enumerate}[label=\textup{(\roman*)}]
\item There exists a nonnegative function $u\in L^p(M)$, smooth and positive
on $M\setminus\Lambda$, with singular set equal to $\Lambda$, such that
$u^{4/(n-2)}g$ is complete on $M\setminus\Lambda$ and $L_gu=\kappa_n u^p$ on
$M$ in the sense of distributions.
\item The complement $M\setminus\Lambda$ admits a complete metric of
constant positive scalar curvature in $[g]|_{M\setminus\Lambda}$.
\item One has
\[
 \begin{gathered}
 Y(M,[g])>0,\qquad C_{1+2/n,n/2}(\Lambda)=0,\\
 (M,[g])\not\simeq(\Sph^n,[g_{\mathrm{rd}}])\quad\text{or}\quad|\Lambda|\ge2.
 \end{gathered}
\]
\end{enumerate}
\end{namedthm}

In \textup{(i)}, having singular set equal to $\Lambda$ means that $u$ is
unbounded on $U\setminus\Lambda$ for every open neighborhood $U$ of every point
of $\Lambda$, and the distributional equation on $M$ reads
\begin{equation}\label{eq:global-weak}
 u\in L^p(M),\qquad
 \int_M u L_g\varphi\dd V_g=\kappa_n\int_M u^p\varphi\dd V_g
 \quad\text{for every }\varphi\in C^\infty(M).
\end{equation}

We remark that the condition $C_{1+2/n,n/2}(\Lambda)=0$ refines the Hausdorff-dimension criterion: it holds when \(\dim_H\Lambda<(n-2)/2\), fails when \(\dim_H\Lambda>(n-2)/2\), and can either hold or fail at equality.
 Furthermore, the necessity argument proves a somewhat stronger statement than is needed for Theorem  1.1~\cite{part4}: any complete conformal metric of nonnegative scalar curvature on the complement forces the same capacity to vanish. Many of the ingredients of the current manuscript are retained in proving Theorem  1.1~\cite{part4}, but there are also substantial new ingredients.  We also obtain several capacity and distributional results, which may be of independent interest. 

When (ii) holds, one can choose a complete metric whose conformal
    factor satisfies \eqref{eq:global-weak} on all of $M$. For the sets
    in Proposition~C.3 of \cite{part4}, there are also complete metrics
    whose conformal factors have a nonzero defect measure.
    
\begin{namedthm}{Proposition C.1~\cite{part4}}
Assume $Y(M,[g])>0$ and let $\Lambda$ be as in Theorem  1.1~\cite{part4}.
The condition $C_{1+2/n,n/2}(\Lambda)=0$ is equivalent to the existence of a
complete scalar-flat metric in $[g]|_{M\setminus\Lambda}$, and also to the
existence of a complete metric in $[g]|_{M\setminus\Lambda}$ whose scalar
curvature is positive at every point of $M\setminus\Lambda$.
\end{namedthm}

\begin{namedthm}{Proposition C.2~\cite{part4}}
For every $\gamma>0$ there is a compact Cantor set $K_\gamma$ in a
coordinate ball of $M$ with $\dim_H K_\gamma=(n-2)/2$.
It carries a Borel probability measure $\sigma$ such that
\begin{equation}\label{eq:gauge}
 c\,h(r)\le\sigma\bigl(B(z,r)\bigr)\le C\,h(r),\qquad
 h(r)=r^{(n-2)/2}\bigl(\log(e/r)\bigr)^{-\gamma},
\end{equation}
for all $z\in K_\gamma$ and $0<r<r_0$. With $p'=p/(p-1)=(n+2)/4$, the capacity
thresholds are
\begin{equation}\label{eq:three-thresholds}
 \begin{aligned}
  C_{1,(n+2)/2}(K_\gamma)>0&\iff\gamma>n/2,\\
  C_{1+2/n,n/2}(K_\gamma)>0&\iff\gamma>(n-2)/2,\\
  C_{2,p'}(K_\gamma)>0&\iff\gamma>(n-2)/4.
 \end{aligned}
\end{equation}
\end{namedthm}

\begin{namedthm}{Proposition C.3~\cite{part4}}
Let $K_\gamma$ and $\sigma$ be as in Proposition C.2~\cite{part4}.
Assume $Y(M,[g])>0$ and $(n-2)/4<\gamma\le(n-2)/2$. For each sufficiently small
$t>0$, there exists
$u_t\in L^p(M)$, smooth and positive on $M\setminus K_\gamma$,  such that
$u_t^{4/(n-2)}g$ is complete and has constant scalar curvature $n(n-1)$ on
$M\setminus K_\gamma$, and
\begin{equation}\label{eq:measure}
 L_gu_t=\kappa_n u_t^p+t\sigma\quad\text{on }M
\end{equation}
in the sense of distributions. Let $G_g=L_g^{-1}$ be the Green operator. Writing $V=G_g\sigma$, one has
$tV\le u_t\le 2tV$ on $M\setminus K_\gamma$. Moreover, any nonnegative global
weak solution $u$ of the homogeneous equation \eqref{eq:global-weak} bounded
above by $CV$ is bounded and smooth on all of $M$.
\end{namedthm}

\begin{namedthm}{Proposition C.4~\cite{part4}}
Let $h>0$ be a smooth solution of $L_gh=\kappa_n h^p$ on $M$.
If $u\in L^p(M)$ is a distributional supersolution satisfying
$L_gu\ge\kappa_n u^p$ on $M$ and $u\ge h$ almost everywhere on $M$, then $u=h$
almost everywhere on $M$.
\end{namedthm}
This rules out constructing a distinct global distributional singular solution simply by keeping it everywhere above a fixed smooth positive Yamabe solution.

\Needspace{16\baselineskip}

\bigskip
\begingroup
\small
\Needspace{7\baselineskip}

\noindent\textsc{Zheng-Chao Han}\\
Department of Mathematics, Rutgers University,\\
Hill Center--Busch Campus, 110 Frelinghuysen Road,\\
Piscataway, NJ 08854, USA\\
\textit{Email address:}
\href{mailto:zchan@math.rutgers.edu}{\nolinkurl{zchan@math.rutgers.edu}}

\medskip
\Needspace{6\baselineskip}

\noindent\textsc{Qinfeng Jiang}\\
School of Mathematical Sciences, Beijing Normal University,\\
Beijing 100875, China\\
\textit{Email address:}
\href{mailto:202531130031@mail.bnu.edu.cn}%
{\nolinkurl{202531130031@mail.bnu.edu.cn}}

\medskip
\Needspace{7\baselineskip}

\noindent\textsc{Hua-Yang Wang}\\
School of Mathematical Sciences, Laboratory of Mathematics and Complex
Systems, MOE,\\
Beijing Normal University, Beijing 100875, China\\
\textit{Email address:}
\href{mailto:wanghuayang@amss.ac.cn}{\nolinkurl{wanghuayang@amss.ac.cn}}

\medskip
\Needspace{10\baselineskip}

\noindent\textsc{Jingang Xiong}\\
School of Mathematical Sciences, Laboratory of Mathematics and Complex
Systems, MOE,\\
Beijing Normal University, Beijing 100875, China;\\
Center for Basic Mathematics, Institute for Advanced Study,\\
Beijing Normal University, Beijing 100875, China\\
\textit{Email address:}
\href{mailto:jx@bnu.edu.cn}{\nolinkurl{jx@bnu.edu.cn}}

\medskip
\Needspace{7\baselineskip}

\noindent\textsc{Lei Zhang}\\
Department of Mathematics, University of Florida,\\
358 Little Hall, P.O. Box 118105,\\
Gainesville, FL 32611-8105, USA\\
\textit{Email address:}
\href{mailto:leizhang@ufl.edu}{\nolinkurl{leizhang@ufl.edu}}
\endgroup

\end{document}